\documentclass[11pt]{amsart}
\usepackage{amsmath,amsthm,epsf,amscd,amssymb,verbatim, eufrak}

\newtheorem{thm}{Theorem}[section]
\newtheorem{lem}[thm]{Lemma}
\newtheorem{cor}[thm]{Corollary}
\newtheorem{prop}[thm]{Proposition}

\newtheorem{prp}[thm]{Proposition}

\newtheorem{lma}[thm]{Lemma}
 
\newtheorem{rem}[thm]{Remark}
\newenvironment{rmk}{\begin{rem}\rm}{\end{rem}}

\newenvironment{pf}{\begin{proof}}{\end{proof}}

\theoremstyle{definition}
\newtheorem{dfn}[thm]{Definition}

\theoremstyle{remark}
\newtheorem{remark}[thm]{Remark}
\newtheorem{ex}[thm]{Example}

\numberwithin{equation}{section}

\newcommand{\bbr}{\begin{remark}}        
\newcommand{\eer}{\end{remark}}

\newcommand{\df}[1]{{\em{#1}}}

\font\bbb=msbm10 scaled 1100

\newcommand{\bea}{\begin{eqnarray}}
\newcommand{\eea}{\end{eqnarray}}
\newcommand{\bmini}{\begin{center}\begin{minipage}{5in}}
\newcommand{\emini}{\end{minipage}\end{center}}
\newcommand{\R}{{\mbox{\bbb R}}}
\newcommand{\C}{{\mbox{\bbb C}}}
\newcommand{\Z}{{\mbox{\bbb Z}}}

\newcommand{\e}{{\mathbf{e}}}
\newcommand{\f}{{\mathbf{f}}}
\newcommand{\g}{{\mathbf{g}}}

\newcommand{\A}{{\mathcal{A}}}
\newcommand{\B}{{\mathcal{B}}}
\newcommand{\U}{{\mathbf{U}}}
\newcommand{\conf}{{\mathcal{C}}}

\newcommand{\func}{{\mathcal{F}}}

\newcommand{\Ordo}{{\mathcal{O}}}
\newcommand{\cand}{{\mathcal{W}}}
\newcommand{\sblv}{{\mathcal{H}}}
\newcommand{\M}{{\mathcal{M}}}
\newcommand{\sgn}{\operatorname{sign}}

\newcommand{\W}{W}
\newcommand{\Disk}{\bar{\Delta}}
\newcommand{\projL}{\Pi_{\C}}
\newcommand{\rar}{\rightarrow}

\newcommand{\dbar}{{\bar{\partial}}}
\newcommand{\ep}{\epsilon}
\newcommand{\de}{\delta}

\newcommand{\hf}{{\tfrac12}}
\newcommand{\GL}{{\mathbf{GL}}}
\newcommand{\action}{{\mathcal{Z}}}
\newcommand{\End}{\operatorname{End}}
\newcommand{\lk}{\operatorname{lk}}

\newcommand{\la}{\langle}
\newcommand{\ra}{\rangle}
\newcommand{\pa}{\partial}
\newcommand{\id}{\operatorname{id}}

\newcommand{\ind}{\operatorname{Index}}
\newcommand{\inr}{\operatorname{int}}
\newcommand{\spa}{\operatorname{Span}}
\newcommand{\krn}{\operatorname{Ker}}
\newcommand{\img}{\operatorname{Im}}
\newcommand{\tb}{\operatorname{tb}}
\newcommand{\diag}{\operatorname{Diag}}
\newcommand{\sign}{\operatorname{Sign}}

\newcommand{\supp}{\operatorname{supp}}
\newcommand{\proj}{\operatorname{pr}}
\newcommand{\ev}{\operatorname{ev}}
\newcommand{\ham}{\operatorname{Ham}}
\newcommand{\pham}{\operatorname{pHam}}
\newcommand{\area}{\operatorname{Area}}
\newcommand{\pih}{\tfrac{\pi}{2}}

\renewcommand{\Im}{\operatorname{Im}}
\renewcommand{\Re}{\operatorname{Re}}

\begin{document}

\title{Legendrian Submanifolds in $R^{2n+1}$ and Contact Homology}

\author{Tobias Ekholm}
\address{Department of Mathematics, Uppsala University, 751 06
   Uppsala, Sweden} 

\author{John Etnyre}
\address{Department of Mathematics, University of Pennsylvania, 209 South 33rd Street, 
	Philadelphia PA 19105-6395}

\author{Michael Sullivan}
\address{Department of Mathematics, University of Michigan,
525 East University Avenue,
Ann Arbor MI 48109-1109}


\begin{abstract}
Contact homology for Legendrian submanifolds 
in standard contact $(2n+1)$-space is rigorously
defined using moduli spaces of holomorphic disks with
Lagrangian boundary conditions in complex $n$-space. 
The homology provides new invariants of Legendrian
isotopy. 
These invariants show that the theory of Legendrian isotopy 
is very rich. 
For example, they detect infinite families of pairwise non-isotopic
Legendrian $n$-spheres, $n$-tori, and surfaces 
which are indistinguishable using previously known invariants.

In a sense, the definition of contact homology presented in this paper
is a high dimensional analog of the work of Chekanov and others on
Legendrian $1$-knots in $3$-space.  
\end{abstract}

\maketitle

\section{Introduction}

A contact manifold is a $(2n+1)$-manifold $N$ equipped with a
completely non-integrable field of hyperplanes $\xi$. An immersion of
an $n$-manifold into $N$ is {\em Legendrian} if it is everywhere
tangent to the hyperplane field $\xi$ and the image of a Legendrian
embedding is a {\em Legendrian submanifold}.  
{\em Standard contact $(2n+1)$-space} is 
Euclidean space $\R^{2n+1}$  
equipped with the hyperplane field $\xi=\krn(\alpha)$, where $\alpha$
is the contact $1$-form $\alpha=dz-\sum_{i=1}^ny_i\,dx_i$ in Euclidean
coordinates $(x_1,y_1,\dots,x_n,y_n,z)$. 

Any closed $n$-manifold $M$ embeds in $\R^{2n+1}$, and it is a consequence
of the h-principle for Legendrian immersions \cite{GromovPDR} that, provided $M$ meets
certain homotopy theoretic conditions (which is the case e.g. if $M$ is
stably parallelizable), any embedding of $M$ into $\R^{2n+1}$ may be
arbitrarily well $C^0$-approximated by Legendrian embeddings. Thus,
Legendrian submanifolds of standard contact $(2n+1)$-space exist in abundance.

Any contact manifold of dimension $2n+1$ is locally contactomorphic (diffeomorphic
through a map which takes contact hyperplanes to contact hyperplanes)
to standard contact $(2n+1)$-space. In this paper we study 
local Legendrian knotting phenomena or, in other
words, the question: {\em When are are two Legendrian
submanifolds of standard 
contact $(2n+1)$-space isotopic through Legendrian submanifolds?} 

For $n=1$, the question above has been extensively studied,
\cite{Chekanov, Eliashberg-Fraser, Etnyre-Honda1, Etnyre-Honda2, ENS}. 
Here,
the {\em classical invariants} of a Legendrian knot are its
topological knot type, its rotation number (the tangential degree of
the curve which arises as the projection of the knot into the
$xy$-plane), and its Thurston-Bennequin invariant (the linking number
of the knot and a copy of the knot shifted slightly in the
$z$-direction). Many examples of Legendrian non-isotopic knots
with the same classical invariants are known. Also, in higher
dimensions, when the ambient contact manifold has more topology (for
example Legendrian knots in $1$-jet spaces of $S^n$) there are interesting
examples of non-trivial Legendrian knots \cite{Eliashberg-Gromov97}.
  
For $n>1$, we define two {\em classical invariants} of an oriented
Legendrian submanifold given by an embedding 
$f\colon L\to\R^{2n+1}$. 
Following \cite{Tabachnikov88}, we first define
its Thurston-Bennequin invariant 
(in the same way as in $\R^3$). 
Second, we note
that the h-principle for Legendrian immersions implies that $f$ is
determined by certain homotopy theoretic invariants, associated 
to its differential $df$, up to regular homotopy through Legendrian
immersions, see Section \ref{1rotclass} for details. We define its 
{\em rotation class} as its Legendrian regular 
homotopy class. The topological embedding invariant in
the $3$-dimensional case disappears in higher dimensions 
since, for $n\ge 2$, any two embeddings of an $n$-manifold into
$\R^{2n+1}$ are isotopic \cite{Hae}. 

Our results indicate that the theory of Legendrian submanifolds of
standard contact $(2n+1)$-space is very rich. For example we show,
generalizing the $3$-dimensional results mentioned above,    

\begin{thm}\label{1intromain}
For any $n>1$ there is an infinite family of Legendrian embeddings 
of the $n$-sphere into $\R^{2n+1}$ that are not Legendrian isotopic even though 
they have the same classical invariants.
\end{thm}

We also prove a similar theorem for Legendrian surfaces and $n$-tori.
We show that for any $N>0$ there exists Legendrian isotopy classes
of $n$-spheres and $n$-tori with fixed Thurston-Bennequin invariants 
and rotation classes which do not admit a representative having
projection into $\R^{2n}$ with less than $N$ double points.
 
Note that Theorem \ref{1intromain} is not known to be true for $n=1$
and is probably false in this case. It is known 
that the number of distinct Legendrian knots with the same classical
invariants can be arbitrarily 
large, but in light of recent work of Colin, Giroux and Honda \cite{CGH}
it seems unlikely that there can be infinitely many.

To show that Legendrian submanifolds are not Legendrian isotopic we
develop the {\em contact homology} of a Legendrian submanifold in
standard contact $(2n+1)$-space. It is defined using
punctured holomorphic disks in $\C^n\approx\R^{2n}$ with boundary on the
projection of the Legendrian submanifold, and which limit to
double points of the projection at the punctures, see below. This is
analogous to 
the approach taken by Chekanov \cite{Chekanov} in dimension $3$;
however, in that 
dimension the entire theory can be reduced to combinatorics \cite{ENS}.  
Our contact homology, realizes in the language of 
Symplectic Field Theory \cite{EGH}, Relative Contact Homology of
standard contact $(2n+1)$-space and in this
framework our main technical theorem can be summarized as

\begin{thm}
The contact homology of Legendrian submanifolds in $\R^{2n+1}$ with
the standard contact form is well defined.
(It is invariant under Legendrian isotopy.)
\end{thm} 

More concretely, if $L\subset\R^{2n+1}$ is a Legendrian submanifold we
associate to $L$ a differential graded algebra $(\A,\pa)$, freely generated
by the double points of the projection of $L$ into $\C^n$. The
differential $\pa$ is defined by counting rigid holomorphic disks with
properties as described above. Thus, contact homology is similar to Floer
homology of Lagrangian intersections and our proof of its invariance
is similar in spirit to Floer's original approach \cite{Floer88a,
Floer88b} in the following way. We analyze bifurcations of moduli
spaces of rigid holomorphic disks under variations of the 
Legendrian submanifold in a generic
$1$-parameter family of Legendrian submanifolds and how these
bifurcations affect the differential graded algebra. Similar
bifurcation analysis is also 
done in \cite{Hutchings, Lee, Sullivan99, Sullivan}. 
Our set-up does not seem well suited to the more popular proof of
Floer theory invariance which uses an elegant 
``homotopy of homotopies'' argument (see, for 
example, \cite{Floer89, Salamon-Zehnder92}).

{\em Acknowledgments:} Much of the collaboration was carried out during
the 2000 Contact Geometry Workshop 
at the American Institute of Mathematics (AIM). The authors are
grateful for AIM's support. 
We also benefited from conversations with F.~Bourgeois, Y.~Chekanov,
K.~Cielebak, Y.~Eliashberg, D.~Salamon, P.~Seidel, and K.~Wysocki.
TE is a Research Fellow of the Royal Swedish Academy of Sciences
supported by the Knut and Alice Wallenberg foundation and
was supported in part by the Swedish research council. 
JE was supported in part by NSF Grant \# DMS-9705949 and \# DMS-0203941 and
MS was supported in part by NSF Grant \# DMS-9729992 and a NSF VIGRE Grant.
MS is also grateful to the Institute for Advanced Study, ETH in
Zurich, and the University of Michigan 
for hosting him while he worked on this project.

\tableofcontents

\part{Contact Homology in $\R^{2n+1}$}

\section{Contact Homology and Differential Graded Algebras}\label{1contacthomology}

In this section we describe how to associate to a Legendrian
submanifold $L$ in standard contact $(2n+1)$-space
a differential graded algebra (DGA) $(\mathcal{A}, \partial).$ Up to a certain equivalence relation
this DGA is an invariant of the Legendrian isotopy class of $L$. In Section~\ref{1First} we recall the 
notion of Lagrangian projection and define the algebra $\mathcal{A}.$
The grading on $\mathcal{A}$ is described in  
Section~\ref{1CZ} after a review of the Maslov index in
Section~\ref{1Maslov}. 

Sections \ref{1modulisection} and \ref{1Differential.section} are devoted to 
the definition of $\partial$ and Section \ref{1Invariance.section}
proves the invariance of the homology of  
$(\mathcal{A}, \partial),$ which we call the contact homology.
The main proofs of these three subsections rely on much analysis, which will be deferred until Part 2.
Finally, in Section \ref{1SFT.section}, we compare contact homology as
defined here with the contact homology sketched in \cite{EGH}.

\subsection{The algebra $\mathcal{A}$}\label{1First}

Throughout this paper we consider the standard contact structure $\xi$ on 
$\R^{2n+1}=\C^n\times\R$ which
is the hyperplane field given as the kernel of the contact 1-form
\begin{equation} \label{1alpha.eqn}
	\alpha=dz-\sum_{j=1}^n y_j dx_j,
\end{equation}
where $x_1,y_1,\ldots, x_n,y_n, z$ are Euclidean coordinates on $\R^{2n+1}.$ 
A \df{Legendrian submanifold} of $\R^{2n+1}$ is an $n$ dimensional
submanifold $L\subset \R^{2n+1}$ everywhere tangent to $\xi$. 
We also recall that the standard symplectic structure on $\C^n$ is given by 
\begin{equation} \notag
	\omega=\sum_{j=1}^n dx_j\wedge dy_j,
\end{equation}
and that an immersion $f\colon L\to\C^n$ of an $n$-dimensional
manifold is {\em Lagrangian} if $f^\ast\omega=0$. 

The \df{Lagrangian projection} projects out the $z$ coordinate:
\begin{equation} \label{1LagProj.eqn}
	\Pi_{\C}\colon \R^{2n+1}\to \C^n;
	\quad (x_1,y_1,\ldots,x_n,y_n,z)\mapsto (x_1,y_1,\ldots,x_n,y_n).
\end{equation}
If $L\subset\C^n\times\R$ is a Legendrian submanifold then 
$\Pi_{\C}\colon L\to\C^n$ is a Lagrangian immersion. Moreover,
for $L$ in an open dense subset of all Legendrian
submanifolds (with $C^\infty$ topology), the self intersection of
$\Pi_{\C}(L)$ consists of a finite number of transverse double
points. We call Legendrian submanifolds with this property {\em chord generic}. 

The Reeb vector field $X$ of a contact form  $\alpha$ is uniquely
defined by the two equations $\alpha(X) = 1$ and $d \alpha(X, \cdot) = 0.$ 
The {\em Reeb chords} of a Legendrian submanifold $L$ are segments of flow
lines of $X$ starting and ending at points of $L.$ 
We see from (\ref{1alpha.eqn}) that in $\R^{2n+1}$,
$X=\frac{\partial}{\partial z}$ and thus  
$\Pi_{\C}$ defines a bijection between Reeb chords of $L$ and double
points of $\Pi_\C(L)$. If $c$ is a Reeb chord we write $c^\ast=\Pi_\C(c)$.


Let $\mathcal{C}=\{c_1,\ldots, c_m\}$ be the set of Reeb chords 
of a chord generic Legendrian submanifold $L\subset\R^{2n+1}.$ 
To such an $L$ we associate an algebra $\mathcal{A}=\mathcal{A}(L)$
which is the free associative unital algebra over the group ring
$\Z_2[H_1(L)]$ generated by $\mathcal{C}$.  We write elements in $\mathcal{A}$ as
\begin{equation}\label{DGAasVSP}
\sum_i  t_1^{n_{1,i}}\ldots t_k^{n_{k,i}}{\mathbf c}_i,
\end{equation}
where the $t_j$'s are formal variables corresponding to a basis 
for $H_1(L)$ thought of multiplicatively and 
${\mathbf c}_i=c_{i_1}\dots c_{i_r}$ is a word in the generators.  
It is also useful to consider the corresponding algebra $\mathcal{A}_{\Z_2}$ over $\Z_2.$ 
The natural map $\Z_2[H_1(L)]\to \Z_2$ induces a reduction of
$\mathcal{A}$ to $\mathcal{A}_{\Z_2}$ (set $t_j=1$, for all $j$).

\subsection{The Maslov index}\label{1Maslov}
Let $\Lambda_n$ be the Grassman manifold of Lagrangian subspaces in
the symplectic vector space $(\C^n,\omega)$  
and recall that $H_1(\Lambda_n)=\pi_1(\Lambda_n)\cong\Z.$ There is a standard isomorphism 
\begin{equation} \notag
	\mu:H_1(\Lambda_n)\to\Z,
\end{equation}
given by intersecting a loop in $\Lambda_n$ with the Maslov cycle $\Sigma.$ To describe $\mu$ more fully
we follow \cite{Robbin-Salamon93} and refer the reader to this paper for proofs of the statements
below. 

Fix a Lagrangian subspace $\Lambda$ in $\C^n$ and let
$\Sigma_k(\Lambda)\subset\Lambda_n$ be the subset of Lagrangian spaces 
that intersects $\Lambda$ in a subspace of $k$ dimensions. The \df{Maslov cycle} is 
\begin{equation} \notag
	\Sigma=\overline{\Sigma_1(\Lambda)}=
	\Sigma_1(\Lambda)\cup\Sigma_2(\Lambda)\cup\dots\cup\Sigma_{n}(\Lambda).  
\end{equation}
This in an algebraic variety of codimension one in
$\Lambda_n.$ If $\Gamma:[0,1]\to \Lambda_n$ is a loop then
$\mu(\Gamma)$ is the intersection number of $\Gamma$ and 
$\Sigma.$  The contribution of an intersection point $t'$ with
$\Gamma(t')\in\Sigma$ to $\mu(\Gamma)$ is calculated as follows.  
Fix a Lagrangian complement $W$ of $\Lambda.$ Then for each $v\in \Gamma(t')\cap \Lambda$
there exists a vector $w(t)\in W$ such that $v+w(t)\in \Gamma(t)$ for $t$ near $t'.$ Define
the quadratic form $Q(v)=\frac{d}{dt}\vert_{t=t'} \omega(v, w(t))$ on $\Gamma(t')\cap \Lambda$ and
observe that it is independent of the complement $W$ chosen. Without loss of generality,
$Q$ can be assumed non-singular and the contribution of the
intersection point to $\mu(\Gamma)$ is the signature 
of $Q.$ Given any loop $\Gamma$ in $\Lambda_n$ we say $\mu(\Gamma)$ is the \df{Maslov index} of the loop.

If $f\colon L\to\C^n$ is a Lagrangian immersion then the tangent
planes of $f(L)$ along any loop $\gamma$ in $L$ gives a loop $\Gamma$ in $\Lambda_n$. 
We define the Maslov index $\mu(\gamma)$ of $\gamma$ as
$\mu(\gamma)=\mu(\Gamma)$ and note that we may view the Maslov index as
a map $\mu: H_1(L)\to \Z.$ Let $m(f)$ be the smallest non-negative
number that is the Maslov index of some 
non-trivial loop in $L.$ We call $m(f)$ the \df{Maslov number} of $f.$
When $L\subset\C^n\times\R$ is a Legendrian submanifold we write
$m(L)$ for the Maslov number of $\Pi_\C\colon L\to\C^n$.

\subsection{The Conley--Zehnder index of a Reeb chord and the grading on $\mathcal{A}$}\label{1CZ}
Let $L\subset\R^{2n+1}$ be a chord generic Legendrian submanifold and let
$c$ be one of its Reeb chords with end points $a,b\in L$,
$z(a)>z(b)$. Choose a path $\gamma:[0,1]\to L$ with $\gamma(0)=a$ and
$\gamma(1)=b$. (We call such path a {\em capping path of $c$}.)  
Then $\Pi_{\C}\circ \gamma$ is a loop in 
$\C^n$ and $\Gamma(t)=d\Pi_\C(T_{\gamma(t)}L)$, $0\le t\le 1$ is a
path of Lagrangian subspaces of $\C^n$. Since $c^\ast=\Pi_\C(c)$ is a transverse double 
point of $\Pi_\C(L)$, $\Gamma$ is not a closed loop.

We close $\Gamma$ in the following way. Let $V_0=\Gamma(0)$
and $V_1=\Gamma(1)$. Choose any complex structure $I$ on $\C^n$
which is compatible with $\omega$ ($\omega(v,Iv)>0$ for all $v$) and
with $I(V_1)=V_0$. (Such an $I$ exists since the 
Lagrangian planes are transverse.) Define the path 
$\lambda(V_0,V_1)(t)=e^{tI}V_1$, $0\leq t\leq \frac{\pi}{2}.$ 
The concatenation, $\Gamma\ast\lambda(V_0,V_1),$ of 
$\Gamma$ and $\lambda(V_0,V_1)$ forms a loop in $\Lambda_n$
and we define the \df{Conley--Zehnder index}, $\nu_\gamma(c),$ of $c$ to be
the Maslov index $\mu(\Gamma\ast\lambda(V_0,V_1))$ of this loop.
It is easy to check that $\nu_\gamma(c)$ is independent of the choice of $I$.
However, $\nu_\gamma(c)$ might depend on the choice of homotopy class of the path $\gamma.$ 
More precisely, if $\gamma_1$ and $\gamma_2$ are two paths with 
properties as $\gamma$ above then 
$$
\nu_{\gamma_1}(c)-\nu_{\gamma_2}(c)=\mu(\gamma_1\ast(-\gamma_2)),
$$
where $(-\gamma_2)$ is the path $\gamma_2$ traversed in the opposite direction.
Thus $\nu_\gamma(c)$ is well defined modulo the Maslov number $m(L).$ 

Let $\mathcal{C}=\{c_1,\ldots, c_m\}$ be the set of Reeb chords of
$L$. Choose a capping path $\gamma_j$ for each $c_j$
and define the \df{grading} of $c_j$ to be 
\[|c_j|=\nu_{\gamma_j}(c_j)-1,\]
and for any $t\in H_1(L)$ define its grading to be $|t|=-\mu(t).$ This makes $\mathcal{A}(L)$ into
a graded ring. Note that the grading depends on the choice of capping
paths but, as we will see below, this choice will be irrelevant.

The above grading on Reeb chords $c_j$ taken modulo $m(L)$ makes   
$\mathcal{A}_{\Z_2}$ a graded algebra with grading in $\Z_{m(L)}.$ (Note that this grading does not
depend on the choice of capping paths.) 
In addition the map from $\mathcal{A}$ to $\mathcal{A}_{\Z_2}$ preserves gradings modulo $m(L).$

\subsection{The moduli spaces}\label{1modulisection}
As mentioned in the introduction, the differential of the algebra
associated to a Legendrian submanifold is defined using spaces of
holomorphic disks. To describe these spaces we need a few preliminary definitions.

Let $D_{m+1}$ be the unit disk
in $\C$ with $m+1$ punctures at the points $p_0,\ldots p_m$ on the
boundary. The orientation of the boundary of the unit disk induces a
cyclic ordering of the punctures. Let $\pa {\hat D}_{m+1}=\pa D_{m+1}\setminus\{p_0,\dots,p_m\}$.

\begin{figure}[ht]
	{\epsfxsize=2.75in\centerline{\epsfbox{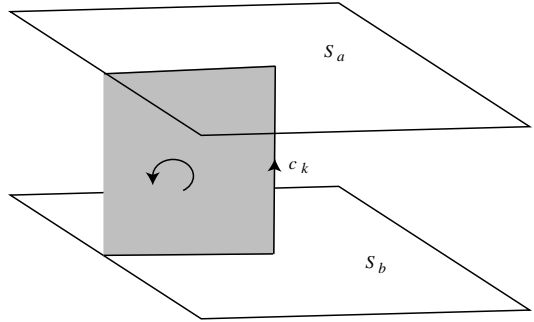}}}
	\caption{Positive puncture lifted to $\R^{2n+1}.$ The gray region is the holomorphic disk and
		the arrows indicate the orientation on the disk and the Reeb chord.}
	\label{1fig:sign}
\end{figure}

Let $L\subset\C^n\times\R$ be a Legendrian submanifold with isolated
Reeb chords. 
If $c$ is a Reeb chord of $L$ with end points $a,b\in L$, $z(a)>z(b)$ then
there are small neighborhoods $S_a\subset L$ of $a$ and $S_b\subset L$ of $b$ 
that are mapped injectively to $\C^n$ by $\Pi_{\C}.$ 
We call $\Pi_{\C}(S_a)$ the \df{upper sheet} of $\Pi_{\C}(L)$ 
at $c^\ast$ and $\Pi_{\C}(S_b)$ the \df{lower sheet}. 
If $u\colon (D_{m+1}, \partial D_{m+1})\to (\C^n,\Pi_{\C}(L))$ is a
continuous map with $u(p_j)=c^\ast$ then we say $p_j$
is \df{positive} (respectively \df{negative}) if $u$ maps points   
clockwise of $p_j$ on $\partial D_{m+1}$ 
to the lower (upper) sheet of $\Pi_{\C}(L)$ and points anti-clockwise
of $p_i$ on $\partial D_{m+1}$ to the upper (lower) sheet of $\Pi_{\C}(L)$ (see Figure~\ref{1fig:sign}).

If $a$ is a Reeb chord of $L$ and if ${\mathbf b}=b_1\ldots b_m$ is an
ordered collection (a word) of Reeb chords then let 
$\mathcal{M}_A(a;{\mathbf b})$ be the space, modulo conformal reparameterization, 
of maps $u:(D_{m+1}, \partial D_{m+1})\to (\C^n,\Pi_{\C}(L))$ which
are continuous on $D_{m+1}$, holomorphic in the interior of
$D_{m+1}$, and which have the following properties  
\begin{itemize}
\item $p_0$ is a positive puncture, $u(p_0)=a^\ast,$
\item $p_j$ are negative punctures for $j>0,$ $u(p_j)=b_j^\ast,$  
\item the restriction $u|\pa {\hat D}_{m+1}$ has a continuous lift 
$\tilde u\colon\pa {\hat D}_{m+1} \to L\subset\C^n\times\R,$ and
\item the homology class of $\tilde u(\partial D_{m+1}^\ast)\cup (\cup_j \gamma_j)$
equals $A\in H_1(L)$,
\end{itemize}
where $\gamma_j$ is the capping path chosen for $c_j$,
$j=1,\dots,m$. Elements in $\M_A(a;{\mathbf b})$ will be called 
{\em holomorphic disks with boundary on $L$} or sometimes simply
holomorphic disks. 

There is a useful fact relating heights of Reeb chords and the area of
a holomorphic disk with punctures mapping
to the corresponding double points. 
The \df{action} (or height) $\action(c)$ of a Reeb chord 
$c$ is simply its length and the action of a word of Reeb chords 
is the sum of the actions of the chords making up the word.
\begin{lem}\label{1height}
	If $u\in\mathcal{M}_A(a;{\mathbf b})$ then
	\begin{equation} \label{1action.eqn}
	\action(a)-\action({\mathbf{b}})=\int_{D_m} u^*\omega = \area(u) \ge 0.
	\end{equation}
\end{lem}

\begin{proof}
        By Stokes theorem,
	$\int_{D_m} u^*\omega=\int_{\partial D_m} u^*(-\sum_j y_jdx_j)=
	\int {\tilde u}^\ast(-dz)=\action(a)-\action({\mathbf b}).$ The second equality
	follows since $u$ is holomorphic and $\omega=\sum_{j=1}^n dx_j\wedge dy_j.$
\end{proof}

Note that the proof of Lemma \ref{1height} implies that any
holomorphic disk with boundary on $L$ must 
have at least one positive puncture. (In contact homology, only disks
with exactly one positive puncture are considered.)

We now proceed to describe the properties of moduli spaces $\mathcal{M}_A(a;{\mathbf b})$ that
are needed to define the differential.  In Part 2 we prove that the
moduli spaces of holomorphic disks with boundary on a
Legendrian submanifold $L$ have these properties provided $L$ is
generic among (belongs to a Baire subset of the space of) 
{\em admissible} Legendrian submanifolds ($L$ is admissible if it is
chord generic and it is real analytic in a neighborhood of all Reeb chord end points).
For more precise definitions of these concepts, see Section \ref{5Isotopies.section},
where it is shown that admissible Legendrian
submanifolds are dense in the space of all Legendrian
submanifolds. In Section \ref{2fas}, we express moduli spaces
$\M_A(a;{\mathbf b})$ as $0$-sets of certain $C^1$-maps between
infinite-dimensional Banach manifolds. We say a moduli space is {\em
transversely cut out} if $0$ is a regular value of the corresponding
map.   

\begin{prop}\label{1manifold}
	For a generic admissible Legendrian submanifold
	$L\subset\C^n\times\R$ the moduli space  
	$\mathcal{M}_A(a;{\mathbf b})$ is a transversely cut out manifold of dimension
	\begin{equation}\label{1dimension}
	d=\mu(A) + |a|-|{\mathbf b}|-1,
	\end{equation}
	provided $d\le 1$. (In particular, if $d<0$ then the moduli space
	is empty.)
\end{prop}
Proposition~\ref{1manifold} is proved in
Subsection~\ref{8gentransv.section}. If $u\in\M_A(a;{\mathbf b})$ we
say 
that $d=\mu(A)+|a|-|{\mathbf b}|$ is the {\em formal dimension of $u$},
and if $v$ is a transversely cut out disk of formal dimension $0$ we
say that $v$ is a {\em rigid disk}. 


The moduli spaces we consider might not be compact, but their lack of
compactness can be understood. It is analogous to ``convergence to broken
trajectories'' in Morse/Floer homology and gives rise to natural 
compactifications of the moduli spaces. This is also called Gromov
compactness, which we cover in more detail in Section
\ref{9compactness.section}.

A {\it broken holomorphic curve}, $u = (u^1, \ldots , u^N)$, is a 
union of holomorphic disks, 
$u^j: (D_{m_j}, \partial D_{m_j}) \rightarrow (\C^n, \projL(L))$, 
where each $u^j$ has exactly one positive puncture $p^j$, with the
following property. To each $p^j$ with $j\ge 2$ is associated a
negative puncture $q^k_j\in D_{m_k}$ for some $k\ne j$ such that
$u^j(p^j)=u^k(q^k_j)$ and $q^{k'}_{j'}\ne q^{k}_{j}$ if $j\ne
j'$, and such that the quotient space 
obtained from $D_{m_1}\cup\dots\cup D_{m_N}$ by identifying $p^j$ and
$q^k_j$ for each $j\ge 2$ is contractible.  
The broken curve can be parameterized by a single smooth
$v : (D_m, \partial D) \rightarrow (\C^n, \projL(L)).$
A sequence $u_\alpha$ of holomorphic disks {\em converges}
to a broken curve $u = (u^1, \ldots, u^N)$ if the
following holds:

\begin{itemize}
\item
For every $j \le N$, there exists
a sequence $\phi_\alpha^j:D_m \rar D_m$ of linear fractional transformations
and a finite set $X^j \subset D_m$ such that $u_\alpha \circ \phi^j_\alpha$
converges to $u^j$ uniformly with all derivatives on compact subsets
of $D_m \setminus X^j$
\item
There exists a sequence of 
orientation-preserving diffeomorphisms $f_\alpha:D_m \rar D_m$ such that
$u_\alpha \circ f_\alpha$ converges in the $C^0$-topology to a parameterization
of $u$.
\end{itemize}

\begin{prop}\label{1compactness}
	Any sequence $u_\alpha$ in $\mathcal{M}_A(a;{\mathbf b})$ has a subsequence	
	converging to a broken holomorphic curve $u = (u^1, \ldots, u^N)$.
	Moreover,  $u^j\in \mathcal{M}_{A_j}(a^j; {\mathbf b}^j)$ with
	$A = \sum_{j=1}^N A_j$ 
and
\begin{equation} \label{1subadditivity.eqn}
	\mu(A) + |a|-|{\mathbf b}| = 	
	\sum_{j=1}^N \left( \mu(A_j) + |a^j|-|{\mathbf b}^j| \right).
\end{equation} 
\end{prop}

Heuristically this is the only type of non-compactness we expect to see in $\M_A(a;{\mathbf b})$: 
since $\pi_2(\C^n) =0$, no holomorphic spheres can ``bubble off'' at an interior
point of the sequence $u_\alpha$, and since $\Pi_\C(L)$ is exact no
disks without positive puncture can form either. Moreover,
since $\Pi_\C(L)$ is compact,
and since $\C^n$ has ``finite geometry at infinity" (see Section \ref{9compactness.section}),
all holomorphic curves with a uniform bound on area must map to a compact set.

\begin{proof}
The main step is to prove convergence to some broken curve, which we defer to Section
\ref{9compactness.section}. The statement about the homology classes follows easily from
the definition of convergence. Equation~\eqref{1subadditivity.eqn} follows from
the definition of broken curves.
\end{proof}

We next show that a broken curve can be glued to form a family of
non-broken curves. For this we need a little notation. 
Let ${\mathbf c^1},\dots,{\mathbf c^r}$ be an ordered collection of
words of Reeb chords. Let the length of (number of
letters in) ${\mathbf c^j}$ be $l(j)$
and let ${\mathbf a}=a_1 \dots a_k$ be a word of
Reeb-chords of length $k>0$. Let $S=\{s_1,\dots,s_r\}$ be $r$ distinct
integers in $\{1,\dots,k\}$. Define the word
${\mathbf a}_S({\mathbf c^1},\dots,{\mathbf c^r})$ of Reeb-chords of
length $k-r+\sum_{j=1}^rl(j)$ as 
follows. For each index $s_j\in S$ remove $a_{s_j}$ from the word
${\mathbf a}$ and insert at its place the word ${\mathbf c^j}$. 
\begin{prop}\label{1gluingI}
        Let $L$ be a generic admissible Legendrian submanifold. Let
        $\M_A(a;{\mathbf{b}})$ and $\M_B(c; {\mathbf d})$ be
        $0$-dimensional transversely cut out moduli spaces and assume
        that the $j$-th Reeb chord in ${\mathbf b}$ is $c$.  
	Then there exist a $\rho>0$ and an embedding
	\[
		G\colon \M_A(a;{\mathbf{b}})\times\M_B(c; {\mathbf d})\times(\rho,\infty)
	\to
		\M_{A+B}(a; {\mathbf b}_{\{j\}}({\mathbf d})).
	\]
	Moreover, if  $u\in\M_A(a;{\mathbf b})$ and 
	$u'\in \M_B(c; {\mathbf d})$ then
	$G(u,u',\rho)$ converges to the broken curve $(u,u')$ as
        $\rho\to \infty,$ and any disk in $\M_A(a;{\mathbf
        b}_{\{j\}}({\mathbf d}))$ with image sufficiently close
	to the image of $(u,u')$ is in the image of $G.$ 
\end{prop}
This follows from Proposition~\ref{glud^2=0} and the definition of
convergence to a broken curve.

\subsection{The differential and contact homology} \label{1Differential.section}
Let $L\subset\C^n\times\R$ be a generic admissible Legendrian submanifold,
let $\mathcal{C}$ be its set of Reeb chords, and
let $\mathcal{A}$ denote its algebra. 
For any generator $a\in \mathcal{C}$ of $\mathcal{A}$ we set 
\begin{equation}
	\partial a
	=\sum_{\hbox{dim } \mathcal{M}_A(a;{\mathbf b})=0} (\#\mathcal{M}_A(a;{\mathbf b}))A {\mathbf b},
\end{equation}
where $\# \mathcal{M}$ is the number of points in $\mathcal{M}$ modulo
$2$, and where the sum ranges over all words ${\mathbf b}$ in the
alphabet $\mathcal{C}$ and $A\in H_1(L)$ for which the above moduli
space has dimension 0. 
We then extend $\partial$ to a map 
$\partial:\mathcal{A}\to \mathcal{A}$ by linearity and the Leibniz rule. 


Since $L$ is generic admissible, it follows from 
Propositions~\ref{1compactness} and \ref{1gluingI} 
that the moduli spaces considered in the definition of $\partial$ are
compact $0$-manifolds and hence consist of a finite number of points. 
Thus $\partial$ is well defined. Moreover,
\begin{lem}\label{1dsquare=0}
	The map $\partial:\mathcal{A}\to\mathcal{A}$ is a differential
	of degree $-1$. That is,
	$\partial\circ \partial=0$
	and $|\partial (a)|=|a|-1$ for any generator $a$ of $\mathcal{A}.$
\end{lem}

\begin{proof}
After Propositions~\ref{1compactness} and \ref{1gluingI} the standard
proof in Morse (or Floer) homology \cite{SalamonNotes} applies. It
follows from (\ref{1dimension}) that $\partial$ lowers degree by 1. 
\end{proof}

The \df{contact homology of} $L$ is 
\begin{equation} \notag
	HC_*(\R^{2n+1},L)= \hbox{Ker } \partial/ \hbox{Im } \partial.
\end{equation}
It is essential to notice that since $\partial$ respects the grading on 
$\mathcal{A}$ the contact homology is a graded algebra.

We note that $\partial$ also defines a differential of degree $-1$ on $\mathcal{A}_{\Z_2}(L).$

\subsection{The invariance of contact homology under Legendrian isotopy} \label{1Invariance.section}
In this section we show 
\begin{prop}\label{1invarianthomology}
        If $L_t\subset\R^{2n+1}$, $0\le t\le 1$ is a Legendrian
	isotopy  between generic admissible Legendrian submanifolds then
	the contact homologies $HC_*(\R^{2n+1},L_0),$ and
	$HC_*(\R^{2n+1},L_1)$ are isomorphic.
\end{prop}

In fact we show something, that at least appears to be, stronger. 
Given a graded algebra $\mathcal{A}=\Z_2[G]\la a_1,\ldots, a_n\ra,$
where $G$ is a finitely generated abelian group, a graded
automorphism $\phi:\mathcal{A}\to\mathcal{A}$ is called 
\df{elementary} if there is some $1\leq j\leq n$ such that 
\begin{equation} \notag
    \phi(a_i) = 
    \begin{cases}
      A_ia_i, & i \neq j \\
      \pm A_ja_j + u, & u \in
      \A(a_1,\ldots,a_{j-1},a_{j+1},\ldots,a_n),\ i=j,
    \end{cases}
\end{equation}
where the $A_i$ are units in $\Z_2[G].$
The composition of elementary automorphisms is called a \df{tame} automorphism. 
An isomorphism from $\mathcal{A}$ to $\mathcal{A}'$ is tame if it is
the composition of a tame automorphism with 
an isomorphism sending the generators of $\mathcal{A}$ to the generators of $\mathcal{A}'.$
An isomorphism of DGA's is called \df{tame} if the isomorphism of the underlying algebras is tame.

Let $(\mathcal{E}_i, \partial_i)$ be a DGA
with generators $\{e_1^i, e_2^i\},$ where 
$|e_1^i|=i, |e_2^i|=i-1$ and $\partial_ie_1^i=e_2^i, \partial_ie_2^i=0.$  
Define the \df{degree $i$ stabilization} 
$S_i(\mathcal{A},\partial)$  of $(\mathcal{A}, \partial)$ 
to be the graded algebra generated by $\{a_1,\ldots, a_n, 
e_1^i, e_2^i\}$ with grading and differential induced from $\mathcal{A}$ and $\mathcal{E}_i.$ 
Two differential graded algebras are called \df{stable tame isomorphic} if they become tame isomorphic
after each is stabilized a suitable number of times. 
\begin{prop}\label{1stabletame}
	If $L_t\subset\R^{2n+1}$, $0\le t\le 1$ is a Legendrian
	isotopy  between generic admissible Legendrian submanifolds then
	the DGA's $(\A(L_0),\pa)$ and $(\A(L_1),\pa)$ are stable tame isomorphic.
\end{prop}

Note that Proposition~\ref{1stabletame} allows us to
associate the stable tame isomorphism class of a DGA to a 
Legendrian isotopy class of Legendrian submanifolds: any Legendrian
isotopy class has a generic admissible 
representative and by Proposition~\ref{1stabletame} the DGA's of any
two generic admissible representatives agree.

It is straightforward to show that two stable tame isomorphic DGA's
have the same homology, see \cite{Chekanov, ENS}. Thus
Proposition~\ref{1invarianthomology} follows from Proposition
\ref{1stabletame}. The proof of the later given below is, in outline,
the same as the proof of invariance of the stable tame isomorphism
class of the DGA of a Legendrian $1$-knot in \cite{Chekanov}. However,
the details in our case require considerably more work. In particular
we must substitute analytic arguments for the purely combinatorial
ones that suffice in dimension three. 

In Section \ref{5Isotopies.section} we show that any two admissible Legendrian
submanifolds of dimension $n>2$ which are Legendrian isotopic are isotopic through a
special kind of Legendrian isotopy: a Legendrian isotopy
$\phi_t\colon L\to\C^n\times\R$, $0\le t\le 1$, is \df{admissible} if $\phi_0(L)$ and
$\phi_1(L)$ are admissible Legendrian submanifolds and if there exist a
finite number of instants $0<t_1<t_2<\dots<t_m<1$ and a $\delta>0$
such that the intervals $[t_j-\delta,t_j+\delta]$ are disjoint subsets
of $(0,1)$ with the following properties.
 
\begin{figure}[ht]
	{\epsfxsize=4in\centerline{\epsfbox{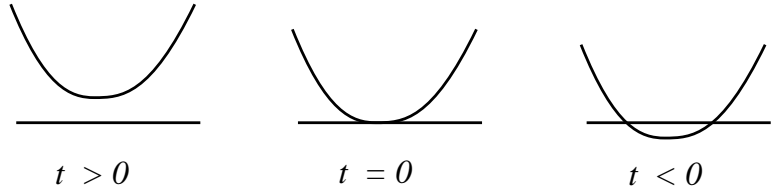}}}
	\caption{Type B double point move.}
	\label{1fig:TypeB}
\end{figure}

\begin{itemize}
\item[(A)]  
	For $t\in[0,t_1-\delta]\cup\Bigl(\bigcup_{j=1}^m[t_j+\delta,
	t_{j+1}-\delta]\Bigr)\cup [t_m+\delta,1]$, $\phi_t(L)$ is an isotopy
	through admissible Legendrian submanifolds.  
\item[(B)] 
    	For $t\in[t_j-\delta,t_j+\delta]$, $j=1,\dots,m$, $\phi_t(L)$
	undergoes a \df{standard self-tangency move}. That is, there
	exists a point $q\in\C^n$ and    
	neighborhoods $N\subset N'$ of $q$ with the following
	properties. The intersection
	$N\cap\Pi_\C(\phi_t(L))$ equals $P_1\cup P_2(t)$ which, up to
	biholomorphism looks like
	$P_1=\gamma_1\times P_1'$ and $P_2=\gamma_2(t)\times P_2'$. Here 
	$\gamma_1$ and $\gamma_2(t)$ are subarcs around $0$ of the curves $y_1=0$
	and $x_1^2+(y_1-1\pm t)^2=1$ in the $z_1$-plane, respectively,
	and $P_1'$ and $P_2'$ are real analytic Lagrangian $(n-1)$-disks in
	$\C^{n-1}=\{z_1=0\}$ intersecting transversely at
	$0$. Outside $N'\times\R$ the isotopy is constant. See
	Figure~\ref{1fig:TypeB}. 
        (The full definition of a standard self tangency move appears
	in Section \ref{5Isotopies.section}. For simplicity, one technical condition
	there has been omitted at this point.)
\end{itemize}
Note that 
two Legendrian isotopic admissible Legendrian submanifolds of
dimension $1$ are in general not isotopic through an admissible
Legendrian isotopy. In this case one must allow also a ``triple point
move'' see \cite{Chekanov, ENS}. 


To prove Proposition~\ref{1stabletame} we need to  
check that the differential graded algebra changes only by stable tame isomorphisms 
under Legendrian isotopies of type (A) and (B).
We start with type (A) isotopies.
\begin{lem}\label{1typeA}
	Let $L_t, t\in [0,1]$ be a type (A) isotopy between generic
	admissible Legendrian submanifolds.
	Then the DGA's associated to $L_0$
	and $L_1$ are tame isomorphic.
\end{lem}

To prove this we use a parameterized version of
Proposition~\ref{1manifold}. If $L_t, t\in I=[0,1]$ is a type (A)
isotopy then the double points of $\Pi_\C(L_t)$ trace
out continuous curves. Thus, when we refer to a Reeb chord
$c$ of $L_{t'}$ for some $t'\in[0,1]$ this unambiguously specifies a
Reeb chord for all $L_t.$ For any $t$ we let 
$\M_A^t(a;{\mathbf b})$ denote the moduli space $\M_A(a;{\mathbf b})$ for $L_t$ and define
\begin{equation}
	\M_A^I(a;{\mathbf b})=\{(u,t)|u\in\M_A^t(a;{\mathbf b})\}.
\end{equation}
As above ``generic'' refers to a member of a Baire subset, see Subsection \ref{8oneparamperturb.section}
for a more precise formulation of this term for $1$-parameter families.   
\begin{prop}\label{1param-manifold}
	For a generic type (A) isotopy $L_t, t\in I=[0,1]$ the following holds.  
	If $a,{\mathbf b}, A$ are such that $\mu(A)+|a|-|{\mathbf
	b}|=d\le 1$ then the moduli space $\mathcal{M}_A^I(a;{\mathbf b})$
	is a transversely cut out $d$-manifold.
       	If $X$ is the union of all these transversely cut out
	manifolds which are $0$-dimensional then
	the components of $X$ are of the form
	$\M_{A_j}^{t_j}(a_j,{\mathbf b_j})$, where 
	$\mu(A_j)+|a_j|-|{\mathbf b_j}|=0$, for a finite number of
	distinct instances $t_1,\dots,t_r\in[0,1]$. Furthermore,
	$t_1,\dots,t_r$ are such that $\M_{B}^{t_j}(c;{\mathbf d})$ is
	a transversely cut out $0$-manifold for every $c,{\mathbf
	d},B$ with $\mu(B)+|c|-|{\mathbf d}|=1$.   
\end{prop}

Proposition~\ref{1param-manifold} is proved in
Subsection~\ref{8degentrans.section}. At an instant $t=t_j$ in 
the above proposition we say a {\em handle slide} occurs, and an element
in $\M_{A_j}^{t_j}(a_j,{\mathbf b}_j)$ will be called a \df{handle slide
disk}. (The term 
handle slide comes form the analogous situation in Morse theory.)

The proof of Lemma~\ref{1typeA} depends, just as the proof of Lemma
\ref{1dsquare=0}, on one compactness- and one gluing result; the following. 
\begin{prop}\label{1param-compactness}
	Any sequence $u_\alpha$ in $\M_A^I(a;{\mathbf b})$ has a subsequence that
	converges to a broken holomorphic curve with the same properties as in
	Proposition \ref{1compactness}.
\end{prop}
The proof of this proposition is identical to that of
Proposition~\ref{1compactness}, see Section~\ref{9compactness.section}.

\begin{prop}\label{1param-gluing}
	Let $\delta>0$ and let $L_t, t \in I=[-\delta,\delta]$ be a small neighborhood
	of a handle slide at $t=0$ in a generic type (A) isotopy. Then
	for $\delta$ sufficiently small, $L_{\pm\delta}$ are generic
	admissible and, with $u\in\M_A^0(a;{\mathbf b})$ denoting the handle
	slide disk, the following holds. 
	\begin{enumerate}
	\item 
Assume that $c$ is the $j$-th letter in 
${\mathbf b}.$
Let $\M_B^0(c;{\mathbf d})$ be a moduli space
of rigid holomorphic disks. 
Then there exist $\rho_0>0$ and an embedding
\[
G\colon \M_B^0(c;{\mathbf d})\times[\rho_0,\infty)
\to
\M_{A+B}^I(a;{\mathbf b}_{\{j\}}({\mathbf d})).
\]
Given $v\in\M_B^0(c;{\mathbf d}),$ $G(v,\rho)$ converges to the broken curve
$(v,u)$ as $\rho\to\infty.$ Moreover, any curve in $\M_{A+B}^I(a;{\mathbf b}_{\{j\}}({\mathbf d}))$
with image sufficiently close to the image of $(v,u)$ is in the image of $G.$ 
	\item 
Let $\M_B^0(c;{\mathbf d})$ be a moduli space
of rigid holomorphic disks, where $S=\{s_1,\dots,s_r\}$,  
and ${\mathbf d}$ has
$a$ at every position of an element in $S.$
Then there exist $\rho_0>0$, and an embedding 
\[
G'\colon\M_B^0(c;{\mathbf d})\times[\rho_0,\infty)\to  
\M_{B+r\cdot A}^I(c;{\mathbf d}_{S}({\mathbf b},\dots,{\mathbf b})).
\]
Given $v\in \M_B^0(c;{\mathbf d}),$ $G'(v,\rho)$ converges to the broken curve
$(v,u, \ldots, u).$  Moreover, any curve in $\M_{B+r\cdot A}^I(c;{\mathbf d}_{S}({\mathbf b},\dots,{\mathbf b}))$
with image sufficiently close to the image of $(v,u, \ldots, u)$ is in the image of $G'.$  
	\end{enumerate}
\end{prop}
\begin{proof}
	This proposition follows from Theorem \ref{gluhandslide1} and 
	Theorem~\ref{gluhandslide2}. 
	We show here why the above are the only kind of broken curves to consider gluing.
	If the broken curve lives in the compactification of the
	one-dimensional $\mathcal{M}^I_B(c_0; {\mathbf c})$, then by 
	(\ref{1subadditivity.eqn}) at least one of its pieces must
	have negative formal dimension. Since the handle slide disk   
	$u$ is the only disk with negative formal dimension, all but
	one of the pieces of the broken disk must be $u$. 
	The requirement that our disks have just one positive puncture and Lemma~\ref{1height}
	reduce all possible configurations of the broken curve to the ones considered above.
\end{proof}

We now prove  Lemma~\ref{1typeA} in two steps. First consider type
(A) isotopies without handle slides.

\begin{lem}\label{1noslide}
	Let $L_t, t\in [0,1]$ be a generic type (A) isotopy of Legendrian submanifolds
	for which no handle slides occur. Then the boundary maps
	$\partial_0$ and $\partial_1$ on $\A=\A(L_0)=\A(L_1)$ satisfies
	$\pa_0=\pa_1$. 
\end{lem}
\begin{proof}
	Propositions~\ref{1param-compactness} and \ref{1param-gluing} imply
	that $\mathcal{M}_A^I(a;{\mathbf B})$ is compact when its
	dimension is one. Since if a sequence in this space converged
	to a broken curve $(u^1,\dots,u^N)$ then at least one $u^j$
	would have negative formal dimension. This contradicts
	the assumptions that no handle slide occurs and that the type
	(A) isotopy is generic. Thus the corresponding 0 dimensional moduli spaces 
	$\mathcal{M}^{0}_A$ and $\mathcal{M}^1_A$ used in the definitions of $\partial_{0}$ and
	$\partial_1,$ respectively, form the boundary of a compact 1-manifold. Hence their
	modulo 2 counts are equal.
\end{proof}

We consider what happens around a handle slide instant. 
Let $L_t$, $t\in[-\delta,\delta]$ and $\mathcal{M}_A^0(a;{\mathbf b})$
be as in Lemma \ref{1param-gluing}. Let $\partial_-$  
denote the differential on $\A=\A(L_{-\delta})$, and $\partial_+$ the one on
$\A=\A(L_{\delta})$.  For generators $c$ in $\A$ define
\begin{equation} \notag
\phi_a(c)=\begin{cases} 
                     c &\text{if $c\ne a$},\\
                     a+A{\mathbf b} &\text{if $c=a$}.
                   \end{cases}
\end{equation}
and extend $\phi_a$ to a tame algebra automorphism of $\A$. 
\begin{lem}\label{1slidegeom}
Let $c$ be a generator of $\A$ then
\begin{equation} \notag
\partial_+ c = 
\begin{cases}
\phi_a(\partial_-c) & \text{if $c\ne a$,}\\
\partial_-(\phi_a(c)) & \text{if $c=a$}.
\end{cases}
\end{equation}
\end{lem}
\begin{proof}
Any $\alpha\in\A$ can be expressed in a unique way as a $\Z_2$-linear
combination of elements $C{\mathbf w}$, where $C\in H_1(L)$ and
${\mathbf w}$ is a word in the generators of $\A$, see \eqref{DGAasVSP}. Let
$\la\alpha,C{\mathbf w}\ra$ denote the coefficient ($0$ or $1$) in such an expansion.   
It follows from Proposition \ref{1param-gluing} that for any generator
$c\ne a$
\begin{equation} \notag
	\langle (\partial_+-\partial_-) c, B {\mathbf w}_1{\mathbf
	b}{\mathbf w}_2\rangle=
	\langle \partial_-c, (BA^{-1}){\mathbf w}_1 a {\mathbf w}_2\rangle.
\end{equation}
From this, the formula for $\partial_+c$ follows when $c\not =a.$
The formula when $c=a$ follows similarly.
\end{proof}

\begin{lem}\label{1slidealg}
	The map $\phi_a:\mathcal{A}\to \A$ is a tame isomorphism from $(\A,\partial_-)$
	to $(\A,\partial_+).$
\end{lem}
\begin{proof}
As $\phi_a$ is clearly a tame isomorphism of algebras we only need to
check that it is also a chain map. If $c\ne a$ is a generator then
$\phi_a\partial_- c=\partial_+c=\partial_+\phi_a c.$
It follows from Lemma~\ref{1height} that $\partial_+a$ contains no
terms which contain an $a$ and that the word ${\mathbf b}$ does not
contain the letter $a$. Thus $\pa_+ A{\mathbf b}=\pa_- A{\mathbf b}$
and hence 
$$
\phi_a\pa_-a=\phi_a\pa_-(\phi_a
(a+A{\mathbf b}))=\phi_a(\pa_+ a+\pa_+A{\mathbf b})=
\pa_+(a+A{\mathbf b})=\pa_+\phi_aa.
$$
\end{proof}

\begin{proof}[Proof of Lemma~\ref{1typeA}]
The lemma follows from Lemmas \ref{1noslide}, \ref{1slidegeom}, and \ref{1slidealg}.
\end{proof}

We consider elementary isotopies of type (B). 
Let $L_t, t\in I=[-\delta,\delta]$ be an isotopy of type (B) where
two Reeb chords $\{a, b\}$ are born as $t$ passes through 0. Let $o$
be the degenerate Reeb chord (double point) at $t=0$ and let
$\mathcal{C}'=\{a_1,\ldots, a_l, b_1,\ldots, b_m\}$ be the other Reeb chords.
Again we note that $c_i\in \mathcal{C}'$ unambiguously defines a
Reeb chord for all $L_t$ and $a$ and $b$ 
unambiguously define two Reeb chords for all $L_t$ when $t>0.$ 
It is easy to see that (with the appropriate choice of capping paths) 
the grading on $a$ and $b$ differ by $1$ so let $|a|=j$ and $|b|=j-1.$
Let $(\A_-, \partial_-)$ and 
$(\A_+, \partial_+)$ be the DGA's associated to
$L_{-\delta}$ and $L_{\delta}$, respectively.
\begin{lem}\label{1typeB}
	The stabilized algebra $S_j(\A_-,\partial_-)$ is tame isomorphic to $(\A_+,\partial_+).$
\end{lem}

\begin{proof}[Proof of Proposition \ref{1stabletame} and \ref{1invarianthomology} ]
The first proposition follows from Lemmas \ref{1typeA} and
\ref{1typeB} and implies in its turn the second.
\end{proof}

We prove Lemma \ref{1typeB} in several steps below. Label the Reeb
chords of $L_t$ so that
\begin{equation} \notag
	\action(b_m)\leq \ldots\leq \action(b_1)\leq
	\action(b)<\action(a)\leq \action(a_1)\leq \ldots \leq
	\action(a_l),
\end{equation}
let $\B=\Z_2[H_1(L)]\la b_1,\dots,b_m\ra$ and note that $\B$ is a
subalgebra of both $\A_-$ and $\A_+$.  
Then
\begin{lem}
	For $\delta>0$ small enough
	\begin{equation} \notag
		\partial_+ a = b + v,
	\end{equation}
	where $v\in\B$.
\end{lem}
\begin{proof}
	Let ${\bf 0}\in H_1(L)$ denote the zero element.
	In the model for the type (B) isotopy 
	there is an obvious disk in $\M^t_{\bf 0}(a;b)$ for $t>0$ small
	which is contained in the $z_1$-plane. We argue that this is
	the only point in the moduli space. We restrict attention to the neighborhood
	$N$ of $o^\ast$ that is biholomorphic to the origin in $\C^n$ as in the description
	of a type (B) move. Let $\pi_i:\C^n\to \C$ be the projection onto the $i^{th}$
	coordinate. If $u:D\to \C^n$ is a holomorphic map in
	$\M^t_{\bf 0}(a;b)$ then $\pi_i\circ u$
	will either be constant or not. If $\pi_i\circ u$ is non-constant for $i>1$ then
	the image of $\pi_1\circ u$ intersected with $N$ has boundary
	on two transverse Lagrangian submanifolds. As such it
	will have a certain area $A_i.$ Since
	$\action(a)-\action(b)\to 0$ as $t\to 0+$ we can choose $t$
	small enough so that $\action(a)-\action(b)< A_i$, for all
	$i>1$. Then $\pi_i\circ u$ must
	be a point for all $i>1$ and for
	$i=1$, it can only be the obvious disk. Lemma~\ref{splittran1} shows 
	that $\M^t_{\bf 0}(a;b)$ is transversely cut out and thus contributes to $\partial_+ a.$ 
	If $u\in\M_A^t(a; b)$, where $A\ne {\mathbf 0}$ then
	the image of $u$ must leave $N$. Thus, the above argument
	shows that $\M_A^t(a; b)=\emptyset$ for $t$ small
	enough. Also, for $t>0$ sufficiently small
	$\action(a)-\action(b)<\action(b_m).$ 
	Hence by Lemma~\ref{1height}, $v\in\B$.
\end{proof}

Define the elementary isomorphism $\Phi_0:\A_+\to S_j(\A_-)$
(on generators) by
\begin{equation} \notag
\Phi_0(c) = 
\begin{cases}
e_1^j& \text{if $c= a$,}\\
e_2^j+v& \text{if $c= b$}\\
c& \text{otherwise.}
\end{cases}
\end{equation}
The map $\Phi_0$ fails to be a tame isomorphism since it is not a
chain map. However, we use it as the first step in an inductive
construction of a tame isomorphism $\Phi_l\colon\A_+\to S_j(\A_-)$. To
this end, for $0\le i\le l$, let $\A_i$ be the subalgebra of $\A_+$ generated by
$\{a_1,\ldots, a_i, a, b, b_1,\ldots, b_m\}$ (note that
$\A_l=\A_+$). Then, with $\tau:S_j(\A_-)\to \A_-$ denoting the natural
projection and with $\partial_-^s$ denoting the differential induced on $S_j(\A_-)$,
we have
\begin{lem}\label{1formula}
	\begin{equation}\label{1eqn:exact} 
		\Phi_0\circ \partial_+ w = \partial_-^s\circ\Phi_0 w 
	\end{equation} 
	for $w\in \A_0$ and
	\begin{equation}\label{1eqn:degen}
		\tau\circ\Phi_0\circ \partial_+=\tau\circ\partial_-^s\circ\Phi_0.
	\end{equation}
\end{lem}

Before proving this lemma, we show how to use it in the inductive
construction which completes the proof of Lemma~\ref{1typeB}.  

\begin{proof}[Proof of Lemma~\ref{1typeB}]
The proof is similar to the proof of Lemmas 6.3 and 6.4 in \cite{ENS} ({\em cf} \cite{Chekanov}).
Define the map $H:S_j(\A_-)\to S_j(\A_-)$ on words ${\mathbf w}$ in the
generators by
\begin{equation} \notag
H({\mathbf w}) = 
\begin{cases}
0& \text{if ${\mathbf w} \in \A_-$,}\\
0& \text{if ${\mathbf w}= \alpha e_1^j\beta$ and $\alpha\in\A_-$}\\
\alpha e_1^j\beta & \text{if ${\mathbf w}=\alpha e_2^j\beta$ and $\alpha\in \A_-$,}
\end{cases}
\end{equation}
and extend it linearly.
Assume inductively that we have defined a graded isomorphism
$\Phi_{i-1}:\A_+\to S_j(\A_-)$ so that it is a 
chain map when restricted to $\A_{i-1}$ and so that
$\Phi_{i-1}(a_k)=a_k$, for $k>i-1$. (Note that $\Phi_0$ has these
properties by Lemma \ref{1formula}.) 

Define the elementary isomorphism $g_i:S_j(\A_-)\to S_j(\A_-)$ on generators by
\begin{equation} \notag
g_i(c) = 
\begin{cases}
c& \text{if $c\not= a_i$,}\\
a_i+H\circ\Phi_{i-1}\circ\partial_+(a_i)& \text{if $c= a_i$}
\end{cases}
\end{equation}
and set $\Phi_i=g_i\circ \Phi_{i-1}.$ Then $\Phi_i$ is a graded isomorphism. 
To see that $\Phi_i$ is a chain map when restricted to
$\A_i$ observe the following facts: 
$\tau\circ H=0, \tau\circ g_i=\tau$, and
$\tau\circ\Phi_i=\tau\circ\Phi_0$ for all $i.$ 
Moreover, $\partial_+a_i\in\A_{i-1}$ and $\tau-\id_{S_j(\A_-)} = \partial_-^s\circ H
+H\circ\partial_-^s,$ where in the last equation we think of
$\tau\colon S_j(\A_-)\to S_j(\A_-)$ as $\tau\colon S_j(\A_-)\to\A_-$
composed with the natural inclusion.

Using these facts we compute 
\begin{align}\notag
\pa_-^s g_i(a_i) =& \pa_-^s (a_i) + (\pa_-^sH)\Phi_{i-1}\pa_+(a_i)=
\pa_-^s(a_i)+(H\pa_-^s+\tau+\id)\Phi_{i-1} \pa_+(a_i)\\\notag
=&\pa_-^s(a_i)+\tau\Phi_0\pa_+(a_i)+\Phi_{i-1}\pa_+(a_i)=
\Phi_{i-1}\pa_+(a_i).
\end{align}
Thus $\Phi_i\circ\partial_+(a_i)=\partial_-^s\circ g_i(a_i)=\partial_-^s
\circ\Phi_i(a_i).$ Since $\Phi_i$ and $\Phi_{i-1}$ agree on $\A_{i-1}$
it follows that $\Phi_i$ is a chain map on $\A_i.$  Continuing we
eventually get a tame chain isomorphism $\Phi_l: \A_+\to S_j(\A_-).$
\end{proof}

The proof of Lemma \ref{1formula} depends on the following two propositions.

\begin{prop}\label{1degenglue}
	Let $L_t, t\in I=[-\delta,\delta]$ be a generic Legendrian
	isotopy of type (B) with notation as above 
	(that is, $o$ is the degenerate Reeb chord of $L_0$ and the
	Reeb chords $a$ and $b$ are 
	born as $t$ increases past 0). 
\begin{enumerate}
\item  Let $\M_A^0(o,{\mathbf c})$ be a moduli space
of rigid holomorphic disks.
Then there exist $\rho>0$ and a local homeomorphism
\[
S\colon\M^0_A(o;{\mathbf c})
\times[\rho,\infty)\to
\M^{(0,\delta]}_A(a;{\mathbf c}),
\]
with the following property. If
$u\in\M^0_A(o;{\mathbf c})$ then any disk
in $\M^{(0,\delta]}_A(a;{\mathbf c})$ 
sufficiently close to the image of $u$  is in the image of $S.$
\item 
Let $\M_A^0(c,{\mathbf d})$ be a moduli 
space of rigid holomorphic disks. 
Let $S\subset\{1,\dots,m\}$ be the subset of positions of ${\mathbf d}$
where the Reeb chord $o$ appears (to avoid trivialities, 
assume $S\ne \emptyset$). Then there exists $\rho>0$ 
and a local homeomorphism
\[
S'\colon\M_A^0(c,{\mathbf d})
\times[\rho,\infty)\to
\M^{(0,\delta]}_A(c,{\mathbf d}_S(b)),
\]
with the following property. If $u\in\M_0(c;{\mathbf d})$ then any disk in
$\M^{(0,\delta]}_A(c;{\mathbf d}_S(b))$
sufficiently close to the image of $u$  is in the image of $S'.$
\end{enumerate}
\end{prop}

This is a rephrasing of Theorem~\ref{shoslft} and the following proposition is a restatement of 
Theorem~\ref{gluslft}

\begin{prop} \label{1degenglue2}
	Let $L_t, t\in I=[-\delta,\delta]$ be a generic 
	isotopy of type (B). Let $\M_{A_1}^0(o;{\mathbf c^1})$, \dots
	, $\M_{A_r}^0(o;{\mathbf c^r})$, and $\M^0_{B}(c;{\mathbf d})$ be moduli spaces of
	rigid holomorphic disks. Let $S\subset\{1,\dots,m\}$ be the
	subset of positions in ${\mathbf d}$ where the Reeb chord $o$
	appears and assume that $S$ contains $r$ elements. Then there exists
	 $\rho>0$ and an embedding 
\[
G\colon\M_B^0(c;{\mathbf d})\times
\Pi_{j=1}^r\M_{A_j}^0(o;{\mathbf d^j})
\times[\rho,\infty)\to
\M^{[-\delta, 0)}_{B+\sum A_j}(c;{\mathbf d}_S({\mathbf c^1},\dots,{\mathbf c^r})),
\]
with the following property. If $v\in\M_0(c;{\mathbf d})$ and $u_j\in\M_0(o;{\mathbf
c^j})$, $j=1,\dots,r$ then any disk in   
$\M^{[-\delta,0)}_{B+\sum A_j}(c;{\mathbf d}_S({\mathbf c^1},\dots,{\mathbf c^r}))$
sufficiently close to the image of $(v, u_1, \ldots, u_r)$ is in the image of $G.$ 
\end{prop}

\begin{proof}[Proof of Lemma~\ref{1formula}]
Equation \eqref{1eqn:exact}
follows from arguments similar to those in Lemma~\ref{1typeA}. Specifically, one
can use these arguments to show that $\partial_+ b_i=\partial_- b_i.$
Then since $\partial_+ b_i\in\B$ and since $\Phi_0$ is the identity on $\B$,
\begin{equation} \notag
	\Phi_0\partial_+ b_i=\partial_+ b_i=\partial_-b_i=\partial_-^s\Phi_0b_i.
\end{equation}
We also compute
\begin{equation} \notag
	\Phi_0\partial_+ a=\Phi_0(b+v)=e^j_2+v+v=e^j_2=\pa_-^s\Phi_0a,
\end{equation}
and, since $\pa_+b$ and $\pa_+ v$ both lie in $\B$, 
\begin{equation} \notag
	\Phi_0\partial_+ b=\pa_+ b,\quad \pa_-^s\Phi_0
	b=\pa_-^s(e^j_1+v)=\pa_-v=\pa_+ v.
\end{equation} 
Since $0=\pa_+\pa_+ a=\pa_+ b+\pa_+v$, we conclude that \eqref{1eqn:exact} holds.

To check (\ref{1eqn:degen}), we write $\pa_+ a_i= W_1+W_2+W_3$, where
$W_1$ lies in the subalgebra generated by
$\{a_1,\dots,a_l,b_1,\dots,b_m\}$, where $W_2$ lies in the ideal
generated by $a$ and where $W_3$ lies in the ideal generated by $b$ in 
the subalgebra generated by $\{a_1,\dots,a_l,b,b_1,\dots,b_m\}$. 

Let $u_{t}$,  be a family of holomorphic disks with
boundary on $L_{t}$. As $t\to 0$, $u_{t}$ converges to a broken disk
$(u^1,\dots,u^N)$ with boundary on $L_0$. This together with the
genericity of the type (B) isotopy implies that for $t\ne 0$ small enough
there are no disks of negative formal dimension with boundary on
$L_t$ since a broken curve 
which is a limit of a sequence of such disks would have at least one
component $u^j$ with negative formal dimension. 

Let $u_{s}\colon D\to\C^n,$ $s\ne 0$ be rigid disks with boundary on $L_s$. If,
the image $u_{-t}(\pa D)$ stays a positive distance away from
$o^\ast$ as $t\to 0+$ then the argument above implies that
$u_{-t}$ converges to a non-broken curve. Hence $\pa_- a_i = W_1 + W_4$ where for each rigid disk
$u_{-t}\colon D\to\C^n$ contributing to a word in $W_4$ there exists points
$q_{-t}\in\pa D$ such that $u_{-t}(q_{-t})\to o^\ast$ as $t\to 0+$. The genericity
assumption on the type (B) isotopy implies that no rigid disk with
boundary on $L_0$ maps any boundary point to $o^\ast$, see
Corollary~\ref{lmagood}. Hence $u_{-t}$ 
must converge to a broken curve $(u^1,\dots,u^N)$ which brakes at
$o^\ast$. Moreover, by genericity and \eqref{1subadditivity.eqn}, every component
$u^j$ of the broken curve must be a rigid disk with boundary on $L_0$.     
Proposition \ref{1degenglue2} shows that any such broken curve may
be glued and Proposition \ref{1degenglue} determines the pieces which
we may glue. It follows that $W_4=\hat W_2$ where $\hat W_2$ is
obtained from $W_2$ by replacing each occurrence of $b$ with
$v$. Therefore, 
\begin{equation}\notag
\tau\Phi_0\pa_+ (a_i)=\tau\Phi_0 (W_1+W_2+W_3) =W_1+\hat W_2=\pa_-
(a_i)=\tau\pa_-^s\Phi_0 (a_i).  
\end{equation}
\end{proof}

\subsection{Relations with the relative contact homology of \cite{EGH}} \label{1SFT.section}

Our description of contact homology is a direct generalization of Chekanov's ideas in \cite{Chekanov}. We
now show how the above theory fits into the more general, though still developing, relative contact
homology of \cite{EGH}.

We start with
a Legendrian submanifold $L$ in a contact manifold $(M,\xi)$ and try to build an invariant
for $L.$ To this end, let $\alpha$ be a contact
form for $\xi$ and $X_\alpha$ its Reeb vector field.
We let $\mathcal{C}$ be the set of all Reeb chords, which under certain non-degeneracy
assumptions is discrete. Let $\mathcal{A}$ be the free associative 
non-commutative unital algebra over $\Z_2[H_1(L)]$  generated by $\mathcal{C}.$ The algebra 
$\mathcal{A}$ can be given a grading using the Conley-Zehnder index (see \cite{EGH}). To do this we must
choose capping paths $\gamma$ in $L$ for each $c\in\mathcal{C}$ which
connects its end points. Note that $c \in \mathcal{C}$, being a piece of a flow
line of a vector field, comes equipped with a 
parameterization $c:[0,T] \rar M$. For later convenience, we
reparameterize $c$ by precomposing it with $\times T\colon[0,1]\to[0,T]$. 

We next wish to define a differential on $\mathcal{A}.$
This is done by counting holomorphic curves in the symplectization of $(M,\xi).$ Recall the 
\df{symplectization} of $(M,\xi)$ is the manifold $W=M\times \R$ with the symplectic form
$\omega=d(e^w\alpha)$ where $w$ is the coordinate in $\R.$ Now choose an almost complex structure 
$J$ on $W$ that is compatible with $\omega$ ($\omega(v,Jv)>0$ if $v\not=0$), leaves $\xi$ invariant and
exchanges $X_\alpha$ and $\frac{\partial}{\partial w}.$ Note that $\overline{L}=L\times \R$ is
a Lagrangian (and hence totally real) submanifold of $(W,\omega).$ Thus we may study holomorphic
curves in $(W,\omega, J)$ with boundary on $\overline{L}.$ Such curves must have punctures. When
the Reeb field has no periodic orbits (as in our case) there can be no internal punctures, so 
all the punctures occur on the boundary. To describe the behavior near the punctures let 
$u:(D_m,\partial D_m)\to (W,\overline{L})$ be a holomorphic curve where $D_m$ is as before.
Each boundary puncture has
a neighborhood that is conformal to a strip $(0,\infty)\times[0,1]$ with coordinates $(s,t)$ such that
approaching $\infty$ in the strip is the same as approaching $p_i$ in the disk. If we write $u$ using
these conformal strip coordinates near $p_i$ then we say $u$ \df{tends asymptotically to a 
Reeb chord} $c(t)$ at $\pm \infty$ if the component of $u(s,t)$ lying
in $M$ limits to $c(t)$ as $s\to \infty$ and the component of $u(s,t)$ lying in $\R$ limits to
$\pm\infty$ as $s\to \infty.$ The map $u$ must tend asymptotically to a Reeb chord at each boundary 
puncture. 
Some cases of this asymptotic analysis were done in \cite{Abbas}.
For $\{a,b_1,\ldots, b_m\}\subset \mathcal{C}$  we consider the 
moduli spaces $\M_A^s(a;b_1,\ldots b_m)$ of holomorphic maps $u$
as above such that: (1) at $p_0,$ $u$ tends asymptotically to $a$ at $+\infty;$ 
(2) at $p_i,$ $u$ tends asymptotically to $b_i$ at $-\infty;$ 
(3) and $\Pi_M(u(\partial D_*))\cup_i \gamma_i$ represents the homology class
$A.$ 
Here the map $\Pi_M:W\to M$ is projection onto the $M$ factor of $W.$  
We may now define a boundary map $\partial$ on the generators $c_i$ of $\mathcal{A}$ 
(and hence on all of $\A$) by 
\begin{equation} \notag
\partial c_i=\sum \#(\M_A^s(c_i; b_1,\ldots, b_m)) A b_1\dots b_m,
\end{equation}
where
the sum is taken over all one dimensional moduli spaces and $\#$ means the modulo two count of the points
in $\M_A^s/\R.$
Here the $\R$-action is induced by a translation in the $w$-direction.

Though this picture of contact homology has been known for some time now, the analysis needed to 
rigorously define it has yet to appear. Moreover, there have been no attempts to make computations
in dimensions above three. Above, by specializing to a nice -- though still rich -- situation, we
gave a rigorous definition of contact homology for Legendrian submanifolds in $\R^{2n+1}.$

Recall that in our setting $(M, \alpha) =(\R^{2n+1}, dz - \sum_{j=1}^n y_j dx_j),$ 
the set of Reeb chords is naturally bijective with the double points of
$\Pi_\C(L).$
Thus, clearly the algebra of this subsection is identical to the one described in Section \ref{1First}.

We now compare the differentials.
We pick the complex structure on the symplectization of $\R^{2n+1}$ as follows.
The projection $\Pi_{\C}: \R^{2n+1}\to \C^n$ gives an isomorphism
$d\Pi_{\C}$ from $\xi_x\subset T_x\R^{2n+1}$ to $T_{\Pi_{\C}(x)}\C^n$ and
thus, via $\Pi_\C$, the standard complex structure on $\C^n$ induces a
complex structure $E:\xi\to \xi$ on $\xi$. Define the complex structure
$J$ on the symplectization $\R^{2n+1}\times\R$ by $J(v)=E(v)$ if $v\in
\xi$ and $J(\frac{\partial}{\partial w})=X.$ Then $J$ is compatible with $\omega=d(e^w\alpha).$ 
Our moduli spaces and the ones used in the standard
definition of contact homology are related as follows. If $u$ in
$\mathcal{M}_A^s(a; b_1,\ldots,b_m)$ then define $p(u)$ to be the map
in $\mathcal{M}_A(a;b_1,\ldots,b_m)$ as $p(u)=\Pi_\C\circ\Pi_M\circ u$ and
$\widetilde{p(u)}=\Pi_M\circ u|\pa \hat D $, where 
$\Pi_M:\R^{2n+1}\times\R\to \R^{2n+1},$ the projection from the symplectization back to the
original contact manifold. 
\begin{lem}
	The map $p:\mathcal{M}^s(a; b_1,\ldots,b_m)/\R\to\mathcal{M}(a;b_1,\ldots,b_m)$ is
	a homeomorphism.
\end{lem} 
\begin{proof}
This was proven in the three dimensional case in \cite{ENS} and the proof here is similar. 
(For details we refer the reader to that paper but we outline the main steps.)
It is clear from the definitions that $p$ is a map between the appropriate spaces (we mod out by
the $\R$ in $\mathcal{M}^s$ since the complex structure on the symplectization is $\R$-invariant and
any two curves that differ by translation in $\R$ will clearly project to the same curve in $\C^n$).
The only non-trivial part of this lemma is that $p$ is invertible. To see this let $u\in \M^s$ be
written $u=(u',z, \tau):D_m\to \C^n\times\R\times \R.$ The fact that $u$ is holomorphic for our
chosen complex structure implies that $z$ is harmonic and hence determined by its boundary 
data. Moreover, the holomorphicity of $u$ also implies that
$\tau$ is determined, up to translation in $w$-direction, by 
$u'$ and $z.$ Thus if we are given a map $u'\in \M$ then we can construct
a $z$ and $\tau$ for which $u=(u',z,\tau)$ will be a holomorphic map $u:D_m\to \R^{2n+1}\times \R.$
If it has the appropriate behavior near the punctures then $u\in\M^s.$
The asymptotic behavior near punctures was studied in \cite{Robbin-Salamon??}.
\end{proof}

\section{Legendrian submanifolds}

In this section, we review the Lagrangian projection and introduce the front projection,
both of which are useful for the calculations of Section \ref{1examples}.
In Sections \ref{1rotclass} and \ref{1tb.section}, we discuss the 
Thurston-Bennequin invariant and the rotation class, which were the only invariants
before contact homology  which could distinguish Legendrian 
isotopy classes.
Finally, we construct in Section \ref{1frontcomp} a useful technique for calculating the
Conley-Zehnder index of Reeb chords.

\subsection{The Lagrangian projection} \label{1lagrangian.section}
Recall that for a Legendrian submanifold $L\subset\C^n\times\R$,
$\Pi_{\C}\colon L\to\C^n$ is a Lagrangian immersion.   
Note that $L\subset\C^n\times\R$ can be recovered, up to rigid translation
in the $z$ direction, from $\Pi_{\C}(L)$: pick a point $p\in \Pi_\C(L)$ and choose any $z$ coordinate
for $p$; the $z$ coordinate of any other point $p'\in L$ is then determined by 
\begin{equation}\label{1recoverz}
	\sum_{j=1}^n \int_{\gamma} y_j dx_j, 
\end{equation}
where $\gamma=\Pi_\C\circ\Gamma$ and $\Gamma$ is any path in $L$ from $p$ to $p'.$ 
Furthermore, given any Lagrangian immersion $f$ into $\C^n$ 
with isolated double points, if the integral in (\ref{1recoverz}) is independent
of the path $\gamma=f\circ\Gamma$ then we obtain a Legendrian
immersion $\tilde f$ into $\R^{2n+1}$ which is an embedding provided
the integral is not zero for paths connecting double points.

A Lagrangian immersion $f: L\to\C^n$ is {\em exact} if
$f^*(\sum_{j=1}^n y_j dx_j)$ is exact and, in this case, \eqref{1recoverz}
is independent of $\gamma$. In particular, if $H^1(L)=0$  
then all Lagrangian immersions of $L$ are exact. Also note that any
Lagrangian regular homotopy $f_t\colon L\to\C^{n}$ of exact Lagrangian
immersions will lift to a Legendrian regular homotopy $\tilde f_t\colon L\to\C^n\times\R$. 

\begin{ex}\label{1basicexample}
Consider $S^n=\{(x,y)\in\R^n\times\R: |x|^2+y^2=1\}$ and define $f\colon S^n\to\C^n$ as 
\begin{equation} \notag
	f(x,y):S^n\to \C^n:(x,y)\mapsto ((1+iy)x).
\end{equation}
Then $f$ is an exact Lagrangian immersion,
with one transverse double point, which
lifts to a Legendrian embedding into $\R^{2n+1}.$ (When $n=1$ the
image of $f$ is a figure eight in the plane with a double point at the origin.)
\end{ex}

\subsection{The front projection} \label{1front.section}
\begin{figure}[ht]
	{\epsfxsize=4in\centerline{\epsfbox{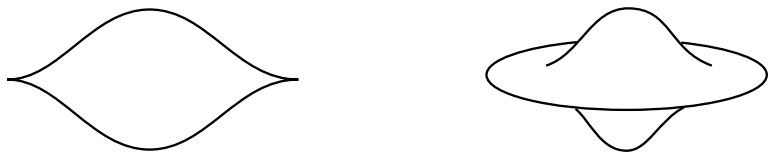}}}
	\caption{Front projection of Example~\ref{1basicexample} 
		in dimension 3, on the left, and 5, on the right.}
	\label{1fig:front}
\end{figure}

The \df{front projection} projects out the $y_j$'s:
\begin{equation} \notag
	\Pi_F:\R^{2n+1}\to \R^{n+1}: (x_1,y_1,\ldots,x_n,y_n,z)\mapsto (x_1,\ldots,x_n, z).
\end{equation}
If $L\subset\R^{2n+1}$ is a Legendrian submanifold then
$\Pi_F(L)\subset \R^{n+1}$ is its \df{front} which is a codimension one subvariety of $\R^{n+1}.$ 
The front has certain singularities. More precisely, for generic $L$, the set of singular
points of $\Pi_F$ is a hypersurface $\Sigma\subset L$ which is smooth outside a
set of codimension $3$ in $L$, and which contains a subset
$\Sigma'\subset\Sigma$ of codimension $2$ in $L$ with the following property. 
If $p$ is a smooth point in
$\Sigma\setminus\Sigma'$ then there are local coordinates $(x_1,\dots,x_n)$ around $p$
in $L$, $(\xi_1,\dots,\xi_n,z)$ around $\Pi_F(p)$ in $\R^{n+1}$, and
constants $\delta= \pm 1$, $\beta,\alpha_2,\dots,\alpha_n$ 
such that
\begin{equation}
\Pi_F(x_1,\dots,x_n)=(x_1^2,x_2,\dots,x_n,\delta x_1^3+\beta x_1^2+\alpha_2x_2+\dots+\alpha_nx_n).
\end{equation}    
For a reference, see \cite{AG} Lemma on page 115. The image under
$\Pi_F$ of the set of smooth points in $\Sigma\setminus\Sigma'$ will
be called the {\em cusp edge} of the front $\Pi_F(L)$. See Figure~\ref{1fig:front}. 

Any map $L\to \R^{n+1}$ with singularities of a generic
front can be lifted (in a unique way) to a Legendrian immersion. 
(The singularities of such a map allow us to solve for the
$y_i$-coordinates from the equation $dz=\sum_{i=1}^n y_idx_i$ and the
solutions give an immersion.) In particular, at a smooth point of the front  
the $y_i$-coordinate equals the slope of the tangent
plane to the front in the $x_iz$-plane. 

The Legendrian immersion of a generic front is an embedding and a
double point of a Legendrian immersion correspond to a double point of
the front with parallel tangent planes. Also note that that 
$\Pi_F(L)$ cannot have tangent planes containing the $z$-direction. 
For a more thorough discussion of singularities occurring in front projections see \cite{AG}.

\subsection{The rotation class}\label{1rotclass}
Let $(M,\xi)$ be a contact $(2n+1)$-manifold with a contact form
$\alpha$. That is, $\alpha$ is a $1$-form on $M$ with
$\xi=\krn(\alpha)$. The complete non-integrability condition on $\xi$
implies $\alpha\wedge (d\alpha)^n\ne 0$ which in turn implies that for
any $p\in M$, $d\alpha_p|\xi_p$ is a symplectic form on 
$\xi_p\subset T_pM$.      

Let $f:L\to (M,\xi),$ be a Legendrian immersion. Then the
image of $df_x:T_xL\to T_{f(x)}M$ is a Lagrangian plane in $\xi_{f(x)}.$ Pick any
complex structure $J$ on $\xi$ which is compatible with its symplectic structure.  
Then the complexification of $df$, $df_{\C}:TL\otimes\C\to \xi$ is a
fiberwise bundle isomorphism. The homotopy
class of $(f,df_{\C})$ in the space of complex fiberwise 
isomorphisms $TL\otimes\C\to\xi$ is called 
the \df{rotation class} of $f$ and is denoted $r(f)$ (or $r(L)$ if
$L\subset M$ is a Legendrian submanifold embedded into $M$ by the inclusion). 
The h-principle for Legendrian immersions \cite{GromovPDR} implies that $r(f)$ is a complete invariant for
$f$ up to regular homotopy through Legendrian immersions.

When the contact manifold under consideration is $\R^{2n+1}$ with the standard contact structure
we may further illuminate the definition of $r(f).$ Let
$(x,y,z)\in\R^n\times\R^n\times\R$ be coordinates on $\R^{2n+1}$ as in
Subsection \ref{1First}. If $J\colon\xi_{(x,y,z)}\to\xi_{(x,y,z)}$ is the complex structure defined by
$J(\pa_{x_j}+y_j\pa_z)=\pa_{y_j}$,
$J(\pa_{y_j})=-(\pa_{x_j}+y_j\pa_{z})$, for $j=1,\dots,n$ then 
the Lagrangian projection 
$\Pi_{\C}:\R^{2n+1} \to \C^n$ gives a complex isomorphism from $(\xi,J)$ to
the trivial bundle with fiber $\C^n$. Thus we may think of $df_{\C}$ as a trivialization
$TL\otimes\C\to \C^n.$ Moreover, we can choose Hermitian metrics on
$TL\otimes \C$ and on $\C^n$ so that $df_{\C}$ is a unitary map. Then
$f$ gives rise to an element in $U(TL\otimes\C, \C^n).$ 
One may check that the group of continuous maps $C(L, U(n))$ acts freely and transitively on
$U(TL\otimes\C, \C^n)$ and thus $\pi_0(U(TL\otimes\C, \C^n))$ is in one to one correspondence with
$[L, U(n)].$ Thus we may think of $r(f)$ as an element in $[L, U(n)].$

We note that when $L=S^n$ then 
\begin{equation} \notag
	r(f)\in \pi_n(U(n))\approx\begin{cases}
      \Z, & n \text{ odd,} \\
      0, & n \text{ even.}
    \end{cases}
\end{equation}
Thus for spheres we will refer to $r(f)$ as the \df{rotation number}.

\subsection{The Thurston-Bennequin invariant} \label{1tb.section}
Given an orientable connected Legendrian submanifold $L$ in an oriented contact
$(2n+1)$-manifold $(M,\xi)$ we define an invariant, 
called the Thurston--Bennequin invariant of $L,$
describing how the contact structure ``twists about $L.$'' The invariant was originally
conceived by Thurston and Bennequin \cite{Bennequin83} when $n=1$ and generalized to higher dimensions by
Tabachnikov \cite{Tabachnikov88}. Here we only define the Thurston--Bennequin invariant when
$L$ is homologically trivial in $M$ (which for $M=\R^{2n+1}$ poses no additional constraints). 

Pick an orientation on $L$. Let $X$ be a Reeb 
vector field for $\xi$  and push $L$ slightly off of itself along $X$
to get another oriented submanifold $L'$ disjoint from $L.$ 
The \df{Thurston--Bennequin invariant} of $L$ is
the linking number
\begin{equation}
\tb(L)=\lk(L,L').
\end{equation}
Note that $\tb(L)$ is independent of the choice of orientation on $L$
since changing it changes also the orientation of $L'$. 
The linking number is computed as follows. Pick any $(n+1)$-chain $C$ in $M$ such that
$\partial C=L$ then $\lk(L,L')$ equals the algebraic intersection
number of $C$ with $L'.$ 

For a chord generic Legendrian submanifold $L\subset \R^{2n+1}$, $\tb(L)$ can be
computed as follows. Let $c$ be a Reeb chord of $L$ with end points $a$
and $b$, $z(a)>z(b)$. Let $V_a=d\Pi_\C(T_aL)$ and $V_b=d\Pi_\C(T_bL).$ Given an
orientation on $L$ these are oriented $n$-dimensional transverse subspaces in $\C^n.$ If the orientation
of $V_a\oplus V_b$ agrees with that of $\C^n$ then we say the sign, $\sgn(c),$ of $c$ is $+1$ otherwise
we say it is $-1.$ Then
\begin{equation}\label{1tbform}
	\tb(L)=\sum_c \sgn(c),
\end{equation}
where the sum is taken over all Reeb chords $c$ of $L.$ 
To verify this formula, use the Reeb-vector field $\pa_z$ to shift $L$
off itself and pick the cycle $C$ as the cone over $L$ through some
point with a very large negative $z$-coordinate. 


Note that the parity of the number of double points of any generic
immersion of an $n$-manifold into $\C^n$ depends only on its regular
homotopy class \cite{Whitney}. Thus the parity of $\tb(L)$ is
determined by the rotation class $r(L).$
Some interesting facts \cite{Eliashberg90} concerning the Thurston-Bennequin invariant are summarized in
\begin{prop}
	Let $L$ be a Legendrian submanifold in standard contact $(2n+1)$-space.
	\begin{enumerate}
	\item If $n>1$ is odd, then for any $k\in\Z$ we can find,
	$C^0$ close to $L$, a Legendrian submanifold $L'$ 
	smoothly isotopic and Legendrian regularly homotopic
	to $L$ with $\tb(L')=2k$. 
	\item If $n$ is even, then $\tb(L)=(-1)^{\frac{n}{2}+1}\frac{1}{2}\chi(L).$ 
	\end{enumerate}
\end{prop}
The ideas associated with (1) are discussed below in Proposition~\ref{1stabcomp}. 
For (2), note that if $n=2k$ is even then the sign of a double point $c$
is independent of the ordering of the subspaces $V_a$ and $V_b$ and in
this case $\tb(L)$ equals Whitney's invariant \cite{Whitney}
for immersions of orientable $2k$-manifolds into oriented $\R^{4k}$ which in
turn equals $-\frac{1}{2}\chi(\nu)$, where $\nu$ is the oriented normal bundle
of the immersion \cite{LashofSmale}. Since the immersion is
Lagrangian into $\C^n$ its normal bundle is isomorphic to the tangent
bundle $TL$ of $L$ (via multiplication with $i$) and as an oriented bundle it is
isomorphic to $TL$ with orientation multiplied by $(-1)^{\frac{n}{2}}$. 

If $n=1$ the situation is much more interesting. In this case there are two types of contact
structures: tight and overtwisted. If the contact structure is overtwisted then the above 
proposition is still true, but if the contact structure is tight (as
is standard contact $3$-space)
then $\tb(L)\leq \chi(L)-| r(L)|.$ There are other interesting bounds on $\tb(L)$ in a tight
contact structure, see \cite{Fuchs-Tabachnikov97, Rudolph95}.

The Thurston-Bennequin invariant of a chord generic Legendrian submanifold can also
be calculated in terms of Conley-Zehnder indices of Reeb chords.
Recall $\mathcal{C}$ is the set of Reeb chords of $L$.
\begin{prop}\label{1proptb}
If $L\subset\R^{2n+1}$ is an orientable chord generic Legendrian submanifold 
then 
\begin{equation} \notag
\tb(L)=(-1)^{\frac{(n-2)(n-1)}{2}}
\sum_{c\in \mathcal{C}}(-1)^{|c|}.
\end{equation}
\end{prop}

\begin{proof}
Recall from (\ref{1tbform}) 
that $\tb(L)$ can be computed by summing $\sgn(c)$ over all Reeb
chords $c$, where $\sgn(c)$ is the oriented intersection between the
upper and lower sheets of $\Pi_\C(L)$ at $c^\ast.$ So to prove the proposition we only need to check
that $\sgn(c)=(-1)^{\frac12(n^2+n+2)}(-1)^{|c|}.$ This will be done after the 
proof of Lemma~\ref{1lmafront}.
\end{proof}

\subsection{Index computations in the front projection}\label{1frontcomp}
Though it was easier to define contact homology using the complex projection 
of a Legendrian submanifold it is frequently easier to construct
Legendrian submanifolds using the front projection.
In preparation for our examples below we discuss Reeb chords and their 
Conley-Zehnder indices in the front projection.

If $L\subset\R^{2n+1}$ is a Legendrian submanifold then the Reeb
chords of $L$ appears in the front projection as vertical line segment
(i.e. a line segment in the $z$-direction) 
connecting two points of $\Pi_F(L)$ with parallel tangent planes.
(See Subsection \ref{1front.section} and note that 
$L$ may be perturbed so that the Reeb chords as seen in $\Pi_F(L)$ do not have
end points lying on singularities of $\Pi_F(L).$)

A generic arc $\gamma$ in $\Pi_F(L)$ connecting two such points $a, b$ 
intersects the cusp edges of $\Pi_F(L)$ transversely and meets no other
singularities of $\Pi_F(L)$ (it might also meet double points of the front projection but
``singularities'' refers to non-immersion parts of $\Pi_F(L)$). Let $p$
be a point on a cusp edge where $\gamma$ intersects it. Note that $\Pi_F(L)$ has
a well defined tangent space at $p$. Choose a line $l$ orthogonal to
this tangent space. Then, since the tangent space does not contain the
vertical direction, orthogonal projection 
to the vertical direction at $p$ gives a linear isomorphism from
$l$ to the $z$-axis through $p$. 
Thus the $z$-axis induces an orientation on $l$. Let $\gamma_p$ be a small part of
$\gamma$ around $p$ and let $h_p\colon\gamma_p\to l$ be orthogonal
projection. The orientation of $\gamma$ induces one on $\gamma_p$ and
we say that the intersection point is an {\em up- (down-) cusp} if $h_p$ is
increasing (decreasing) around $p$.  

Let $c$ be a Reeb chord of $L$ with end points $a$ and $b$, $z(a)>z(b)$. Let $q$ be the
intersection point of the vertical line containing $c\subset\R^{n+1}$ and
$\{z=0\}\subset\R^{n+1}$. Small parts of $\Pi_F(L)$ 
around $a$ and $b$, respectively, can be viewed as the graphs of
functions $h_a$ and $h_b$ from a neighborhood of $q$ in $\R^n$ to $\R$
(the $z$-axis). Let $h_{ab}=h_a-h_b$. Since the tangent planes of $\Pi_F(L)$ at $a$ and
$b$ are parallel, the differential of $h_{ab}$ vanishes at
$q$. If the double point $c^\ast$ of $\Pi_\C(L)$ corresponding to $c$ is
transverse then the Hessian $d^2 h_{ab}$ is a non-degenerate 
quadratic form (see the proof below). Let $\ind(d^2 h_{ab})$ denote
its number of negative eigenvalues. 

\begin{lem}\label{1lmafront}
If $\gamma$ is a generic path in $\Pi_F(L)$ connecting $a$ to $b$ then
\begin{equation} \notag
\nu_\gamma(c)= D(\gamma) - U(\gamma) + \ind(d^2h_{ab}),
\end{equation}
where $D(\gamma)$ and $U(\gamma)$ is the number of down- and up-cusps
of $\gamma$, respectively.
\end{lem}
\begin{proof}
To compute the Maslov index as described in Section~\ref{1CZ} we use
the Lagrangian reference space $x=0$ in
$\R^{2n}$ (that is, the subspace $\spa(\pa_{y_1},\dots,\pa_{y_n})$) with
Lagrangian complement $y=0$ ($\spa(\pa_{x_1},\dots,\pa_{x_n})$).

We must compute the Maslov index of the loop
$\Gamma\ast\lambda(V_b,V_a)$ where $V_b$ and $V_a$ are the Lagrangian
subspaces $d\Pi_\C(T_b L)$ and $d\Pi_\C(T_aL)$ and $\Gamma(t)$ is the
path of Lagrangian subspaces induced from $\gamma.$  
We first note that
$\Gamma(t)$ intersects our reference space transversely (in $0$) if
$\gamma(t)$ is a smooth point of $\Pi_F(L)$, since near such points $\Pi_F(L)$ can be thought
of as a graph of a function over some open set in
$x$-space (i.e. $\{z=0\}\subset\R^{n+1}$). Thus, for generic $\gamma$,
the only contributions  to the Maslov index come from cusp-edge
intersections and the path $\lambda(V_b,V_a)$. 

We first consider the contribution from $\lambda(V_b,V_a)$.
There exists orthonormal coordinates $u=(u_1,\dots,u_n)$ in $x$-space so
that in these coordinates
\begin{equation} \notag
d^2h_{ab}=\diag(\lambda_1,\dots,\lambda_n).
\end{equation}

We use coordinates $(u,v)$ on
$\R^{2n}=\C^n$ where $u$ is as
above, $\pa_j=\pa_{u_j}$, and $\pa_{v_j}=i\pa_j$ ($i=\sqrt{-1}$). In these coordinates our
symplectic form is simply $\omega=\sum_{j=1}^n du_j\wedge dv_j$, and
our two Lagrangian spaces are given by 
$V_a=\spa_{j=1}^n(\pa_j+id^2h_a\pa_j),$ 
$V_b=\spa_{j=1}^n(\pa_j+id^2h_b\pa_j).$ One easily computes  
\begin{equation} \notag
\omega(\pa_j+id^2h_b\pa_j,\pa_j+id^2h_a\pa_j)=
\omega(\pa_j,id^2h_{ab}\pa_j)=\lambda_j.
\end{equation}
Moreover, let
\begin{equation} \notag
W_j=\spa(\pa_j+id^2h_a\pa_j,\pa_j+id^2h_b\pa_j)=
    \spa(\pa_j+id^2h_{a}\pa_j,i\pa_j)=\spa(\pa_j+id^2h_b\pa_j,i\pa_j),
\end{equation}
then $W_j$ and $W_k$ are symplectically orthogonal, $\omega(W_j,W_k)=0,$ for $j\ne k$. 

Let $v_a(j)$ be a unit vector in direction $\pa_j+id^2h_a\pa_j$ and
similarly for $v_b(j)$. Define
the almost complex structure $I$ 
as follows
\begin{equation} \notag
I(v_b(j))=\sign(\lambda_j)v_a(j), 
\end{equation}
and note that it is compatible with $\omega$. 
Then $e^{sI}v_b(j)$, $0\le s\le\frac{\pi}{2}$, intersects the line in 
direction $i\pa_j$ if and only if $\lambda_j<0$ and does so in the
positive direction.

It follows that the contribution of  $e^{sI}V_b$, $0\le s\le\frac{\pi}{2}$
(i.e. $\lambda(V_b,V_a)$) to the Conley-Zehnder index is $\ind(d^2 h_{ab})$.

Second we consider cusp-edge intersections:
at a cusp-edge intersection $p$ (which we take to be the origin) there
are coordinates $u=(u_1,\dots,u_n)$ such that the front locally around 
$p=0$ is given by $u\mapsto(x(u),z(u))$, where
\begin{equation} \notag
x(u)=(u_1^2,u_2,\dots,u_n),\quad z(u)=\delta u_1^3+\beta
u_1^2+\alpha_2u_2+\dots+\alpha_n u_n,
\end{equation}
where $\delta$ is $\pm 1$, and  $\beta$ and $\alpha_j$ 
are real constants. 
We can assume the oriented curve $\gamma$ is given by 
$u(t)=(\epsilon t,0,\dots,0)$, where $\epsilon=\pm 1$. If we
take the coorienting line $l$ to be in the direction of the vector
\begin{equation} \notag
v(p)=(-\beta,-\alpha_2,\dots,-\alpha_n,1).
\end{equation}
then the function
$h_p$ is   
\begin{equation} \notag
h_p(t)=\delta\epsilon^3t^3,
\end{equation}
and we have an up-cusp if $\delta\epsilon>0$ and a down cusp if
$\delta\epsilon<0$. 

The curve $\Gamma(t)$ of Lagrangian tangent planes of $\Pi_\C(L^n)$ along 
$\gamma$ is given by
\begin{equation} \notag
\Gamma(t)=\spa(2\epsilon
t\pa_1+i\frac{3\delta}{2}\pa_1,\pa_2,\dots,\pa_n). 
\end{equation}
The plane $\Gamma(0)$ intersects our reference plane at $t=0$ along
the line in direction $i\pa_1$. As described in Section~\ref{1Maslov} the sign of the intersection is given
by the sign of 
\begin{equation} \notag
\frac{d}{dt}\omega(i\tfrac{3\delta}{2\sigma}
\pa_{1},2\sigma\epsilon\pa_{1})=-3\delta\epsilon
\end{equation}
Thus, we get negative signs at up-cusps and positive at down-cusps.
The lemma follows.
\end{proof}

\begin{proof}[Completion of Proof of Proposition~\ref{1proptb}]
We say the orientations on two hyperplanes transverse the the $z$-axis in $\R^{n+1}$ {\em agree} if
their projection to $\{z=0\}\subset\R^{n+1}$ induce the same
orientation on this $n$-dimensional subspace. Let $c$ be a Reeb chord
of $L$ and let $a$ and $b$ denote its end points on $\Pi_F(L)$. 
If the orientations on $T_a\Pi_F(L)$ and $T_b\Pi_F(L)$ agree then the above proof shows
that the bases 
\begin{align}\notag
&\Bigl(\pa_1+id^2h_a\pa_1,\dots,\pa_n+id^2h_a\pa_n,\pa_1+id^2h_b\pa_1,\dots,\pa_n+id^2h_b\pa_n\Bigr)\simeq\\
&\Bigl(id^2h_{ab}\pa_1,\dots,id^2h_{ab}\pa_n,\pa_1,\dots,\pa_n\Bigr),
\end{align}
provide oriented bases for $d\Pi_\C(T_a L)\oplus d\Pi_\C(T_bL).$
Note that the standard orientation of $\C^n$ is given by the positive
basis $(\pa_1,i\pa_1,\dots,\pa_n,i\pa_n)$ which after multiplication
with $(-1)^{\frac{n(n+1)}{2}}$ agrees with the orientation given by
the basis $(i\pa_1,\dots,i\pa_n,\pa_1,\dots,\pa_n)$.
Thus
$$\sgn{c}= (-1)^{\frac{n(n+1)}{2}}(-1)^{\ind(d^2(h_{ab}))}.$$
However, the orientations of $T_a\Pi_F(L)$ and $T_b\Pi_F(L)$ do not always agree.
Let $\gamma$ be the path in $L$ connecting $a$ to $b.$ The orientations on $T_{\gamma(t)}(\Pi_F(L))$
do not change as long as $\gamma$ does not pass a cusp edge. 
It follows from the local model for a cusp edge that each time $\gamma$ transversely
crosses a cusp edge the orientation on $T_{\gamma(t)}(\Pi_F(L))$ changes. Thus 
$\sgn{c}= (-1)^{\frac{n(n+1)}{2}}(-1)^{D(\gamma)+U(\gamma)}(-1)^{\ind(d^2h_{ab})} 
= (-1)^{\frac{1}{2}(n^2+n+2)}(-1)^{|c|}$ as we needed to show.
\end{proof}

\section{Examples and constructions}\label{1examples}
Before describing our examples, we discuss the linearized contact homology in Section
\ref{1linhom}. 
This is an invariant of Legendrian submanifolds derived from the DGA. Its main advantage over contact
homology it that it is easier to compute. In Section \ref{1examples.section}, we
do several simple computations of contact homology. In Sections \ref{1stabilization.section}
and \ref{1frontspinning.section}, we describe two 
constructions: {\em stabilization} and {\em front spinning}.
In these subsections, we construct infinite families of pairwise
non-isotopic Legendrian $n$-spheres, $n$-tori 
and surfaces which are indistinguishable by the classical invariants.  
In Section \ref{1remarksontheexamples.section}, we
summarize the information gleaned from our examples and discuss what this says about
Legendrian submanifolds in general.

\subsection{Linearized homology}\label{1linhom}
To distinguish Legendrian submanifolds using contact homology one must
find computable invariants of stable tame isomorphism classes of DGA's.  
We use an idea of Chekanov \cite{Chekanov} to ``linearize'' the
homology of such algebras. To keep the discussion simple we will only consider algebras 
generated over $\Z_2$ and not $\Z_2[H_1(L)].$


Let $\A$ be an algebra generated by $\{c_1,\dots,c_m\}$. For
$j=0,1,2,\dots$ let $\A_j$ denote the ideal of $\A$ generated by all
words ${\mathbf c}$ in the generators with $l({\mathbf c})\ge j$. 
A differential $\partial\colon\A\to\A$ is called \df{augmented} if 
$\pa(A_1)\subset A_1$ (in other words if $\pa c_j$ does not contain 1
for any $j$). If $(\A,\pa)$ is augmented then $\pa(A_j)\subset
A_j$ for all $j$. A DGA $(\A,\partial)$ is called \df{good} if its
differential is augmented. 

Let $(\A,\pa)$ be a DGA with generators $\{c_1,\dots,c_m\}$ and consider
the vector space $\mathcal{V}=\A_1/\A_2$ over $\Z_2$. 
If $(\A,\pa)$ is good then $\pa\colon\A\to\A$ induces a differential
$\pa_1\colon{\mathcal V}\to{\mathcal V}$. Note that
$\{c_1,\dots,c_m\}$ gives a basis in ${\mathcal V}$ and that in this basis $\pa_1 c_j$ 
equals the part of $\pa c_j$ which is linear in the generators. We
define the {\em linearized homology} of a $(\A,\pa)$ as
$$
\krn(\pa_1)/\img(\pa_1),
$$
which is a graded vector space over $\Z_2$.


We want to apply this construction to DGA's associated to Legendrian
isotopy classes. Let $L\subset \R^{2n+1}$ be an
admissible Legendrian submanifold with algebra $(\A(L),\pa)$ generated
by $\{c_1,\dots,c_m\}$. Let $G$ be the set of tame isomorphisms of
$\A(L)$ and for $g\in G$ let $\pa^g\colon\A(L)\to\A(L)$ be $\pa^g=g\pa
g^{-1}$. We define the {\em linearized contact homology of $L$},
$HLC_\ast(\R^{2n+1},L)$ to be the set of isomorphism classes of
linearized homologies of $(\A,\pa^g)$, where $g\in G$ is such that
$(\A,\pa^g)$ is good. (Note that this set may be empty.) Define
$G_0\subset G$ to be the subgroup of tame 
isomorphisms $g_0$ such that $g_0(c_j)=c_j+a_j$ for all $j$, where
$a_j=0$ or $a_j=1$. Note that $a_j=0$ if $|c_j|\ne 0$ since $g_0$ is
graded and that $G_0\approx \Z_2^k$, where $k$ is the number of
generators of $\A$ of degree $0$.   

\begin{lma}
If $L_t\subset\R^{2n+1}$ is a Legendrian isotopy between admissible
Legendrian submanifolds then $HLC_\ast(\R^{2n+1},L_0)$ is isomorphic
to $HLC_\ast(\R^{2n+1},L_1)$. Moreover, if $L\subset\R^{2n+1}$ is an
admissible Legendrian submanifold then $HLC_\ast(\R^{2n+1},L)$ is
equal to the set of isomorphism classes of linearized homologies of
$(\A,d^{g_0})$, where $g_0\in G_0$ is such that $(\A,\pa^{g_0})$ is
good.  
\end{lma}

\begin{proof}
The first statement follows from the observation that the
stabilization $(S_j(\A),\pa)$ of a good DGA $(\A,\pa)$ is good and
that the linearized homologies of $(S_j(\A),\pa)$ and $(\A,\pa)$ are
isomorphic. The second statement is proved in \cite{Chekanov}. 
\end{proof}

Let $L\subset\R^{2n+1}$ be an admissible Legendrian submanifold. Note
that if $\A(L)$ has no generator of degree $1$ then $(\A,\pa)$ is
automatically good and if $\A$ has no generator of degree $0$ then
$G_0$ contains only the identity element. If the set
$HLC_\ast(\R^{2n+1},L)$ contains only one element we will sometimes
below identify this set with its only element.

\subsection{Examples} \label{1examples.section}
In this subsection we describe several relatively simple examples in which the contact homology is
easy to compute 
and defer more complicated computations to the following subsections.

\begin{ex}\label{1basicex}{ 
The simplest example in all dimensions is $L_0$ described in
Example~\ref{1basicexample}, with a  single Reeb chord $c.$ Using Lemma~\ref{1lmafront} and
the fact that the difference of the $z$-coordinates at the end points
of $c$ is a local maximum we find $|c|=n.$ 
So  $\A(L_0)= \langle c \rangle$ and the differential is $\partial c=0,$ showing (if $n>1$) the 
contact homology is 
\begin{equation} \notag
	HC_k(\R^{2n+1},L_0)=
	\begin{cases}
      	0, & k \not\equiv 0 \mod n, \mbox{ or } k<0, \\
      	\Z_2, & \hbox{otherwise.}
    \end{cases}
\end{equation}
If $n=1$ this is still true but $\M(c;\emptyset)$ is not empty 
(it contains two elements \cite{Chekanov}).}
\end{ex}

\begin{ex}{
Generalizing Example~\ref{1basicex} 
above we can consider the Legendrian sphere $L'$ in $\R^{2n+1}$ with
$3$ cusp edges in its front projection. See Figure~\ref{1fig:3cusp}.
\begin{figure}[ht]
	{\epsfxsize=1.5in\centerline{\epsfbox{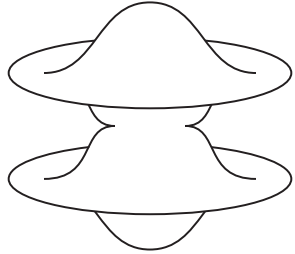}}}
	\caption{The sphere $L'$ with 3 cusps.}
	\label{1fig:3cusp}
\end{figure} 
If one draws the pictures with an $SO(n)$ symmetry about the
$z$-axis then there will be one Reeb chord running from the top of the sphere to the bottom, call it
$c$ and a $(n-1)$-spheres worth of Reeb chords. Perturbing the symmetric picture
slightly yields two Reeb chords $a, b$ in place of the spheres worth
in the symmetric picture. The gradings are 
\[
|c|= n+2,\quad
|a|=1,\quad
|b|=n.
\]
It is clear that whatever the contact homology of $L'$ is, it is different from that in the example above.
When $n>2$ (respectively $n=2$) there are two (respectively three) possibilities for the boundary map. 
Thus for $n$ even
we have examples of non-isotopic Legendrian spheres with the same classical invariants. 
}\end{ex}

\begin{figure}[ht]
	{\epsfxsize=4.5in\centerline{\epsfbox{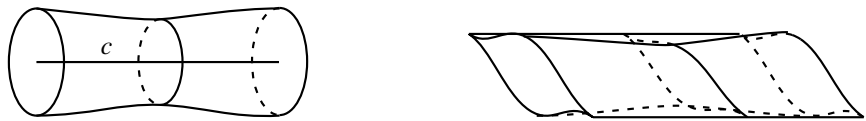}}}
	\caption{The neighborhood of $c$ and its deformation into the front projection of a Legendrian
		tube.}
	\label{1fig:tube}
\end{figure}

Given two Legendrian submanifolds $K$ and $K'$ we can form their ``(cusp) connected sum'' as follows:
isotope $K$ and $K'$ so that their fronts are separated by a
hyperplane in $\R^{n+1}$ containing the $z$-direction and let $c$ be an 
arc beginning at a cusp edge of $K$ and ending at a cusp edge of $K'$ 
and parameterized by
$s\in[-1,1].$ Take a neighborhood $N$ of $c$ whose vertical cross sections
consist of round balls  whose radii vary with $s$ and have exactly one minimum at $s=0$ and no
other critical points. Introducing cusps along $N$ as indicated in Figure~\ref{1fig:tube}.
Define the ``connected sum'' $K\# K'$ to be the Legendrian
submanifold obtained from the joining of  $K\setminus(K\cap N), K'\setminus(K'\cap N)$ 
and  $\partial N.$ Note this operation might 
depend on the cusp edges one chooses on $K$ and $K'$ but we will make this choice explicit
in our examples. In dimension 3 it can be shown that the connected sum of two knots is
well defined \cite{Chekanov, Etnyre-Honda2}. It would be interesting to understand this
operation better in higher dimensions. See Remark~\ref{1connectsum.remark} below.
\begin{lem}\label{1sumcomp}
Let ${\mathcal C}$ and ${\mathcal C}'$ be the sets of Reeb chords of
$K$ and $K'$, respectively, and let $|\cdot|_K$, $|\cdot|_{K'}$, and
$|\cdot|$ denote grading in $\A(K)$, $\A(K')$, and $\A(K\# K')$,
respectively. It is possible to perform the connected 
sum so that the set of Reeb chords of $K\# K'$ is ${\mathcal
C}\cup{\mathcal C}'\cup\{h\}$ and so that the following holds.
\begin{enumerate}
\item If $c\in {\mathcal C}$ then $|c|_K=|c|$, if
$c'\in{\mathcal C}'$ then $|c|_{K'}=|c|$, and $|h|=n-1.$
\item $\partial h=0.$
\item If $\A_K$ and $\A_{K'}$ denote the subalgebras of $\A(K\# K')$
generated by ${\mathcal C}\cup\{h\}$ and ${\mathcal C}'\cup\{h\}$,
respectively, then $\pa(\A_K)\subset\A_K$ and $\pa(\A_{K'})\subset\A_{K'}$.  
\item If $c\in{\mathcal C}$ then $\pa c\in \A(K\# K')_1$ if and only
if $\pa_K c\in\A(K)_1$ and similarly for $c'\in{\mathcal C}'$. (In
other words, the constant part of $\pa c$ ($\pa c'$) does not change after the
connected summation.) 
\end{enumerate}
\end{lem}
\begin{proof}
We may assume that $K$ and $K'$ are on opposite sides of the hyperplane $\{x_1=0\}$ and there is
a unique point $p,$ respectively $p',$ on a cusp edge of $K,$ respectively $K',$ that is closest to
$K',$ respectively $K.$ We may further assume that all the coordinates
but the $x_1$ coordinate of $p$ and  
$p'$ agree. Define $K\# K'$ using  $c,$ the obvious horizontal arc
connecting $p$ and $p.$  It is now clear  
that all the Reeb chords in $K$ and $K'$ are in $K\# K'$ and there is
exactly one extra chord $h,$ coming 
from the minimum in the neighborhood $N$ of $c.$  It is also clear that the gradings of the
inherited chords are unchanged and that $|h|=n-1.$ 

Denote $z_j=x_j+iy_j.$ The image of a holomorphic disk $u:D\to \C^n$
with positive puncture at $h^\ast$ must lie in the complex hyperplane $\{z_1=0\}.$ 
To see this notice that the projection of $K\# K'$ onto the
$z_1$-plane is as shown in Figure~\ref{1fig:sum}.
\begin{figure}[ht]
	{\epsfxsize=3.5in\centerline{\epsfbox{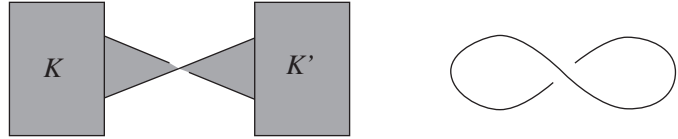}}}
	\caption{The projection onto the $z_1$-plane (left). The intersection with the $z_2$-plane (right).} 
	\label{1fig:sum}
\end{figure} 
Let $u_1$ be the composition of $u$ with this projection. If
$u_1$ is not constant then $u(\partial D)$ 
must lie in the shaded region in Figure~\ref{1fig:sum}. 
Thus the corner
at $h^\ast$ (note $h^\ast$ projects to $0$ in this figure) must be a
negative puncture. Since any holomorphic disk with positive puncture
at $h^\ast$ must lie entirely in the 
hyperplane $\{z_1=0\}$ it cannot have any negative punctures. Thus
$\partial h$ has only a constant part. 
For $n>2$ this implies $\partial h=0$ immediately. If $n=2$ then $\partial h=0$ since 
in this case,  there are exactly two holomorphic disks in the
$z_2$-plane, see Figure~\ref{1fig:sum}, and  
Lemma~\ref{splittran2} implies both these disks contribute to the
boundary map.

To see (3) consider the projection of $K\# K'$ onto the $z_1$-plane, see Figure~\ref{1fig:sum}.
If a holomorphic disk $D$ intersected the projection of $K$ and $K'$ then it would
intersect the $y_1$-axis in a closed interval, with non-trivial interior,  containing the origin.
This contradicts the maximum principle since the intersection of the boundary of $D$
with the $y_1$-axis can contain only the origin.

For the last statement, consider Reeb chords in $K$, those in $K'$
can be handled in exactly the same way. We use Corollary~\ref{lmagood} which implies that we may
choose the points $p$ so that no rigid holomorphic disk $u\colon
D\to\C^n$ with boundary on $K$ maps any point in $\pa D$ to $p$.
Since the
space of rigid disks is a compact $0$ manifold there are only finitely
many rigid disks, $u_1,\dots,u_r$, say. Since each $u_k$ is continuous on the
boundary $\pa D$ we find that $u_1(\pa D)\cup\dots\cup u_r(\pa D)$
stays a positive distance $d$ away from $p$. Consider the ball
$B(p,\frac12 d)$ and use a tube attached inside $B(p,\frac14 d)$ for
the connected sum. If $c\in{\mathcal C}$ and $v$ is a rigid
holomorphic disk with
boundary on $K\# K'$, with positive puncture at $c$, no other
punctures, and such that the image $v(\pa D)$ is disjoint from 
$\pa B(p,\frac12 d)$ then $v$ is also a disk with boundary on $K$ and
hence $v=u_j$ for some $j$. 

Since no holomorphic disk with boundary on $K\#K'$ which touches a
point in $K$ can pass the hyperplane
$\{x_1=0\}\subset\C^n$ it will also represent a disk on the connected
sum $K\# L_0$, where $L_0$ is a small standard sphere. Pick a generic
Legendrian isotopy $K_t$, $0\le t\le 1$ of $K\#L_0$ to $K$ which is
supported in 
$(K\cap B(p,\frac14 d))\# L_0$. Then either there exists $t<1$
such that all rigid disks $v$ on $K_t$ for $t>0$ satisfies 
$v(\pa  D)\cap(B(p,\frac 12 d)\setminus B(p,\frac14 d))=\emptyset$ or there exists a sequence
of rigid disks $v_j$ with boundary on $L_{t_j}$, $t_j\to 1$ as
$j\to\infty$ such that $v_j(\pa D)\cap(B(p,\frac12 d)\setminus
B(p,\frac14 d))\ne
\emptyset$. In the first case the lemma follows from the
observation above. 
We show the second case cannot appear: by Gromov compactness, the
sequence $v_j$ has a subsequence which converges to a broken
disk $(v^1,\dots,v^N)$ with boundary on $K$. Since $K$ is generic
there are no disks with negative formal dimension and all
components of $(v^1,\dots,v^N)$ must be rigid. But since $v_j$ is rigid the broken
disk must in fact be unbroken by \eqref{1subadditivity.eqn}. Thus we find a rigid
disk $v^1$ with boundary on $K$ such that 
$v^1(\pa D)\cap B(p,d)\ne\emptyset$ contradicting our choice of
$p$. The lemma follows. 
\end{proof} 


\subsection{Stabilization and the proof of Theorem \ref{1intromain}} \label{1stabilization.section}

In this subsection we describe a general construction that can be applied to Legendrian submanifolds called 
\df{stabilization}. Using the stabilization technique, we prove Theorem \ref{1intromain}.
We begin with a model situation.

In $\R^{n+1}$ consider two unit balls $F$ and $E$ in the hyperplanes 
$\{z=0\}$ and $\{z=1\},$ respectively. Let $M$ be a $k$-manifold
embedded in $F.$ Let $N$ be a regular $\epsilon$-neighborhood of $M$
in $F$ for some positive $\epsilon\ll 1.$  
Deform $F$ to $F'$ by pushing $M$ up to $z=\epsilon$ and deform 
$N$ so that the $z$-coordinate of 
$p\in N$ is $\epsilon-\text{dist}(p, M).$ Note that are many Reeb chords, one for each point in $M$ and
$F\setminus N.$ To deform this into a generic picture choose a Morse function $f:M\to [0,1]$ and
$g: \overline{(F\setminus N)}\to [0,1]$ such that $g^{-1}(1)=\partial F$ and 
$g^{-1}(0)=\partial N.$ (It is important to notice that we may, if we wish, modify the boundary conditions
on $g|_{\partial F}$ depending on our circumstances.) 
Take a positive $\delta\ll\epsilon$ and further deform $F'$ by adding $\delta f(p)$
to the $z$ coordinate of points in $M$ and subtracting $\delta g(p)$ from the $z$ coordinate of
points in $F\setminus N.$ The result is a generic pair of Lagrangian disks $F'$ and $E$ 
with one Reeb chord for each critical point of $f$ and $g.$ Define $F''$ as we defined $F'$ but begin by
dragging $M$ up to $z=1+\epsilon$ (instead of $z=\epsilon$ as we did for $F'$).
 
Now if $\Pi_F(L)$ is the front projection of a Legendrian submanifold $L$ and there are two horizontal 
disks in $\Pi_F(L)$ we can identify them with $F$ and $E$ above. (Note we can always assume there
are horizontal disks by either looking near a cusp and flattening out a region,
or letting $F$ and $E$ be the disks obtained by flattening out the regions around the top and
bottom of a Reeb chord.) 
Legendrian isotope $L$ so that $F$ becomes $F'.$ 
Replacing $F'$ in $\Pi_F(L)$ by $F''$ 
will result in the front of a 
Legendrian submanifolds $L'$ which is called the \df{stabilization of $L$ along $M$}. 

\begin{figure}[ht]
	{\epsfxsize=4.5in\centerline{\epsfbox{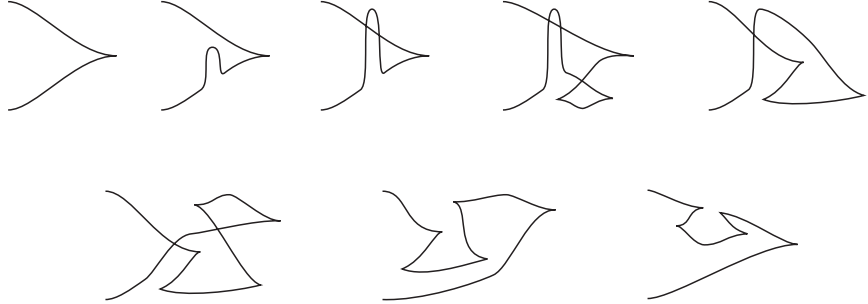}}}
	\caption{Our stabilization in dimension 3 is equivalent to two normal 3 dimensional
		stabilizations \cite{Etnyre-Honda1}.} 
	\label{1doubstab}
\end{figure}

\begin{prop}\label{1stabcomp}
If $L'$ is the stabilization of $L$ with notation as above, then 
\begin{enumerate}
\item The rotation class of $L'$ is the same as that of $L.$
\item The invariant $\tb$ is given by
	\[\tb(L')=
	\begin{cases}
      	\tb(L), &  n \mbox{  even}, \\
      	\tb(L)+ (-1)^{(D-U)}2\chi(M), & n \hbox{ odd,}
    \end{cases}\]
	where $D, U,$ is the number of down, up, cusps along a generic path from $E$ to $F$ in
	$\Pi_F(L).$ 
\item The Reeb chords of $L$ and $L'$ are naturally identified. The grading of any chord not associated
	with $M, F$ and $E$ is the same for both $L$ and $L'.$ Let $c$ be a chord 
	associated to $M, F$ and $E$ and let $|c|_L$ be its grading in $L$ and $|c|_{L'}$ its
	grading in $L'.$ Then 
	\[|c|_{L'}=n-2-|c|_L\]
\end{enumerate}
\end{prop}
This theorem may seem a little strange if one is used to Legendrian knots in $\R^3.$ In particular, it
is well known that in 3 dimensions there are two different stabilizations and both change the 
rotation number. 
What is called a stabilization in dimension 3 is really a ``half stabilization,''
as defined here. 
(Recall such a ``half stabilization'' corresponds to adding zig-zags to the front projection and
looks like a ``Reidemeister Type 1'' move in the Lagrangian projection \cite{Etnyre-Honda1}.) 
In particular if one does the above described stabilization near a cusp in dimension
3 it will be equivalent to doing both types of half stabilizations.  See Figure~\ref{1doubstab}. 
\begin{rmk}
The
stabilization procedure will typically produce non-topologically isotopic knots when done in dimension
3.
\end{rmk}

\begin{proof}
Recall the rotation class is defined as the Legendrian regular
homotopy class.  Now (1) is easy to see since
the straight line homotopy from $\Pi_F(L)$ to $\Pi_F(L')$ will give a regular Legendrian homotopy between 
$L$ and $L'.$ Statement (2) follows from (3). As for (3) let
$c$ be a chord corresponding to a critical point of the Morse function $f$ then Lemma~\ref{1lmafront}
implies that $|c|_L= (k-\text{Morse Index}_c(f))+D-U-1$ and $|c|_{L'}=\text{Morse Index}_c(f)+(n-k)-D+U-1$ where $k$ is the dimension of $M.$
\end{proof}

We now consider some examples to see what effect stabilization has on contact homology.

\begin{ex}
Let $L$ be a Legendrian submanifold in $\R^{2n+1}$ and $p$ a point on
a cusp of $\Pi_F(L).$ Consider a small 
ball $B$ around $p$ in $\R^{n+1}.$ We can isotope the front projection so as to create two
new Reeb chords $c_1$ and $c_2$ in $B$, see Figure~\ref{1doubstab}, such that $|c_1|=0$ and $|c_2|=1$. 
Let $F'$ be the front obtained by pushing
the lower end point of $c_1$ past the 
upper sheet of $\Pi_F(L)$ in $B$ and let $L'$ be the corresponding Legendrian submanifold.
\begin{prop}\label{1stabvanish} 
	The contact homology of $L'$ is 
	\[HC_k(\R^{2n+1}, L')=0.\]
\end{prop}
\begin{proof}
We can assume that $p$ is at the origin in $\R^{n+1}.$ 
For any $\epsilon$ define $B_\epsilon$ to be the product of the ball of radius $\epsilon$ about
$p$ in the $x_1z$-plane times $[-\epsilon, \epsilon]^{n-1}$ (in $x_2\ldots x_n$-space). We may now 
assume that $\Pi_F(L)\cap B_\epsilon$ is the cusp shown in Figure~\ref{1doubstab} times 
$[-\epsilon, \epsilon]^{n-1}$ and that the stabilization is done in $B_{\frac{\epsilon}{2}}.$
A monotonicity argument as in the proof of Lemma~\ref{9hbar.lma} shows that  
any disk with a positive puncture at $c_2$ (or $c_1$) and leaving
$B_\epsilon$ has area bounded below. However, the action of $c_2$ can
be made arbitrarily small. Therefore any disk with a positive corner
at $c_2$ must stay in the ball $B_\epsilon.$  The projection of
$\Pi_{\C}(L')\cap B_\epsilon$ to a $z_j$-plane $j\ne 1$     
is shown on the left hand side of Figure~\ref{1stabb}. (The reason for
the appearance of this picture is that we can choose the front so that
$\frac{\pa z}{\pa x_j}\cdot x_j\le 0$, for all $j>1$.)
\begin{figure}[ht]
	{\epsfxsize=4.5in\centerline{\epsfbox{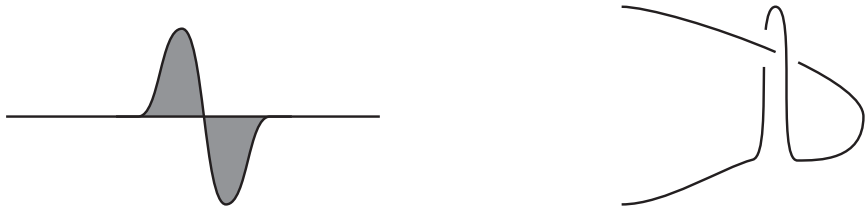}}}
	\caption{$\Pi_{\C}(L')$ projected onto a $z_j$-line, $j\ne 1$
	(left) and intersected with the $z_1$-plane  
	(right).} 
	\label{1stabb} 
\end{figure} 
The boundary of a projection of a holomorphic curve must 
lie in the shaded region of the figure; moreover, the corner at $c_2$ of such a disk is negative.
Thus any holomorphic curve with positive puncture at $c_2$ must lie
entirely in the $z_1$-plane. The right hand side of 
Figure~\ref{1stabb} shows $\Pi_{\C}(L')\cap B_\epsilon\cap\{z_1-\text{plane}\}.$ We see there
one disk which, by Lemma~\ref{splittran2} contributes to the boundary
of $c_2.$ Thus $\partial c_2=1$, and one may easily check this 
implies $HC_k(\R^{2n+1}, L')=0.$
\end{proof}
\end{ex}
This last example is not particularly surprising given the analogous theorem, long
known in dimension $3$ \cite{Chekanov}, that stabilizations (or actually ``half stabilizations'' even) kill
the contact homology. With this in mind the following examples might be a little surprising. It 
shows that in higher dimensions stabilization does not always kill the
contact homology. The main
difference with dimension $3$ is the stabilizations we do below would, in dimension $3$, change the
knot type.

\begin{ex} \label{1BS}
When $n=2$ we define the
sphere $L_1$ via its front projection, which is described in Figure~\ref{1fig:l1}. 
For $n>2$ there is an analogous front projection: take two copies, $L_0, L'_0,$ of 
the Legendrian sphere $L_0$ from Example~\ref{1basicex} 
and arrange them as shown in the figure. Deform $L_0$ as shown in Figure~\ref{1fig:l1}. 
Take a curve $c,$ parameterized by $s\in[-1,1],$  from
the cusp edge on $L_0$ to the cusp edge on $L'_0.$ By taking this curve to be 
very large we can 
assume the rate of change in its $z$-coordinate is very small. Moreover we will assume that
by the time it passes under $L_0$ its $z$-coordinate is less than the $z$-coordinates of $L_0'$ and
thus has to ``slope up'' to connect with $L_0'.$ (These choices will minimize the number of Reeb chords.)  
Take a neighborhood $N$ of $c$ whose vertical cross sections
consist of round balls whose radii vary with $s$ and have exactly one minimum at $s=0$ and no
other critical points. Introducing cusps along $N$ as indicated in Figure~\ref{1fig:tube}
\begin{figure}[ht]
	{\epsfxsize=4.5in\centerline{\epsfbox{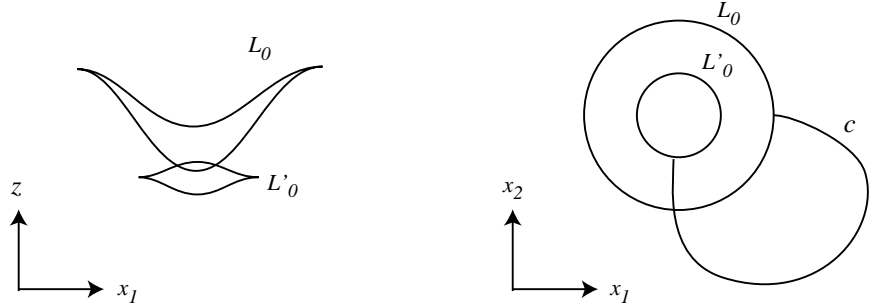}}}
	\caption{On the left hand side the $x_1z$-slice of part of $L_1$ is show. To see this portion
		in $\R^3$ rotate the figure about its center axis. On the right hand side we indicate
		the arc $c$ connecting the two copies of $L_0.$}
	\label{1fig:l1}
\end{figure} 
we can join $L_0,L'_0$ and $\partial N$ together to form a front projection for a Legendrian sphere
in $\R^{2n+1}.$
There are exactly six Reeb chords involving only $L_0$ and $L'_0$ which we label $a_1,\ldots, a_6.$ 
There is also a Reeb chord $b$ that occurs in $N$ where the radii of the
cross sectional balls have a minimum.  
Using Lemma~\ref{1lmafront} we compute:
	\begin{eqnarray*}
	 	|a_1|&=&|a_2|=|a_5|= n,\\
		|b|&=&n-1,\\
		|a_4|&=&0,\\
		|a_3|&=&|a_6|=-1.\\
	\end{eqnarray*}

\begin{prop}\label{1propl1}
The following are true
\begin{enumerate}
\item $L_1$ is a stabilization of $L_0.$
\item For all $n$ the rotation classes of $L_1$ and $L_0$ agree.
\item When $n$ is even $\tb(L_1)=\tb(L_0)$ and when $n$ is odd $\tb(L_1)=\tb(L_0)-2$.
\item The linearized contact homology of $L_1$ in homology grading $-1$ is
	\[HLC_{-1}(\R^{2n+1}, L_1)=\Z_2.\]
\item $L_0$ and $L_1$ are not Legendrian isotopic.
\end{enumerate}	 
\end{prop}
\begin{proof}
Let $L_1'$ be the Legendrian sphere whose front is the same as the
front of $L_1$ except that $L_0'$ has 
been moved down so as to make $L_0$ and $L_0'$ disjoint. Then 
$L_1'$ is clearly Legendrian isotopic  
to $L_0$ and stabilizing $L_1'$ (using $M$ a point) results in $L_1.$
Thus Statement~(1) holds.   
Statements~(2) and (3) follow from (1) and
Proposition~\ref{1stabcomp}. Statement~(5) follows from~(4).  

The Reeb chords for $L_1'$ and $L_1$ are easily identified and their gradings are the same except
for $|a_5|_{L_1'}=-2.$ 
At this point it is clear that $HLC_{-1}(\R^{2n+1}, L_1)=\Z_2$ or $\Z_2\oplus \Z_2.$ (This
is good enough to distinguish $L_0$ and $L_1$.)
Since $L_1'$ and $L_0$ are Legendrian isotopic their linearized
contact homologies must agree. Furthermore, the linearized contact homology
of $L_0$ is a one element set, 
$$
HLC_n=HLC_n(\R^{2n+1},L_0)=\Z_2,\quad 
HLC_j=HLC_j(\R^{2n+1},L_0)=0,\,\, j\ne n. 
$$
Thus, if $\pa_1'$ denotes the (linearized) differential on $\A(L_1')_1/\A(L_1')_2$ we conclude
the following.
\begin{itemize}
\item[(a)]
$\pa_1' a_5=0$ since $a_5$ is the generator of lowest grading.
\item[(b)]
$\img(\pa_1'|\spa(a_3,a_6))=\spa(a_5)$ since
$HCL_{-2}=0$ and thus $\krn(\pa_1'|\spa(a_3,a_6))$ is
$1$-dimensional.
\item[(c)]
If $n>2$ then $\pa_1' a_4$ spans $\krn(\pa_1'|\spa(a_3,a_4))$ since
$HLC_{-1}=0$, if $n=2$ then
$\img(\pa_1'|\spa(a_4,b))=\krn(\pa_1'|\spa(a_3,a_4))$ and therefore
$\krn(\pa_1'|\spa(a_4,b))$ is $1$-dimensional.
\item[(d)]
If $n>2$ then $\pa_1' b=0$. Also, $\img(\pa_1'|\spa(a_1,a_2))=\krn(\pa_1'|\spa(a_4,b))$.
\end{itemize}

Let $\hat L_1$ be the Legendrian immersion ``between'' $L_1'$ and
$L_1$ with one double point which arises as the length of the Reeb
chord $a_5$ shrinks to $0$. Take $\hat L_1$ to be generic admissible. 
Moreover, by Corollary~\ref{lmagood} we may assume that no rigid holomorphic
disk with boundary on $\hat L$ and without puncture at $a_5^\ast$ maps
any boundary point to 
$a_5^\ast$. As in the proof of Lemma \ref{1sumcomp}, we find a ball
$B(a_5^\ast,d)$ such that no rigid disk without puncture
at $a_5^\ast$ maps a boundary point into $B(a_5^\ast, d)$. 

Let $K_t, t\in [-\delta,\delta]$ be a small Legendrian regular homotopy
such that $K_0=\hat L_1$, $K_\delta$ is Legendrian isotopic to $L_1'$
and $K_{-\delta}$ is Legendrian isotopic to $L_1$. Moreover, we take
$K_t$ supported inside a small neighborhood of $a_5$ which map into
$B(a_5^\ast,\frac14 d)$ by $\Pi_\C$. Now if $u\colon D\to K_0$ is a
disk on $K_0$ which maps no boundary 
point into $B(a_5^\ast,\frac12 d)$ then $u$ can be viewed as a disk with boundary on $K_t$
and vice versa.  

We show that there exists $\epsilon>0$ such that for $|t|<\epsilon$
there exist no rigid disk with boundary on $K_t$ and without puncture at
$a_5$ which map a boundary point to $B(a_5^\ast,\frac12 d)$. If this is not the case we
extract a subsequence $v_j$ of such maps which, by Gromov compactness,
converges to a broken disk $(v^1,\dots,v^N)$ with boundary on
$K_0$. If $N>1$ then, by \eqref{1subadditivity.eqn} at least one of the disks $v^j$
must have negative formal dimension but since $K_0$ is generic
admissible no such disks exists and the limiting disk $v^1$ is
unbroken. Now $v^1$ is a rigid disk with boundary on $K_0$ and without
puncture at $a_5^\ast$ which maps boundary points to
$B(a_5^\ast,d)$. This contradicts the choice of $K_0$ and hence proves
the existence of $\epsilon>0$ with properties as claimed. Thus, (c)
above implies, with $\pa_1$ the differential on 
$\A(L_1)_1/\A(L_1)_2$, $\pa_1(\spa(a_4))$ ($\pa_1(\spa(a_4,b))$ if $n=2$) is
$1$-dimensional and hence (4) holds.
\end{proof}

Let $L_2$ be the Legendrian sphere obtained by connect summing two copies of $L_1.$ Note $L_1$
only has one cusp edge so there is no ambiguity in the construction; thus, we choose any arc
which is disjoint from the fronts of the two spheres that are being connect summed. We similarly
define $L_k$ by connect summing $L_{k-1}$ with $L_1.$
\begin{thm}
	The Legendrian spheres $L_k$ are all non Legendrian isotopic
	and, for $n$ even, have the same classical invariants. 
\end{thm}
\begin{proof}
This follows since $HLC_{-1}(\R^{2n+1}, L_k)$ equals $\Z_2^k.$
\end{proof}
\end{ex} 

In order to construct examples in dimensions $2n+1$ where $n$ is odd we consider a variant of 
this example.
\begin{ex}
Let $L_1'$ be constructed as $L_1$ is in Example~\ref{1BS} except start with $L_0$ and $L_0'$ 
as shown in Figure~\ref{1fig:l1prime}.
\begin{figure}[ht]
	{\epsfxsize=1.5in\centerline{\epsfbox{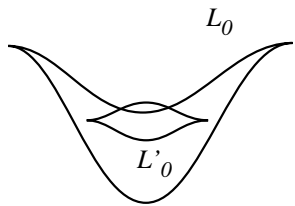}}}
	\caption{The position of $L_0$ and $L_0'$ to construct $L_1'.$}
	\label{1fig:l1prime}
\end{figure} 
Like $L_1$, $L_1'$ will have seven Reeb chords which we label in a similar manner. 
Here the grading on the chords is
\begin{eqnarray*}
	 	|a_1|&=&|a_2|=|a_5|= n,\\
		|a_3|&=&|a_6|=|b|=n-1,\\
		|a_4|&=&0.\\
\end{eqnarray*}
\begin{prop}
\label{4L1prime.prop}
The following are true
\begin{enumerate}
\item $L_1'$ is a stabilization of $L_0.$
\item For all $n$ the rotation class of $L_1'$ and $L_0$ agree.
\item When $n$ is even $\tb(L_1')=\tb(L_0)$ and when $n$ is odd $\tb(L_1')=\tb(L_0)+2$.
\item The linearized contact homology of $L_1'$ has only one element and in homology grading $0$ is
	\[HLC_{0}(\R^{2n+1}, L_1')=\Z_2.\]
\item $L_0$ and $L_1'$ are not Legendrian isotopic.
\end{enumerate}	 
\end{prop}
The proof of this proposition is identical to the proof of Proposition~\ref{1propl1}. 
To obtain interesting
examples when $n$ is odd we let $K_1$ be the connected sum of $L_1$ and $L_1'$ and let $K_k$ be the
connected sum of $K_{k-1}$ with $K_1.$
\begin{thm}\label{1mainex1}
The classical invariants of $K_k$ agree with those of $L_0$ but $K_k$
and $K_j$ are not Legendrian isotopic if $k\ne j$. 
\end{thm} 
This follows from Propositions \ref{1proptb}, \ref{1propl1}, 
\ref{4L1prime.prop},
Lemma~\ref{1sumcomp}, and the computations of the
linearized contact homology for $L_1$ and $L_1'.$
\end{ex}

\begin{ex}
Let $F_g$ be the Legendrian surface of genus $g$ with front obtained by ``connect
summing'' several standard $2$-spheres as shown in Figure \ref{1fig:fg}
\begin{figure}[ht]
	{\epsfxsize=3.5in\centerline{\epsfbox{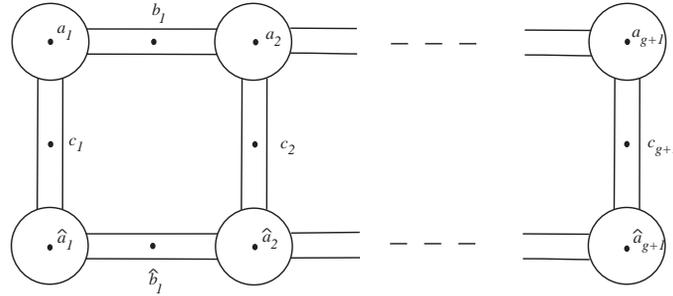}}}
	\caption{Top view of $F_g.$} 
	\label{1fig:fg}
\end{figure} .
Then $\A(F_g)$ is generated by
$$
\{a_j,\hat a_j,b_k,\hat b_k, c_j\}_{1\le j\le 1+g,\,\, 1\le k\le g},
$$ 
where $|a_j|=|\hat a_j|=2$ and
$|b_k|=|\hat b_k|=|c_j|=1$. Using projection to and slicing with the
$z_1$- and $z_2$-planes as above we find
\begin{align}\notag
\pa a_1 &= b_1 + c_1,\\\notag
\pa \hat a_1 &= \hat b_1 + c_1,\\\notag
\pa a_j &= b_{j-1}+b_j+c_j,\text{ for }j\ne 1, 1+g,\\\notag
\pa \hat a_j &= \hat b_{j-1}+\hat b_j+c_j,\text{ for }1<j<1+g,\\\notag
\pa a_{1+g} &= b_g + c_{1+g},\\\notag
\pa \hat a_{1+g} &= \hat b_g + c_{1+g},\\\notag
\pa b_k &=\pa \hat b_k=\pa c_j=0,\text{ for all }j,k.
\end{align} 
We find $HC_\ast(\R^5,F_g)=\Z_2\la a, b_1,\dots,b_g\ra$ where $|a|=2$
and $|b_i|=1$. Let $L_1$ be as in Example \ref{1BS} and
define $F_g^0 = F_g$ and $F_g^k=F^{k-1}_g\#L_1$. Then
the subspace of elements of grading $-1$ in $HLC_\ast(\R^5,F^k_g)$ is
$k$-dimensional. Thus, $F^k_g$ and $F^j_g$ are not Legendrian isotopic
if $j\ne k$.   
Clearly $tb(F^k_g) = tb(F^j_g).$
To see that $r(F^k_g) = r(F^j_g),$ it suffices to check, via front projections,
that the Maslov classes are the same on the generators of $H_1(F^j_g) = 
H_1(F^j_g).$
\end{ex}

\subsection{Front spinning} \label{1frontspinning.section}
\begin{figure}[ht]
	{\epsfxsize=3in\centerline{\epsfbox{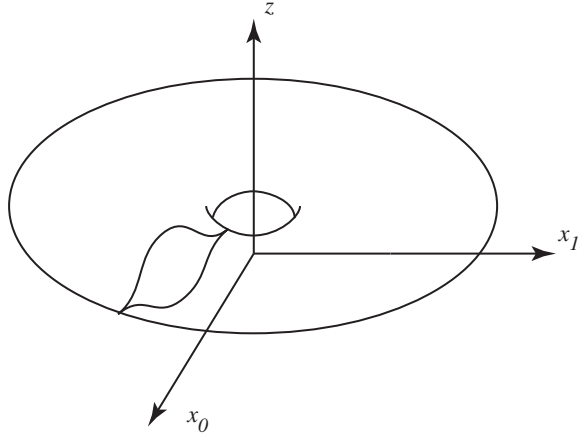}}}
	\caption{The front of $\Sigma L.$}
	\label{1fig:1spin}
\end{figure}

Given a Legendrian manifold $L\subset \R^{2n+1}$ we construct the \df{suspension} of $L,$
denoted $\Sigma L$ as follows: let $f\colon M\to \R^{2n+1}$ be a parameterization of $L,$ and
write 
\begin{equation} \notag
f(p)=(x_1(p),y_1(p), \ldots, x_n(p),y_n(p),z(p)),  
\end{equation}
The front projection $\Pi_F(L)$ of $L$ is the subvariety
of $\R^{n+1}$ parameterized by $\Pi_F\circ f(p)=(x_1(p),\ldots, x_n(p),z(p)).$ 
We may assume that $L$ has been translated so that $\Pi_F(L)\subset \{x_1>0\}.$
If we embed $\R^{n+1}$ into $\R^{n+2}$ via $(x_1, \ldots, x_n, z)\mapsto (x_0, x_1,\ldots x_n,z)$
then $\Pi_F(\Sigma L)$ is obtained from $\Pi_F(L)\subset \R^{n+1}$ by
rotating it around the subspace 
$\{x_0=x_1=0\}.$ See Figure~\ref{1fig:1spin}. 
We can parameterize $\Pi_F(\Sigma L)$ by 
$(\sin\theta x_1(p), \cos\theta x_1(p), x_2(p), \ldots, x_n(p)),$
$\theta\in S^1$. Thus, $\Pi_F(\Sigma L)$ is the front for a
Legendrian embedding $L\times S^1\to\R^{2n+3}$ we denote the
corresponding Legendrian submanifold $\Sigma L.$ We have the following simple lemma.
\begin{lem}\label{1spinclass}
	The Legendrian submanifold $\Sigma L\subset\R^{2n+3}$ has
	\begin{enumerate}
		\item the topological type or $L\times S^1,$
		\item the Thurston--Bennequin invariant $\tb(\Sigma L)=0,$
		\item Maslov class determined by
			\[\mu_{\Sigma L}(g)= \begin{cases}
      				\mu_L(h), &  \text{if $g=\iota h$ where $\iota:\pi_1(L)\to \pi_1(\Sigma L)$
							is the natural inclusion} \\
     				 0, & \text{if $g=[\text{point}\times S^1]$},
    			\end{cases}\]
		\item the same Maslov number as $L,$ $m(\Sigma L)=m(L),$ and
		\item the rotation class of $\Sigma L$ is determined by the rotation class of $L.$
	\end{enumerate}
\end{lem}

Though it seems difficult to compute the full contact homology of
$\Sigma L$ we can extract useful information about its linear
part. To this end we introduce the following notation. Let
$L\subset\R^{2n+1}$ be a Legendrian submanifold and let
$\A=\A(L)=\Z_2[H_1(L)]\la c_1,\dots,c_m\ra$ be the graded algebra generated by
its Reeb chords. We associate auxiliary algebras to $L$ which are
free unital algebras over $\Z_2[H_1(\Sigma L)]$. For any integer $N$,
let $\Z_{2N}^0\subset\Z_{2N}$ denote the subgroup of even
elements and let $\Z_{2N}^1=\Z_{2N}\setminus\Z_{2N}^0$.  
\begin{itemize}
\item Let  
$$
\A_\Sigma^N(L)=\Z_2[H_1(\Sigma L)]\Bigl\la c_j[\alpha],\hat
c_j[\beta]\Bigr\ra_{1\le j\le m,\,\,\alpha\in\Z_{2N}^0,\,\,\beta\in\Z_{2N}^1},
$$ 
where
$|c_j[\alpha]|=|c_j|$, $\alpha\in\Z_{2N}^0$ and $|\hat c_j[\beta]|=|c_j|+1$, $\beta\in\Z_{2N}^1$.
\item
For $\beta\in\Z_{2N}$, define the subalgebra
$\A_\Sigma^N[\beta]\subset\A_\Sigma^N=\A_\Sigma^N(L)$ as 
$$
\A_\Sigma^N[\beta]=\Z_2[H_1(\Sigma L)]\la
c_j[\beta-1],c_j[\beta+1],\hat c_j[\beta]\ra_{1\le j\le m}.
$$ 
\item
Define the algebra
$$
\A_\sigma(L)=\Z_2[H_1(\Sigma L)]\la c_j,\hat c_j\ra_{1\le j\le m},
$$
where $|\hat c_j|=|c_j|+1$.
\end{itemize}
We note that there is a natural homomorphism
$\pi\colon\A_\Sigma^N\to\A_\sigma$ defined on generators by
$\pi(c_j[\alpha])=c_j$, and $\pi(\hat c_j[\beta])=\hat c_j$. Also note
that for each $\alpha\in\Z_{2N}^0$ there is a natural inclusion
$\Delta[\alpha]\colon\A\to\A_\Sigma^N$ defined on generators by 
$\Delta[\alpha](c_i)=c_i[\alpha]$, and using the
natural inclusion $H_1(L)\to H_1(\Sigma L)$ on coefficients.   

Viewing $\A(L)$ and $\A_\sigma(L)$ as a vector space
over $\Z_2$, see \eqref{DGAasVSP}, and again using 
$H_1(L)\to H_1(\Sigma L)$ we define the linear map 
$\Gamma\colon \A(L)\to\A_\sigma(L)$ by
$$
\Gamma(1)=0,\quad \Gamma(t_1^{n_1}\dots t_s^{n_s}c_{i_1}\dots
c_{i_r})=t_1^{n_1}\dots t_r^{n_r}\Bigl(
\sum_{j=1}^r c_{i_1}\dots c_{i_{j-1}} \hat c_{i_j} c_{i_{j+1}}\dots
c_{i_r}
\Bigr).
$$

\begin{prop}\label{1spinhom}
Let $c_1,\ldots, c_m$
be the Reeb chords of $L$ and let $(\A,\partial)$ denote its DGA.  
Then there exits an even integer $N$ and a representative $X$ of the
Legendrian  isotopy class of $\Sigma L$ with associated DGA $(\A(X),\pa_\Sigma)$ satisfying
\begin{align}
\A(X)&=\A_\Sigma^N,\\
\pa_\Sigma c_i[\alpha]&=\Delta[\alpha](\pa c_i),\text{ for all }\alpha\in\Z_{2N}^0,\\
\pa_\Sigma \hat c_i[\beta]&=c_i[\beta-1]+c_i[\beta+1]+\gamma_i^{1}[\beta]+\gamma_i^2[\beta],
\text{ for all }\beta\in\Z_{2N}^1
\end{align} 
where, $\gamma_i^2[\beta]$ lies
in the ideal of $\A_\Sigma^N[\beta]$ generated by all monomials
which are quadratic in the variables $\hat c_1[\beta],\ldots, \hat
c_m[\beta],$, where $\gamma_i^1[\beta]\in\A_\Sigma^N[\beta]$ is linear in
the generators $\hat c_i[\beta]$ and satisfies
\begin{equation}
\pi(\gamma_i^1[\beta])=\Gamma(\pa c_i).
\end{equation}
Moreover, $(\A(X),\pa_\Sigma)$ is stable tame isomorphic to $(\A_\Sigma^2,\pa_\Sigma)$.
\end{prop}

We will prove this proposition in the next subsection but first we consider its consequences.
To simplify notation, we consider the algebra generated over $\Z_2$ instead of $\Z_2[H_1(\Sigma L)].$ 


\begin{ex}\label{torusknot.eg}
Let $T_k$ be the Legendrian torus knot in Figure~\ref{1fig:knotex}
\begin{figure}[ht]
	{\epsfxsize=3in\centerline{\epsfbox{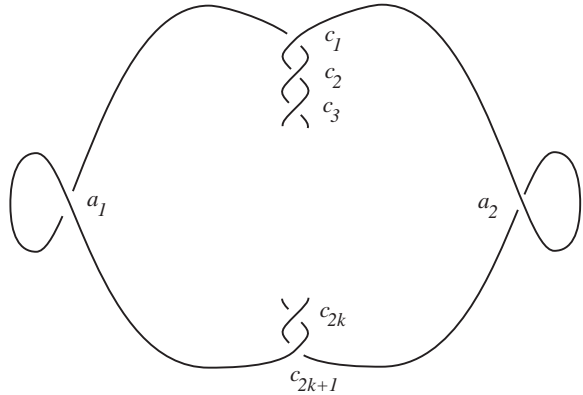}}}
	\caption{The knots $T_k.$}
	\label{1fig:knotex}
\end{figure} 
with rotation number $r(T_k)=0$. The algebra for
$T_k$ is $\A(T_k)=\Z_2\la a_1, a_2, c_1,\dots, c_{2k+1}\ra$ with $|a_1|=|a_2|=1$ and
$|c_j|=0$ for all $j$. With Greek letters running over the 
integers in $[1,2k+1]$ we have
\begin{align*}
	 	\partial a_1&= 1+\sum_{\alpha}c_{\alpha}
		+\sum_{\alpha>\beta>\gamma} c_\alpha c_\beta c_\gamma+\dots
		+c_{2k+1}c_{2k}\dots c_{1},\\
		\partial a_2&= 1+\sum_{\alpha}c_{\alpha}+
		\sum_{\alpha<\beta<\gamma} c_\alpha c_\beta
		c_\gamma+\dots
		+c_1c_2\dots c_{2k+1},\\
		\partial c_j&= 0, \quad \text{all j}.\\
\end{align*}

We note that $\pa^g$, where $g$ is the elementary automorphism with
$g(c_1)=c_1+1$ and which fixes all other generators, is augmented and
that the linearized homology of $(\A,\pa^g)$ is (as a vector space
without grading) $\Z_2^{2k+1}$. Applying the suspension operation $n$ times we get a
Legendrian $n$-tori $\Sigma^n T_k$ with $\tb(\Sigma^n T_k)=0$ for all
$n>0$, with rotation classes independent of $k$ (see
Lemma~\ref{1spinclass}), and with Maslov number equal to $0$. The algebras of 
$\Sigma^n T_k$ admit an elementary isomorphism (add $1$ to each
$c_1[\alpha_1][\alpha_2]\dots[\alpha_n]$ with $|c_1[\alpha_1][\alpha_2]\dots[\alpha_n]|=0$)
making them good and such that the corresponding linearized
homology is isomorphic to $\Z_2^{2^n(2k+1)}$. This implies that every chord generic Legendrian
representative of $\Sigma^n T_k$ has at least $2^n(2k+1)$ Reeb
chords. Moreover, since $\Sigma^n T_j$ has a representative with
$2^n(2j+3)$ Reeb chords it is easy to extract an infinite family of
pairwise distinct Legendrian $n$-tori from the above.  
\end{ex}

Using Example~\ref{torusknot.eg} we find
    
\begin{thm}\label{manytori.thm}
There are infinitely many Legendrian 
$n$-tori in $\R^{2n+1}$ that are pairwise not Legendrian isotopic even
though their classical invariants agree.
\end{thm}


\begin{ex} \label{1Whdouble.eg}
As a final family of examples we consider the Whitehead doubles of the unknot $W_s$ shown in 
Figure~\ref{1fig:wd}. Note that $r(W_s) = 0.$
\begin{figure}[ht]
	{\epsfxsize=4in\centerline{\epsfbox{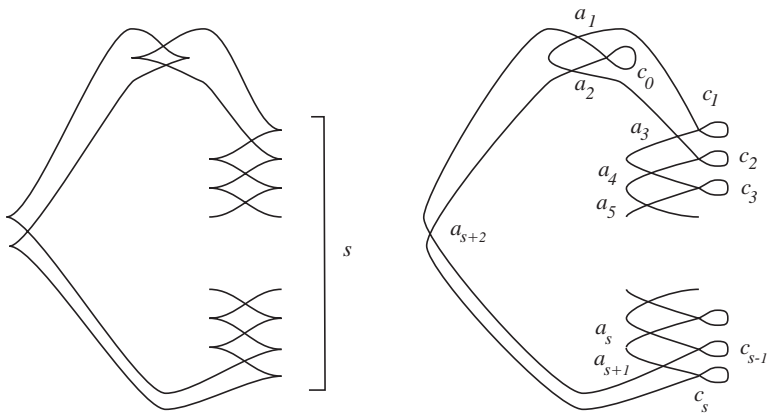}}}
	\caption{The front projection (left) and Lagrangian projection (right) of the
		knots $W_s.$}
	\label{1fig:wd} 
\end{figure} 
The algebra for $W_s$ is 
\[\A(\R^3, W_s)=\Z_2 \langle c_0, \ldots,c_s, a_1,\ldots, a_{s+2} \rangle\] 
with
$|c_i|=1, |a_1|=-|a_2|=s-2$ and $|a_i|=0$ for $i>2.$ Moreover,
\begin{eqnarray*}
	 	\partial c_0&=& 1+ a_1a_2+ a_{s+2},\\
		\partial c_1&=&  1+ a_3,\\
		\partial c_i&=& 1+a_{i+1}a_{i+2} \quad \text{for $i>1$},\\
		\partial a_i&=&0\quad \text{for all $i$}.
\end{eqnarray*}
The differential is clearly not augmented, but there will be a unique
tame graded automorphism making it augmented. 
The only feature of the linearization we use is that 
\begin{eqnarray*}
	 	HLC_{2-s}(\R^5, \Sigma W_s)&=& \Z_2,\\
		HLC_i(\R^5, \Sigma W_s)&=& 0 \quad \text{for $i<2-s$}.
\end{eqnarray*}
The computation of these groups for $\Sigma^n W_s$ yields the same
answer. Thus, they are all distinct and we get another proof of
Theorem \ref{manytori.thm}.
\end{ex}

\subsection{Proof of the front spinning proposition}
To prove Proposition~\ref{1spinhom} we first analyze another
Legendrian submanifold. Let $\psi:\R\to \R_{\geq0}$ 
be a smooth small perturbation of the constant function 1 that is
a 2-periodic and has non-degenerate local maxima at even integers and local minima at odd
integers. Given a Legendrian submanifold $L\subset\R^{2n+1}$ parameterized as above we define
the front  $\Pi_F(L\times \R)$ in $\R^{n+1}$ by 
\begin{equation} \label{1G.eqn}
	G(p,t)=(\psi(t)z(p),t, x_1(p),\ldots, x_n(p)).
\end{equation}
Denote the resulting Legendrian submanifold of $\R^{2(n+1)+1}$ by $L\times \R.$ Heuristically, $L\times \R$
is a kind of ``cover'' of $\Sigma L$ and the boundary map of $\Sigma L$ shall be determined by 
studying the boundary map of $L\times \R.$ We begin with a simple lemma.
\begin{lem}
	For each Reeb chord $c_j$ of $L$ there are $\Z$ Reeb chords, $c_j[n],$ for $L\times \R$; moreover,
	$|c_j[2n]|=|c_j|+1$  and $|c_j[2n+1]|=|c_j|.$
\end{lem}
\begin{lem}\label{1splitR}
A holomorphic disk in $\C^{n+1}$ with boundary on $\Pi_{\C}(L\times \R)$ cannot intersect the
hyperplane $z_0=k, k\in \Z.$ In addition any holomorphic disk with a negative corner at $c_j[2k]$
or a positive corner at $c_j[2k+1]$ must lie entirely in the plane $z_0=2k, z_0=2k+1,$ respectively.
\end{lem}
The proof of this lemma is identical to the proof of Lemma~\ref{1sumcomp} once one has
drawn the projection of $\Pi_{\C}(L\times \R)$ onto the
$z_0$-plane. Also, an argument similar to that in the proof of
Lemma~\ref{1sumcomp} in combination with Lemma~\ref{splittran1} shows:
\begin{lem}\label{1twodisks}
There is a unique holomorphic disk with positive corner at $c_j[2n]$
and negative corner at $c_j[2n\pm 1]$ and this disk is transversely cut out.
\end{lem}
We now discuss perturbations of $L$ necessary to ensure the appropriate moduli spaces are  manifolds.
The reader might want to skip this until Part 2 has been read, in particular the proof of 
Lemma~\ref{splittran1}.
To ensure all our moduli spaces are cut out transversely we might have to perturb our Legendrian near
positive corners of non-transversely cut out holomorphic disks. Note that due to Lemma~\ref{1splitR}
we see disks with positive corners at a Reeb chord in $x_0=2k+1$ lie in $z_0=2k+1.$ Thus the linear problem
splits for these disks and an argument similar to the proof of
Lemma~\ref{splittran2} shows they are all transversely cut out. In
perturbing $L\times\R$ to be generic we can assume the perturbation 
is near the hyperplanes $x_0=2k, k\in \Z$ and none of the Reeb chords move.

Now let $\B_{2k}=\Z_2[H_1(L)]\langle c_j[2k-1],c_j[2k], c_j[2k+1]\rangle_{j=1}^m$ and $\B_{2k+1}=
\Z_2[H_1(L)]\langle c_j[2k+1]\rangle_{j=1}^m.$ These are all sub-algebras of the algebra $\B$ generated
by all the Reeb chords for $L\times\R.$ Let $\partial_{\R}$ be the boundary map for $L\times\R.$
From the above Lemmas we clearly have 
\[\partial_{\R}(\B_{2k+1})\subset \B_{2k+1}\]
and
\[\partial_{\R}(\B_{2k})\subset \B_{2k}.\]
Moreover Lemma~\ref{1splitR} and our discussion of the generic perturbation above give
\begin{lem}
	Let $\Gamma_{\pm 1}:\A\to \B_{\pm 1}$ be given by $\Gamma_{\pm 1}(c_j)=c_j[\pm 1].$ Then
	\[\partial_{\R} c_j[\pm 1]= \Gamma_{\pm 1}(\partial c_j).\]
\end{lem}

To understand $\partial_{\R}$ on $\B_0$ we begin with
\begin{lem}
	$\partial^2_{\R}=0.$
\end{lem}
\begin{proof}
Let $F$ be the part of the front of $L\times\R$ between $x_0=-(2k+\frac32)$ and $x_0=2k+\frac32,$ say. 
Let $F'$ be $F$ translated $4k+10$ units in the $x_1$-direction. For sufficiently large $k$ $F\cap F'=
\emptyset.$ For such a $k$ let $G\cup G'$ be $F\cup F'$ rotated by $\frac\pi2$ around the affine subspace 
$\{x_0=0, x_1=2k+5\}.$  Now make $f\cup F'\cup G\cup G'$ into a closed Legendrian submanifold $L'$ by
adding ``round corners.''

Using the lemmas above and a monotonicity argument as in the proof of Lemma~\ref{1stabvanish}
it is easy to see that the boundary map for $L'$ agrees with $\partial_{\R}$ on $\B_0$ and $B_1.$ 
Since we know the square of the boundary map for a closed compact Legendrian is 0 the lemma follows.
\end{proof}
\begin{lem}\label{spindiff}
	We can choose $L\times \R$ so that 
	the part of $\partial_{\R}(c_j[0])$ that is constant in the generators $c_0[0],\ldots, c_m[0],$ 
	is 
	\[c_j[-1]+c_j[1].\]
\end{lem}
\begin{proof}
This term in present in $\partial_{\R}(c_j[0])$ by Lemma~\ref{1twodisks}. 

To see there are no
disks $D$ with just one corner (which of course is positive at $c_j[0]$) assume we have such a
disk $D$. Then consider $\psi_s:\R\to \R$ where $\psi_0=\psi$ and $\psi_1(s)=1$ and
the corresponding Legendrian submanifolds $(L\times \R)_s$ whose fronts are defined using $\psi_s$ 
just as the front of $L\times \R$ used $\psi$ in (\ref{1G.eqn}). 
As we isotope $L\times\R=(L\times\R)_0$ to $(L\times \R)_1$ we see that $D$
will have to converge to a (possibly broken) disk for $(L\times\R)_1.$ But arguing as in the 
lemmas above we can see that any such disk will have to have $z_0$ constant and thus corresponds to
a disk for $L.$ However there can be no rigid holomorphic disk for $L$ with corner at $c_j$ since
$|c_j|=|c_j[0]|-1=1-1=0.$ Moreover, if we have a broken holomorphic disk one can similarly see that
one of the pieces of the broken disk will have negative formal
dimension and thus cannot exist since we took $L$ to be generic. 

One may similarly argue that
there are no holomorphic disks with one positive corner at $c_j[0]$
and all negative corners at Reeb chords $c_k[\pm 1]$, where $k\ne j$ for
some $j$.
\end{proof}

Lemma \ref{spindiff} implies
\[\partial_{\R} (c_j[0])= c_j[-1]+c_j[1]+ \eta_j + r_j,\]
where $\eta_j$ is the part of $\partial_{\R} (c_j[0])$ linear in $c_1[0],\ldots, c_m[0]$ and
$r_j$ is the remainder (terms which are at least quadratic in the $c_j[0]$'s). Since $\partial_{\R}^2=0$
we see that 
\begin{equation}\label{1step1}
\partial_{\R}(c_j[-1]+c_j[1])=\sigma(\eta_j),
\end{equation} 
where $\sigma$ is the algebra homomorphism defined by $\sigma c_j[0]= c_j[-1]+c_j[1]$ and
$\sigma c_j[\pm 1]=0.$
A straightforward calculation shows that 
\begin{equation} \notag
\eta_j=\Gamma_0(\partial c_j)
\end{equation}
is a solution to (\ref{1step1}) where $\Gamma_0:\A\to \B_0$ is the linear map defined on monomials
by
\begin{align*}\Gamma_0(c_{j_1}\ldots c_{j_r})= &c_{j_1}[0]c_{j_2}[-1]\ldots 
c_{j_r}[-1] + c_{j_1}[1]c_{j_2}[0]
c_{j_3}[-1]\ldots c_{j_r}[-1] +\dots\\ + &c_{j_1}[1]\ldots c_{j_{r-1}}[1] c_{j_r}[0].\end{align*}
While this is not the only solution to (\ref{1step1}), it is unique in the following sense:
let $\B'=\Z_2[H_1(L)]\langle c_1,\ldots, c_m, c_1[0],\ldots, c_m[0]\rangle$ and define 
$\pi: \B_0\to \B'$ by $\pi(c_j[\pm 1])=c_j$ and $\pi(c_j[0])=c_j[0].$
\begin{lem}\label{1uniquesoln}
	If $\alpha$ is linear in the $c_j[0]$'s and $\sigma(\alpha)=0$ then $\pi(\alpha)=0.$
\end{lem}
\begin{proof}
Assume $\alpha$ satisfies the hypothesis of the lemma. Let $a c_j[0] b$ 
be a monomial in $\alpha$
 with
$l(a)=l_1$ and $l(b)=l_2$ where $l(\cdot)= \text{length of monomial}.$  
We claim that since $\sigma(\alpha)=0,$ there must be another monomial
in $\alpha$ of the form $a' c_j[0] b'$ where $\pi(a)=\pi(a')$ and
$\pi(b)=\pi(b').$  The lemma will follow. To see this note that 
\[\sigma(a c_j[0] b)=a c_j[-1] b+a c_j[1] b.\]
Since the $(l_1+1)$th letter in these monomials is different they can only be canceled by a monomial 
of the form  $\sigma (a' c_j[0] b')$ or by terms coming from two separate monomials $w_1$ and $w_2$ in $\alpha.$
So either we are done or two of the four terms in $\sigma(w_1+w_2)$ are canceled by $\sigma(a c_j[0] b)$
leaving two other terms. Note the two leftover terms still have different $(l_1+1)$th letters.
So once again either there is a term of the form $a' c_j[0] b'$ to cancel these two terms or there are 
two more monomials in $\alpha.$ We clearly will eventually find the desired  term $a' c_j[0] b'$ 
(induction on the number of terms in $\alpha$). One should observe that $a'$ and $b'$ have the properties 
stated above because in each step in the above cancellation process we are replacing $c_j[\pm1]$'s
with $c_j[\mp1]$'s.
\end{proof}

We are now ready to prove our main result of this subsection.

\begin{proof}[Proof of Proposition~\ref{1spinhom}]
Consider $\Sigma L$ represented by rotating the front of $L$ around a circle $C$ with radius 
$\frac{1}{\pi}N$ for some even integer $N.$  Perturb this non-generic
front with a function $\phi$ on $C$ similar to (\ref{1G.eqn}) so that
$\phi$ approximates the constant function 1, has local minima at
angles $2m\cdot\frac{\pi}{N}$ and local maxima at
$(2m+1)\cdot\frac{\pi}{N}$, $m=1,\dots,N-1$. Let $X_N$ denote the
corresponding Legendrian submanifold. Then there is a natural $1-1$
correspondence between the generators 
of the algebra $\A(X_N)$ and the generators of $\A_\Sigma^N(L)$. 

Considering the projections of $X_N$ to the complex lines in $\C^{n+1}$
which intersect $\R^{n+1}$ in lines through antipodal local minima
of $\phi$ we see that the differential $\pa_\Sigma$ of
$\A(X_N)=\A_\Sigma^N$ preserves the subalgebras
$\A^N_\Sigma[\beta]$ for every $\beta\in\Z_{2N}^1$.
Moreover, a finite part of the front of $X_N$ over an arc on $C$ between two
minima is for $N$ sufficiently large an arbitrarily good approximation
of the part of the front of $L\times\R$ between $-1$ and $1$. In fact, since all spaces of rigid
disks on $L\times\R$ are transversely cut out there is a neighborhood
of $L\times\R$ in the space of (admissible) Legendrian submanifolds
such that the moduli spaces of rigid disks on any Legendrian
submanifold in this neighborhood is canonically isomorphic to those of
$L\times\R$. This can be seen as follows: by Gromov compactness, there exists
a neighborhood of $L\times\R$ in the space of admissible Legendrian
submanifolds such that for any $Y$ in this neighborhood there are no 
holomorphic disks with boundary on $Y$ and with negative formal
dimension, pick a generic type (A) isotopy from $L\times\R$ to $Y$ and
apply Lemma \ref{1typeA}. Thus,  for sufficiently large $N$ the subalgebras
$(\A_\Sigma^N[\beta],\pa_\Sigma)$ are all isomorphic to the algebras
$(\B_0,\pa)$. The first part of the proposition now follows from
Lemmas \ref{spindiff} and \ref{1uniquesoln}.

For the statement of stable tame isomorphism class, note that the
subalgebras of $\A_\Sigma^N$ generated by all Reeb chords corresponding
to maxima and minima over the circle in the closed upper (lower) half planes
are both isomorphic to the subalgebra of $\A(L\times\R)$ generated by
Reeb chords between $1$ and $N+1$. Change the front of $L\times\R$
by shrinking the minima of $\psi$ over $1$ and $N+1$ until the
corresponding Reeb chords are shorter than all other Reeb chords. We
can then find a Legendrian isotopy which cancels pairs of maxima and
minima of $\psi$ between $1$ and $N+1$, leaving one maximum. Thus the
subalgebra generated by Reeb chords between $1$ and $N+1$ is stable
tame isomorphic to $\B_0$. Moreover, if the stabilizations and tame
isomorphisms corresponding to self tangencies and handle slides in the
canceling process are
constructed as in the proofs of Lemmas \ref{1typeA} and \ref{1typeB}, respectively, then 
Reeb chords which are shorter than the positive chord of a handle
slide disk and shorter than both chords canceling in a self tangency
are left unchanged by the tame isomorphisms. Thus the subalgebra
$\B_1$ and $\B_{N+1}$ are left unchanged by this chain of
stabilizations and tame isomorphisms and hence it induces a chain of
stabilizations and tame isomorphisms connecting $\A_\Sigma^2(L)$ to
$\A_\Sigma^N(L)$.      
\end{proof}

\subsection{Remarks on the examples} \label{1remarksontheexamples.section}
\begin{rmk}
\label{1connectsum.remark}
In dimension $3$ the connected sum of Legendrian knots is well defined \cite{Etnyre-Honda2}. 
However in higher dimensions there are several ways to make a Legendrian version of the connected sum.
Lemma~\ref{1sumcomp} discusses one such way. 
However there are other direct generalization of the $3$ dimensional connected sum. 
Thus the correct definition of connected sum is not clear.
Even if we just consider the ``cusp connected sum''
(from Lemma~\ref{1sumcomp}) it is still not clear if it is well defined. So we ask
\bmini{\em
Is the connected sum well defined?}
\emini
And more specifically
\bmini{\em
Does the cusp connected sum depend on the cusps chosen in the construction?}
\emini
\end{rmk}

\begin{rmk}
Colin, Giroux and Honda have announced the following result: in dimension $3$ if you fix the 
$\tb, r$ and a knot type there are only finitely many Legendrian knots realizing this data. 
If one considers this 
question in higher dimensions most of our examples described above provide counterexamples to the
corresponding assertion. In particular consider Theorem~\ref{1mainex1}.
\end{rmk}

\begin{rmk}
Given a Legendrian submanifold $L$ we can define an invariant $N(L)$ to be the minimal number of
Reeb chords associated to a generic representative of $L.$ 
One can ask
\bmini{\em
Is there an effective bound on $N(L)$ in terms of the Thurston-Bennequin invariant and rotation class?}
\emini
The result of Colin, Giroux and Honda mentioned above indicate a positive answer to this question 
in dimension $3.$ While our examples above show that these answer is a resounding NO in dimensions
above $3.$
\end{rmk}

\part{Analytic Theory}

In this part of the paper we present the details of the analytic parts
of the proofs from Part 1. 

\section{Admissible Legendrian submanifolds and isotopies}
\label{5Isotopies.section}
\subsection{Chord genericity}
Recall that a Legendrian submanifold $L\subset\R^{2n+1}$ is chord
generic if all its Reeb chords correspond to transverse double points
of the Lagrangian projection $\Pi_\C$.
For a dense open set in the space of paths of Legendrian embeddings, 
the corresponding $1$-parameter families $L_t$, $0\le t\le 1$,  are
chord generic except for a finite number of parameter values
$t_1,\dots,t_k$ where $\Pi_\C(L_{t_j})$ has one double point with
self-tangency, and where for some $\delta>0$ $\Pi_\C(L_t)$,
$(t_j-\delta,t_j+\delta)$, is a versal deformation of $\Pi_\C(L_{t_j})$, 
for $j=1,\dots,k$.  We call $1$-parameter families
$L_t$ with this property {\em chord generic 1-parameter families}. 

\subsection{Local real analyticity}
For technical reasons, we require our Legendrian
submanifolds to be real analytic in a neighborhood of
the endpoints of their Reeb chords and that self-tangency instants in
$1$-parameter families have a very special form.

\begin{dfn}\label{dfnadm}
A chord generic Legendrian submanifold $L\subset\C^n\times\R$ is 
{\em admissible} if for any Reeb chord $c$ of $L$ with endpoints $q_1$
and $q_2$ there are neighborhoods  
$U_1\subset L$ and $U_2\subset L$ of $q_1$ and $q_2$, respectively, such
that $\Pi_\C(U_1)$ and $\Pi_\C(U_2)$ are real analytic submanifolds of 
$\C^n$.
\end{dfn}  

We will require that self-tangency instants in $1$-parameter families
have the following special form. Consider
$0\in\C^n$ and coordinates $(z_1,\dots,z_n)$ on $\C^n$. Let $P_1$ and
$P_2$ be Lagrangian submanifolds of $\C^n$  
passing through $0$. Let $x=(x_1,\dots,x_n)\in\R^n$ and
$y=(y_1,\dots,y_n)\in\R^n$  be coordinates on $P_1$ and $P_2$,
respectively. Let $R_1\subset P_1$ and $R_2\subset P_2$ be the boxes
$|x_j|\le 1$ and $|y_j|\le 1$, $j=1,\dots,n$. Let $B_j(2)$ and
$B_j(2+\epsilon)$ for some small $\epsilon>0$ be the balls of radii
$2$ and $2+\epsilon$ around $0\in P_j$, $j=1,2$. We require that in
$R_1$, $P_1$ has the form
\begin{equation}\label{P_1}
\gamma_1\times \hat P_1
\end{equation}  
where $\gamma_1$ is an arc around $0$ in the real line in the
$z_1$-plane and where $\hat P_1$ is a Lagrangian submanifold of
$\C^{n-1}\approx \{z_1=0\}$. We require that in $R_2$, $P_2$ has the
form 
\begin{equation}\label{P_2}
\gamma_2(t)\times \hat P_2,
\end{equation}
where $\gamma_2$ is an arc around $0$ in the unit-radius circle
centered at $i$ in $z_1$-plane and where $\hat P_2$ is a Lagrangian
submanifold of $\C^{n-1}\approx \{z_1=0\}$, which meets $\hat P_1$
transversally at $0$.

If $q\in\C^n$ let $\lambda_q$ denote the complex line in $T_q\C^n$
parallel to the $z_1$-line. We also require that for every point  
$p\in B_j(2+\epsilon)\setminus B_j(2)$ the tangent plane $T_p P_j$
satisfies 
\begin{equation}\label{B(2)}
T_p P_j\cap \lambda_p=0,\quad j=1,2.
\end{equation}

\begin{dfn}\label{dfnadmst}
Let $L_t$ be a chord generic $1$-parameter family of Legendrian
submanifolds such that $L_0$ has a self-tangency.
We say that the self-tangency instant $L_0$ is {\em standard} if
there is some neighborhood $U$ of the self-tangency point and a
biholomorphism $\phi\colon U\to V\subset\C^n$ such that 
\begin{equation}\label{2LocFormofST}
\phi(L_t\cap U)=P_1\cup P_2(t)\cap N,
\end{equation} 
where $N$ is some neighborhood of $0\in\C^n$, and where $P_2(t)$ is
$P_2$ transalted $t$ units in the $y_1$-direction.
\end{dfn}

\begin{dfn}\label{5admissile_param.dfn}
Let $L_t$, $0\le t\le 1$ be a chord generic $1$-parameter family of
Legendrian submanifolds. Let $t_1,\dots,t_k$ be its self tangency
instants. We 
say that $L_t$ is an {\em admissible} $1$-parameter family if 
$L_t$ is admissible for all $t\ne t_k$, if there exists small disjoint
intervals $(t_j-\delta,t_j+\delta)$ where the
$1$-parameter family is constant outside some small neighborhood $W$ 
of the self-tangency point, and if all self-tangency instants
are standard. 
\end{dfn}

\begin{dfn}
A Legendrian submanifold $L\subset\R\times\C^n$ which is a
self-tangency instant of 
an admissible $1$-parameter family will be called 
{\em semi-admissible}.
\end{dfn}

\subsection{Reducing the Legendrian isotopy problem}\label{2Ham.iso}
We prove a sequence of lemmas which  reduce the
classification of Legendrian submanifolds up to Legendrian isotopy to
the classification of admissible Legendrian submanifolds
up to admissible Legendrian isotopy.

We start with a general remark concerning lifts of Hamiltonian
isotopies in $\C^n$.
If $h$ is a smooth function with compact support in $\C^n$ then the
Hamiltonian vector field    
$$
X_h=-\frac{\pa h}{\pa y_i}\pa_{x_i}+\frac{\pa h}{\pa x_i}\pa_{y_i}
$$ 
associated to $h$ generates a $1$-parameter family of diffeomorphisms
$\Phi_h^t$ of $\C^n$. Moreover, the vector field $X_h$ lifts to a
contact vector field 
$$
\tilde X_h=-\frac{\pa h}{\pa y_i}\pa_{x_i}+
\frac{\pa h}{\pa x_i}\pa_{y_i} +\left(h-y_i\frac{\pa h}{\pa y_i}\right)\pa_z
$$ 
on $\C^n\times\R$, which generates a $1$-parameter family
$\tilde\Phi_h^t$ of contact diffeomorphisms of $\C^n\times\R$ 
which is a lift of $\Phi_h^t$. We write $\Phi_h=\Phi^1_h$ and
similarly $\tilde\Phi_h=\tilde\Phi_h^1$.

We note for future reference that in case the preimage of the support
of $h$ in $L$ has more than one connected component we may define a
Legendrian isotopy of $L$ by moving only one of these components
(for a short time) using $\tilde X_h$ and keeping the rest of them
fixed.    

An {\em $\epsilon$-isotopy} is an isotopy during which no point moves a 
distance larger than $\epsilon>0$.

\begin{lma}\label{admstat}
Let $L$ be a Legendrian submanifold. Then, for any $\epsilon>0$, there
is an admissible Legendrian submanifold $L_\epsilon$ which is
Legendrian $\epsilon$-isotopic to $L$. 
\end{lma}

\begin{pf}
As mentioned, we may after arbitrarily small Legendrian isotopy assume 
that $L$ is chord generic. Thus,
it is enough to consider one transverse double point. We may assume
that one of 
the sheets of $L$ close to this double point is given by
$x\mapsto (x,df(x),f(x))$ for some smooth function $f$. Let $g$ be a
real analytic function approximating  
$f$ (e.g. its Taylor polynomial of some degree). Consider a
Hamiltonian $h$ which is $h(x,y)=g(x)-f(x)$ in this neighborhood and
$0$ outside some slightly larger neighborhood. It is clear that the
corresponding Hamiltonian vector field can be made arbitrarily
small. Its flow map at time $1$ is given by
$\Phi_h^1(x,y)=\left(x,y+dg(x)-df(x)\right)$. 
Using this and suitable cut-off functions for the lifted Legendrian
isotopies the lemma follows.
\end{pf}

\begin{lma}
Let $L_t$ be any chord generic Legendrian isotopy from an admissible
Legendrian submanifold $L_0$ to another one $L_1$. Then for any $\epsilon>0$,
there is an admissible Legendrian isotopy
$\epsilon$-close to $L_t$ connecting $L_0$ to $L_1$.
\end{lma}

\begin{pf}
Let $t_1,\dots,t_k$ be the self tangency instants of the
isotopy. First change the isotopy so that there exists small disjoint
intervals $(t_j-\delta,t_j+\delta)$ where the
$1$-parameter family is constant outside some small neighborhood $W$ 
of the self-tangency point. Consider the restriction of the isotopy to
the self-tangency free regions. 
The $1$-parametric version of the proof of Lemma \ref{admstat} clearly 
applies to transform this part of the isotopy into one consisting of
admissible Legendrian submanifolds. Then change the isotopy in the
neighborhoods of the self tangency instants, 
using essentially the same argument as above, to a  
self-tangency of the from given in \eqref{P_1} and \eqref{P_2}. 

It remains to show how to fulfill the condition \eqref{B(2)}. To this
end, consider a
Lagrangian submanifold of the form \eqref{P_1}. Locally it is given by
$(x,df(\hat x))$, where $\hat x=(x_2,\dots,x_n)$. Let $\phi(x)$ be a
function which equals $0$ in $B(2-\frac12\epsilon)$ and outside 
$B(2+2\epsilon)$ and has $\frac{\pa^2\phi}{\pa x_1 x_j}\ne 0$ for
some $j$ at all points in $B(2+\epsilon)\setminus B(2)$. (For example
if $K$ is a small constant a suitable cut-off of the function
$Kx_1(x_2+\dots x_n)$ has this property). We see as above that our
original Legendrian is Legendrian isotopic to
$(x,df(x)+d\phi(x))$. The tangent space of the latter submanifold is
spanned by the vectors
\begin{align}
&\pa_{x_1}+\sum_j\frac{\pa^2\phi}{\pa x_j\pa x_1}\pa_{y_j},\\
&\pa_{x_r}+\sum_j\frac{\pa^2\phi}{\pa x_j\pa x_r}\pa_{y_j}+
\frac{\pa^2 f}{\pa x_j\pa x_r}\pa_{y_j},\quad 2\le r\le n.
\end{align}    
Any non-trivial linear combination of the last $n-1$ vectors projects
non-trivially to the subspace $dx_1=dy_1=dy_2=\dots=dy_n=0$. The
first vector lies in the subspace $dx_2=\dots=dx_n=0;$ thus,
since the first vector does not lie in the $z_1$-line because
$\frac{\pa^2\phi}{\pa x_1x_j}\ne 0$ for some $j\ne 1$, no linear
combination of the vectors does either. $P_2$ can be deformed in a
similar manner.

After the self-tangency moment is passed it is easy to Legendrian
isotope back to the original family through admissible
Legendrian submanifolds.
\end{pf}

\section{Holomorphic disks}\label{2hodi}
In this section we establish notation and ideas that will be used throughout
the rest of the paper.

\subsection{Reeb chord notation}
\label{hodic}
\begin{figure}[ht]
	{\epsfxsize=2.75in\centerline{\epsfbox{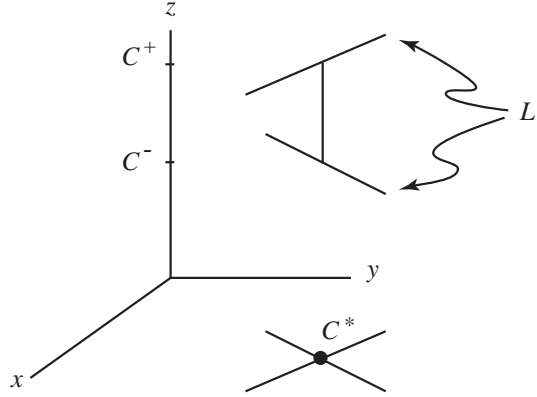}}}
	\caption{A Reeb chord in $\R^3.$ }
	\label{1fig:cpic}
\end{figure}

Let $L\subset\C^n\times\R$ be a Legendrian submanifold and let $c$ be
a Reeb chord of $L$. The $z$-coordinate of the upper and lower end
points of $c$ will be denoted by 
$c^+$ and $c^-,$ respectively. See Figure~\ref{1fig:cpic}. So as a
point set $c=c^\ast\times[c^-,c^+]$ and the action of $c$ is simply
$\action(c)=c^+-c^-.$  

If $r>0$ is small enough so that $\Pi_{\C}^{-1}(B(c^\ast, r))$
intersects $L$ is  
exactly two disk about the upper and lower end points of $c,$ then we
define $U(c^\pm,r)$ to 
be the component of $\Pi_{\C}^{-1}(B(c^\ast, r))\cap L$ containing
$c^\ast\times c^\pm.$  

\subsection{Definition of holomorphic disks}
If $M$ is a smooth manifold then let $\sblv^{\rm loc}_k(M,\C^n)$ denote the
Frechet space of all functions which agree locally with a function
with $k$ derivatives in $L^2$. Let $\Delta_m\subset\C$ denote the unit disk
with $m$ punctures on the boundary, let $L\subset\C^n\times\R$ be a
(semi-)admissible Legendrian submanifold.

\begin{dfn} \label{6hodi.dfn}
A {\em holomorphic disk with boundary on $L$} consists of two functions 
$u\in\sblv^{\rm loc}_2(\Delta_m,\C^n)$ and 
$h\in\sblv_{\frac32}(\pa \Delta_m,\R)$
such that 
\begin{align}\label{eqhol}
\bar\pa u(\zeta) = 0, &\text{ for $\zeta\in\inr(\Delta_m)$,}\\
(u(\zeta),h(\zeta))\in L, &\text{ for $\zeta\in\pa \Delta_m$,}\label{eqbdry}
\end{align}
and such that for every puncture $p$ on $\pa \Delta_m$ there exists a
Reeb chord $c$ of $L$ such that 
\begin{equation}\label{eqaspt}
\lim_{\zeta\to p}u(\zeta)=c^\ast.
\end{equation}
\end{dfn}   
When \eqref{eqaspt} holds we say that $(u,h)$ {\em maps the puncture
$p$ to the Reeb chord $c$}.

\begin{rmk}
Since $u\in\sblv^{\rm loc}_2(\Delta_m,\C^n)$, the restriction of $u$ to the
boundary lies in $\sblv^{\rm loc}_{\frac32}(\pa\Delta_m,\C^n)$. 
Therefore both $u$ and its restriction to the boundary are
continuous. Hence \eqref{eqbdry} and \eqref{eqaspt} make
sense.   
\end{rmk}

\begin{rmk}
If $u\in\sblv^{\rm loc}_2(\Delta_m,\C^n)$ then 
$\bar\pa u\in\sblv_1^{\rm loc}(\Delta_m, {T^\ast}^{0,1}D_m\otimes\C^n)$ and
hence the trace of $\bar\pa u$ (its restriction to the boundary $\pa
\Delta_m$) lies in $\sblv^{\rm loc}_{\frac12}(\pa
\Delta_m,{T^\ast}^{0,1}D_m\otimes\C^n)$. If 
$u$ is a holomorphic disk then $\bar\pa u=0$ and hence its trace 
$\bar\pa u|\pa \Delta_m$ is also $0.$ 
\end{rmk}

\begin{rmk}
It turns out, see Section \ref{9strongI.section}, that if $(u,f)$ is a
holomorphic 
disk with boundary on a smooth $L$, then the function $u$ is in
fact smooth up to and including the boundary and thus $f$ is also
smooth. Hence, it is possible to rephrase Definition \ref{6hodi.dfn} in
terms of smooth functions. (Also, it follows that the definition above
agrees with that given in Section~\ref{1modulisection}.)
The advantage of the present definition is
that it allows for Legendrian submanifolds of lower regularity.
(The Legendrian condition applies to submanifolds $L$ which are merely
$C^1$-smooth.) 
\end{rmk}

\subsection{Conformal structures}
We describe the space of conformal structures on $\Delta_m$ as
follows. If $m\le 3$, then the conformal structure is unique. 
Let $m>4$ and let the punctures of $\Delta_m$ be $p_1,\dots,p_m$. Then
fixing the positions of the punctures $p_1,p_2,p_3$ the conformal
structure on $\Delta_m$ is determined by the position of the remaining
$m-3$ punctures. In this way we identify the space of conformal
structures $\conf_m$ on $\Delta_m$ with an open simplex of dimension
$m-3$.   


\subsection{A family of metrics}\label{2hodi.metric}
Let $\Delta$ denote the unit disk in the complex plane.
Consider $\Delta_m$ with $m$ punctures 
$p_1,\dots,p_m$ on the boundary and conformal structure $\kappa$.
Let $d$ be the smallest distance along $\pa\Delta$ between two
punctures and take  
$$
\delta=\min\left\{\frac{d}{100},\frac{\pi}{100}\right\}.
$$
Define $D(p,\delta)$ to be a disk such that $\pa \Delta(p,\delta)$
intersects $\pa\Delta$ orthogonally at two points $a_+$ and $a_-$ of
distance $\delta$ (in $\pa\Delta$) from $p$.  

Let $L_p$ be the oriented tangent-line of $\pa\Delta$ at $p$ and let $g_p$ be
the unique M{\"o}bius transformation which fixes $p$, maps $a_+$ to
the point of distance $\delta$ from $p$ along $L_p$, maps $a_-$ to 
the point of distance $-\delta$ from $p$ along $L_p$, and such that
the image of $g_p(\Delta)$ intersects the component of
$\C-L_p$ which intersects $\Delta$.   

The function 
$h_p\colon D(p,\delta)\cap\Delta_m\to [0,\infty)\times[0,1]$ defined by
\begin{equation} \notag
h_p(\zeta)=-\frac{1}{\pi}\Bigl(            
\log\left(-i{\bar p}(g_p(\zeta)-p)\right)
-\log(\delta)
\Bigr),
\end{equation}
is a conformal equivalence. Let $g_0$ denote the Euclidean metric on $\C$. 
Then there exists a function 
$s\colon [0,\frac12]\times[0,1]\to\R$ such that
 ${h_p^{-1}}^\ast g_0=s(\zeta)g_0$ on $[0,\frac12]\times[0,1].$ 
Let $\phi\colon[0,\infty)\to[0,1]$ be a smooth
function which is $0$ in a neighborhood of $0$ and $1$ in 
a neighborhood of $\frac12$ for $\tau>\frac12.$ Let $g_p$ be the metric 
\begin{equation} \notag
g_p(\tau+it)=\Bigl(\phi(\tau)+(1-\phi(\tau))s(t+it)\Bigr)g_0,
\end{equation}
on $[0,\infty)\times[0,1]$.

Now consider $\Delta_m$ with the metric $g(\kappa)$ which agrees with 
$h^\ast_{p_j} g_{p_j}$ on $h_{p_j}^{-1}([\frac{1}{2\pi},\infty)\times[0,1])$ 
for each puncture $p_j$, and with $g_0$ on $\Delta_m-(D(p_1,\delta)\cup\dots,
D(p_{m},\delta))$. Then $(\Delta_m,g(\kappa))$ is conformally equivalent to 
$(\Delta_m,g_0)$. 

We denote by $D_m(\kappa)$ the disk $\Delta_m$ with the metric
$g(\kappa).$ If the 
specific $\kappa$ is unimportant or clear from context we will 
simply write $D_m$. Also $E_{p_j}\subset D_m$ will denote the
Euclidean neighborhood $[1,\infty)\times[0,1]$ of the $j^{\rm th}$
puncture $p_j$ of $D_m$. We use coordinates $\zeta=\tau+it$ on
$E_{p_j}$ and let $E_{p_j}[M]$ denote the subset of
$\tau+it\in E_{p_j}$ with $|\tau|\ge M$.

\subsection{Sobolev spaces}
Consider $D_m$ with metric $g(\kappa)$ for some $\kappa\in\conf_m$. 
Let $\hat D_m$ denote the open Riemannian manifold which is obtained
by adding an open collar to $\pa D_m$ and extending the metric in a
smooth and bounded way to $\hat g(c)$.

The Sobolev spaces $\sblv_k^{\rm loc}(\hat D_m,\C^n)$ are now defined in
the standard way as the space of $\C^n$-valued
functions (distributions) the restrictions of which to any open ball $B$
in any relatively compact coordinate chart $\approx\R^2$ lies in
the usual Sobolev space $\sblv_k(B,\C^n)$. 

Using the metric $\hat g(\kappa)$ and the finite cover 
$$
\bigcup_j \inr(\hat E_{p_j}[1])\cup (\hat D_m-\cup_j {\hat E}_{p_j}[2]),
$$
where $\hat E_{p_j}$ is the union of $E_{p_j}$ and the corresponding
part of the collar of $\hat D_m$,  we define, for each integer $k$,
the space  $\sblv_k(\hat D_m,\C^n)$ as the subspace of all  
$f\in\sblv_k^{\rm loc}(\hat D_m,\C^n)$ with 
$\|f\|_k<\infty$.   

We consider $\sblv_k(\hat D_m,\C^n)$ as a space of distributions
acting on $C^\infty_0(\hat D_m,\C^n)$. We write 
\begin{itemize}
\item
$\sblv_k(D_m,\C^n)$ for the space of restrictions to 
$\inr(D_m)\subset\hat D_m$ of elements in $\sblv_k(\hat D_m,\C^n)$, and
\item
$\dot\sblv_k(A,\C^n)$ for the set of distributions in 
$\sblv_k(\hat D_m,\C^n)$ supported in $A \subset \hat D_m$. 
\end{itemize}

Then $\dot\sblv_k(D_m,\C^n)$ is a closed subspace of 
$\sblv_k(\hat D_m,\C^n)$ and if $K_m=\hat D_m-\inr(D_m)$ then  
$$
\sblv_k(D_m,\C^n)=
\sblv_k(\hat D_m,\C^n)/\dot\sblv_k(K_m,\C^n).
$$

We endow $\sblv_k(D_m,\C^n)$ and $\dot\sblv_k(D_m,\C^n)$ with 
the quotient- and induced topology, respectively. Let 
$C^\infty_0(D_m,\C^n)$ denote the space of restrictions of
elements in $C^\infty_0(\hat D_m,\C^n)$ to $D_m$.

\begin{lma}
$C^\infty_0(D_m,\C^n)$ is dense in 
$\sblv_k(D_m,\C^n)$,
$C^\infty_0(\inr(D_m),\C^n)$ is dense in $\dot\sblv_k(D_m,\C^n)$,
and the spaces 
$\sblv_{k}(D_m,\C^n)$ and $\dot\sblv_{-k}(D_m,\C^n)$ are dual
with respect to the extension of the bilinear form
$$
\int_{D_m}\la u,v\ra\,dA
$$
where $u\in C^\infty_0(D_m,\C^n), 
v\in C^\infty_0(\inr(D_m),\C^n)$
and $\la\, ,\,\ra$ denotes the standard Riemannian inner product on
$\C^n\approx\R^{2n}$. 
\end{lma} 
This is essentially Theorem B.2.1 p. 479 in \cite{ho}.

We will also use weighted Sobolev spaces: for $a\in\R$, let
$e^j_a\colon D_m\to\R$ be a smooth function such that
$e^j_a(\tau+it)=e^{a\tau}$ for $\tau+it\in E_{p_j}[3]$ 
and $e_a(\zeta)=1$ for $\zeta\in D_m-E_{p_j}[2]$. 
For $\mu=(\mu_1,\dots,\mu_m)\in\R^{m}$, let 
$\e_{\mu}\colon D_m\to\R\otimes\id\subset\GL(\C^n)$ be 
$$
\e_\mu(\zeta) = \Pi_{j=1}^m e^j_{\mu_j}(\zeta) \id,
$$
Note that $\e_\mu(\zeta)$ preserves Lagrangian subspaces.
We can now define $\sblv_{k,\mu}(D_m,\C^n)=\{u\in \sblv_k^{\rm
loc}(D_m,\C^n) : \e_{\mu} u\in \sblv_k(D_m,\C^n)\},$ with norm 
$\|u\|_{k,\mu}=\|\e_\mu u\|_k$.

\subsection{Asymptotics}\label{2asymp}
Let $\Lambda_0$ and $\Lambda_1$ be (ordered) Lagrangian subspaces of
$\C^n$. Define the {\em complex angle}
$\theta(\Lambda_0,\Lambda_1)\in[0,\pi)^n$ inductively as follows:

If $\dim(\Lambda_0\cap\Lambda_1)=r\ge 0$ let
$\theta_1=\dots=\theta_r=0$ and let $\C^{n-r}$ denote the Hermitian
complement of $\C\otimes \Lambda_0\cap\Lambda_1$ and let 
$\Lambda_i'=\Lambda_i\cap\C^{n-r}$ for $i=0,1$. 
If $\dim(\Lambda_0\cap\Lambda_1)=0$ then let
$\Lambda_i'=\Lambda_i$, $i=0,1$ and let $r=0$. Then $\Lambda'_0$
and $\Lambda'_1$ are Lagrangian subspaces. Let $\alpha$ be
smallest angle such that
$\dim(e^{i\alpha}\Lambda_0\cap\Lambda_1)=r'>0$. Let 
$\theta_{r+1}=\dots=\theta_{r+r'}=\alpha$. Now repeat the construction
until $\theta_n$ has been defined.    
Note that $\theta(A\Lambda_0,A\Lambda_1)=\theta(\Lambda_0,\Lambda_1)$
for every $A\in\U(n)$ since multiplication with $e^{i\alpha}$ commutes
with everything in $\U(n)$. 

\begin{prp}\label{SalRob}
Let $(u,h)$ be a holomorphic disk with boundary on a
(semi-)admissible Legendrian submanifold $L$. Let $p$ be
a puncture on $D_m$ such that $p$ maps to the Reeb chord $c.$ For  
$M>0$ sufficiently large the following is true:  

If $\Pi_\C(L)$ self-intersects transversely at $c^\ast$ then 
\begin{equation}\label{transv} 
|u(\tau+it)|=\Ordo(e^{-\theta\tau}),\quad \tau+it\in E_p[M],
\end{equation}
where $\theta>0$ is the smallest complex angle of $c$.

If $\Pi_\C(L)$ has a self-tangency at $c^\ast$ then either there exists
a real number $c_0$ such that
\begin{equation}\label{selftan1}
u(\tau+it)=\Bigl(\frac{\pm2}{c_0+\tau+it},0,\dots,0\Bigr)
+\Ordo(e^{-\theta\tau}) 
\quad \tau+it\in E_p[M], 
\end{equation}
or 
\begin{equation}\label{selftan2}
|u(\tau+it)|=\Ordo(e^{-\theta\tau}),\quad \tau+it\in E_p[M],
\end{equation}
where $\theta$ is the smallest non-zero complex angle of $L$ at $p$.

In particular, if the punctures $p_1,\dots, p_m$ on $D_m$
map to Reeb chords $c_1,\dots,c_m$ and if $f\colon D_m\to\C^n$ is any  
smooth function which is constantly equal to
$c^\ast_1,\dots,c_m^\ast$ in neighborhoods of
$p_1,\dots,p_m$, then $u-f\in \sblv_2(D_m,\C^n)$.
\end{prp}

\begin{pf}
Equation \eqref{transv} is a consequence of
Theorem B in \cite{Robbin-Salamon??}. To prove the corresponding statement for a
self-tangency double point we may assume that the self-tangency point
is $0\in\C^n$ and that around $0,$ $\Pi_\C(L)$ agrees with the local
model in Definition~\ref{5admissile_param.dfn}. Elementary complex analysis (see 
Lemma~\ref{lmaralocsol} below) shows that for a 
standard self tangency the first component $u_1$ of a holomorphic disk
is given by  
\begin{equation}
u_1(\zeta)=\frac{\pm2}{\zeta-c_0+\sum_{n\in\Z} c_n\exp(n\pi\zeta)},
\end{equation}
where $c_j$ are real constants, in $E_{p}[M]$. Applying
\cite{Robbin-Salamon??} again to the remaining 
components $u'$ of $u$ gives the claim.
The last statement follows immediately from the asymptotics at
punctures.  
\end{pf}

\section{Functional analytic setup}\label{2fas}
As explained in Section~\ref{1contacthomology}, contact homology is built using
moduli-spaces of holomorphic disks. In this section we construct
Banach manifolds of maps of punctured disks into $\C^n$ which satisfy
certain boundary conditions. In this setting, moduli-spaces will
appear as the zero-sets of bundle maps.  

In Section~\ref{fcanA} we define our Banach manifolds as
submanifolds in a natural bundle of Banach spaces. 
To find atlases for our Banach manifolds we proceed in the standard
way: construct an ``exponential map'' from the proposed tangent space 
and show it is a diffeomorphism near the origin. To do this,
in Section~\ref{fcanB}, 
we turn our attention to a special metric on the
tangent bundle of the Legendrian submanifold. 
From this we construct a family of metrics on $\C^n$ in Section~\ref{fcanC} and
use it to define a preliminary version of the ``exponential map'' for the 
Banach manifold. Section~\ref{fcanE} contains some technical results  needed
to deal with families of Legendrian submanifolds.
In Section \ref{fcanF} we show how to construct the
atlas. Section \ref{fcanG} discusses how to invoke variations of
the conformal structure of the source space into the present setup. In
Section \ref{fcanH} we linearize the bundle map, the zero set of which
is the moduli-space. Section~\ref{fcanI} discusses some issues involving
the semi-admissible case.

\subsection{Bundles of affine Banach spaces}\label{fcanA}
Let $L_\lambda\subset\C^n\times\R$, $\lambda\in\Lambda$, where
$\Lambda$ is an open subset of a Banach space, denote
a smooth family of chord generic admissible Legendrian
submanifolds. That is, $\Lambda$ is smoothly mapped into
the space of admissible Legendrian embeddings of $L$
endowed with the $C^\infty$-topology. 

We also study the semi-admissible case. To this end we also let
$L_\lambda$, $\lambda\in\Lambda,$ be a smooth family of  
semi-admissible Legendrian submanifolds. 
For simplicity, and
since it will suffice for our applications, we assume that in this case,
the self tangency point of $\Pi_\C(L_\lambda)$ remain fixed as
$\lambda$ varies and that in a neighborhood of this point the product
structure of $\Pi_\C(L_\lambda)$ is preserved and the first components
$\gamma_1$ and $\gamma_2,$ shown in Figure~\ref{1fig:TypeB} remain
fixed as $\lambda$ varies.

Let ${\bf a}(\lambda)=(a_1(\lambda),\dots,a_m(\lambda))$,
$\lambda\in\Lambda$  
be an ordered collection of Reeb chords of $L_\lambda$
depending continuously on $\lambda.$ 
Consider $D_m$ with  punctures $p_1,\dots,p_{m}$, and a
conformal structure $\kappa\in\conf_m.$ 

Fix families, smoothly depending on
$(\lambda,\kappa)\in\Lambda\times\conf_{m}$, of smooth reference functions  
\begin{equation} \notag
u_{\rm ref}[{\bf a}(\lambda),\kappa]\colon D_{m}\to\C^n
\end{equation}
such that $u_{\rm ref}[{\bf a}(\lambda),\kappa]$ is 
constantly equal to $a_k^\ast$ in $E_{p_k}$, and  
\begin{equation} \notag
h_{\rm ref}[{\bf a}(\lambda),\kappa]\colon\pa D_{m}\to\R
\end{equation}
such that $h_{\rm ref}[{\bf a}(\lambda),\kappa]$ is constantly equal to 
$a_1^-(\lambda)$ and $a_1^+(\lambda)$ on $[1,\infty)\subset E_{p_1}$,
and $[1,\infty)+i\subset E_{p_1}$, respectively, and, for 
$k\ge 2$, constantly equal to $a_k^+(\lambda)$ and $a_k^-(\lambda)$ on   
$[1,\infty)\subset E_{p_k}$,
and $[1,\infty)+i\subset E_{p_k}$, respectively.

Let $\epsilon=(\epsilon_1,\dots,\epsilon_m)\in[0,\infty)^m$.
For $u\colon D_{m}\to\C^n$ and $h\colon\pa D_{m}\to\R$ consider the
conditions 
\begin{align}\label{ucnd}
u-u_{\rm ref}[{\bf a}(\lambda),\kappa]
&\in\sblv_{2,\epsilon}(D_{m},\C^n),\\ \label{hcnd} 
h-h_{\rm ref}[{\bf a}(\lambda),\kappa]
&\in\sblv_{\frac32,\epsilon}(\pa D_{m+1},\R).  
\end{align}
(Note that the $\kappa$-dependence of the right hand sides in
\eqref{ucnd} and \eqref{hcnd} has been dropped from the notation.)
Define the affine Banach space 
\begin{equation} \notag
\func_{2,\epsilon}({\bf a}(\lambda),\kappa)=\Bigl\{
(u,h)\colon D_m\to\C^n\times\R\,\, \colon
u \text{ satisfies \eqref{ucnd}, } h \text{ satisfies \eqref{hcnd}} 
\Bigr\},
\end{equation}
endowed with the norm which is the sum of the norms of the components.
Let   
\begin{equation} \notag
\func_{2,\epsilon,\Lambda}({\bf a},\kappa)=
\bigcup_{\lambda\in\Lambda}
\func_{2,\epsilon}({\bf a}(\lambda),\kappa) 
\end{equation}
be the metric space with distance function
\begin{align}\notag
d((v,f,\lambda),(w.g,\mu))=&
\|(v-u_{\rm ref}[{\bf a}(\lambda),\kappa])-
(w-u_{\rm ref}[{\bf a}(\mu),\kappa])\|_{2,\epsilon}\\\notag
+&\|(f-h_{\rm ref}[{\bf a}(\lambda),\kappa])-
(g-h_{\rm ref}[{\bf a}(\mu),\kappa])\|_{\frac{3}{2},\epsilon}\\
+&|\lambda-\mu|.
\end{align}

We give $\func_{2,\epsilon,\Lambda}({\bf a},\kappa)$ the structure of
a Banach manifold by producing an atlas as follows. Let
$(w,f,\lambda)\in\func_{2,\epsilon,\Lambda}({\bf a},\kappa)$. Let 
$(w_\mu,f_\mu,\mu)$ be any family such that
$(w_\lambda,f_\lambda,\lambda)=(w,f,\lambda)$ and such that 
\begin{equation} \notag
\mu\mapsto
(w_\mu-u_{\rm ref}[{\mathbf a}(\mu),\kappa],
f_\mu-h_{\rm ref}[{\mathbf a}(\mu),\kappa])
\end{equation}
is a smooth map into 
$\sblv_{2,\epsilon}(D_m,\C^n)\times\sblv_{\frac32,\epsilon}(\pa D_m,\R)$.
Then a chart is given by 
\begin{align}\notag
&\sblv_{2,\epsilon}(D_m,\C^n)\times\sblv_{\frac32,\epsilon}(\pa D_m,\R)
\times\Lambda\to
\func_{2,\epsilon,\Lambda}({\mathbf a},\kappa);\\ 
&(g,r,\mu)\mapsto(w_\mu+g,f_\mu+r,\mu). 
\end{align}

If $(u,h,\lambda)\in\func_{2,\epsilon,\Lambda}({\bf
a}(\lambda),\kappa)$ then   
$\bar\pa u \in\sblv_{1,\epsilon}(D_m,{T^\ast}^{0,1}D^m\otimes\C^n)$ and its
trace $\bar\pa u|\pa D_m$  lies in
$\sblv_{\frac12}(D_m,{T^\ast}^{0,1}D_m\otimes\C^n))$.

\begin{dfn}
Let $\cand_{2,\epsilon,\Lambda}({\bf
a},\kappa)\subset\func_{2,\epsilon,\Lambda}({\bf a},\kappa)$  
denote the  
subset of elements $(u,h,\lambda)$ which satisfy
\begin{align}\label{W1}
&(u,h)(\zeta)\in L_\lambda\text{ for all }\zeta\in\pa D_m,\\ \label{W2}
&\int_{\pa D_m}\la\bar\pa u,v\ra\,ds=0,\text{ for every }
v\in C^\infty_0(\pa D_m,{T^{\ast}}^{0,1}D_m\otimes\C^n),
\end{align}
where $\la\,,\ra$ denotes the inner product on
${T^\ast}^{0,1}\otimes\C^n$ induced from the standard (Riemannian)
inner product on $\C^n$.  
\end{dfn}
\begin{lma}
$\cand_{2,\epsilon,\Lambda}({\bf a},\kappa)$ is a closed subset.
\end{lma}
\begin{pf}
If $(u_k,h_k,\lambda_k)$ is a sequence in 
$\cand_{2,\epsilon,\Lambda}({\bf a},\kappa)$
which converges  
in $\func_{2,\epsilon,\Lambda}({\bf a},\kappa)$ then $\lambda_k\to\lambda$
and the sequence $(u_k|\pa D_m,h_k)$ converges in $\sup$-norm. 
Hence \eqref{W1} is a closed condition. 
Also, $\bar\pa$ is continuous as is the trace map. It
follows that \eqref{W2} is a closed condition as well. 
\end{pf}

\subsection{The normal bundle of a Lagrangian immersion with a special
metric}\label{fcanB}   
Let $L\subset\C^n\times\R$ be a an instant of a chord generic
$1$-parameter family of Legendrian submanifolds. Then $\Pi_\C\colon
L\to\C^n$ is a Lagrangian 
immersion and the normal bundle of $\Pi_\C$ is isomorphic to the
tangent bundle $TL$ of $L$. On the restriction $T_L(TL)$ of the
tangent bundle $T(TL)$ of $TL$ to the zero-section $L$ there is a
natural endomorphism $J\colon T_L(TL)\to T_L(TL)$ such that $J^2=-1$.
It is defined as follows.
If $p\in L$ then $T_{(p,0)}(TL)$ is a direct sum of the space
of horizontal vectors tangent to $L$ at $p$ and the space of vertical
vectors tangent to the fiber of $\pi\colon TL\to L$ at $p$.
If $v\in T_L(TL)$ is tangent to $L$ at $p\in L$ then 
$Jv$ is the vector $v$ viewed as a tangent vector to the fiber $T_p L$ of
$\pi\colon TL\to L$ at $(p,0)$, and if $w$ is a vector tangent to the
fiber of $\pi$ at $(p,0)$ then $Jw=-w$, where $-w$ is viewed as a
tangent vector in $T_pL$. This defines $J$ on the two direct summands.
Extend it linearly.

The immersion $\Pi_\C\colon L\to\C^n$ extends to an immersion $P$ of a
neighborhood of the zero-section in $TL$ and $P$ can be
chosen so that along $L$, $i\circ dP=dP\circ J$.

From a Riemannian metric  $g$ on $L$, we construct a metric $\hat g$
on a neighborhood of the zero section in $TL$ in the following way.  
Let $v\in TL$ with $\pi(v)=p$. Let $X$ be a tangent
vector of $TL$ at $v$. The Levi-Civita connection of $g$ 
gives the decomposition $X=X^H+X^V$, where $X^V$ 
is a  vertical vector, tangent to the fiber, and $X^H$ lies
in the horizontal subspace at $v$ determined by the connection. 
Thus $X^V$ is a vector in $T_p L$ with its endpoint at $v$. It can be
translated linearly to the origin $0\in T_p L$. We use the same
symbol $X^V$ to denote this vector translated to $0\in T_p L$. 
Write $\pi X\in T_p L$ for the
image of $X$ under the differential of the projection $\pi$ and let
$R$ denote the curvature tensor of $g$.

Define the field of quadratic forms $\hat g$ on $TL$ as
\begin{equation}\label{TLmetric}
\hat g(v)(X,Y) = g(p)(\pi X,\pi Y)+g(p)(X^V,Y^V)+g(p)(R(\pi X,v)\pi Y,v),
\end{equation}
where $v\in TL$, $\pi(v)=p$, and $X,Y\in T_v(TL)$.

\begin{prp}\label{prpgTL}
There exists $\rho>0$ such that $\hat g$ is a Riemannian metric on 
\begin{equation} \notag
\{v\in TL\colon g(v,v)<\rho\}.
\end{equation}
In this metric, the zero section $L$ is totally geodesic and the
geodesics in $L$ are exactly those in the metric $g$. Moreover, if
$\gamma$ is a geodesic in $L$ and $X$ is a vector field in $T(TL)$
along $\gamma$ then $X$ satisfies the Jacobi equation if and only if
$JX$ does.
\end{prp}

\begin{pf}
Since $g(R(\pi X,v)\pi Y,v)=g(R(\pi Y,v)\pi X,v)$, $\hat g$ is
symmetric. When restricted to the $0$-section $\hat g$ is
non-degenerate. The first statement 
follows from the compactness of $L$. 

In Lemmas \ref{Chsthatg} and \ref{Jachatg} below we show $L$ is
totally geodesic and the statement about Jacobi-fields, respectively.  
\end{pf}

Let $x=(x_1,\dots,x_n)$ be local coordinates around $p\in L$ and let
$(x,\xi)\in\R^{2n}$ be the corresponding coordinates on $TM$, where 
$\xi=\xi_s\pa_s$ (here, and in the rest of this section, we use the
Einstein summation convention, repeated indices are summed over) where 
$\pa_j$ is the tangent vector of $TL$ in the $x_j$-direction.
We write $\pa_{j^\ast}$ for the tangent vector of $TL$ in the
$\xi_j$-direction. Let $\nabla$, $\hat\nabla$ denote the Levi-Civita
connections of $g$ and $\hat g$, respectively. Let 
Roman and Greek indices run over the sets $\{1,\dots,n\}$ and
$\{1,1^\ast,2,2^\ast,\dots,n,n^\ast\}$, respectively
and recall the following standard notation: 
\begin{align*} \notag
&g_{ij}=g(\pa_i,\pa_j),\quad \hat g_{\alpha\beta}=\hat
g(\pa_\alpha,\pa_\beta),\\ \notag
&\nabla_{\pa_i}\pa_j=\Gamma_{ij}^k\pa_k,\quad
\hat\nabla_{\pa_\alpha}\pa_\beta=\hat\Gamma_{
\alpha\beta}^\gamma\pa_{\gamma}\\ \notag
& R(\pa_i,\pa_j)\pa_k=R_{ijk}^l\pa_l,\quad
g(R(\pa_i,\pa_j))\pa_k,\pa_r)=R_{ijkr}.
 \end{align*}

\begin{lma}\label{comphatg}
The components of the metric $\hat g$ satisfies
\begin{align}\label{comphatgi}
\hat g_{ij}(x,\xi) &=
g_{ij}(x)+\xi_s\xi_t\Bigl(g_{kr}(x)\Gamma^k_{is}(x)\Gamma^r_{jt}(x)
+R_{isjt}(x)\Bigr)
,\\\label{comphatgii}
\hat g_{i^\ast j^\ast}(x,\xi) &= g_{ij}(x),\\\label{comphatgiii}
\hat g_{ij^\ast}(x,\xi) &= \xi_sg_{jk}(x)\Gamma_{is}^k(x).
\end{align}
\end{lma}

\begin{pf}
Since $\pa_{j^\ast}$ is vertical, \eqref{comphatgii} holds. Note that
the horizontal space at $(x,\xi)$ 
is spanned by the velocity vectors of the curves obtained by
parallel translating $\xi$ along the coordinate directions through
$x$. Let  $V(t)$ be a parallel vector field through $x$ in the
$\pa_j$-direction with $V(0)=\xi$ and $\dot V(0)=a_k\pa_k$. Then
\begin{equation} \notag
V(t)=(\xi_k+ta_k+\Ordo(t^2))\pa_k
\end{equation}
and applying $\nabla_{\pa_j}$ to $V(t)$ we get
\begin{equation} \notag
0=\nabla_{\pa_j} V(t)=\xi_s\nabla_{\pa_j}\pa_s+a_k\pa_k+\Ordo(t).
\end{equation}
Taking the limit as $t\to 0$ we find
$a_k\pa_k=-\xi_s\Gamma_{js}^k(x)\pa_k$.
Hence the horizontal space at $(x,\xi)$ is
spanned by the vectors $\pa_j-\xi_s\Gamma_{js}^k\pa_{k^\ast}$,
$j=1,\dots,n$ and therefore,
\begin{equation} \notag
\pa_j^V = \xi_s\Gamma_{js}^k(x)\pa_{k^\ast}.
\end{equation}
Straightforward calculation gives \eqref{comphatgi} and
\eqref{comphatgiii}. 
\end{pf}

\begin{lma}\label{Chsthatg}
The Christofel symbols of the metric $\hat g$ at $(x,0)$ satisfies
\begin{align}\label{Chsthatgi}
&\hat\Gamma_{ij}^k(x,0)=\hat\Gamma_{ij^\ast}^{k^\ast}(x,0)
=\Gamma_{ij}^k(x),\\ \label{Chsthatgii}
&\hat\Gamma_{ij}^{k^\ast}(x,0)=\hat\Gamma_{ij^\ast}^{k}(x,0)=
\Gamma_{i^\ast j^\ast}^{k^\ast}(x,0)=0.
\end{align}
Hence if $\gamma$ is a geodesic in $(L,g)$ then it is also a geodesic
in $(TL,\hat g)$.
\end{lma}

\begin{pf}
The equations
\begin{equation} \notag
\hat\Gamma_{\alpha\beta}^{\gamma}=
\frac12 \hat g^{\gamma\delta}
(\hat g_{\alpha\delta,\,\beta} +\hat g_{\beta\delta,\,\alpha}
-\hat g_{\alpha\beta,\,\delta}),
\end{equation}
where $\hat g^{\alpha\beta}$ denotes the components of the inverse
matrix of $\hat g$ and Lemma \ref{comphatg} together imply
\eqref{Chsthatgi} and \eqref{Chsthatgii}. 

Let $x(t)$ be a geodesic in 
$(L,g)$. Then $(x,x^\ast)=(x(t),0)$ satisfies
\begin{align} \notag
&\ddot x_k + \hat\Gamma_{ij}^k\dot x_i\dot x_j
+\hat \Gamma_{i^\ast j}^k\dot x_{i^\ast}\dot x_j
+\hat \Gamma_{i j^\ast}^k\dot x_i \dot x_{j^\ast}
+\hat \Gamma_{i^\ast j^\ast}^k\dot x_{i^\ast}\dot x_{j^\ast}
=\ddot x_k + \Gamma_{ij}^k\dot x_i\dot x_j=0,\\ \notag
&\ddot x_{k^\ast} + \hat\Gamma_{ij}^{k^\ast}\dot x_i\dot x_j
+\hat \Gamma_{i^\ast j}^{k^\ast}\dot x_{i^\ast}\dot x_j
+\hat \Gamma_{i j^\ast}^{k^\ast}\dot x_i \dot x_{j^\ast}
+\hat \Gamma_{i^\ast j^\ast}^{k^\ast}\dot x_{i^\ast}\dot x_{j^\ast}=0.
\end{align}
This proves the second statement.
\end{pf}

\begin{lma}\label{Jachatg}
If $\gamma$ is a geodesic in $(TL,\hat g)$ which lies in $L$ then 
$X$ is a Jacobi-field along $\gamma$ if and only if $JX$ is.
\end{lma}

\begin{pf} 
We establish the following two properties of the metric $\hat g$ and
the endomorphism $J$. If $\gamma$ is a curve
in $L$ with tangent vector $T$ and $X$ is 
any vector field in $T(TL)$ along $\gamma$ then
\begin{equation}\label{x''part}
\hat\nabla_T JX=J\hat\nabla_T X.
\end{equation}
If $X$, $Y$, and $Z$ are tangent vectors to $TL$ at $(p,0)\in L$
such that $Y$ and $Z$ are horizontal (i.e. tangent to $L$) and if
$\hat R$ denotes the curvature tensor of $\hat g$ at $(p,0)$ then
\begin{equation}\label{Jcomm}
\hat R(JX,Y)Z=J\hat R(X,Y)Z.
\end{equation} 

For \eqref{x''part}, use local coordinates and write, for $\gamma(t)=x(t)$,
$T(x)=a_k(x)\pa_k$,
$X(x)=b_j(x)\pa_j+b_{j^\ast}(x)\pa_{j^\ast}$. By Lemma \ref{Chsthatg}, 
\begin{align}\notag
\hat\nabla_T JX &=
a_k\hat\nabla_{\pa_k}(-b_{j^\ast}\pa_j+b_{j}\pa_{j^\ast})\\ \notag
& = a_k\Bigl[-(\pa_kb_{j^\ast})\pa_j+(\pa_kb_j)\pa_{j^\ast}\\ \notag
&\quad\quad\quad  -b_{j^\ast}(\hat\Gamma_{kj}^r\pa_r
+\hat\Gamma_{kj}^{r^\ast }\pa_{r^\ast})
+b_{j}(\hat\Gamma_{kj^\ast}^r\pa_r
+\hat\Gamma_{kj^\ast}^{r^\ast }\pa_{r^\ast})\Bigr]\\    \notag
&=
a_k\left[-(\pa_kb_{j^\ast})\pa_j+(\pa_kb_j)\pa_{j^\ast}
-b_{j^\ast}\hat\Gamma_{kj}^r\pa_r
+b_{j}\hat\Gamma_{kj^\ast}^{r^\ast }\pa_{r^\ast}\right]\\  \notag
&=Ja_k\left[(\pa_kb_{j^\ast})\pa_{j^\ast}+(\pa_kb_j)\pa_j
+b_{j^\ast}\hat\Gamma_{kj^\ast}^{r^\ast}\pa_{r^\ast}
+b_{j}\hat\Gamma_{kj}^{r}\pa_{r}\right]
=J\hat\nabla_T X.
\end{align}

For \eqref{Jcomm}, introduce normal coordinates $x$ around $p$. Then
\begin{equation}\label{normal}
g_{ij}(0)=\delta_{ij},\quad \Gamma_{ij}^k(0)=0
\end{equation} 
for all $i,j,k$, and hence Lemma \ref{comphatg} implies, 
\begin{equation} \notag
\hat g_{ij}(0,\xi)=\delta_{ij}+\Ordo(\xi^2),\quad
\hat g_{i^\ast j}(0,\xi)=0,\quad
\hat g_{i^\ast j^\ast}(0,\xi)=\delta_{ij}.
\end{equation}  
Therefore,
\begin{equation}\label{invcomp}
\hat g^{ij}(0,\xi)=\delta^{ij}+\Ordo(\xi^2),\quad
\hat g^{i^\ast j}(0,\xi)=0,\quad
\hat g^{i^\ast j^\ast}(0,\xi)=\delta^{ij}.
\end{equation}

We show that, in these normal coordinates,  
\begin{equation}\label{normalJcomm}
\hat R(\pa_{i^\ast},\pa_j)\pa_k=J \hat R(\pa_i,\pa_j)\pa_k
\end{equation} 
at $(0,0)$. Since $\hat R$ is a tensor field, \eqref{normalJcomm}
implies \eqref{Jcomm}.  

Lemma \ref{Chsthatg} implies that all Christofel symbols of $\hat g$
vanishes at $(x,\xi)=(0,0)$ and also that
$\pa_i\Gamma_{jk}^{r^\ast}(x,0)=0$ all  
$i,j,k,r^\ast$. Hence,
\begin{align}\notag
\hat R(\pa_i,\pa_j)\pa_k &=
\hat\nabla_{\pa_i}\hat\nabla_{\pa_j}\pa_k-
\hat\nabla_{\pa_j}\hat\nabla_{\pa_i}\pa_k\\ \notag
&=
(\pa_i\hat\Gamma_{jk}^r)\pa_r+(\pa_i\hat\Gamma_{jk}^{r^\ast})\pa_{r^\ast}
-(\pa_j\hat\Gamma_{ik}^r)\pa_r-(\pa_j\hat\Gamma_{ik}^{r^\ast})\pa_{r^\ast}\\ \notag
&=(\pa_i\Gamma_{jk}^r)\pa_r-(\pa_j\Gamma_{ik}^r)\pa_r,
\end{align}
and thus
\begin{equation}\label{JRcalc}
J\hat R(\pa_i,\pa_j)\pa_k=
(\pa_i\Gamma_{jk}^r)\pa_{r^\ast}-(\pa_j\Gamma_{ik}^r)\pa_{r^\ast}.
\end{equation}

We compute the left hand side of \eqref{normalJcomm}:
\begin{align}\notag
\hat R(\pa_{i^\ast},\pa_j)\pa_k&=
\hat\nabla_{\pa_{i^\ast}}\hat\nabla_{\pa_j}\pa_k-
\hat\nabla_{\pa_j}\hat\nabla_{\pa_{i^\ast}}\pa_k\\ \label{way} 
&=
\hat\nabla_{\pa_{i^\ast}}\left(\hat\Gamma^{r}_{jk}\pa_{r}+
\hat\Gamma^{r^\ast}_{jk}\pa_{r^\ast}\right)
-\hat\nabla_{\pa_j}\left(\hat\Gamma^{r}_{i^\ast k}\pa_{r}+
\hat\Gamma^{r^\ast}_{i^\ast k}\pa_{r^\ast}
\right).
\end{align}
Lemma \ref{Chsthatg} gives
$\pa_{j}\hat\Gamma^{r}_{i^\ast k}=0$, and Lemma \ref{comphatg} in
combination with \eqref{invcomp} give
$\pa_{i^\ast}\hat\Gamma^{r}_{jk}=0$. Hence, 
\begin{equation}\label{oway}
\hat R(\pa_{i^\ast},\pa_j)\pa_k=
(\pa_{i^\ast}\hat\Gamma^{r^\ast}_{jk})\pa_{r^\ast}
-(\pa_j\hat\Gamma_{ik^\ast}^{r^\ast})\pa_{r^\ast}
=
(\pa_{i^\ast}\hat\Gamma^{r^\ast}_{jk})\pa_{r^\ast}
-(\pa_j\Gamma_{ik}^{r})\pa_{r^\ast}.
\end{equation}
It thus remains to compute $\pa_{i^\ast}\hat\Gamma^{r^\ast}_{jk}$.
\begin{align}\notag
\pa_{i^\ast}\hat\Gamma^{r^\ast}_{jk} &=
\frac12\pa_{i^\ast}
\left(\hat g^{r^\ast l^\ast}
(\hat g_{jl^\ast,k}+ \hat g_{kl^\ast,j}-\hat g_{jk,l^\ast}) 
+ \hat g^{r^\ast l}
(\hat g_{jl,k}+\hat g_{kl,j}-\hat g_{jk,l})\right)\\ \notag
&=\frac12
\hat g^{rl}\pa_{i^\ast}
(\hat g_{jl^\ast,k}+\hat g_{kl^\ast,j}-\hat g_{jk,l^\ast})
\quad\text{[by \eqref{invcomp}]}\\ \notag
&=\frac12 
g^{rl}\Bigl((\pa_k\Gamma_{ji}^m)g_{ml}+
(\pa_j\Gamma_{ki}^m)g_{ml}
-(R_{jikl}+R_{jlki})
\Bigr)\quad\text{[Lemma \ref{comphatg}, \eqref{normal}]}\\ 
&=\frac12(\pa_k\Gamma_{ji}^r+\pa_j\Gamma_{ki}^r
-(R_{jikr}+R_{jrki}))\quad\text{[\eqref{normal}]}.
\end{align}
But
\begin{equation} \notag
R_{jikr}=g(\nabla_{\pa_j}\nabla_{\pa_i}\pa_k
-\nabla_{\pa_i}\nabla_{\pa_j}\pa_k,\pa_r)
=\pa_j\Gamma^r_{ik}-\pa_i\Gamma^r_{jk}
=\pa_j\Gamma^r_{ki}-\pa_i\Gamma^r_{jk},
\end{equation}
and 
\begin{equation} \notag
R_{jrki}=R_{kijr}=\pa_k\Gamma^r_{ij}-\pa_i\Gamma^r_{kj}=
\pa_k\Gamma^r_{ji}-\pa_i\Gamma^r_{jk}.
\end{equation}
Hence 
\begin{equation} \notag
\pa_{i^\ast}\hat\Gamma^{r^\ast}_{jk}=\pa_i\Gamma^r_{jk},
\end{equation}
which together with \eqref{way} and \eqref{oway} imply
\eqref{normalJcomm}.
 
Consider a geodesic of $(TL,\hat g)$ in $L$ with tangent
vector $T$. By \eqref{x''part} and \eqref{Jcomm},
\begin{equation}
\hat\nabla_T\hat\nabla_T JX + \hat R(JX,T)T=
J(\hat\nabla_T\hat\nabla_T X + \hat R(X,T)T).
\end{equation}
Thus $X$ is a Jacobi field if and only if $JX$ is.
\end{pf}

\subsection{A family of metrics on $\C^n$}\label{fcanC}
Let $L\subset\C^n\times\R$ be an instant of a chord generic
$1$-parameter family of Legendrian
submanifolds and fix a Riemannian metric $g$ on $L$. Using the  
metric $\hat g$ on $TL$ (see Section \ref{fcanB}), we construct
a $1$-parameter family of metrics $g(L,\sigma)$, $0\le\sigma\le1$, on
$\C^n$ with good properties with respect to $\Pi_\C(L)$. 

Let $c_1,\dots,c_m$ be the Reeb chords of $L$. Fix $\delta>0$ such
that all the $6\delta$-balls 
$B(c_j^\ast,6\delta)$ are disjoint and such that the intersections
$B(c_j^\ast,6\delta)\cap\Pi_\C(L)$ are homeomorphic to two $n$-disks
intersecting at a point.     

Identify the normal bundle of the immersion $\Pi_\C$ with the tangent
bundle $TL$. Consider the metric $\hat g$ on a $\rho$-neighborhood of 
the $0$-section in $TL$ ($\rho>0$ as in Proposition \ref{prpgTL}). Let
$P\colon W\to\C^n$ be an immersion of a $\rho'$-neighborhood
$N(\rho')$ of the $0$-section $\rho'\le\rho$ such that $i\circ
dP=dP\circ J$ along the $0$-section.  


Consider the $P$-push-forward of the metric $\hat g$ to the image of
$N(\rho')$ restricted to $L\setminus\bigcup_jU(c_j^-,\delta)$. Note that if $\rho'>0$
is small enough this restriction of $P$ is an embedding and
the push-forward metric is defined in a neighborhood of
$\Pi_\C(L\setminus\bigcup_j U(c_j^-,2\delta))$. Extended it to a metric $g^1$ on
all of $\C^n$, which agrees with the standard metric outside a
neighborhood of $\Pi_\C(L)$.  

Consider the $P$-push-forward of the metric $\hat g$ to the image
of the $\rho'$-neighborhood of the $0$-section 
restricted to $L\setminus\bigcup_jU(c_j^+,\delta)$. This metric is defined in
a neighborhood of $\Pi_\C(L\setminus\bigcup_j U(c_j^+,2\delta))$ and can be
extended to a metric $g^0$ on all of $\C^n$, which agrees with the
standard metric outside a neighborhood of $\Pi_\C(L)$.

Choose the metrics $g^0$ and $g^1$ so that they agree outside 
$\cup_j B(c_j^\ast,3\delta)$ and let $g^\sigma$, $0\le\sigma\le 1$ be a
smooth $1$-parameter family of metrics on $\C^n$ with the following
properties:
\begin{itemize}
\item $g^\sigma=g^0$ in a neighborhood of $\sigma=0$, 
\item $g^\sigma=g^1$ in a neighborhood of $\sigma=1$, 
\item $g^\sigma$ is constant in $\sigma$ outside 
$\cup_j B(c_j^\ast,4\delta)$.
\end{itemize}
We take  $g(L,\sigma)=g^\sigma$.

\begin{rmk}
If $L_\lambda$, $\lambda\in\Lambda$ is a smooth  family of
chord generic Legendrian submanifolds then, as is easily seen, the above
construction can be carried out in such a way that the 
family of $1$-parameter families of metrics $g(L_\lambda,\sigma)$
becomes smooth in $\lambda$.
\end{rmk}

Given a vector field $v$ along a disk $u\colon D_m\to\C^n$ with boundary
on $L$, we would like to be able to exponentiate $v$ to get a
variation of $u$ among disks with boundaries on $L$.
We will not be able to use a fixed metric $g^\sigma$ to do this. 
To solve this problem
let $\sigma\colon\C^n\times\R\to[0,1]$ be a smooth function which 
equals $0$ on  
\begin{equation} \notag
\C^n\times\R-\bigcup_j
B(c_j^\ast,5\delta)\times 
\left[c_j^+ -\frac12\action(c_j),c_j^++1 \right]
\end{equation} 
and equals $1$ on
\begin{equation} \notag
\bigcup_j B(c_j^\ast,4\delta)\times
\left[c_j^+-\frac14 \action(c_j),c_j^++\frac12 \right].
\end{equation}   

Let $\exp^g_p$ denote the exponential map of the metric $g$ at the
point $p$. 
If $p\in L_\lambda$ and $v$ is tangent to $L_\lambda$ at $p$, then 
write $x(p)=\Pi_\C(p)$ and $\xi(v)=\Pi_\C(v)$. One may now easily prove the following
lemma.
\begin{lma}\label{stayin}
Let $L_\lambda$, $\lambda\in\Lambda$ be a family of
(semi-)admissible Legendrian submanifolds.
Let $0\in\Lambda$ and let $\sigma\colon\C^n\times\R\to[0,1]$ be the
function constructed from $L_0$  as above.
There exists $\rho>0$ and a neighborhood $W\subset\Lambda$ of $0$ such
that if $p$ is any point in $L_\lambda$, $\lambda\in W$ and $v$ any  
vector tangent to $L_\lambda$ at $p$ with $|\xi(v)|<\rho$ then 
\begin{equation} \notag
\exp_{x}^{g(L_\lambda,\sigma(p))}t\xi\in\Pi_\C(L_\lambda)
\text{ for } 0\le t\le 1.
\end{equation}
\end{lma} 



\subsection{Extending families of Legendrian embeddings and their
differentials}\label{fcanE} 
In the next subsection we will need to exponentiate 
vector fields along a disk whose boundary is
in $L_0$ $(0\in\Lambda)$ to get a disk with boundary in $L_\lambda$ for
$\lambda$ near $0.$ To accomplish this 
we construct diffeomorphisms of $\C^n.$

Consider $L_\lambda\subset\C^n\times\R$, $\lambda\in\Lambda$ and let
$0\in\Lambda$.  
There exists a smooth family of Legendrian embeddings  
\begin{equation} \notag
k_\lambda\colon L_0\to \C^n\times\R,
\end{equation}
such that $k_0$ is the inclusion,  $k_\lambda(L_0)=L_\lambda$,
and $k_\lambda(c_j^\pm(0))=c_j^\pm(\lambda)$ for each $j$.

As in Section \ref{fcanC}, fix $\delta>0$ such
that all the $6\delta$-balls 
$B(c_j^\ast(0),6\delta)$ are disjoint and such that the intersections
$B(c_j^\ast(0),6\delta)\cap\Pi_\C(L_0)$ are homeomorphic to two
$n$-disks intersecting at a point. 

Let $W\subset\Lambda$ be a neighborhood of
$0$ such that $c_j^\ast(\lambda)\in B(c_j^\ast(0),\delta)$ for
$\lambda\in W$. We construct a smooth $\Lambda$-family $(\lambda\in W)$ of
$1$-parameter families of diffeomorphisms
$\psi_\lambda^\sigma\colon\C^n\to\C^n$, $0\le\sigma\le 1$, 
$\lambda\in W$. 
Note that 
\begin{align}\notag
K_\lambda^1&=\Pi_\C\circ k_\lambda \colon L_0^1=L_0\setminus \bigcup_j
U(c_j^+,3\delta)\to\C^n,\\
K^0_\lambda&=\Pi_\C\circ k_\lambda \colon L_0^0=L_0\setminus \bigcup_j
U(c_j^-,3\delta)\to\C^n 
\end{align}
are Lagrangian embeddings and that
$K_\lambda^1(c_j^\ast(0))=K_\lambda^0(c_j^\ast(0))=c_j^\ast(\lambda)$, for each
Reeb chord $c_j(0)$ of $L_0$. 

Identify tubular neighborhoods of $L_0^1$ and $L_0^0$ with their
respective tangent bundles so that $J$ along the $0$-section of the
tangent bundles corresponds to $i$ in $\C^n$ (see Section
\ref{fcanC}). Define for $(p,v)\in T L_0^\alpha\subset\C^n$, $\alpha=0,1$, 
\begin{equation}\label{hatK}
\hat K_\lambda^\alpha(p,v)=K^\alpha_\lambda(p)+idK^\alpha_\lambda(v). 
\end{equation}
Then $\hat K_\lambda^\alpha$ is a diffeomorphism on some neighborhood of
$L^\alpha_0\subset\C^n$, $\alpha=0,1$. Note that the diffeomorphisms
$\hat K_\lambda^0$ and $\hat K_\lambda^1$ agree
outside $\bigcup_j B(c_j^\ast,4\delta)$.

Extend $\hat K_\lambda^0$ and $\hat K_\lambda^1$ to 
diffeomorphisms on all of $\C^n$ in such a way that their extensions
agree outside $\bigcup_j B(c_j^\ast,4\delta)$. Call these extensions
$\psi^\alpha_\lambda$, $\alpha=0,1$. 

Let $\psi^\sigma_\lambda$, $0\le\sigma\le 1$ be a 
$\Lambda$-family of $1$-parameter families
of diffeomorphisms which are constant in $\sigma$ near $\sigma=0$ and
$\sigma=1$ and with the following properties.
First, $\psi_\lambda^\sigma$, $0\le\sigma\le 1$ connects
$\psi^0_\lambda$ to $\psi^1_\lambda$. Second, $\psi_\lambda^\sigma$ is
constant in $\sigma$ outside $\cup_jB(c_j^\ast,5\delta)$ and in 
$\cup_j (B(c_j^\ast,5\delta)\setminus B(c_j^\ast,4\delta))\cap L_0$. Third
$\psi^\sigma_\lambda(c_j^\ast(0))=c_j^\ast(\lambda)$, $0\le \sigma\le 1$.   

For future reference we let $Y^\sigma_\lambda$ denote the
$1$-parameter family of $1$-forms on $\Lambda$ with coefficients in
smooth vector fields on $\C^n$ defined by
\begin{equation}\label{defYY}
Y^\sigma_\lambda(x,\mu)=D_\lambda\psi^\sigma_\lambda(x)\cdot\mu,\quad
\lambda\in\Lambda,\mu\in T_\lambda\Lambda, x\in\C^n,\sigma\in[0,1].
\end{equation}

By \eqref{hatK}, $d\psi^\alpha_\lambda$, $\alpha=0,1$ are complex
linear maps when restricted to the restriction of the tangent bundle
of $\C^n$ to $L_0^\alpha$. Moreover, these maps fit together to a smooth
$\Lambda$-family of maps $\hat A_\lambda\colon L_0\to \GL(\C^n)$ which is
obtained as follows. Pick a smooth function $\beta$ on $L_0$ with
values in $[0,1]$ which is $0$ outside $U(c_j^+,5\delta)$ and $1$
inside $U(c_j^+,4\delta)$ define
\begin{equation} \notag
\hat A_\lambda(p)=d\psi^{\beta(p)}(d\Pi_\C(T_pL)).
\end{equation} 
Let $A_\lambda^\sigma\colon\C^n\to \GL(\C^n)$ be an $s$-parameter
family of $1$-parameter families of maps with the following properties. 
\begin{itemize}
\item
$A_\lambda^\sigma=\hat A_\lambda$ on 
$\Pi_\C(L_0)\setminus \Pi_\C(U(c_j^+,5\delta))$ 
\item 
$A_\lambda^1=\hat A_\lambda$ on $\Pi_\C(U(c_j^+,4\delta))$ 
\item
$A_\lambda^\sigma$ is constant in $\sigma$ on
$B(c_j^\ast,5\delta)\setminus B(c_j^\ast,4\delta)\cap L_0$ 
\item 
$\bar\pa A_\lambda^0=0$  
along $\Pi_\C(L_0)\setminus\Pi_\C(U(c_j^+,4\delta))$ and 
$\bar\pa A_\lambda^1=0$  
along $\Pi_\C(L_0)\setminus\Pi_\C(U(c_j^-,4\delta))$.
\item
$\|A_\lambda^\sigma-\id\|_{C^\infty}\le2\|\hat A_\lambda-\id\|_{C^\infty}$. 
\end{itemize}


\subsection{Local coordinates}\label{fcanF}
We consider first the chord generic case.
Let $L_\lambda\subset\C^n\times\R$, $\lambda\in\Lambda$ be a family of
chord generic Legendrian submanifolds. We 
construct local coordinates on $\cand_{2,\epsilon,\Lambda}({\mathbf a},\kappa)$.  

Let $\sigma\colon\C^n\times\R\to[0,1]$
be the function constructed from $L_0$, $0\in\Lambda$.
For $p\in L_\lambda$ and $v$ a tangent vector of 
$\Pi_\C(L_\lambda)$ at $q=\Pi_\C(p)$, write
\begin{equation} \notag
\exp^{g(L_\lambda,\sigma(p))}_q v=\exp^{\lambda,\sigma}_q v.
\end{equation} 
Moreover, if $\rho>0$ is as in Lemma \ref{stayin} and $|v|\le\rho$ we
write $z(p,v)$ for the $z$-coordinate of the endpoint of the unique
continuous lift of the path  
$\exp^{\lambda,\sigma(p)}_qtv$, $0\le t\le1$, to
$L\subset\C^n\times\R$.

Let $(w,f)\in\cand_{2,\epsilon,\Lambda}({\mathbf a},\kappa)$. 
Let $F\colon D_m\to\R$ be an extension of $f$ such that 
$F\in\sblv_{2,\epsilon}(D_m,\R)$ (in particular $F$ is continuous)
and such that $F$ is smooth with all  
derivatives uniformly bounded outside a small neighborhood of $\pa D_m$.
Then $w\times F \colon D_m\to\C^n \times\R$. In the case that $w$ and
$f$ are smooth we take $F$ to be smooth. Furthermore, in the case that 
$w$ and $f$ are constant close to each puncture we take $F$ to be an
affine parameterization of the corresponding Reeb-chord in a
neighborhood of each puncture where $w$ and $f$ are constant.
The purpose of this choice of $F$ is that when we exponentiate a vector field
at the disk $(w,f)$, we need $(w,F)$ to determine the metric.

For $r>0$, define
$$
\B_{2,\epsilon}((w,f),r)\subset\sblv_{2,\epsilon}(D_m,\C^n) 
$$
as the intersection of the closed subspace of
$v\in\sblv_{2,\epsilon}(D_m,\C^n)$ which satisfies 
\begin{align}
&v(\zeta)\in \Pi_\C\Bigr(T_{(f(\zeta),w(\zeta))}L\Bigl),\text{ for }
\zeta\in\pa D_m,\\
& \int_{\pa D_m}\la\bar\pa v,a\ra\, ds=0, \text{ for every }
a\in C^\infty_0(\pa D_m,\C^n)
\end{align}
and the ball $\{u\colon \|u\|_{2,\epsilon}<r\}$. 

When the parameter space $\Lambda$ is $0$-dimensional 
we can define a coordinate chart around 
$(f,w,0)\in \cand_{2,\epsilon,\Lambda}({\mathbf a},\kappa)$
($0\in\Lambda$) by 
\begin{align} \notag
&\Phi[(w,f,0)]\colon 
\B_{2,\epsilon}((w,f),r)\times \Lambda
\to\func_{2,\epsilon,\Lambda}({\mathbf a},\kappa);\\ \notag
&\Phi[(w,f,0)](v,\lambda)=(u,l,\lambda)
\end{align}
where 
\begin{align} \notag
u(\zeta)&=
\exp^{\lambda, \sigma(\zeta)}_{w(\zeta)}
\Bigl(v(\zeta)\Bigr),\\ \notag
l(\zeta)&=
z\Bigl((w(\zeta), f(\zeta)),v(\zeta)\Bigr),\quad \zeta\in\pa D_m. 
\end{align}

When $\Lambda$ is not $0$-dimensional we will need to use the maps $A_\lambda^{\sigma}$ to
move the ``vector field'' $v$ from $L_0$ to $L_\lambda.$ Moreover, to ensure our new maps are
in the appropriate space of functions we will also need to cut off the original
map $w.$ To this end 
let $(w,f,\lambda)\in\cand_{2,\epsilon,\Lambda}({\mathbf a},\kappa)$. Then
there exists $M>0$ and vector-valued functions $\xi_j$, $j=1,\dots,m$
such that 
\begin{equation} \notag
w(\tau+it)=\exp_{a_j^\ast}^{\lambda, \omega(t)}\xi_j(\tau+it),\quad
\text{for } \tau+it\in E_{p_j}[M],
\end{equation}
where $\omega\colon[0,1]\to[0,1]$ is a smooth approximation of the
identity, which is constant in neighborhoods of the endpoints of the
interval. Define $(w[M],f[M])$ as follows. Let 
\begin{equation} \notag
w[M](\zeta)=
\begin{cases}
w(\zeta), &\text{ for }\zeta\notin \cup_j E_{p_j}[M],\\
\exp_{a_j^\ast}^{\lambda, \omega(t)}(\alpha\xi_j), &\text{ for
}\zeta=\tau+it\in E_{p_j}[M],
\end{cases}
\end{equation}
where $\alpha\colon E_{p_j}\to\C$ is a smooth function which is $1$ on
$E_{p_j}\setminus E_{p_j}[M+1]$, $0$ on $E_{p_j}[2M]$, and holomorphic 
on the boundary. Let $f[M]$ be the natural lift of the boundary values 
of $w[M]$. It is clear that $(w[M],f[M])\to (w,f)$ as
$M\to\infty$. For convenience we use the notation
$(w[\infty],f[\infty])$ to denote this limit. We write $F[M]$ for the
extension of $f[M]$ to $D_m$.

Let $(w,f,0)\in\cand_{2,\epsilon,\Lambda}({\mathbf a},\kappa)$
($0\in\Lambda$). For 
large $M>0$ consider $(w[M],F[M])$. To simplify notation, write 
$\sigma[M](\zeta)=\sigma(w[M](\zeta),F[M](\zeta))$ and
$w[M]_\lambda(\zeta)=\psi^{\sigma[M](\zeta)}_\lambda(w[M](\zeta))$.
Define
\begin{align} \notag
&\Phi[(w,f,0);M]\colon 
\B_{2,\epsilon}((w[M],f[M]),r)\times \Lambda
\to\func_{2,\epsilon,\Lambda}({\mathbf a},\kappa);\\ \notag
&\Phi[(w,f,0);M](v,\lambda)=(u,l,\lambda)
\end{align}
where 
\begin{align} \notag
u(\zeta)&=
\exp^{\lambda, \sigma[M](\zeta)}_{w[M]_\lambda(\zeta)}
\Bigl(A_\lambda^{\sigma[M](\zeta)}v(\zeta)\Bigr),\\ \notag
l(\zeta)&=
z\Bigl((w[M]_\lambda(\zeta), f[M]_\lambda(\zeta)),
A_\lambda^{\sigma[M](\zeta)}v(\zeta)\Bigr),\quad \zeta\in\pa D_m. 
\end{align}

In the semi-admissible case we use the above construction close to
all Reeb chords except the chord $c_0$ at the self-tangency point. At
$c_0^\ast$ we utilize the fact that we have a local product structure
of $\Pi_\C(L_\lambda)$ which is assumed to be preserved in a rather
strong sense under $\lambda\in\Lambda$, see Section \ref{fcanA}. 
This allows us to construct the family of
metrics $g^\sigma_\lambda$ as product metrics close to
$c_0^\ast$. Once we have metrics with this property, we can apply the
cut-off procedure above to the last $(n-1)$ coordinates of an element
$(w,f,0)\in\cand_{2,\epsilon,\Lambda}$ and just keep the first
coordinate of $w$ in a neighborhood of $c_0^\ast$ as it is. We use the
same notation $(w[M],f[M])$ for the map which results from this
modified cut-off procedure from $(w,f)$ in the semi-admissible case.

\begin{prp}\label{Localcoord}
Let $\epsilon\in[0,\infty)^m$. Then there exists $r>0$, $M>0$, and a
neighborhood $W\subset\Lambda$ of $0$ 
such that the map 
\begin{equation} \notag
\Phi[(w,f,0)]\colon\B_{2,\epsilon}((w[M],f[M]),r)\times W
\to\func_{2,\epsilon,\Lambda}({\mathbf a},\kappa)
\end{equation} 
is $C^1$ and gives local 
coordinates on some open subset of
$\cand_{2,\epsilon,\Lambda}({\mathbf a},\kappa)$ 
containing $(w,f,0)$. Moreover, if $\Lambda$ is $0$-dimensional then we
may take $M=\infty$.
\end{prp}

\begin{pf}
Fix some small $r>0$. Consider the auxiliary map  
\begin{align} \notag
&\Psi\colon
\B_{2,\epsilon}((w[M],f[M]),r)\times\sblv_{\frac32,\epsilon}(\pa
D_m,\R)\times \Lambda
\to
\func_{2,\epsilon,\Lambda}({\mathbf a}),\\ \notag
&\Psi(v,r,\lambda)=
\Phi[(w[M],f[M]),0](v,\lambda)+(0,0,r)
\end{align} 
where $(u,h, \mu)+(0,0, r)=(u,h, \mu+r)$. 

We show in Lemma \ref{funcoord} that
$\Psi$ is $C^1$ with differential in a neighborhood of $(0,0,0)$ 
which maps injectively into the tangent space of the target and has
closed images. These closed images have direct complements and hence
the implicit function theorem applies and shows that the image is a
submanifold. Moreover, for $M$ large enough $(w,f,0)$ is in the image.

We finally prove in Lemma \ref{lmafinish} that
$\cand_{2,\epsilon}({\mathbf a},\kappa)$ lies
inside the image and that it corresponds exactly to $r=0$ in the given
coordinates.
\end{pf}

Lemma~\ref{funcoord} is a consequence of the following technical
lemma. 

\begin{lma}\label{Tech}
Let $\Lambda$ be an open neighborhood of $0$ in a Banach space. 
Let $(w,f,\lambda)\in\func_{2,\epsilon,\Lambda}({\mathbf a},\kappa)$,
and $v,u,q\in\B_{2,\epsilon}((w,f),r)$. Let $\zeta$ be a coordinate on
$D_m$ and let $\epsilon\in[0,\infty)^m$.   
\begin{itemize}
\item[{\rm (a)}]
Let 
\begin{equation} \notag
G\colon \C^n\times\C^n\times\C^n\times\C^n\times[0,1]\times\Lambda\to\C^n
\end{equation}
be a smooth function with all derivatives uniformly bounded and let
$\sigma\colon\C^n\times\R\to[0,1]$ be a smooth function with the same
property. If 
\begin{align}\label{a)G=0i}
&G(x,0,0,\theta,\sigma,\lambda)=0,\\ \label{a)G=0ii}
&G(x,\xi,0,\theta,\sigma,0)=0,
\end{align} 
then there exists a constant $C$ (depending
on $\|Dw\|_{1,\epsilon}$, $\|DF\|_{1,\epsilon}$ and $r$) such that 
$G(\zeta,\lambda)=G(w(\zeta),v(\zeta),u(\zeta),q(\zeta),
\sigma(F(\zeta),w(\zeta)),\lambda)$ satisfies
\begin{equation}\label{a)est}
\|G(\zeta,\lambda)\|_{2,\epsilon}\le  
C(\|u\|_{2,\epsilon}+\|v\|_{2,\epsilon}+|\lambda|).
\end{equation} 
\item[{\rm (b)}]
Let 
\begin{equation} \notag
G\colon\C^n\times\C^n\times\C^n\times[0,1]\times\Lambda\to\C^n
\end{equation}
be a smooth function with all derivatives uniformly bounded.
If 
\begin{align}\label{b)G=0i}
&G(x,0,0,\sigma,\lambda)=0,\\\label{b)G=0ii}
&G(x,\xi,0,\sigma,0)=0,\\\label{b)G=0iii}
&D_3G(x,\xi,0,\sigma,0)=0,\text{ and}\\\label{b)G=0iv}
&D_5G(x,\xi,0,\sigma,0)=0
\end{align}
then there exists a constant $C$ (depending on
$\|Dw\|_{1,\epsilon}$, $\|DF\|_{1,\epsilon}$ and $r$) such that
$G(\zeta,\lambda)=G(w(\zeta),v(\zeta),u(\zeta)
,\sigma(F(\zeta),w(\zeta)),\lambda)$ satisfies
\begin{equation} \notag
\|G(\zeta,\lambda)\|_{2,\epsilon} 
\le C(\|u\|_{2,\epsilon}^2+|\lambda|^2).
\end{equation}
\end{itemize}
\end{lma}

\begin{pf}
For simplicity, we suppress intermediate functions in the
notation, e.g., we write $\sigma(\zeta)$ for $\sigma(w(\zeta),F(\zeta))$.
Consider (a). Assume that $w,v,u,q,F$ are smooth functions. By \eqref{a)G=0i}
\begin{equation}\label{a)0th}
|G(\zeta,\lambda)|\le C(|v|+|u|),
\end{equation}
since the derivatives of $G$ are uniformly bounded. 

(For simplicity, we will use the letter $C$ to denote many different
constants in this proof. This (constant!) change of notation will not
be pointed out each time.)

Let $\hat{G}(\zeta) = G(\zeta, \lambda).$
We write $(w,v,u,q,\sigma)=(x_1,x_2,x_3,x_4,x_5)$ and use 
the Einstein summation convention. The derivative of $\hat{G}(\zeta)$ is
\begin{equation} \notag
D\hat{G}(\zeta)=D_j\hat{G}\cdot Dx_j,
\end{equation}
where $D$ without subscript refers to derivatives with respect to
$\zeta$, and $D_j\hat{G}$ refers to the derivative of $\hat{G}$ with respect to
its $j$-th argument. We use the following notation for functions
$(y_1,\dots,y_l)$,
\begin{equation} \notag
|D^{j_1}y|^{k_1}\dots|D^{j_m}y|^{k_m}=
\sum_{\alpha\in A}
\Pi_{k=1}^l|D^{j_1}y_k|^{\alpha_k^1}\dots|D^{j_m}y_k|^{\alpha_k^m}
\end{equation}
where
$A=\{\alpha\in(\Z_{\ge 0})^{lm}\colon\alpha_1^r+\dots+\alpha_l^r=k_r\}$.

Let $(w,F,q)=(y_1,y_2,y_3)$ and $(v,u)=(z_1,z_2)$ then, by \eqref{a)G=0i}
\begin{align} \notag
&|D_j\hat{G}|\le C|z|,\quad j\in\{1,4,5\}\\ \notag
&|D_j\hat{G}|\le C,\quad j \in \{2,3\}\\ \notag
& D\sigma=D_1\sigma\cdot DF+D_2\sigma\cdot Dw,\text{ hence }
|D\sigma|\le C|Dy|.
\end{align}
Then
\begin{equation}\label{a)1st}
|D\hat{G}(\zeta)|^2\le 
C\Bigl(|z|^2|Dy|^2+|z||Dz||Dy|+|Dz|^2\Bigr).
\end{equation}

The second derivative of $\hat{G}(\zeta)$ is
\begin{equation} \notag
D^2\hat{G}(\zeta)=D_iD_j\hat{G}\cdot Dx_i\cdot Dx_j + D_j\hat{G}\cdot D^2x_j
\end{equation}
By \eqref{a)G=0i},
\begin{align} \notag
|D_iD_j\hat{G}|&\le C|z|,\quad i,j\in\{1,4,5\}\\ \notag
|D_iD_j\hat{G}|&\le C,\quad  j\in\{2,3\}\\ \notag
D^2\sigma&=D_1^2\sigma\cdot DF\cdot DF+
2D_2D_1\sigma\cdot Dw\cdot DF\\ \notag
&+D^2_2\sigma Dw\cdot Dw+
D_1\sigma\cdot D^2 F+D_2\sigma\cdot D^2w,\\ \notag
\text{hence }&
|D^2\sigma|\le C(|Dy|^2+|D^2y|).
\end{align}
Thus  
\begin{align}\notag
|D^2 \hat{G}(\zeta)|^2\le
C\Bigl(&|z^2|(|Dy|^4 +|Dy|^2|D^2 y|)+|z||Dz||Dy||D^2y|\Bigr.\\ \label{a)2nd}
&\Bigl.+|Dz|^4+|Dz||D^2z||Dy|+|Dz|^2|D^2 z|+|D^2 z|^2
\Bigr).
\end{align}

Note that by \eqref{b)G=0i} and \eqref{b)G=0iii}, 
$r,$  which the constant $C$ absorbs,
controls the $q$ (or $y_3$) norms.
Moreover, the remaining $y_1$ and $y_2$ norms are also absorbed by $C.$
Thus, using \eqref{a)0th}, \eqref{a)1st}, and \eqref{a)2nd} we derive the
estimate 
\begin{equation}\label{a)pre}
\|\hat{G}(\zeta)\|_{2,\epsilon}\le
C(\|u\|_{2,\epsilon}+\|v\|_{2,\epsilon})
\end{equation}
as follows. The Sobolev-Gagliardo-Nirenberg theorem implies 
$\|Dy\|_{L^4}\le C\|Dy\|_{1,2}$ (and the corresponding statement for
$u$ and $v$). Morrey's theorem implies that $\|u\|_{2,\epsilon}$
controls the $\sup$-norm of 
$u$ (and the corresponding statement for $v$). 
These facts together with H{\"o}lder's inequality gives \eqref{a)pre}.

It is now straightforward to prove (a). Let
$\Omega=(x,\xi,\eta,\theta,\sigma)$ then
\begin{equation}\label{a)lexpan}
G(\Omega,\lambda)=G(\Omega,0)+D_6G(\Omega,0)\cdot\lambda
+R(\Omega,\lambda)\cdot\lambda\cdot\lambda.
\end{equation}
Differentiating \eqref{a)lexpan} twice with respect to $\lambda$ of and
applying \eqref{a)G=0i} we find $R(x,0,0,\theta,\sigma,\lambda)=0$.
Applying the argument above to $D_6G$ and $R$, and to $G(\Omega,0)$ but
using \eqref{a)G=0ii} and $u$ instead of \eqref{a)G=0i} and $(u,v)$,
\eqref{a)est} follows.

The proof of (b) is similar. We first use \eqref{b)G=0ii} and
\eqref{b)G=0iii} to conclude
\begin{equation}\label{b)0th}
\hat G=G(w,v,u,\sigma,0)\le C|u|^2.
\end{equation}
The derivative of $\hat G(\zeta)$ is
\begin{equation} \notag
D\hat G(\zeta)=D_j\hat G\cdot Dx_j,
\end{equation}
and with $(w,F,v)=(y_1,y_2,y_3)$
\begin{align} \notag
&|D_j\hat G|\le C|u|^2,\quad j\in\{1,2,4\}\\ \notag
&|D_3\hat G|\le C|u|,\\ \notag
& D\sigma=D_1\sigma\cdot DF+D_2\sigma\cdot Dw,\text{ hence }
|D\sigma|\le C|Dy|.
\end{align}
Thus
\begin{equation}\label{b)1st}
|D\hat G(\zeta)|^2\le 
C\Bigl(|u|^4|Dy|^2+|u|^3|Du||Dy|+|u|^2|Du|^2 \Bigr).
\end{equation} 
The second derivative of $\hat G(\zeta)$ is
\begin{equation} \notag
D^2\hat G(\zeta)=D_iD_j\hat{G}\cdot Dx_i\cdot Dx_j + D_j\hat{G}\cdot D^2x_j
\end{equation}
We have
\begin{align} \notag
&|D_iD_j\hat{G}|\le C|u|^2,\quad i,j\in\{1,2,4\}\\ \notag
&|D_iD_3\hat{G}|\le C|u|,\quad i\in\{1,2,4\},\\ \notag
&|D_3^2\hat{G}|\le C. 
\end{align}
This implies
\begin{align}\notag
|D^2\hat{G}(\zeta)|^2\le C\Bigl(&|u|^4(|Dy|^4+|Dy|^2|D^2y|+|D^2y|^2)\\\notag 
&+|u|^3(|Du||Dy|^3 +|Du||Dy||D^2y|)+|u|^2|Du|^2|Dy|^2\\\label{b)2nd}
&+|u|^2|Du|^4+|Du||Dy||D^2y|+ |D^2u||D^2y|+ |u|^2|D^2u|^2\Bigr).
\end{align} 
In the same way as above we derive from \eqref{b)0th}, \eqref{b)1st},
and \eqref{b)2nd} the estimate
\begin{equation}\label{b)pre}
\|\hat G(\zeta)\|_{2,\epsilon}\le C\|u\|_{2,\epsilon}^2.
\end{equation} 
The proof of (b) can now be completed as follows. Write
$\Omega=(x,\xi,\eta,\sigma)$ then
\begin{equation} \notag
G(\Omega,\lambda)=G(\Omega,0)+D_5G(\Omega,0)\cdot\lambda+
D_5^2G(\Omega,0)\cdot\lambda\cdot\lambda+
R(\Omega,\lambda)\cdot\lambda\cdot\lambda\cdot\lambda,
\end{equation}
and differentiation gives $R(x,0,0,\sigma,\lambda)=0$. For $G(\Omega,0)$,
we use \eqref{b)pre}. The term $D_5\hat G(\zeta)$ can be estimated
as in (a) by $C\|u\|_{2,\epsilon}$. The two remaining terms are also
estimated as in (a) by $C(\|u\|_{2,\epsilon}+\|v\|_{2,\epsilon})$.  
\end{pf}

In order to express the derivative of $\Psi$ we will use the function
$K\colon\C^n\times\C^n\times[0,1]\times\Lambda\to\C^n$ defined by 
\begin{equation}\label{2mapKdef} 
K(x,\xi,\sigma,\lambda)=
\exp_{\psi^\sigma_\lambda(x)}^{\lambda,\sigma}
A^\sigma_\lambda\xi -\psi^\sigma_\lambda(x). 
\end{equation}
We will  need to lift $K$ (at least on part of its domain) so that it maps to 
$\C^n\times \R.$ We describe this lift.

Consider $L_\lambda\subset\C^n\times\R$, $\lambda\in\Lambda$.
Let $K_\lambda\colon TL_0\to\C^n\times\R$ be an embedding
extension of $k_\lambda$ (see Section \ref{fcanE}). Consider the
immersion $P_\lambda\colon V\subset TL_0\to\C^n$ 
which extends $\Pi_\C\circ k_\lambda$, where $V$
is a neighborhood of the $0$-section in $TL_0$. 
Choose $V$ and a neighborhood $W\subset\Lambda$ of $0$, so small that
the self-intersection of $P_\lambda$ is contained inside
$\bigcup_jB(c_j^\ast(0),2\delta)$. Consider the following subset $N$ of
the product $\C^n\times[0,1]$.
\begin{align}\notag
N=& P(V)\setminus\bigcup_jB(c_j^\ast(0),3\delta)\times[0,1]
\,\,\cup\,\, \bigcup_j P(V|L_\lambda\cap
U(c_j^+,4\delta))\times[1-\epsilon,1]\\ \notag
& \cup\,\, \bigcup_j P(V|L_\lambda\cap
U(c_j^-,4\delta))\times[0,\epsilon].
\end{align}    
We define a map $\psi_\lambda\colon N\to\C^n\times\R$ in the natural
way, $\psi_\lambda(q,\sigma)=K_\lambda(p_\sigma,v_\sigma)$ where
$(p_\sigma,v_\sigma)$ is the preimage of $q$ under $P$ with $p\in
U(c_j^{\pm},4\delta)$ where the sign is determined by $\sigma$. 

Using this construction we may do the following. If
$W\subset\C^n\times\C^n\times[0,1]\times \Lambda$ and
$G\colon W\to\C^n$ is a function such that $(G,\sigma)(W)\in N$ then
we may define a lift $\tilde G\colon W\to\C^n\times\R$.    


We now use this construction to lift the function $K$ defined in \eqref{2mapKdef}.
For $x$ sufficiently close to $L_0$, $\xi$ sufficiently small and
$\sigma$ sufficiently close to $0$ or $1$ when $x$ is close to double
points of $\Pi_\C(L_0)$ the lift $\tilde K$ of $K$ can be defined. Let
$K_\R$ denote the $\R$-coordinate of $\tilde K$.

\begin{lma}\label{funcoord}
If $\dim \Lambda>0$, let $M<\infty$. If $\dim \Lambda=0$, let
$M=\infty$. The map 
\begin{equation} \notag
\Psi\colon\B_{2,\epsilon}((w[M],f[M]),r)\times
\sblv_{\frac32,\epsilon}(\pa D_m,\R)\times\Lambda\to
\func_{2,\epsilon,\Lambda}({\bf a},\kappa)
\end{equation}
is $C^1$. Its derivative at $(v,h,\mu)$ is the map
\begin{equation} \notag
(u,l,\lambda)\mapsto 
\Bigl(D_2K\cdot v+D_4K\cdot\lambda,
\la D_2K_\R\cdot v+D_4K_\R\cdot\lambda\ra+l,\lambda\Bigr),
\end{equation}
where all derivatives of $K$ and $\tilde K$ are evaluated at 
$(w[M]_\lambda,v,\sigma[M],\lambda)$ and where $\la u\ra$
denotes restriction of $u\colon D_m\to\R$ to the boundary.
\end{lma}

\begin{pf}
Using local coordinates on $\func_{2,\epsilon,\Lambda}$ as described
in Section~\ref{fcanA}, we write $\Psi=(\Psi_1,\Psi_2,\Psi_3)$
Statements concerning $\Psi_3$ are trivial.
Note that $\Psi_1 = K + \Psi_\lambda^\sigma(x).$
So to see that $\Psi_1$ is continuous we note that
$K(x,0,\sigma,\lambda)=0$ and apply Lemma \ref{Tech} (a)
to get that $K$ is Lipschitz in $v$ and $\lambda$ and hence continuous.
To see that
$\Psi_1$ is differentiable we note that if
\begin{align} \notag
G(x,\xi,\eta,\sigma,\lambda)&=
K(x,\xi+\eta,\sigma,\mu+\lambda)-K(x,\xi,\sigma,\mu)\\ \notag
&-\Bigl(D_2K(x,\xi,\sigma,\mu)\cdot\eta
+D_4K(x,\xi,\sigma,\mu)\cdot\lambda\Bigr),
\end{align}
then the conditions \eqref{b)G=0i}--\eqref{b)G=0iv} are fulfilled and
Lemma \ref{Tech} (b) implies $\Psi_1$ is differentiable and has 
differential as claimed. Finally, applying Lemma \ref{Tech} (a) to the
map
\begin{equation} \notag
G(x,\xi,\eta,\sigma,\lambda)=D_2K(x,\xi,\sigma,\mu)\cdot\eta
+D_4K(x,\xi,\sigma,\mu)\cdot\lambda
\end{equation} 
shows $\Psi_1$ is $C^1$. 

Using $\tilde K$, we can extend the $\R$-valued function
$z((w[M]_\lambda,f[M]_\lambda),A_\lambda^\sigma v)$ to a small
neighborhood of $\pa D_m$ in $D_m$. With this done the (non-trivial
part) of the derivative of $\Psi_2$ can be handled exactly as above.   
\end{pf}

Let $\zeta=x_1+ix_2$ be a complex local coordinate in $D_m$. Then, if
$u\colon D_m\to\C^n$, we may view $\bar\pa u$ as $\pa_1u+i\pa_2u$. 
As in the proof of Lemma \ref{funcoord} we use local coordinates on
$\func_{2,\epsilon, \Lambda}$ and write $\Psi=(\Psi_1,\Psi_2,\Psi_3)$.  

\begin{lma}\label{Bdryval}
Assume that $w\colon D_m\to \C^n$ and $v\colon D_m\to \C^n$ are smooth functions and
let $g$ be any metric on $\C^n$.
If $u(\zeta)=\exp_{w(\zeta)}(v(\zeta))$ then 
$\bar\pa u=X_1(1)+iX_2(1)$, where $X_j$, $j=1,2$ are the 
Jacobi-fields along the geodesic $\exp_{w(\zeta)}(tv(\zeta))$, 
$0\le t\le 1$, with $X(0)=\pa_j w(\zeta)$ and $X'(0)=\pa_j v(\zeta)$.

In particular, there exists $r>0$ such that if $(v(\zeta),\lambda)\in
\sblv_{2,\epsilon}(D_m,\C^n)\times\Lambda$, $\|v\|_{2,\epsilon}\le r$
then the restriction of 
$\bar\pa\Psi_1(v,\lambda)$ to $\pa D_m$ equals $0$ if and only if
the restriction of $\bar\pa v=0.$ 
\end{lma}

\begin{pf}
Consider
\begin{equation} \notag
\alpha(s,t)=\exp_{w(\zeta+s)}(tv(\zeta+s)),\quad  0\le t\le 1,
-\epsilon\le s\le\epsilon. 
\end{equation} 
Since for fixed $s$, $t\mapsto\alpha(s,t)$ is a geodesic we find that
$\pa_s\alpha(0,t)=X_1(t)$ is a Jacobi field along the geodesic
$t\mapsto\exp_{w(\zeta)}(tv(\zeta))$ with initial conditions
\begin{align} \notag
X_1(0)&=\pa_s\exp_{w(\zeta+s)}(0\cdot v)=\pa_1w(\zeta),\\ \notag
X'_1(0)&=\pa_t\pa_s\alpha(0,0)=\pa_s\pa_t\alpha(0,0)=\pa_1v(\zeta).
\end{align}
Moreover, 
\begin{equation} \notag
\exp_{w(\zeta+s)}(v(\zeta+s))=\alpha(s,1)
\end{equation}
and hence
\begin{equation} \notag
\pa_1\exp_{w(\zeta)}(v(\zeta))=\pa_s\alpha(0,1)=X_1(1).
\end{equation}
A similar analysis shows that
\begin{equation} \notag
\pa_2\exp_{w(\zeta)}(v(\zeta))=X_2(1).
\end{equation}
This proves the first statement.

Consider the second statement. Note that the metrics
$g(L_\lambda,\sigma(w[M](\zeta),F[M](\zeta)))$ are constant in $\zeta$
for $\zeta$ in a neighborhood of $\pa D_m$. Consider first the case
that $w[M]$ and $v$ are smooth. Then the above result together with the
Jacobi-field property of the metric $\hat g$ (see Lemma
\ref{Jachatg}), from which $g(L_\lambda,\sigma)$ is constructed
implies that for $\zeta\in\pa D_m$,
$\bar\pa(\Psi_1(v,\lambda))=X_1(1)+iX_2(1)$ equals the 
value of the Jacobi-field $X_1+JX_2=Y$ with initial condition
$Y(0)=0$, $Y'(0)=\bar\pa A^\sigma_\lambda v.$ (Note $Y(0)=0$ since 
along the boundary $\bar\pa w=0.$) Hence $\bar\pa(\Psi_1(v,\lambda))=0$ for
$\zeta\in\pa D_m$ if and only if the same is true for $v$ provided $v$
is shorter than the minimum of injectivity radii of
$g(L_\lambda,\sigma)$.      
An approximation argument together with the continuity of $\Psi_1$
(also in $w[M]$, see the proof of Lemma~\ref{Tech} (a)),  $\bar\pa$,
and of restriction to the boundary gives the second statement in full
generality. 
\end{pf}

\begin{lma}\label{lmafinish}
For $r>0$ small enough the image of $\Psi$ is a
submanifold of $\func_{2,\epsilon,\Lambda}$. Moreover, there exists
$M>0$, $r>0$, and a neighborhood $U$ of $(w[M],f[M],0)$ in
$\func_{2,\epsilon,\Lambda}$ such that $(w,f)\in U$ and
$U\cap\cand_{2,\epsilon,\Lambda}$ is contained in the image of $\Psi$
and corresponds to the subset $h=0$ in the coordinates
\begin{equation} \notag
(v,h,\lambda)\in 
B_{2,\epsilon}((w[M],f[M]),r)\times\sblv_{\frac32,\epsilon}(\pa
D_m,\R)\times\Lambda. 
\end{equation}
\end{lma}
\begin{pf}
Let $(w,f)\in\cand_{2,\epsilon}({\bf a},\kappa)$.
Let $K$ be as in Lemma~\ref{funcoord}. Then 
$D_2K(x,0,\sigma,0)\cdot\eta=\eta$ and
$D_4K(x,0,\sigma,0)=0$. Hence the differential of $\Psi$ at $(0,0,0)$ is 
\begin{equation} \notag
d\Psi(0,0,0)=
\left(\begin{matrix}
\iota          & 0   & 0    \\
\la\iota_\R\ra & \id & 0   \\
0              &  0  & \id
\end{matrix}\right),
\end{equation}
where $\iota$ denotes the inclusion of the tangent space of
$B_{2,\epsilon}((f,w),r)$ into $\sblv_{2,\epsilon}(D_m,\C^n)$ and
$\la\iota_\R\ra v$ denotes the $\R$-component of the vector field
$\tilde v$ which maps to $v$ under $\Pi_\C$ and is tangent to $L_0$. 
Note that the tangent space of $B_{2,\epsilon}((w,f),r)$ is a closed
subspace of $\sblv_{2,\epsilon}(D_m,\C^n)$. 
 
Thus, $d\Psi(0,0,0)$ is an injective map with closed image. Since the
first component of $\func_{2,\epsilon}$ is modeled on a Banach space
which allow a Hilbert-space structure we see that the image of the
differential admits a direct complement. Moreover, applying  
Lemma~\ref{Tech} to the explicit differential in Lemma~\ref{funcoord} 
we conclude
that the norm of the differential of $\Psi$ is Lipschitz in $v$
and $\lambda$ with Lipschitz constant depending only on
$\|Dw[M]\|_{1,\epsilon}$ and 
$\|DF[M]\|_{1,\epsilon}$. Hence the the implicit function theorem
shows that there exists $r>0$ and $W\subset\Lambda$ (independent of
$M$) such that the image of $B((w[M],f[M]),r)\times W$ is a
submanifold. From the norm-estimates on the differential it follows
that for $M$ large enough $(w,f)$ lies in this image.

The statement about surjectivity onto
$U\cap\cand_{2,\epsilon,\Lambda}$ follows from the fact that
$\Pi_\C(L_\lambda)$ is totally geodesic in the metric
$g(L_\lambda,\sigma)$ and Lemma \ref{Bdryval}. The statement on
coordinates is trivial.    
\end{pf}

\subsection{Bundle over conformal structures}\label{fcanG}
The constructions above all depend on the conformal structure $\kappa$ on
$D_m$. This conformal structure is unique if $m\le 3$. Assume that
$m>3$ and recall that we identified the space of conformal structures
$\conf_m$ on $D_m$ with an open simplex of dimension $m-3$. 

The space 
\begin{equation} \notag
\cand_{2,\epsilon,\Lambda}({\bf a})=\bigcup_{\kappa\in\conf_m}
\cand_{2,\epsilon,\Lambda}({\bf a},\kappa),
\end{equation}
has a natural structure of a locally trivial Banach manifold bundle
over $\conf_m$. To see this we must present local trivializations.

Let $\Delta$ denote the unit disk in the complex plane and let
$\Delta_m$ denote the same disk with $m$
punctures $p_1,\dots,p_m$ on the boundary and conformal structure
$\kappa$. Fixing the positions of $p_1,p_2,p_3$, this structure is
determined by the positions of the remaining $m-3$ punctures. 
We coordinatize a neighborhood of the conformal structure $\kappa$ in
$\conf_m$ as follows. Pick $m-3$ vector 
fields $v_1,\dots,v_{m-3}$, with 
$v_k$ supported in a neighborhood of $p_{k+3}$, $k=1,\dots,m-3$ in
such a way that $v_k$ generate a $1$-parameter family of
diffeomorphism $\phi_{p_{k+3}}^{\tau_k}\colon\Delta\to\Delta$,
$\tau_k\in\R$ which is a rigid rotation around $p_{k+3}$ and which is
holomorphic on the boundary. Let the supports of $v_k$ be
sufficiently small so that the supports of $\phi_{p_{k+3}}^{\tau_k}$,
$k=1,\dots,m-3$ are disjoint. Then the diffeomorphisms
$\phi_{p_3}^{\tau_1},\dots,\phi_{p_m}^{\tau_{m-3}}$ all commute. 
Define, for $\tau=(\tau_1,\dots,\tau_{m-3})\in\R^{m-3}$,
$\phi^\tau=\phi_{p_3}^{\tau_1}\circ\dots\circ\phi_{p_m}^{\tau_{m-3}}$
and a local 
coordinate system around $\kappa$ in $\conf_m$ by
\begin{equation}\notag
\tau\mapsto\Bigl(d\phi^\tau\Bigr)^{-1}\circ j_\kappa\circ d\phi^\tau. 
\end{equation}
These local coordinate systems give an atlas on $\conf_m$. 

Using this family we define the trivialization over 
$\R^{m-3}\approx U\subset\conf_m$ by composition with
$\phi^{-\tau}$. That is, a local 
trivialization over $U$ is given by
\begin{align}\notag
&\Phi\colon \cand_{2,\epsilon,\Lambda}({\mathbf a},\kappa)\times U
\to\cand_{2,\epsilon,\Lambda}({\mathbf a});\\ \notag
&\Phi(w,f,\lambda,\tau)=
\Bigl(w\circ\phi^{-\tau},
f\circ\phi^{-\tau},\lambda,\theta\Bigr) 
\end{align}

In a similar way we endow the space 
\begin{equation} \notag
\sblv_{1,\epsilon}(D_m,T^\ast D_m)
=\bigcup_{\kappa\in\conf_m}
\sblv_{1,\epsilon}(D_m,T^\ast D_m,g(\kappa)),
\end{equation}
with its natural structure as a locally trivial Banach space bundle
over $\conf_m$.

Representing the space of conformal structures $\conf_m$ in this way we
are led to consider its tangent space $T_\kappa\conf_m$ as generated by
$\gamma_1,\dots,\gamma_{m-3}$, 
where $\gamma_k=i\cdot\bar\pa v_k$, in the following sense. If $\gamma$ is
any section of $\End(T D_m)$ which anti-commutes with $j_\kappa$ and which
vanishes on the boundary then there exists unique numbers
$a_1,\dots,a_{m-3}$ and a unique vector field $v$ on $\Delta_m$ which
is holomorphic on the boundary and which vanish at $p_k$,
$k=1,\dots,m$ such that
\begin{equation}\label{solRH}
\gamma=\sum_{k} a_k \gamma_k+i\bar\pa v.
\end{equation}
The existence of such $v$ is a consequence of the fact that the
classical Riemann-Hilbert problem for the $\bar\pa$-operator on the
unit disk with tangential boundary conditions has index $3$ and is
surjective (the kernel being spanned by the vector fields 
$z\mapsto iz,z\mapsto i(z^2+1),z\mapsto z^2-1$). 
 
Going from the punctured disk $\Delta_m$ to $D_m$ with our
standard metric,  the behavior of the vector fields $v_j$ close to 
punctures where they are supported is easily described. In fact
the vector fields can be taken as $\pa_x$ in coordinates
$z=x+iy\in(\C_+,\R,0)$ in a neighborhood of the puncture $p$ on 
$\pa \Delta_m$. 
The change of coordinates taking
us to the standard end $[0,\infty)\times[0,1]$ is
$\tau+it=\zeta=-\frac{1}{\pi}\log z$ and we see the corresponding
vector field on $[0,\infty)\times[0,1]$ is $\frac{1}{\pi}e^{\pi\zeta}$ (where
we identify vector fields with complex valued functions).   
As in Proposition~\ref{R-Hcompactification}, we see that equation \eqref{solRH} holds on
$D_m$ with $v$ in a Sobolev space with (small) negative exponential weights
at the punctures.


\subsection{The $\bar\pa$-map and its linearization}\label{fcanH}
Consider the space $\sblv_{1,\epsilon}(D_m,{T^\ast}^{0,1}D_m\otimes\C^n)$ and
the closed subspace 
$\sblv_{1,\epsilon}[0](D_m,{T^\ast}^{0,1} D_m\otimes\C^n)$ consisting of elements whose 
trace (restriction to the boundary) is $0.$
Note that this space depends on the metric on $D_m$. In our case it
thus depends on the conformal structure.  For simplicity we keep the
notation and consider 
$\sblv_{1,\epsilon}[0](D_m,{T^\ast}^{0,1} D_m\otimes\C^n)$ as a bundle over
$\conf_m$. We extend this bundle to a bundle over $\Lambda$ making it
trivial in the $\Lambda$ directions and denote the result 
$\sblv_{1,\epsilon,\Lambda}[0](D_m,{T^\ast}^{0,1}D_m\otimes\C^n)$.

The $\bar\pa$-map is the map
\begin{align}\notag
&\hat\Gamma\colon\cand_{2,\epsilon,\Lambda}({\mathbf a})
\to\sblv_{1,\epsilon,\Lambda}[0](D_m,{T^\ast}^{0,1} D_m\otimes \C^n);\\ \notag
&\hat\Gamma(w,f,\kappa,\lambda)= \Bigl(dw + i\circ dw\circ
j_\kappa,\kappa,\lambda\Bigl).
\end{align}  
We will denote the first component of this map simply $\Gamma$.
An element $(w,f,\kappa,\lambda)$ is thus holomorphic with respect to
the complex structure $j_\kappa$ if and only if
$\hat\Gamma(w,f,\kappa,\lambda)=(0,\kappa,\lambda)$. Hence, if 
$L_\lambda$, $\lambda\in\Lambda$ is a family of chord
generic Legendrian submanifolds then the (parameterized)
moduli-space of holomorphic disks with 
boundary on $L_\lambda$, positive puncture at $a_1$, and negative
punctures at $a_2,\dots,a_m$ is naturally identified with the preimage
under $\hat \Gamma$ of the $0$-section in
$\sblv_{1,\epsilon,\Lambda}[0](D_m,{T^\ast}^{0,1}D_m\otimes\C^n)$ for
sufficiently small $\epsilon\in[0,\infty)$. 

We compute the linearization of the $\bar\pa$-map. As in Section~\ref{fcanG} 
we think of tangent vectors $\gamma$ to $\conf_m$ at $\kappa$ as
sections of $\End(T D_m)$. For $\kappa\in\conf_m$ and $u\colon
D_m\to\C^n$, let  $\bar\pa_\kappa u=du+i\circ du\circ j_\kappa$.
 
Let $(w,f,\kappa,0)\in\cand_{2,\epsilon,\Lambda}({\mathbf
a})$. Identify the tangent space of $\cand_{2,\epsilon,\lambda}({\bf
a})$ at $(w,f,\kappa,0)$ with 
$T \B_{2,\epsilon}((w,f),r)\times T_\kappa\conf_m\times T_0\Lambda$. 

\begin{lma}\label{lem:fcan.gamma}
The differential of $\Gamma$ at $(w,f,\kappa,0)$ is the map
\begin{equation}\label{DGamma}
d\Gamma(v,\gamma,\lambda)=\bar\pa_\kappa v 
+\bar\pa_\kappa\bigl(Y_0^\sigma(w,\lambda) \bigr)
+i\circ dw\circ \gamma.
\end{equation}
Recall $Y_0^\sigma$ was defined in \eqref{defYY}.
\end{lma} 

\begin{pf}
Assume first $w$ and $f$ are constant close to punctures.
Let $\B_{2,\epsilon}((w,f),r)\times\conf_m\times\Lambda$ be
a local coordinates around $(w,f,\kappa,0)$.

Let $K(x,\xi,\sigma,\lambda)=\psi_\lambda^\sigma(x)+\xi$.
Then 
\begin{equation} \notag
R(x,\xi,\sigma,\lambda)=
\exp_{\psi_\lambda^\sigma(x)}^{\lambda,\sigma}A_\lambda^\sigma\xi-
K(x,\xi,\sigma,\lambda) 
\end{equation}
satisfies
\begin{equation} \notag
R(x,0,\sigma,\lambda)=0,\quad D_2R(x,0,\sigma,0)=0, \quad D_4R(x,0,\sigma,0)=0;
\end{equation} 
thus, Lemma \ref{Tech} (b) implies that 
\begin{equation} \notag
\|R(w,v,\sigma,\lambda)\|_{2,\epsilon}\le
C(\|v\|_{2,\epsilon}^2+|\lambda|^2). 
\end{equation}

Continuity of the linear operators  
\begin{equation} \notag
\bar\pa_{\kappa+\gamma}\colon
\sblv_{2,\epsilon}(D_m,\C^n)\to\sblv_{1,\epsilon}(D_m,T^\ast
D_m\otimes\C^n), 
\end{equation}
where we use local coordinates $\R^{m-3}$ on $\conf_m$ and
$\kappa+\gamma\in\R^{m-3}\subset\conf_m$,  
shows that
\begin{equation}\label{remest}
\|\bar\pa_{\kappa+\gamma} R(w,v,\sigma,\lambda)\|_{1,\epsilon}\le
C(\|v\|_{2,\epsilon}^2+|\lambda|^2). 
\end{equation}
It is straightforward to check that  
\begin{align}\notag
\Bigl\|&\bar\pa_{\kappa+\gamma} K(w,v,\kappa+\mu,\lambda)
  -\bar\pa_{\kappa}w
-\Bigl(\bar\pa_{\kappa}v+ 
\bar\pa_{\kappa}\bigl(Y_0(w,\lambda) \bigr)
+i\circ dw\circ \gamma\Bigr)\Bigr\|_{1,\epsilon}\\\label{Deasy}
&\le C(\|v\|_{2,\epsilon}^2+|\lambda|^2+|\gamma|^2)
\end{align}
Equations \eqref{Deasy} and \eqref{remest} imply the lemma in the special
case when $(w,f)$ is constant close to punctures (and in the general case if
$\dim(\Lambda)=0$). 

If $(w,f)$ is not constant close to punctures consider the maps 
$(w[M],f[M])$ which are constant close to punctures. We have
$(w[M],f[M])\to (w,f)$ as $M\to\infty$. Since the local coordinates
are $C^1$ a limiting argument proves \eqref{DGamma} in the general case. 
\end{pf}

\subsection{Auxiliary spaces in the semi-admissible case}\label{fcanI}
In Section~\ref{8degentrans.section} we show that for a dense open set of 
semi-admissible Legendrian submanifolds $L$ no rigid holomorphic
disks with boundary on $L$ have exponential decay at their degenerate
corners. Once this has been shown we know that if $0$ is the degenerate
corner and $L$ has the form \eqref{2LocFormofST} around $0$ then for any rigid
holomorphic disk $u\colon D_m\to\C^n$ with puncture $p$ mapping to $0$
there exists 
$M>0$ and $c\in\R$ such that 
\begin{equation} \notag
u(\zeta)=\Bigl(-2(\zeta+c)^{-1},0,\dots,0\Bigr)+\Ordo(e^{-\theta|\zeta|}),
\text{ for }\zeta\in E_{p}[\pm M], 
\end{equation}
where $\theta>0$ is the smallest non-zero complex angle of the Reeb
chord at $0$. (Here we implicitly assume that $P_2$ in our standard
self tangency model lies above $P_1$ in the $z$-direction, and that
neighborhoods of positive (negative) punctures are parameterized by
$[1,\infty)\times[0,1]$ ($(-\infty,-1]\times[0,1]$).)
To study disks of this type we use the following construction. 

Let $a_0$ denote the Reeb-chord at $0$. Assume that ${\mathbf a}$ has the
Reeb-chord $a_0$ in $k$ positions. For $C=(c_1,\dots,c_k)\in\R^k$ fix 
a smooth reference function which equals
\begin{equation} \notag
u_{\rm ref}^C(\zeta)=
\Bigl(-2(\zeta+c_j)^{-1},0,\dots,0\Bigr), 
\end{equation}
in a neighborhood of the $j^{\rm th}$ puncture mapping to $a_0$. 

Let $L_\lambda$, $\lambda\in\Lambda$ be a family of semi-admissible
Legendrian submanifolds. 
We construct 
for $\epsilon\in[0,\infty)^m$, with those components $\epsilon_j$ of
$\epsilon$ which correspond to punctures mapping to the degenerate
corner satisfying $0<\epsilon_j<\theta$ and for fixed $C\in\R^k$, the spaces
\begin{equation} \notag
\func_{2,\epsilon}^C({\mathbf a})
\end{equation}
by using reference functions looking like $u_{\rm ref}^C$ for
$C\in\R^k$ in neighborhoods of punctures mapping to
$a_0^\ast=0\in\C^n$. We construct local coordinates as in Section~\ref{fcanA}  
taking advantage of the fact that $\lambda\in\Lambda$ fixes
$a_0^\ast$. Also we consider the space 
\begin{equation} \notag
\cand^C_{2,\epsilon,\Lambda}({\mathbf a}),
\end{equation}
which is defined in the same way as before. We note that the
construction giving local coordinates on this space in Section~\ref{fcanF} 
can still be used since in the semi-admissible case we need
not cut-off the first component of $w$ in $(w,f)$ close to punctures
mapping to $c_0$ since $\lambda\in\Lambda$ are assumed to preserve the
product structure and $\gamma_1$ and $\gamma_2$.

With this done we consider the bundle
\begin{equation}\label{tildecand}
\widetilde\cand_{2,\epsilon,\Lambda}=
\bigcup_{C\in\R^k}\cand_{2,\epsilon,\Lambda}^C({\mathbf  a}), 
\end{equation}
which we present as a locally trivial bundle over $\R^k$ as follows.

In the case that ${\mathbf a}$ has $\ge 3$ elements we fix for $C\in\R^k$ the
diffeomorphism $\phi^C\colon D_m\to D_m$ which equals to 
$\zeta\to\zeta+c_j$ in $E_{p_j}[M]$ for any puncture $p_j$ mapping to
$c_0$, equals the identity on $D_m\setminus \bigcup_j E_{p_j}[M-2]$,
and is holomorphic on the boundary.

(Since we often reduce the few punctured cases to the many punctured
case, see Section~\ref{ssecmarked}, the following two constructions will not be
used in the sequel, we add them here for completeness.)
In case ${\mathbf a}$ has length $1$ we think of $D_1$ as of the upper
half-plane $\C_+$ with the puncture at $\infty$. The map
$z\mapsto-\frac{1}{\pi}\log z$ identifies the region 
$\{z\in\C_+\colon |z|>R\}$ with the strip 
$[\frac{1}{\pi}\log R,\infty)\times[0,1]$ where we think of the latter
space as a part of $E_p$, where $p$ is the puncture of $D_1$. Also,
this map takes the conformal reparameterization $z\mapsto e^{\pi C}z$ to
$\phi^C\colon \zeta \mapsto \zeta+C$ in $E_p$ and we identify $\R$ with this set of
conformal reparameterizations $\{\phi^C\}_{C\in\R}$.
In case ${\mathbf a}$ has length $2$ we think of $D_2$ as the strip
$\R\times [0,1]$ and identify $\R$ with the conformal
reparameterizations $\phi^C(\zeta)=\zeta\mapsto\zeta+ C$.      
 
Using composition with the maps $\phi^C$ we construct local
trivializations of the bundle in \eqref{tildecand}. We then find local
coordinates 
$B_{2,\epsilon}(0,r)\times\conf_m\times\R^k\times\Lambda$  
on 
$\widetilde\cand_{2,\epsilon,\Lambda}$ and the linearization of the
$\bar\pa$-map $\Gamma$ at $(w,f,\kappa,0,0)$ is
\begin{align}\notag
d\Gamma(v,\gamma,c,\lambda)&=\bar\pa_\kappa v 
+\bar\pa_\kappa\bigl((Y_0(w,\lambda))\bigr)\\\label{Cdiff}
&+i\circ dw\circ \gamma+ 
\bar\pa_\kappa\left(dw\cdot \left(\frac{\pa\phi^C}{\pa C}|_{C=0}\right)\cdot c\right).
\end{align}
Here $c=(c_1,\dots,c_k)$ is a tangent vector to $\R^k$ written in the
basis $\{\hat C_1,\dots,\hat C_k\}$ where $\hat C_j$ is a unit vector
in the tangent space to $C_j\in\R$. 
We notice that the second term in \eqref{Cdiff} lies in  
$\sblv_{1,\epsilon}[0](D_m, {T^\ast}^{0,1} D_m)$ because of the special
assumptions on $L_\lambda$ in a neighborhood of $c_0^\ast$ and that
the last term does as well since the difference of $w$ and the holomorphic
function $u_{\rm ref}^C$ lies in a Sobolev space weighted by
$e^{\epsilon\tau}$ in $E_{p_j}$ and for a holomorphic function the
last term vanishes in the region where $\phi^C$ is just a translation.

\subsection{Homology decomposition}\label{fcanJ}
Let $L\subset\C^n\times\R$ be a (semi-)admissible
Legendrian submanifold. Let ${\mathbf c}=c_0c_1\dots c_m$ be a word of
Reeb chords of $L$. If
$(u,f)\in\cand_{2,\epsilon}({\mathbf c})$ then the homotopy classes of
the paths induced by $(u|\pa D_m,f)$ in $L$ connecting the Reeb chord endpoints
determines the path component of $(u,f)\in\cand_{2,\epsilon}({\mathbf c})$. 

Let $A\in H_1(L)$ and let 
$\cand_{2,\epsilon}({\mathbf c};A)\subset\cand_{2,\epsilon}({\mathbf
c})$ be the union of those path components of
$\cand_{2,\epsilon}({\mathbf c})$ such that the homology class of the loop
$f(\pa D_m)\cup (\bigcup_j\gamma_j)$  equals $A$, where $\gamma_j$ is
the capping 
path chosen for the Reeb chord $c_j$ endowed with the appropriate
orientation, see Section \ref{1modulisection}. For fixed conformal
structure $\kappa$ we write $\cand_{2,\epsilon}({\mathbf c},\kappa;A)$
and in the chord semi generic case 
${\widetilde\cand}_{2,\epsilon}({\mathbf c};A)$ and interpret these
notions in a similar way.

\section{Fredholm properties of the linearized equation}
\label{7index.section}

In this section we study properties of the linearized
$\bar\pa$-equation. In particular we determine the index of the
$\bar\pa$-operator with Legendrian boundary conditions. It
will be essential for our  geometric applications to use weighted
Sobolev spaces and to understand how constants in certain
elliptic estimates depend on the weights. 

Our presentation has two parts: the ``model'' case where the domain
is a strip or half-plane; and the harder case where the domain is $D_m.$

In Section \ref{7cokerregularity.section}, we show that an element of
the cokernel has a smooth representative.
In Section \ref{kercokelts}, we derive expansions for
the kernel and cokernel elements.
We use these two subsections
in Sections \ref{model} through \ref{7modelproblem3.section}, 
to prove the elliptic estimate for the model problem, as well as
derive a formula for the index.
In Sections~\ref{7boundaryconditions.section},
we set up the boundary conditions for the linearized problem with
domain $D_m.$
In Sections \ref{7deform.section} through \ref{7Jtrick.section},
we prove the Fredholm properties for the $D_m$ case.
In Sections \ref{7Jtrick.section} and \ref{7CZST.section},
we connect the index formula to the Conley-Zehnder index of
Section \ref{1frontcomp}.

\subsection{Cokernel regularity}
\label{7cokerregularity.section}

To control the cokernels of the operators studied below we use the
following regularity lemma.

For this subsection only, we use coordinates $(x,y)$ for
the half-plane $\R^2_+ = \{(x,y)\colon y \ge 0\}.$
Let $A\colon\R\to \GL(\C^n)$ be a smooth map with $\det(A)$ uniformly
bounded away from $0$ and all derivatives uniformly bounded. 
We also simplify notation for this subsection only and
define the following Sobolev spaces:
let $\sblv_k = \sblv_k(\R^2,\C^n)$;
let $\bar\sblv_k$
denote the space of restrictions of elements in $\sblv_k$ to $\inr(\R^2_+)$;
let $\dot\sblv_k$ denote the subspace
of elements in $\sblv_k$ with support in $\R^2_+$;
let $\bar\sblv_1[0]$ denote the subspace of all elements in
$\bar\sblv_1$ which vanish on the boundary; 
and let
$\bar\sblv_2[A]$ denote the subspace of 
elements $u$ in $\bar\sblv_2$ such that $u(x,0)\in A(x)\R^n$ and such
that the trace of $\bar\pa u$ (its restriction to the boundary) equals
$0$ in $\sblv_{\frac32}(\R,\C^n)$. 

An element $\xi$ in the cokernel of $\bar\pa$ will be in the dual
space of $\sblv_{1}[0]$. The dual of $\sblv_{1}$ is
$\sblv_{-1}$ and thus the dual of $\sblv_{1}[0]$ is the
quotient space 
\begin{equation} \label{7coker.eqn}
\sblv_{-1}/\sblv_{1}[0]^\perp,
\end{equation}
where $\sblv_{1}[0]^\perp$ denotes the annihilator of
$\sblv_{1}[0]$ in $\sblv_{-1}$. As usual, let $\la\,,\,\ra$ denote the
standard Riemannian inner product on $\C^n\approx\R^n$.

\begin{lma}\label{lmacokerreg} 
Fix $h>0$ and
assume that $v\in\dot\sblv_{-1}$ satisfies
\begin{equation}
\int_{\R\times[0,h)}\la\bar\pa u,v\ra\,dx\wedge dy=0,
\end{equation}
for all $u\in\bar\sblv_2[A]$ with compact support in
$\R\times[0,h)$. Then, for 
every $\epsilon$ with $0<\epsilon<h$ and every $k>0$, the class
$[v]\in\dot\sblv_{-1}/\bar\sblv_1[0]^\perp$ of 
$v$ contains an element $v_0$ which is $C^k$ in
$\R\times[0,\epsilon)$, up to and including the boundary.  
\end{lma}

\begin{pf}
Let $m>2$ be such that $\sblv_m\subset C^k$ and $\epsilon<h$ be
given. Extend $A$ to a smooth map  
$\hat A$ on $\R^2_+$ with uniformly bounded derivatives as
follows. Let $0<\delta<\frac{\epsilon}{4}$ and let $\eta\colon\R_+\to\R_+$
be a smooth non-decreasing function such that $\eta(y)=y$ for $0\le
y\le\frac{\delta}{2}$ and $\eta(y)=\delta$ for $y\ge \delta$. Define 
\begin{equation} \notag
\hat A(x,y)=A(x)+\sum_{k=1}^{m+1}\frac{i^{k}\eta(y)^k}{k!}\pa^k_xA(x), 
\end{equation}  
and choose $\delta>0$ so small that 
$|\det(\hat A(x,y))|>\rho>0$. Note that $\bar\pa \hat A(x,0)$ vanishes to
order $m$. Therefore multiplication with $\hat A$ and $\hat A^{-1}$
gives the following commutative diagram where horizontal arrows are
isomorphisms:
\begin{equation} \notag
\begin{CD}
\bar\sblv_2[\id] @>{\times \hat A}>> \bar\sblv_2[A]\\
@V{\bar\pa+\hat A^{-1}\bar\pa\hat A}VV @VV{\bar\pa}V\\
\bar\sblv_1[0] @<<{\times \hat A^{-1}}< \bar\sblv_1[0]\\
\end{CD}.
\end{equation}  
The lemma follows once its analogue with the vertical operator on the
left in the diagram replacing $\bar\pa$ is proved. Let
$B=\hat A^{-1}\bar\pa \hat A$.

Consider smooth functions $u$ with compact support in
$\R\times(0,h)$. For such $u$ we have (recall
$\dot\sblv_{k}\subset\sblv_{k}$, all $k$)  
\begin{equation}
0=\int\la(\bar\pa+B) u,v\ra\,dx\wedge dy
=\int\la u,(\pa+B^\ast) v\ra\,dx\wedge dy.
\end{equation}
Hence $(\pa + B^\ast) v=0$ in $\R\times(0,h)$.
Let $\alpha$ and $\hat\alpha$ be smooth functions with all derivatives
uniformly bounded, with support in
$\R\times(0,h)$, and such that $\hat\alpha=1$ on the support of
$\alpha$. Then the elliptic estimate for $\pa$ 
on $\R^2$ implies (with
$K>0$ a sufficiently large constant and $k$ any integer)
\begin{align}\notag
\|\alpha v\|_k &\le C(\|\alpha v\|_{k-1}+\|\pa(\alpha v)\|_{k-1})\\\notag
&\le C(\|\alpha v\|_{k-1} + \|(\pa+B^\ast)\alpha v\|_{k-1} 
+\|B^\ast \alpha v\|_{k-1})\\\notag
&\le C'(\|\alpha v\|_{k-1} + \|(\pa+B^\ast)\alpha v\|_{k-1})\\\notag
&\le 
C'(\|\alpha v\|_{k-1}+
\|((\pa+B^\ast)\alpha)v\|_{k-1}+
\|\alpha(\pa+B^\ast)v\|_{k-1})\\\label{vint}
&\le C''\|\hat \alpha v\|_{k-1}
\end{align} 
since $(\pa+B^\ast)v=0$ in $\R\times(0,h)$. 
It follows from \eqref{vint} 
that $v$ is smooth in $\R\times(0,h)$. Thus, $(\pa+B^\ast)v\in\dot\sblv_{-2}$
is a distribution such that 
$\supp((\pa +B^\ast) v)\cap\R\times[0,h)\subset \R\times\{0\}$ and hence
there exist distributions $f$ and $g$ on $\R$ such that   
\begin{equation}
(\pa+B^\ast)v=\delta(y)\otimes f(x) + \delta'(y)\otimes g(x)
\end{equation}
in $\R\times[0,h)$, where $\delta$ denotes the Dirac distribution and
$\delta'$ its derivative. 
Since any function 
$\phi\in\sblv_2(\R,\C^n)$ can be extended constantly
in the $y$-direction in a neighborhood of  
the real axis so that it lies in $\sblv_2$ we find that
$f(x)\in\sblv_{-2}(\R,\C^n)$. Thus, 
$\delta(y)\otimes f(x)$ lies in $\dot\sblv_{-2}$ and
therefore so does $\delta'(y)\otimes g(x)$. 

Let $u\in\bar\sblv_2[\id]$ have support in $\R\times[0,h)$ and let
$\hat u$ be an extension of $u$ to some neighborhood of $\R^2_+$ in
$\R^2$. Let $\tilde{B}$ and $\tilde{B}^\ast$ denote the 
extensions of $B$ and $B^\ast$
to $\R^2$ by defining them to be $0$ on $\inr(\R^2_-)$. Note that by
the definition of $\hat A$, $\tilde{B}$  and $\tilde{B}^\ast$ are 
$C^m$-functions with
uniformly bounded derivatives. Let $\gamma$ be a smooth 
function which equals $1$ on $\R^2_+$ and equals $0$ outside a
neighborhood of $\R^2_+$ in $\R^2$. Then
\begin{align}\notag
&\int\la u,(\pa+\tilde{B}^\ast) v\ra\,dx\wedge dy=
\int\la \gamma \hat u,(\pa+ \tilde{B}^\ast) v\ra\,dx\wedge dy\\\label{genint}
&=\int\la (\bar\pa\gamma)\hat u,v\ra\,dx\wedge dy+
\int\la (\gamma)(\bar\pa+\tilde{B})\hat u,v\ra\,dx\wedge dy=0+0.
\end{align} 
Let $\phi\colon\R\to\R^n$ be any smooth compactly supported function. 
Let $\beta$ be a function with $\beta(0)=0$ and $\beta'(0)=1$ and
$\beta=0$ outside a small neighborhood of $0$. Then
$$
u(x,y)=\phi(x)+i\beta(y)\phi'(x),
$$ 
lies in $\sblv_2[\id]$ and by \eqref{genint} 
$$
0=\int\la u,(\pa+\tilde{B}^\ast) v\ra\,dx\wedge dy
=\int (\Re f-\Im g')\cdot\phi\,dx.
$$
Hence $\Re f=\Im g'$. Since 
$\Re f\in\sblv_{-2}(\R,\C^n)$ we find that  
if $v^{\rm I}=\delta(y)\otimes\Im g$ then 
$$
\bar\pa v^{\rm I}=\delta(y)\otimes\Re f
+i\delta'(y)\otimes\Im g\in\sblv_{-2}
$$  
and hence
$$
\|v^{\rm I}\|_{-1}\le C(\|v^{\rm I}\|_{-2}+\|\bar\pa v^{\rm
I}\|_{-2})<\infty. 
$$
Let $a_j\in C^\infty_0(\R^2_+)$ be a sequence such that $a_j\to v^{\rm
I}$ in $\sblv_{-1}$ as $j\to\infty$. Define 
$a_j^{-}(x,y)=a_j(x,-y)$. Then $a_j^-$ are supported
in the lower half plane and
$\bar\pa a_j(x,-y)=\pa a_j^{-}(x,y)$. 
Hence $a^-_j$
approaches a distribution $v^{\rm II}$ with support on the boundary, 
$\pa v^{\rm II}=\bar\pa v^{\rm I}$, and since 
$v^{\rm I}\in\bar\sblv_{-1}[0]$ also $v^{\rm II}\in\bar\sblv_1[0]^\perp$. 
Let $v_{\rm I}=v-v^{\rm II}$ then   
$$
(\pa+\tilde{B}^\ast) v_{\rm I}= i\delta(y)\otimes \Im f(x) + \delta'(y)\otimes
\Re g(x), 
$$   
in $\R\times[0,h)$.

Let $b_j\in C^\infty_0(\R^2_+)$ be a sequence such that 
$b_j\to v_{\rm I}$ as $j\to\infty$   
and define $b_j^\ast(x,y)=\overline{b_j(x,-y)}$. 
Then $\pa b^\ast(x,y)=\overline{\pa b_j(x,y)}$. 
Define 
$$
B^d(x,y)=\begin{cases}
         B(x,y), &\text{ if }y\ge 0\\ 
         \overline{B(x,-y)}, &\text{ if } y\le0
         \end{cases}
$$	
Then $B^d$ is a $C^m$ function. If $b^d_j=\frac12(b_j+ b_j^\ast)$ then $b^d_j$
approaches a distribution $v^d_{\rm I}$ with 
$$
(\pa + (B^d)^\ast)v^d_{\rm I}=\delta'(y)\otimes \Re g(x),
$$
in $\R\times[0,h)$. Again,
let $\phi\colon\R\to\R$ be any smooth function with compact support.  
Let $\theta\colon\R\to\R$ be an odd smooth function with
$\theta(0)=0$, $\theta'(0)=1$ and $\theta=0$ outside a small
neighborhood of $0$. Then with $u(x,y)=\phi(x)\theta(y)$  
$$
\int\la(\bar\pa+B^d) u,v^d_{\rm I}\ra\,dx\wedge dy=
-\int\la u,(\pa + (B^d)^\ast)v^d_{\rm I}\ra\,dx\wedge dy
=\int_\R\phi(x)\Re g(x)\,dx. 
$$
But, writing $B^d=B^d_{\Re}+iB^d_{\Im}$, 
\begin{align}\notag
&\int\la(\bar\pa+B^d) u,v^d_{\rm I}\ra\,dx\wedge dy=\\\notag
&\int\Bigl(\Re v^d_{\rm I}\cdot(B^d_{\Re}(x,y)+\phi'(x)\theta(y)-
B^d_{\Im}(x,y)-\phi(x)\theta'(y))\\\notag
&+\Im v^d_{\rm I}\cdot(B^d_{\Re}(x,y)+\phi(x)\theta'(y)+
B^d_{\Im}(x,y)+\phi'(x)\theta(y))\Bigr)\,dx\wedge dy=\\\notag
\lim_{j\to\infty}&\int\Bigl(\Re
b_j^d(x,y)\cdot(B^d_{\Re}(x,y)+\phi'(x)\theta(y)- 
B^d_{\Im}(x,y)-\phi(x)\theta'(y))\\
&+\Im b_j^d(x,y)\cdot(B^d_{\Re}(x,y)+\phi(x)\theta'(y)+
B^d_{\Im}(x,y)+\phi'(x)\theta(y))\Bigr)\,dx\wedge dy=0,
\end{align}
since both summands in the last integral are odd in $y$. Hence 
$\Re g(x)=0$ and $(\pa+(B^d)^\ast) v^d_{\rm I}=0$ in $\R\times[0,h)$ and
therefore by the elliptic estimate for $\pa$, $v^d_{\rm I}$ lies in  
$\sblv_{m+1}$ in $\R\times[0,h)$.  
Let $\eta\colon \R^2\to[0,1]$ be a smooth function which is $1$ on
$\R\times [0,\epsilon)$ and $0$ outside $[0,h)$. Let $v_{\rm II}=\eta
v^d_{\rm I}|\R^2_+$.  
Then $v_{\rm II}\in\dot\sblv_0\subset\dot\sblv_{-1}$ and 
$v^{\rm III}=v_{\rm II}-\eta v_{\rm I}\in\dot\sblv_{-1}$ is a 
distribution with support on the boundary. Hence
$v^{\rm III}=\delta(y)\otimes h(x)$ and thus $v^{\rm III}\in
\sblv_1[0]^\perp$.   
Since 
\begin{equation}
v=v_{\rm I}+v^{\rm II}=\eta v_{\rm I}+ (1-\eta)v_{\rm I}+ v^{\rm II}=
v_{\rm II}+(1-\eta)v_{\rm I}-v^{\rm III}+ v^{\rm II},
\end{equation}
where $v^{\rm III}+ v^{\rm II}\in\bar\sblv_1[0]^\perp$ 
we find that $v_0=v_{\rm II}+(1-\eta)v_{\rm I}$ is a representative of 
$[v]$ which is as smooth as required in $\R\times[0,\epsilon)$.
\end{pf}

\subsection{Kernel and cokernel elements}\label{kercokelts}

Consider the strip $\R\times[0,1]\subset\C$ endowed with the standard
flat metric, the corresponding complex structure and coordinates 
$\zeta=\tau+it$.    
For $k\ge 0$, let
$$
\sblv_k=\sblv_k(\R\times[0,1],\C^n),
$$ 
and for $k\le 0$, let $\sblv_k$ denote the $L^2$-dual of $\sblv_{-k}.$ 
We also use the notions $\sblv_k^{\rm loc}$ which are
to be understood in the corresponding way.

If $u\in\sblv_k^{\rm loc}$ then the restriction
of $u$ to $\pa(\R\times[0,1])=\R\cup\R+i$ lies in  
$\sblv_{k-\frac12}^{\rm loc}(\R\cup\R+i,\C^n)$.
For $u\in\sblv_1^{\rm loc}$ consider the boundary conditions 
\begin{align}
\int_\R\la u,v\ra\,d\tau=0\quad & \label{rabdr0}
\text{ for all } v\in C^\infty_0(\R,i\R^n),\\
\int_{\R+i}\la u,v\ra\,d\tau=0\quad & \label{rabdr1}
\text{ for all }v\in C^\infty_0(\R+i,\R^n), 
\end{align}

Let $f\colon\R\times[0,1]\to\C^n$ be a smooth function satisfying
\eqref{rabdr0} and \eqref{rabdr1}. Define the function 
$f^d\colon\R\times[0,2]\to\C^n$ as
$$
f^d(\tau+it)=
\begin{cases}
f(\tau+it) &\text{ for } 0\le t\le 1,\\ 
-{\overline f}(\tau + i(2-t)) &\text{ for } 1< t\le 2,
\end{cases}
$$ 
where ${\overline w}$ denotes the complex conjugate of $w\in\C^n$.
Then $f^d$ and $\pa_\tau f^d$ are continuous, $\pa_t f^d$ may have a
jump discontinuity over the line $\R+i$, $f^d(\tau+0i)=-f^d(\tau+2i)$,
and $\|f^d\|_1= 2\|f\|_1$. Hence we can define the double
$u^d\in\sblv_1^{\rm loc}(\R\times[0,2])$ of any
$u\in\sblv_1^{\rm loc}$ which satisfies \eqref{rabdr0} and
\eqref{rabdr1}. For $u\in\sblv_k^{\rm loc}$, let 
$\bar\pa u=(\pa_\tau+i\pa_t)u$ and $\pa u=(\pa_\tau-i\pa_t)u$.    

\begin{lma}\label{lmaralocsol}
If $u\in\sblv_1^{\rm loc}$ satisfies \eqref{rabdr0} and \eqref{rabdr1}
and 
\begin{itemize}
\item[{\rm (a)}] $\bar\pa u=0$ in the interior of $\R\times[0,1]$ then
$$
u(\zeta)=\sum_{n\in\Z}C_n\exp\left((\pih+n\pi)\zeta\right),
$$
where $C_n\in\R$.
\item[{\rm (b)}] $\pa u=0$ in the interior of $\R\times[0,1]$ then 
$$
u(\zeta)=\sum_{n\in\Z}C_n\exp\left((\pih+n\pi)\bar\zeta\right),
$$
where $C_n\in\R$.
\end{itemize}
Moreover, if $u$ satisfies {\rm (a)} or {\rm (b)} and  $u\in\sblv_k$ for some 
$k\in\Z$ then $u=0$.  
\end{lma}

\begin{pf}
We prove (a), (b) is proved in the same way. Clearly it is enough to
consider one coordinate at a time. So assume the target is $\C$
and let $u$ be as in the statement. 

Consider $u^d$, then $\bar\pa u^d$ is an element of 
$\sblv_{0}^{\rm loc}(\R\times[0,2],\C)$ 
with support on $\R + i\cup\pa(\R\times[0,2])$. Such a distribution is
a three-term linear combination of tensor products of a Dirac-delta in
the $t$-variable and a distribution on $\R$ and hence lies
in $\sblv_0(\R\times[0,2],\C)$ only if it is
zero. Thus $\bar\pa u=0$ and we may use elliptic regularity to
conclude that $u$ is smooth in the interior of $\R\times[0,2]$. 
(In fact, doubling again and using the same argument, we find that $u$
is smooth also on the boundary.)   

We may now
Fourier expand $u^d(\tau,\cdot)$ in the eigenfunctions $\phi$ of the
operator $i\pa_t$ which satisfy the boundary condition 
$\phi(0)=-\phi(2)$. These eigenfunctions are 
$$
t\mapsto\exp\left(i(\pih+n\pi)t\right), \text{ for }n\in\Z.
$$
We find
$$
u^d=\sum_n c_n(\tau)\exp\left(i(\pih+n\pi)t\right)
$$
where, by the definition of $u^d$, $c_n(\tau)$ are real valued
functions and
$$
\bar\pa u^d=\sum_n
\left(c'_n(\tau)-(\pih+n\pi)c_n(\tau)\right) 
\exp\left(i(\pih+n\pi)t\right). 
$$
Hence,
$$
u(\zeta)=\sum_n C_n\exp\left((\pih+n\pi)\zeta\right).
$$

Assume that $u\in\sblv_k$ for some $k\in\Z$. Then, since for 
$j\ge 0$ the restriction of any $v\in\sblv_j(\R\times[0,2],\C)$ to
$\R\times[0,1]$ lies in $\sblv_j$,    
$$
\lambda_u(v)=\int_{\R\times[0,2]}\la v, u^d\ra\,d\tau\wedge dt,
$$ 
is a continuous linear functional on $\sblv_j(\R\times[0,2],\C)$
for $j=k$ if $k\ge 0$ or $j=-k$ if $k<0$. 

Let $\psi\colon\R\to[0,1]$ be a smooth function equal to $1$ on
$[0,1]$ and $0$ outside $[-1,2]$. For $n,r\in\Z$ let 
$$
\alpha_{n,r}(\tau+it)=\psi(\tau+r)\exp(i(\pih+n\pi)t).
$$   
Then $\alpha_{n,r}\in\sblv_j(\R\times[0,2],\C)$ and $\|\alpha_{n,r}\|_j=K(n)$ 
for some constant $K(n)$ and all $r$. It is straightforward to see that 
$$
\lambda_u(\alpha_{n,r})=
2C_n\int_{r-1}^{r+2}\psi(\tau+r)\exp((\pih+n\pi)\tau)\,d\tau=l_{n,r}.
$$
The set $\{l_{n,r}\}_{r\in\Z}$ is unbounded unless $C_n=0$. 
Hence $\lambda_u$ is continuous only if each $C_n=0$.
\end{pf}

\subsection{The right angle model problem}\label{model}
As mentioned we will use weighted Sobolev spaces. The weight functions
are functions on $\R\times[0,1]$ which are independent of $t$ and
have the following properties.

For $a=(a^+,a^-)\in\R^2$ and $\theta\in [0,\pi)$, let
\begin{equation}\label{m(a)}
m(\theta,a)=\min
\Bigl\{|n\pi +\theta + a^+|, |n\pi +\theta + a^-|\Bigr\}_{n\in\Z}.
\end{equation}

For $a\in\R^2$ with $m(\pih,a)>0$, let
$e_a\colon\R\to\R$ be a smooth positive function 
with the following properties:
\begin{itemize}
\item[{\bf P1}]
There exists $M>0$ such that $e_a(\tau)=e^{a^+\tau}$ for $\tau\ge M$
and $e_a(\tau)=e^{a^-\tau}$ for $\tau\le -M$.
\item[{\bf P2}]
The logarithmic derivative of $e_a$,
$\alpha(\tau)=\frac{e_a'(\tau)}{e_a(\tau)},$ is (weakly)
monotone and $\alpha'(\tau)=0$ if and only if $\alpha(\tau)$ equals
the global maximum or minimum of $\alpha$. 
\item[{\bf P3}]
The derivative of $\alpha$ satisfies
$|\alpha'(\tau)|<\frac15 m(\pih,a)^2$ for all $\tau\in\R$.
\end{itemize}  
Let 
$$
\mu=(\mu_1,\dots,\mu_n)=
(\mu_1^+,\mu_1^-,\dots,\mu_n^+,\mu_n^-)\in\R^{2n}, 
$$
be such that $m(\pih,\mu_j)>0$, for $j=1,\dots,n$. 
Define the $(n\times n)$-matrix valued function ${\e_\mu}$ on $\R$ 
as 
$$
\e_\mu(\tau)=\diag(e_{\mu_1}(\tau),\dots,e_{\mu_n}(\tau)).
$$
Define the weighted Sobolev spaces  
\begin{equation*}
\sblv_{k,\mu} =\left\{
u\in\sblv_{k}^{\rm loc}\colon \e_\mu u\in\sblv_k
\right\},\text{ with norm }\|u\|_{k,\mu}=\|\e_\mu u\|_k.
\end{equation*}
To make the doubling operation used in Section \ref{kercokelts} work on
$\sblv_2$, we  impose further boundary conditions. If
$u\in\sblv_1^{\rm loc}$ then its trace lies in
$\sblv_{\frac12}^{\rm loc}(\R\cup\R+i,\C^n)$. We say that $u$ {\em
vanishes on the boundary if}
\begin{equation}\label{rabdr}
\int_{\R\cup\R+i}\la u,v\ra\,d \tau =0\text{ for every
}v\in\C^\infty_0(\R\cup\R+i,\C^n). 
\end{equation} 
Define
\begin{align*}
\sblv_{2,\mu}(\underbrace{\pih,\dots,\pih}_n) &= \left\{
u\in\sblv_{2,\mu}\colon
u \text{ satisfies \eqref{rabdr0}, \eqref{rabdr1}, and $\bar\pa u$
satisfies \eqref{rabdr}}
\right\},\\
\sblv_{1,\mu}[0] &= \left\{
u\in\sblv_{1,\mu}\colon
u \text{ satisfies \eqref{rabdr}}
\right\}.
\end{align*}

\begin{prp}\label{modelwght}
If $m(\pih, \mu_j)>0$ for $j=1,\dots,n$ then the operator 
$$
\bar\pa\colon\sblv_{2,\mu}(\pih,\dots,\pih)\to
\sblv_{1,\mu}[0] 
$$
is Fredholm with index 
$$
\sum_{j=1}^n \sharp\left(-\tfrac{\mu_j^-}{\pi}-\tfrac12,
-\tfrac{\mu_j^+}{\pi}-\tfrac12\right)
-\sharp\left(
\tfrac{\mu_j^-}{\pi}-\tfrac12,\tfrac{\mu_j^+}{\pi}-\tfrac12\right)
$$
where $\sharp(a,b)$ denotes the number of integers in the interval
$(a,b)$. 

Moreover, if $\mu_j^+=\mu_j^-$ for all $j$ and
$M(\mu)=\min\{m(\pih, \mu_1),\dots,m(\pih, \mu_n)\}$ then
$u\in\sblv_{2,\mu}(\pih,\dots,\pih)$ satisfies  
\begin{equation}\label{blowup1}
\|u\|_{2,\mu}\le C(\mu)\|\bar\pa u\|_{1,\mu},
\end{equation} 
where $C(\mu)\le\frac{K}{M(\mu)}$, for some constant $K$.
\end{prp}

\begin{pf}
The problem studied is split and it is clearly sufficient to consider
the case $n=1$. We first determine the dimensions of the kernel and
cokernel. It is immediate from Lemma
\ref{lmaralocsol} 
that the kernel of $\bar\pa$ is finite dimensional on
$\sblv_{2,\mu}(\pih)$ and that the number of linearly independent
solutions is exactly $\sharp\left(-\frac{\mu^-}{\pi}-\frac12,
-\frac{\mu^+}{\pi}-\frac12\right)$. 

Recall that an element in the cokernel of $\bar\pa$ is an element $\xi$ in the dual
space of $\sblv_{1,\mu}[0]$. 
The dual of $\sblv_{1,\mu}$ is
$\sblv_{-1,-\mu}$ and thus, as in (\ref{7coker.eqn}), 
the dual of $\sblv_{1,\mu}[0]$ is the quotient space 
\begin{equation} \notag
\sblv_{-1,-\mu}/\sblv_{1,\mu}[0]^\perp.
\end{equation}
Lemma \ref{lmacokerreg} implies that any element in the
cokernel has a smooth representative. Let $v$ be a smooth
representative. 
Then 
\begin{equation} \notag
\int_{\R\times[0,1]}\la\bar\pa u, v\ra \,d\tau\wedge dt=0,
\end{equation} 
for any smooth compactly supported function $u$ which meets the
boundary conditions \eqref{rabdr0}, \eqref{rabdr1}, and \eqref{rabdr}.
Using partial integration we conclude  
\begin{equation}\label{intreg}
\int_{\R\times[0,1]}\la u,\pa v\ra \,d\tau\wedge dt=0.
\end{equation}
Thus $\pa v=0$ in the interior. Noting that for any two functions
$\phi_0\in C^\infty_0(\R,\R)$ and $\phi_1\in C^\infty_0(\R,i\R)$ there 
exists a function $u\in\C^\infty_0(\R\times[0,1],\C)$ such that 
$\bar\pa u|\pa(\R\times[0,1])=0$, $u|\R=\phi_0$, and $u|\R+i=\phi_1$
we find that $iv$ satisfies \eqref{rabdr0} and \eqref{rabdr1}. Lemma
\ref{lmaralocsol} then implies that the cokernel has dimension 
$\sharp\left(\frac{\mu_j^-}{\pi}-\frac12,\frac{\mu_j^+}{\pi}-\frac12\right)$. 

We now prove that the image of $\bar\pa$ is closed, and in doing
so also establish \eqref{blowup1}. Let 
$$
A(\tau)=\exp\left(\int_0^\tau\alpha(\sigma)\,d\sigma\right).
$$ 
Then multiplication with $A$ defines a Banach space 
isomorphism $A\colon\sblv_{k,\mu}\to\sblv_k$. 
The inverse $A^{-1}$ of $A$ is multiplication with $A(\tau)^{-1}$.
These isomorphisms gives the following commutative diagram
$$
\begin{CD}
\sblv_{2,\mu}(\pih) @<{A^{-1}}<< \sblv_2(\pih^\ast)\\
@V{\bar\pa}VV @VV{\bar\pa-\alpha}V\\
\sblv_{1,\mu}[0] @>>A> \sblv_1[0],
\end{CD}
$$
where $\sblv_2(\pih^\ast)$ is defined as   
$\sblv_2(\pih)$ except that instead of requiring that $\bar\pa u$
vanishes on the boundary we require that $(\bar\pa-\alpha)u$ does.
We prove that the operator $\bar\pa-\alpha$ on the right in the above
diagram has closed range and conclude the corresponding statement for
the operator on the left. Note that if $u\in\sblv_2(\pih^\ast)$ then
both $\pa_\tau u$ and $\pa_t u$ satisfy \eqref{rabdr0} and
\eqref{rabdr1}. Hence the doubling operation described in Section
\ref{kercokelts} induces a map $\sblv_2(\pih^\ast)\to\sblv_2(\R\times [0,2])$
with $\|u^d\|_2=2\|u\|_2$.

Let 
\begin{equation}\label{Smu1}
S(\mu)=\{n\in\Z\colon -\tfrac{\mu^-}{\pi}-\hf<n<-\tfrac{\mu^+}{\pi}-\hf\}.
\end{equation}
(Note that $S(\mu)=\emptyset$ if $\mu^+\ge\mu^-$.) The map
$\gamma_n\colon\sblv_{2}(\pih^\ast)\to\sblv_2(\R,\R)$, 
\begin{equation}\label{gamma_n}
u\mapsto c_n(\tau)=
\int_0^2 u^d(\tau,t)\exp\left(-i(\pih+n\pi)t\right)\, dt
\end{equation}
is continuous. Let $W_2\subset\sblv_{2}(\pih^\ast)$ be the closed subspace
\begin{equation}\label{defW}
W_2=\bigcap_{n\in S(\mu)}\ker(\gamma_n).
\end{equation}
Using the Fourier expansion of $u^d$ we see that $W_2$ has a
direct complement  
\begin{equation}\label{defV}
V_2=\bigcap_{n\notin S(\mu)}\ker(\gamma_n).
\end{equation}
(Note that if $\mu^+>\mu^-$ then $W_2=\sblv_2(\pih^\ast)$ and $V_2=\emptyset$.)

Similarly, we view the maps $\gamma_n$ defined by $\eqref{gamma_n}$ as
maps $\sblv_1[0]\to\sblv_1(\R,\R)$ and get the corresponding direct
sum decomposition $\sblv_1[0]=W_1\oplus V_1$. If
$u\in\sblv_2(\pih^\ast)$ then the Fourier expansion of $u^d$ is
\begin{equation} \notag
u^{d}(\tau+it)=\sum_n c_n(\tau)e^{i(\pih+n\pi)t}.
\end{equation}
Hence 
\begin{equation}\label{barpau^d}
(\bar\pa-\alpha)u^d(\tau+it)=
\sum_n (c_n'(\tau)-(\alpha(\tau)+\pih+n))e^{i(\pih+n\pi)t}.
\end{equation}
It follows that $\bar\pa (W_2)\subset W_1$ 
and $\bar\pa (V_2)\subset V_1$.

Let $w\in W_2$. Fourier expansion of $w^d$ gives 
\begin{align}\notag
&2\|(\bar\pa-\alpha)w\|_0^2=\\ \notag
&=\sum_{n\notin S(\mu)}
\int_\R \left(|c'_n|^2+\left((\pih+n\pi+\alpha(\tau))^2+\alpha'\right)
|c_n|^2\right)\,d\tau\\\label{Fexpans}  
&\ge 2C\|w\|_1^2,
\end{align} 
where the constant $C$ is obtained as follows.  
If $\mu^+>\mu^-$ then {\bf P2} implies that the coefficients of
$|c_n|^2$ are strictly positive, and 
if $\mu^+<\mu^-$ then  
{\bf P3} implies that the coefficients in front of  
$|c_n|^2$ are larger than
$\frac{4}{5}m(\pih,\mu)^2$ 
since $n\notin S(\mu)$. Finally, if $\mu^-=\mu^+$
then $\alpha'=0$ and the coefficients in front of $|c_n|^2$ are larger
than $m(\pih, \mu)^2$ for all $n.$ 

If $w\in\sblv_2(\pih^\ast),$ then $\pa_\tau w$ and
$i\pa_t w$ satisfies \eqref{rabdr0} and \eqref{rabdr1} and the Fourier
coefficients $c_n(\tau)$ of their doubles vanish for 
$n\in S(\mu)$. Thus, the same argument applies to these functions and 
the following estimates are obtained
\begin{align*}
\|(\bar\pa-\alpha)\pa_\tau w\|_0  &\ge  C\|\pa_\tau w\|_1,\\
\|(\bar\pa-\alpha)\pa_t w\|_0  &\ge  C\|\pa_t w\|_1.
\end{align*} 
If $\mu_+=\mu_-$ then $\alpha'=0$ and  $\bar\pa-\alpha$ commutes with
both $\pa_t$ and $\pa_\tau$. Hence
\begin{align*}
\|(\bar\pa-\alpha) w\|_1&\ge\frac{1}{2}(
\|(\bar\pa-\alpha) w\|_0+\|\pa_\tau(\bar\pa-\alpha) w\|_0+
\|\pa_t(\bar\pa-\alpha) w\|_0)\\
&\ge C(\|w\|_1+\|\pa_\tau w\|_1+\|\pa_t w\|_1)\ge C\|w\|_2,
\end{align*}
where $C=Km(\mu,\pih)$. This proves \eqref{blowup1}.

If $\mu_+\ne\mu_-$ then 
$\pa_\tau(\bar\pa-\alpha)w=(\bar\pa-\alpha)\pa_\tau w-\alpha' w$,
and with $K>0$ we conclude from the triangle inequality
\begin{align*}
&K\|(\bar\pa-\alpha) w\|_0+\|\pa_\tau(\bar\pa-\alpha) w\|_0+
\|\pa_t(\bar\pa-\alpha) w\|_0\\
&\ge KC\|w\|_1+C\|\pa_\tau w\|_1-\|\alpha' w\|_0+C\|\pa_t w\|_1\\
&\ge\left(KC-\frac{m(\pih,\mu)^2}{5}\right)
\|w\|_1+C\|\pa_\tau w\|_1+C\|\pa_t w\|_1,
\end{align*}
since $|\alpha'|<\frac{m(\pih,\mu)^2}{5}$. Thus choosing $K$ sufficiently
large we find that there exist a constant $K_1$ such that for $w\in W$
\begin{equation}\label{estW}
\|w\|_2\le K_1\|(\bar\pa-\alpha) w\|_1. 
\end{equation}
Thus, if $\mu^+>\mu^-$ we conclude that the range of $\bar\pa-\alpha$
is closed. If $\mu^+<\mu^-$ we need to consider also $V_2$.

For $v\in V_2$ we have
$$
v^d(\tau,t)=\sum_{n\in S(\mu)} c_n(\tau)\exp\left(i(\pih+n\pi)t\right).
$$
Let $V^\perp_2$ be the space of functions in $V_2$ which, under doubling,
map to the orthogonal complement of the doubles $\phi_n^d$ of the
functions 
$\phi_n(\zeta)= \exp((\pih+n\pi)\zeta+\int\alpha\,d\tau)$, $n\in S(\mu)$
with respect to the $L^2$-pairing
on $\sblv_2(\R\times[0,2],\C)$. Then $V^\perp_2$ is a closed subspace of
finite codimension in $V_2$.
 
We claim there exists a constant $K_2$ such that for all
$v^\perp\in V^\perp_2$ 
\begin{equation}\label{estV}
\|v^\perp\|_2\le K_2\|(\bar\pa-\alpha) v^\perp\|_1.
\end{equation}
Assume that this is not the case. Then there exists a sequence $v^\perp_j$
of elements in $V^\perp_2$ such that 
\begin{align}\label{bd}
&\|v^\perp_j\|_2=1, \\ 
&\|(\bar\pa-\alpha) v^\perp_j\|_1\to 0. \label{sol}
\end{align}

Let $P>M$ be an integer (see condition {\bf P1}) and let 
$v^\perp\in V^\perp$. Consider the 
restriction of $v^\perp$ and $\bar\pa v^\perp$ to
$\Theta_P=\{\tau+it\colon |\tau|\ge P\}$. 
Using Fourier expansion as in \eqref{Fexpans}, partial integration, 
and the fact that $\alpha'(\tau)=0$ for $|\tau|>M$ we find
\begin{align}\notag
&2\|(\bar\pa-\alpha)v^\perp|\Theta_P\|_1\ge\\ \label{frstV}
& C\left(\|v^\perp|\Theta_P\|_2+ 
\sum_{n\in S(\mu)} 
\mu^+(|c_n(P)|^2+|c_n'(P)|^2) -\mu^-(|c_n(-P)|^2+|c_n'(-P)|^2) \right).
\end{align} 

By a compact Sobolev embedding we find for each positive integer $P$ a
subsequence $\{v^\perp_{j(P)}\}$ which  
converges in $\sblv_1([-P,P]\times[0,1],\C)$. Moreover, we may assume
that these subsequences satisfies
$\{v^\perp_{j(P)}\}\supset\{v^\perp_{j(Q)}\}$ if $P<Q$.  

Let $(c_n)_j$ be the sequence of Fourier coefficient functions
associated to the sequence $v^\perp_j$. The estimates
\begin{equation}\label{ODEest}
\|c\|_k\le
C(\|c\|_{k-1}+\|(\frac{d}{d\tau}-(\pih+n\pi+\alpha))c\|_{k-1}),
\end{equation}
and \eqref{sol} implies that $(c_n)_{j(P)}$ converges to a smooth
solution of the equation $(\frac{d}{d\tau}-(\pih+n\pi+\alpha))c=0$ on
$ [-P,P]$. Hence, 
$v^\perp_{j(P)}$ converges to a smooth solution of 
$(\bar\pa-\alpha)u=0$ on $\Theta_P$ satisfying the boundary conditions
\eqref{rabdr0} and \eqref{rabdr1}. Such a solution has the form 
$$
\sum_{n\in S(\mu)}k_n\phi_n(\zeta),
$$ 
where $k_n$ are real constants.

We next show that in fact all $k_n$ must be zero. Note that by
Morrey's theorem and \eqref{bd} we get a uniform $C^0$-bound 
$|v_j^\perp|\le K$. Therefore, $|(c_n)_j|\le 2K$ and hence
\begin{align*}
&\int\la (v^\perp_j)^d,(\phi_n)^d\ra\, d\tau\wedge dt=\\
&\int_\R (c_n)_j\exp\left((\pih+n\pi)\tau+
\int\alpha\,d\tau\right)\,d\tau= 
\int_{-P}^P (c_n)_j\exp\left((\pih+n\pi)\tau+
\int\alpha\,d\tau\right)\,d\tau\\
&+\int_{P}^\infty (c_n)_j\exp((\pih+n\pi+\mu^+)\tau)\,d\tau+
\int_{-\infty}^{-P} (c_n)_j\exp((\pih+n\pi+\mu^+)\tau)\,d\tau.
\end{align*} 
But 
\begin{align*}
&\left|\int_{P}^\infty (c_n)_j\exp((\pih+n\pi+\mu^+)\tau)\,d\tau\right|+
\left|\int_{-\infty}^{-P}
(c_n)_j\exp((\pih+n\pi+\mu^-)\tau)\,d\tau\right|\le \\
&\frac{2K}{m(\pih,\mu)}\left(\exp((\pih+n\pi+\mu^+)P)
+\exp(-(\pih+n\pi+\mu^-)P)\right)\to 0\text{ as }P\to\infty.
\end{align*}
We conclude from this that unless $k_n=0$, $v^\perp_{j(P)}$
violates the orthogonality conditions for $P$ and $j(P)$ sufficiently
large.  

Consider \eqref{frstV} applied to elements in the sequence
$\{v^\perp_j\}$. As $j\to\infty$ the term on the left hand side and
the sum in the right hand side tends to $0$. Hence
$\|v^\perp_{j}|\Theta_P\|_2\to 0$. Applying 
\eqref{ODEest} to $(c_n)_{j}$ and noting that both terms on the right
hand side goes to $0$ we conclude that also 
$\|v^\perp_{j}|[-P\times P]\times[0,1]\|_2\to 0$. This contradicts
\eqref{bd} and hence \eqref{estV} holds.

The estimates \eqref{estW} and \eqref{estV}
together with the direct sum decompositions   
$\sblv_2(\pih^\ast)=W_2\oplus V_2$ and $\sblv_1[0]=W_1\oplus V_1$, and
the fact that $\bar\pa-\alpha$ respects this decomposition 
shows that the image of $\bar\pa-\alpha$ is closed also in the case
$\mu^+<\mu^-$.  
\end{pf}

\begin{rmk}
In many cases, the first statement in Proposition \ref{modelwght}
still holds with weaker assumptions on the weight function than
{\bf P1}--{\bf P3}. For example, if $\mu_+<\mu_-$ then we need only
know that $\max\{\alpha',0\}$ is sufficiently small compared to 
$(\pih+n\pi+\alpha)^2$ for $n\notin S(\mu)$ to derive
\eqref{estW} and the derivation of 
\eqref{estV} is quite independent of $\alpha'$ as long as $\alpha$
eventually becomes constant.  
\end{rmk}

\subsection{The model problem with angles}
\label{7modelproblem2.section}

We study more general boundary conditions than those
in Section \ref{model}. 
Recall $(x_1+iy_1,\dots,x_n+iy_n)$ are coordinates on
$\C^n.$
Let $\pa_j$ denote the unit tangent vector in the
$x_j$-direction,  
for $j=1,\dots,n$. For $\theta=(\theta_1,\dots,\theta_n)\in [0,\pi)^n$,
let $\Lambda(\theta)$ be the Lagrangian subspace of $\C^n$ spanned by
the vectors $e^{i\theta_1}\pa_1,\dots,e^{i\theta_n}\pa_n$. Consider
the following boundary conditions for $u\in\sblv_1^{\rm loc}$.
\begin{align}
\int_\R\la u,v\ra\,d\tau=0 & \label{nglbdr0}
\text{ for all } v\in C^\infty_0(\R,i\R^n),\\
\int_{\R+i}\la u,v\ra\,d\tau=0 & \label{nglbdr1}
\text{ for all } v\in C^\infty_0(\R+i,i\Lambda(\theta)). 
\end{align}
If $m(\theta_j,\mu_j)>0$ (see \eqref{m(a)}) for all $j$ then define
\begin{align*}
\sblv_{2,\mu}(\theta) &= \left\{
u\in \sblv_{2,\mu}\colon
u \text{ satisfies \eqref{nglbdr0}, \eqref{nglbdr1}, and $\bar\pa u$
satisfies \eqref{rabdr}}
\right\},\\
\sblv_{1,\mu}[0] &= \left\{
u\in\sblv_{1,\mu}\colon
u \text{ satisfies \eqref{rabdr}}
\right\}.
\end{align*}

\begin{prp}\label{prpngl}
If $m(\theta_j,\mu_j)>0$ for $j=1,\dots,n$ then the operator 
$$
\bar\pa\colon
\sblv_{2,\mu}(\theta)\to\sblv_{1,\mu}[0]
$$
is Fredholm of index 
\begin{equation}\label{nglindex}
\sum_{j=1}^n \sharp\left(-\tfrac{\mu_j^-+\theta_j}{\pi},
-\tfrac{\mu_j^++\theta_j}{\pi}\right)
-\sharp\left(\tfrac{\mu_j^-+\theta_j}{\pi}-1,
\tfrac{\mu_j^++\theta_j}{\pi}-1\right).
\end{equation}
Moreover, if $\mu_j^+=\mu_j^-$ for all $j$ and
$M(\mu)=\min\{m(\mu_1,\theta_1),\dots,m(\mu_n,\theta_n)\}$ then
$u\in\sblv_{2,\mu}(\theta_1,\dots,\theta_n)$ satisfies  
\begin{equation}\label{blowup2}
\|u\|_{2,\mu}\le C(\mu)\|\bar\pa u\|_{1,\mu},
\end{equation} 
where $C(\mu)\le\frac{K}{M(\mu)}$, for some constant $K$.
\end{prp}
\begin{pf}
Consider the holomorphic $(n\times n)$-matrix
$$
\g_\theta(\zeta)=\diag\left((\exp(\pih-\theta_1)\zeta),\dots,
(\exp(\pih-\theta_n)\zeta)\right).
$$
Multiplication with $\g_\theta$ defines isomorphisms 
\begin{align*}
\sblv_{2,\mu}(\theta)&\to\sblv_{2,\lambda}(\pih,\dots,\pih) \quad \mbox{and}\\
\sblv_{1,\mu}[0]&\to\sblv_{1,\lambda}[0]
\end{align*} 
where $\lambda=(\lambda_1,\dots,\lambda_n)$ and
$\lambda_j^\pm=\mu_j^\pm-\pih+\theta_j$. Since $\g_\theta$ is
holomorphic it commutes with $\bar\pa$. The proposition now follows
from Proposition \ref{modelwght}.
\end{pf}

\subsection{Smooth perturbations of the model problem with angles}
\label{7modelproblem3.section}

Let $B\colon \R\times[0,1]\to\U(n)$ be a smooth map such that 
\begin{equation}\label{Bcond}
\bar\pa B|\pa \R\times[0,1]=0.
\end{equation}
Let $\theta\in[0,\pi)$ and consider the following boundary conditions
for $u\in\sblv_1^{\rm loc}$:
\begin{align}\notag
& \int_{\R}\la u,v\ra\,d\tau\wedge dt=0,\\ \label{pert0}
&\text{ for all }v\in
C^\infty_0(\R,\C^n) \text{ such that } v(\tau)\in iB(\tau)\R^n,\\ \notag
& \int_{\R+i}\la u,v\ra\,d\tau\wedge dt=0,\\ \label{pert1}
& \text{ for all }v\in 
C^\infty_0(\R+i,\C^n) \text{ such that } v(\tau+i)\in
iB(\tau)\Lambda(\theta). 
\end{align}

For $\mu=(\mu^+,\mu^-)\in\R^2$ let 
$\lambda(\mu)=(\mu^+,\mu^-,\mu^+,\mu^-,\dots,\mu^+,\mu^-)\in\R^{2n}$ define
$$
\sblv_{2,\mu}(\theta,B)=\{u\in\sblv_2^{\rm loc}\colon u \text{ satisfies
\eqref{pert0} and \eqref{pert1},}\bar\pa u\text{ satisfies
\eqref{rabdr} and } \e_{\lambda(\mu)} u\in\sblv_2 \}.   
$$
\begin{prp}\label{prppert}
If $m(\theta_j,\mu)>0$ for $j=1,\dots,n$ then there exists
$\delta>0$ such that for all $B$ satisfying \eqref{Bcond} with
$\|B-\id\|_{C^2}<\delta$, the operator 
$$
\bar\pa\colon
\sblv_{2,\mu}(\theta,B)\to\sblv_{1,\mu}[0]
$$
is Fredholm of index 
\begin{equation}
\sum_{j=1}^n \sharp\left(-\tfrac{\mu^-+\theta_j}{\pi},
-\tfrac{\mu^++\theta_j}{\pi}\right)
-\sharp\left(\tfrac{\mu^-+\theta_j}{\pi}-1,
\tfrac{\mu^++\theta_j}{\pi}-1\right).
\end{equation}
\end{prp}
\begin{pf}
Multiplication with $B$ and $B^{-1}$ defines Banach space isomorphisms 
$$
\sblv_{2,\mu}(\theta)\stackrel{\times
B}{\longrightarrow}\sblv_{2,\mu}(\theta,B), 
$$
and 
$$
\sblv_{1,\mu}[0]\stackrel{\times
B^{-1}}{\longrightarrow}\sblv_{1,\mu}[0]. 
$$
Thus up to conjugation the operator considered is the same as 
$$
\bar\pa+B^{-1}\bar\pa
B\colon\sblv_{2,\mu}(\theta)\to\sblv_{1,\mu}[0].
$$
The theorem now follows from Proposition \ref{prpngl}, and the fact that
the subspace of Fredholm operators is open and that the index is 
constant on path components of this subspace. 
\end{pf}

\subsection{Boundary conditions}
\label{7boundaryconditions.section}

In the upcoming subsections we study the linearized $\bar\pa$-problem on a disk
$D_m$ with $m$ punctures. Refer back to Section~\ref{2hodi.metric} for
notation concerning $D_m.$ 
%

\begin{dfn} \label{7atinfinity.dfn}
A smooth map $A\colon\pa D_m\to\U(n)$ will be called {\em small at
infinity} if there exists 
$M>1$ such that for each $j=1,\dots,m$ the restriction of $A$ to   
$\pa E_{p_j}[M]$ approaches a constant map in the
$C^2$-norm on each component of $\pa E_{p_j}[M']$ as $M'\to\infty$. 
It will be called {\em constant at infinity} if there exists $M>1$
such that for each $j=1,\dots,m$ the restriction of $A$ to each
component component of $\pa E_{p_j}[M]$ is constant. 
\end{dfn}

Let $A\colon\pa D_m\to\U(n)$ be small at infinity. 
For $u\in\sblv_1^{\rm loc}(D_m,\C^n)$, 
consider the boundary condition:   
\begin{align}\notag
&\int_{\pa D_m}\la u,v\ra\,ds=0,\\ \label{linwkbdr}
&\text{for all $v\in C^0_0(\pa D_m,\C^n)$ such that
$v(\zeta)\in iA(\zeta)\R^n$ for all $\zeta\in\pa D_m$}.
\end{align}

In previous subsections coordinates $\zeta=\tau+it$ on $\R\times[0,1]$
were used and we implicitly considered the bundle
${T^\ast}^{0,1}\R\times[0,1]$ as trivialized by the form $d\bar\zeta$, and
sections in this bundle as $\C^n$-valued functions. We do not want to
specify any trivialization of ${T^\ast}^{0,1}D_m$ and so we view the
$\bar\pa$-operator as a map from
$\sblv_2$-functions into $\sblv_1$-sections of 
${T^\ast}^{0,1}D_m\otimes\C^n$. Consider, for 
$u\in\sblv_1^{\rm loc}(D_m,{T^\ast}^{0,1}D_m\otimes\C^n)$, the 
boundary condition
\begin{equation}\label{barpa0}
\int_{\pa D_m}\la u,v\ra\,ds=0, 
\text{ for all $v\in C^0_0(\pa D_m,T^{0,1}D_m\otimes\C^n)$.}
\end{equation}

Henceforth, to simplify notation, if the source space $X$ in a
Sobolev space $\sblv_k(X,Y)$ is $D_m$ we will drop it from the notation.
If $u\in\sblv_2^{\rm loc}(\C^n)$ then 
$\bar\pa u\in\sblv_1^{\rm loc}({T^\ast}^{0,1}D_m\otimes\C^n)$. 
Define
$$
\sblv_2(\C^n;A)=\left\{u\in\sblv_2(\C^n)\colon 
u\text{ satisfies \eqref{linwkbdr} and }
\bar\pa u \text{ satisfies \eqref{barpa0}}\right\},
$$
and 
$$
\sblv_1({T^\ast}^{0,1}D_m\otimes\C^n;[0])=\left\{
u\in\sblv_1({T^\ast}^{0,1}D_m\otimes\C^n)\colon 
u \text{ satisfies \eqref{barpa0}}\right\},
$$

Define 
\begin{align*}
&\sblv_{2,\mu}(\C^n;A)=\\
&\left\{u\in\sblv_2^{\rm loc}(\C^n)\colon
u\text{ satisfies \eqref{linwkbdr}, }
\bar\pa u \text{ satisfies \eqref{barpa0}, and }
\e_\mu u\in\sblv_2(\C^n)\right\}.
\end{align*}
and
\begin{align*}
&\sblv_{1,\mu}({T^\ast}^{0,1}D_m\otimes\C^n;[0])=\\
&\left\{ u\in\sblv_1^{\rm loc}({T^\ast}^{0,1}D_m\otimes\C^n)\colon 
u \text{ satisfies \eqref{barpa0} and }
\e_\mu u\in\sblv_1({T^\ast}^{0,1}D_m\otimes\C^n)
 \right\}.
\end{align*}

Let $p_j$ be a puncture of $D_m$. The orientation of $D_m$ induces an
orientation of $\pa D_m$. Let $A^{0}_j$ and $A^{1}_j$ denote the
constant maps to which $A$ converges on the component of $\pa E_{p_j}$ 
close to $p_j$ corresponding to $\R$ and $\R+i$, respectively. Define
$$
\theta(j)=\theta(A^{0}_j\R^n,A^{1}_j\R^n).
$$ 
Then there are unique unitary complex coordinates 
$$
z(j)=(x(j)_1+iy(j)_1,\dots,x(j)_n+iy(j)_n)
$$
in $\C^n$ such that 
\begin{align*}
&A^0_j\R^n=\spa\left\la\pa(j)_1,\dots,\pa(j)_n\right\ra,\\
&A^1_j\R^n=\spa\left\la
e^{i\theta(j)_1}\pa(j)_1,\dots,e^{i\theta(j)_n}\pa(j)_n\right\ra.  
\end{align*}

\begin{prp}\label{linfred}
Let $A\colon\pa D_m\to\U(n)$ be small at infinity.
If $\mu$ satisfies $\mu_j\ne -\theta(j)_r+k\pi$ for $j=1,\dots,m$,
$r=1,\dots,n$, and every $k\in\Z,$ 
then the operator
\begin{equation}\label{riktig}
\bar\pa\colon\sblv_{2,\mu}(\C^n;A)\to
\sblv_{1,\mu}({T^\ast}^{0,1}D_m\otimes\C^n;[0])
\end{equation}
is Fredholm.
\end{prp}

\begin{pf}
Assume that for $M>0$, $A|\pa E_{p_j}[M-1]$ is sufficiently close to a
constant map (see Proposition \ref{prppert}). Choose smooth 
complex-valued functions $\alpha_0,\alpha_1,\dots,\alpha_m$ with the
following properties: $\alpha_j$ is constantly $1$ on $E_{p_j}[M+2]$;
the sum $\sum_j\alpha_j$ is close to the constant function $1$, 
$\bar\pa\alpha_j=0$ on $\pa D_m$; and $\alpha_j$ is constantly equal
to $0$ on $D_m-E_{p_j}[M+1]$, for $j=1,\dots,m$.  

Glue to each $E_{p_j}[M]$ a
half-infinite strip $(-\infty,M]\times[0,1]$ and denote the result
$\bar E_{p_j}$. Extend the boundary conditions from $E_{p_j}[M]$ to  
$\bar E_{p_j}$ keeping them close to constant.
Let the weight in
the weight function remain constant. Glue to $D_m-\bigcup_jE_{p_j}[M+2]$,
$m$ half disks and extend the boundary conditions smoothly. Denote the
result $\bar D_m$. Note 
that the boundary value problem on $\bar D_m$ is the
vector-Riemann-Hilbert problem, which is known to be Fredholm, and that
the weighted norm on this compact disk is equivalent to the standard
norm.
 
Now let $u\in\sblv_{2,\mu}(\C^n;A)$. Then $\alpha_j u$ is in the
appropriate Sobolev space for the extended boundary value problem on
$\bar E_{p_j}$ ($\bar D_m$ if $j=0$) and because the elliptic estimate
holds for all of these problems and since all of them except possibly
the one on $\bar D_m$ has no kernel, there exists a constant $C$ such
that 
\begin{align} \notag
\|u\|_{2,\mu}\le & 
\|\alpha_0 u\|_{2,\mu}+\sum_{j=1}^n\|\alpha_j u\|_{2,\mu}\\ \notag
\le & C\left(\|\alpha_0 u\|_{1,\mu}
+\sum_{j=0}^n\|\bar\pa(\alpha_j u)\|_{1,\mu}\right)\\ \label{linest}
\le & C\left(\sum_{j=0}^n\|\bar\pa\alpha_j u\|_{1,\mu}+ \|\alpha_0 u\|_{1,\mu}
+\sum_{j=0}^n\|\alpha_j\bar\pa u\|_{1,\mu}
\right). 
\end{align}
 
We shall show that \eqref{linest} implies that every bounded sequence
$u_r$ such that $\bar\pa u_r$ converges has a convergent
subsequence. This implies that $\bar\pa$ has a closed image and a finite
dimensional kernel (\cite{ho} Proposition 19.1.3). Clearly it
is sufficient to consider the case $\bar\pa u_r\to 0$. Consider the
restrictions of $u_r$ to a compact subset $K$ of $D_m$ such that
$$
\supp(\alpha_0)\cup\supp(\bar\pa\alpha_0)\cup
\dots\cup\supp(\bar\pa\alpha_m)\subset K. 
$$
A compact Sobolev embedding argument gives a subsequence $\{u_{r'}\}$
which converges in $\sblv_{1}(K,\C^n)$. Thus,
\eqref{linest} implies that $\{u_{r'}\}$ is a Cauchy sequence in
$\sblv_{2,\mu}(A;\C^n)$ and hence it converges.  

It remains to prove that the cokernel is finite dimensional. Lemma
\ref{lmacokerreg} shows that any element in the 
cokernel of $\bar\pa$ can be represented by a smooth function $v$ on
$D_m$. Partial integration implies this function  satisfies $\pa v=0$
with boundary conditions given by the matrix function $iA$. Assume
first that $A$ is constant at 
infinity. Then, Lemma \ref{lmaralocsol} and conjugation with the
holomorphic $(n\times n)$-matrix $\g_\theta$ as in the proof of
Proposition \ref{prpngl}  
gives  explicit formulas for the restrictions of these smooth
functions to $E_{p_j}[M]$, for each $j$. It is straightforward to check
from these local formulas that $v$ lies in
$\sblv_{2,-\mu}(\C^n,iA)$. Thus, 
repeating the argument above with $\pa$ replacing $\bar\pa$ shows that
the cokernel is finite dimensional. The lemma follows in the case when
$A$ is constant at infinity. The general case then follows by an
approximation argument as in the proof of Proposition \ref{prppert}. 
\end{pf}

\subsection{Index-preserving deformations}
\label{7deform.section}

We compute the index of the operator in
\eqref{riktig}. Using approximations it is easy to see that it is sufficient to
consider the case when $A\colon\pa D_m\to\U(n)$ is constant 
at infinity. Thus, let $A$ be such a map which is constant on $\pa
E_{p_j}[M]$ for every $j$ and consider the Fredholm operator 
\begin{equation}\label{calcoprtr}
\bar\pa\colon\sblv_{2,\mu}(\C^n;A)\to
\sblv_{1,\mu}({T^\ast}^{0,1}D_m\otimes\C^n;[0]),
\end{equation}
where $\mu=(\mu_1,\dots,\mu_m)\in\R^{m}$ satisfies
\begin{equation}\label{wghtcond}
\mu_j\ne-\theta(j)_r+ n\pi \text{ for every } j,r,n.  
\end{equation}

\begin{lma}\label{lmadef}
Let $B_s\colon D_m\to\U(n)$, $s\in[0,1]$, be a continuous family of
smooth maps such that 
\begin{align}\notag
&B_s \text{ is bounded in the $C^2$-norm,}\\ \notag
&B_s|\pa E_{p_j}[M] \text{ is constant in } \tau +it,\\ \notag 
&\bar\pa B_s|\pa D_m=0,\text{ and}\\ \label{defcond}
&B_0\equiv\id.
\end{align}
Let $\lambda\colon[0,1]\to\R^m$ be a continuous map such that
$\lambda(0)=\mu$ and $\lambda(s)$ satisfies
\eqref{wghtcond} for every $s\in[0,1]$.
Then the operator
\begin{equation*}
\bar\pa\colon\sblv_{2,\lambda(1)}(\C^n;B_1A)\to
\sblv_{1,\lambda(1)}({T^\ast}^{0,1}D_m\otimes\C^n;[0])
\end{equation*}
has the same Fredholm index as the operator in \eqref{calcoprtr}.
\end{lma}

\begin{pf}
The Fredholm operator 
\begin{equation*}
\bar\pa\colon\sblv_{2,\lambda(s)}(\C^n;B_sA)\to
\sblv_{1,\lambda(s)}({T^\ast}^{0,1}D_m\otimes\C^n;[0])
\end{equation*}
is conjugate to 
\begin{equation*}
\bar\pa-B_s\bar\pa B_s^{-1}\colon\sblv_{2,\mu}(\C^n;A)\to
\sblv_{1,\mu}({T^\ast}^{0,1}D_m\otimes\C^n;[0]).
\end{equation*}
The family $\bar\pa-B_s\bar\pa B_s^{-1}$ is then a continuous family of 
Fredholm operators.
\end{pf}

In order to apply Lemma \ref{lmadef} we shall show how to deform given
weights and boundary conditions into other boundary conditions and
weights keeping the Fredholm index constant using the conditions in
Lemma \ref{lmadef}. We accomplish this in two steps: first deform the
problem so that the boundary value matrix is diagonal; then change
the weights and angles at the ends into a special form where
compactification is possible. 

\begin{lma}\label{def1}
Let $A\colon\pa D_m\to\U(n)$ be constant at infinity. 
Then there exists a continuous family 
$B_s\colon D_m\to\U(n)$, $0\le s\le 1$, 
of maps satisfying \eqref{defcond} such that
$$ 
B_1(\zeta)A(\zeta)=\diag(b_1(\zeta),\dots,b_n(\zeta)), \zeta\in\pa D_m.
$$
\end{lma}

\begin{pf}
We first make $A$ diagonal on the ends where it is constant. Note
that in canonical coordinates $z(j)$ on the end $E_{p_j}[M]$ the
matrix $A$ is diagonal. Let $B_j\in\U(n)$ be the matrix which
transforms the complex basis $\pa(j)_1,\dots,\pa(j)_n$ to the standard 
basis. Let $B_j(s)$ be a smooth path in $\U(n)$, starting at $\id$ and
ending at $B_j$. Define $B_s=B_j(s)$ on $E_{p_j}[M]$ for each $j$. 

We
need to extend this map to all of $D_m$.
To this end consider the loop on the boundary of $S=D_m- E_{p_j}[M]$. There
exists a $1$-parameter family of functions $B_s\colon S\to \U(n)$ such 
that $B_0=\id$ and $B_1 A$ is diagonal, since any loop is homotopic to 
a loop of diagonal matrices. The loops $B_s$ can be smoothly extended
to all of $D_m$ 

Finally, we need that $\bar\pa B_s=0$ on the
boundary. We get this as follows: let $C$ be a collar on the boundary
with coordinates $\tau$ along the boundary and $t$ orthogonal to the
boundary, $0\le t\le\epsilon$ and let $\phi\colon[0,\epsilon]\to\R$ be 
a smooth function which equals the identity on
$[0,\frac{\epsilon}{4}]$ and $0$ for $t\ge
\frac{\epsilon}{2}$. Redefine $B_s$ on the collar as
$$
\tilde B_s=B_s(\zeta)\exp(i\phi(t)B^{-1}_s(\zeta)\bar\pa B_s(\zeta)).
$$   
Then $\tilde B_s$ satisfies the boundary conditions and equals $B_s$ on
the boundary and in the complement of the collar.
 
Consider the loop on the boundary of $S=D_m-\pa E_{p_j}[M]$. There
exists a $1$-parameter family of functions $B_s\colon S\to \U(n)$ such 
that $B_0=\id$ and $B_1 A$ is diagonal, since any loop is homotopic to 
a loop of diagonal matrices. The loops $B_s$ can be smoothly extended
to all of $D_m$ and the above trick makes $B_s$ satisfy the boundary
conditions.
\end{pf}

Now let $A\colon D_m\to\U(n)$ take values in diagonal matrices.
Assume that $A$ is constant near the punctures and that 
$\mu=(\mu_1,\dots,\mu_m)\in\R^m$ satisfies \eqref{wghtcond}.

\begin{lma}\label{def2}
There are continuous families of smooth maps
$B_s\colon D_m\to\diag\subset\U(n)$ and $\lambda\colon[0,1]\to\R$
which satisfy \eqref{defcond} and \eqref{wghtcond} (where
the $\theta(j)$ are computed w.r.t. $B_s$)
respectively such that
$$
B_1A=\id
$$
in a neighborhood of each puncture.
\end{lma}

\begin{pf}
Let $M>0$ be such that $A$ is constant in $E_{p_j}[M]$ for each
$j$. Let $\phi\colon[0,1]\to[0,1]$ be an approximation of the identity
which is constant near the endpoints of the interval. Let
$\psi\colon[M,\infty)\to[0,1]$ be a smooth increasing function which is
identically $0$ on $[M,M+1]$ and identically $1$ on $[M+2,\infty)$.
For $\alpha=(\alpha_1,\dots,\alpha_m)\in(-\pi,\pi)^m$ let 
$$
\tilde g_\alpha(\zeta)=
\begin{cases}
1 &\text{ for }\zeta\in D_m-\bigcup_j E_{p_j}[M],\\
e^{i\psi(\tau)\alpha_j\phi(t)} &
\text{ for }\zeta=\tau+it\in E_{p_j}[M],
\end{cases}
$$
and let $g_\alpha$ be a function which agrees with $\tilde g_\alpha$
except on $E_{p_j}[M]-E_{p_j}[M+2]$ and which satisfies $\bar\pa
g_\alpha|\pa D_m=0$. 

Consider the complex angle $\theta(j)\in[0,\pi)^n$ and the weight
$\mu_j$. Assume first that $\mu_j\ne k\pi$ for all $k\in\Z$ and $j=1,\dots,m$. 
Let $m_j$ be the unique number $0\le m_j\le\pi$ such that
$m_j=k\pi-\mu_j$ for some $k\in\Z$. By
\eqref{wghtcond} $\theta(j)_r\ne m_j$ for all $r$. If
$\theta(j)_r>m_j$ define $\alpha_r=\pi-\theta(j)_r$, and if
$\theta(j)_r<m^\ast_j$ define $\alpha(j)_r=-\theta(j)_r$.
Define
$$
B_s=\diag\left(g_{s\alpha_1},\dots,g_{s\alpha_n}\right)
$$
and let $\lambda(s)\equiv\mu$.

Assume now that $\mu_j=k\pi$ for some $j$. For $0\le s\le \frac12$,
let $B_s=\id$ and take
$\lambda_j=\mu_j-\epsilon s$  for some sufficiently small
$\epsilon>0$. Repeat the above construction to construct $B_s$ for
$s\le\frac12\le 1$.   
\end{pf}

\subsection{The Fredholm index of the standardized problem}
\label{7RiemannHilbert.section}

Consider $D_m$ with $m$ punctures on the boundary, conformal structure
$\kappa$ and metric $g(\kappa)$ as above and neighborhoods $E_{p_j}$ of the
punctures $p_1,\dots,p_m$. 

Let $\Delta_m$ denote the representative of the conformal structure
$\kappa$ on $D_m$ which is the unit disk in $\C$ with $m$ punctures at
$1,i,-1,q_3,\dots,q_m$ with the flat metric. 
Then there exists a
conformal and therefore holomorphic map $\Gamma\colon D_m\to\Delta_m$. We
study the behavior of $\Gamma$ on $E_{p_j}$. Let $p=p_j$ and let $q$ be
the puncture on $\Delta_m$ to which $p$ maps. After translation and
rotation in $\C$ we may assume that the point $q=0$
and that $\Delta_m$ is the disk of radius $1$ centered at $i$. We may
then find a holomorphic function on a neighborhood $U\subset\Delta_m$
of $q=0$ which fixes $0$ and maps $\pa\Delta_m\cap U$ to the real
line. Composing with this map we find that $\Gamma$ maps $\infty$ to
$0$, $\tau+0i$ to the negative real axis and $\tau+i$ to the positive
real axis for $\tau>M$ for some $M$. Thus this composition equals
$C\exp(-\pi\zeta)$ where $C<0$ is some negative real constant.       
Thus, up to a bounded holomorphic change of coordinates on a
neighborhood of $q$ the map $\Gamma$ on $E_{p_j}$ looks like
$\Gamma(\zeta)=\exp(-\pi\zeta)$ and its inverse $\Gamma^{-1}$ in these
coordinates satisfies $\Gamma^{-1}(z)=-\tfrac{1}{\pi}\log(z)$.    

Let $A\colon\pa D_m\to\diag\subset\U(n)$ be a smooth function which is
constantly equal to $\id$ close to each puncture. We may now think of
$A$ as being defined on $\pa \Delta_m$. We extend $A$ smoothly to $\pa\Delta$
by defining its extension $\hat A$ at the punctures as $\hat
A(p_j)=\id$ for each $j$.

Consider the following boundary condition for $u\in\sblv_2(\Delta,\C^n)$:
\begin{equation}\label{RH1}
\int_{\pa\Delta}\la u,v\ra\, ds=0\text{ for all }v\in 
C^\infty_0(\pa\Delta,\C^n) \text{ with } v(z)\in i\hat A(z)\R^n \text{
for all } z\in\pa\Delta.
\end{equation}
For $u\in\sblv_1(\Delta, {T^\ast}^{0,1}\Delta\otimes\C^n)$
consider the boundary conditions
\begin{equation}\label{RH2}
\int_{\pa\Delta}\la u,v\ra\, ds=0\text{ for all }v\in 
C^\infty_0(\pa\Delta,T^{0,1}\Delta\otimes\C^n).
\end{equation}
Define
\begin{align} \notag
\sblv_2(\Delta,\C^n;\hat A)&=\left\{
u\in\sblv_2(\Delta,\C)\colon u \text{ satisfies \eqref{RH1} and }
\bar\pa u \text{ satisfies \eqref{RH2}}
\right\}\\
\notag
\sblv_1(\Delta,{T^\ast}^{0,1}\Delta;[0])&=
\left\{
u\in\sblv_2(\Delta,\C)\colon u \text{ satisfies \eqref{RH2}}
\right\}.
\end{align}

\begin{lma}\label{R-H}
The operator 
$$
\bar\pa\colon\sblv_2(\Delta, \C^n;\hat A)\to
\sblv_1(\Delta,{T^\ast}^{0,1}\Delta;[0])
$$
is Fredholm of index $n+\mu(\hat A)$, where $\mu(\hat A)$ denotes the
Maslov index 
of the loop $z\mapsto A(z)\R^n$, $z\in\pa\Delta$, of Lagrangian
subspaces in $\C^n$.     
\end{lma}

\begin{pf}
This is (a direct sum of) classical Riemann-Hilbert problems.
\end{pf}

Let $\lambda(a)=(a,\dots,a)\in\R^m$.
\begin{prp}\label{R-Hcompactification}
For $-\pi<a<0$ the Fredholm index of the operator
$$
\bar\pa\colon\sblv_{2,\lambda(a)}(\C^n;A)\to
\sblv_{1,\lambda(a)}({T^\ast}^{0,1}D_m\otimes\C^n;[0])  
$$
equals $n+\mu(A)$.
\end{prp}

\begin{pf}
The holomorphic map $\Gamma\colon D_m\to\Delta_m$ and its holomorphic
inverse commute with the 
$\bar\pa$ operator. 
Any solution on $D_m$ must look like 
$\sum_{n\le 0} c_n e^{\pi n\zeta}$ (the negative weights allows for
$c_0\ne 0$) in canonical coordinates close to
each puncture. Thus $\Gamma^{-1}$ pulls back solutions on $D_m$ to
solutions on $\Delta_m$. Using also $\Gamma$ we see that the kernels
are isomorphic.

Elements in the cokernel on $D_m$ are of the form
$(\sum_{n<0}c_ne^{n\pi\bar\zeta})d\bar\zeta$ (the positive weight
implies $c_0=0$). Pulling back with $\Gamma^{-1}$ gives elements of the
form $(\sum_{n>0} \bar z^n)\tfrac{d\bar z}{\bar z}$ which are in the
cokernel of the $\bar\pa$ on $\Delta$. So the cokernels are also isomorphic.
\end{pf}

\subsection{The index of the linearized problem}
\label{7indexlinear.section}

In this subsection we determine the Fredholm indices of the problems
which are important in our applications to contact geometry.

Let $A\colon D_m\to \U(n))$ be a map which is small at
infinity.  Assume that $A_j^0\R^n$ and $A_j^1\R^n$ are transverse for
all $j$. For $0\le s\le 1$, let $\f_j(s)\in\U(n)$ be the matrix which
in the canonical coordinates $z(j)$ is represented by the matrix
$$
\diag(e^{-i(\pi-\theta(j)_1)s},\dots,e^{-i(\pi-\theta(j)_n)s}).
$$  
If $p$ and $q$ are consecutive punctures on $\pa D_m$
then let $I(a,b)$ denote the (oriented) path in $\pa D_m$ which
connects them.
Define the loop $\Gamma_A$ of Lagrangian subspaces in $\C^n$ by
letting the loop
$$
\left(A|I(p_1,p_2)\right)\ast\f_2\ast\left(A|I(p_2,p_3)\right)
\ast\f_3\ast\dots\ast\left(A|I(p_m,p_1)\right)\ast\f_1
$$
of elements of $\U(n)$ act on $\R^n\subset\C^n$.
\begin{prp}\label{indexw/oJtrick}
For $A$ as above the index of the operator
$$
\bar\pa\colon\sblv_{2}(\C^n;A)\to\sblv_1({T^\ast}^{0,1}D_m\otimes\C^n;[0])
$$
equals $n+\mu(\Gamma_A)$ where $\mu$ is the Maslov index.
\end{prp}

\begin{pf}
Using Lemmas \ref{def1} and \ref{def2} we deform $A$ to put the problem into
standardized form with weight $-\epsilon$ at each corner, without
changing the index. Call the new matrix $B$. We need to consider how
$B$ is constructed from $A$. The key step to understand is the point
where we make $B$ equal the identity on the ends. This is achieved
by first introducing a small {\em negative} weight and then rotating the 
space 
$$
A_j^1\R^n=\spa\left\la
e^{i\theta(j)_1}\pa_1,\dots,e^{i\theta(j)_n}\pa_n
\right\ra
$$ 
to $A_j^0$ according to 
\begin{equation}\label{rotation}
\spa\left\la e^{i(\theta(j)_1+s\phi(\tau)(\pi-\theta(j)_1)}\pa_1,\dots,
e^{i(\theta(j)_n+s\phi(\tau)(\pi-\theta(j)_n)}\pa_n\right\ra,
\end{equation}
where $0\le s\le 1$ and $\phi\colon[M,\infty)\to [0,1]$ equals $1$ on
$[M+2,\infty)$ and $0$ on $[M,M+1]$.

We now calculate the Maslov-index $\mu(B)$. 
Since as we follow $\R+i$ along the negative $\tau$-direction from
$M+1$ to $M$, $B$ experiences the inverse of the rotation
\eqref{rotation}, the proposition follows. 
\end{pf}

We now consider the simplest degeneration at a corner.
Compare this with Theorem 4.A of \cite{Floer88a} or 
the appendix of \cite{Sullivan99}. 
Let $\epsilon>0$ be a small number.
Let 
$A_s\colon D_m\to\U(n)$, $0\le s\le 1$ be a family matrices which
are small at 
infinity and constant in $s$ near each puncture in
$S\subset\{1,\dots,m\}$, where each component of the complex angle is
assumed to be positive. At $p_r$, $r\not\in S$, assume that
$\theta(r)_s=(\pi-s,\theta_2(r),\dots,\theta_n(r))$, where
$\theta_j(r)\ne 0$, $j=2,\dots,n.$. Let $\beta(\epsilon)\in\R^m$
satisfy $\beta(\epsilon)_r=0$ if $r\in S$ and
$\beta(\epsilon)_r=-\epsilon$ if $r\notin S.$ 

\begin{prp} \label{7degenindex.prop}
The index of the operators
$$
\bar\pa\colon\sblv_2(\C^n;A_s)\to\sblv_1(T^{0,1}D_m;[0])
$$
for $s>0$ and of the operator
$$
\bar\pa\colon\sblv_{2,\beta(\epsilon)}(\C^n;A_0)
\to\sblv_{1,(-\epsilon,0,\dots,0)} (T^{0,1}D_m;[0]), 
$$
are the same.
\end{prp}

\begin{pf}
This is a consequence of Lemma \ref{lmadef}.
\end{pf}

Finally, we show how the index is affected if the weight is changed. 

\begin{prp}\label{7degenindex2.prop}
Let $A\colon D_m\to\U(n)$ be constant at infinity and suppose that the
complex angle at each puncture except possibly $p_1$ has positive
components. Assume that
$0\le\pi-\theta(1)_1<\pi-\theta(1)_2<\dots<\pi-\theta(1)_n$.
Let $\epsilon>0$ be smaller than $\min_{r}(\pi-\theta(r)_r)$, 
and let $\pi-\theta(1)_{j}<\delta<\pi-\theta(1)_{j-1}$  
Then the index of the problem
$$
\bar\pa\colon\sblv_{2,(-\epsilon,0,\dots,0)}(\C^n;A)\to
\sblv_{1,(-\epsilon,0,\dots,0)}({T^\ast}^{0,1}D_m;[0]) 
$$ 
is $j$ larger than that of
$$
\bar\pa\colon\sblv_{2,(\delta,0,\dots,0)}(\C^n;A)
\to\sblv_{1,(\delta,0,\dots,0)}({T^\ast}^{0,1}D_m;[0])
$$ 
\end{prp}

\begin{pf}
First deform the matrix into diagonal form without changing the
weights. If $n>1$ this can be done in such a way that the index
corresponding to the first component is positive. Then put the first
component in standardized form. We must consider the index difference
arising from the first component as the weight changes from negative
to positive. The condition that a solution lies in $\sblv_{2,\delta}$
means that the corresponding solution on $\Delta$ vanishes at
$p_1$. Thus the dimension of the kernel increases by $1$. The cokernel
remains zero-dimensional. This
argument can then be repeated for other components. To handle the
$1$-dimensional case one may either use similar arguments for
cokernels or reduce to the higher dimensional case by adding extra
dimensions.   
\end{pf}

\subsection{The index and the Conley-Zehnder index}
\label{7Jtrick.section}

We translate Proposition \ref{indexw/oJtrick} into a more invariant language.
Recall from Section \ref{1CZ} that
we denote by $\nu_\gamma(c)$ the Conley-Zehnder index of Reeb chord
$c$ with capping path $\gamma$. In the following proposition we
suppress $\gamma$ from the notation.

\begin{prp}\label{indexJtrick}
Let $(u,f)\in\cand_2({\mathbf c},\kappa;B)$ be a holomorphic disk with
boundary on an admissible $L$, and with $j$ 
positive punctures at Reeb chords $a_1,\dots,a_j$ and $k$
negative punctures at Reeb chords $b_1,\dots,b_k$.
Then the index of $d\Gamma_{(u,f)}$ equals
\begin{equation}\label{indexJtrick.eqn}
\mu(B) + (1-j)n +\sum_{r=1}^j\nu(a_r) -\sum_{r=1}^k\nu(b_r).
\end{equation}
\end{prp}
\begin{rmk}
Note that \eqref{indexJtrick.eqn} is independent of the choices of
capping paths.
\end{rmk}

\begin{pf}
We simply translate the result of Proposition \ref{indexw/oJtrick}.
At a positive puncture $p$, the tangent space corresponding to $\R+0i$ 
($\R+i$) in $\pa E_{p}$ is the lower (upper) one and at a negative 
puncture the situation is reversed. We must compare the rotation path
$\lambda(V_1,V_0)$ used in the definition of the Conley--Zehnder index
with the rotation used in the construction of the arcs ${\mathbf f}_i$
in Proposition~\ref{indexw/oJtrick}. At a negative puncture, the path
${\mathbf f_i}$ is the inverse path of $\lambda(V_1,V_0)$. Hence the
contribution to the Maslov index of ${\mathbf f}_i$ at a negative
corner equals minus the contribution from $\lambda$. Consider the
situation at positive puncture mapping to $a^\ast$. Let
$\lambda(V_1,V_0)$ be the path 
used in the definition of the Conley--Zehnder index. Then
$\lambda(V_1,V_0)$ rotates the lower tangent space $V_1$ of $\Pi_\C(L)$ at
$a^\ast$ to the upper $V_0$ according to $e^{sI}$, $0\le
s\le\frac{\pi}{2}$, where $I$ is a complex structure compatible with
$\omega$. Let $\lambda(V_0,V_1)$ be the path which rotates $V_0$ to
$V_1$ in the same fashion. Then the path ${\mathbf f}_j$ is the inverse path of
$\lambda(V_0,V_1)$ and hence the contribution to the Maslov index of
${\mathbf f}_j$ equals the contribution from $\lambda(V_1,V_0)$ minus $n$.   

To get the loop $B$ from $\Gamma_A$ (see Proposition
\ref{indexw/oJtrick})   
the arcs ${\mathbf f}_i,$ must be removed and replaced by the arcs
$\Gamma_i,$ induced from the capping paths of the Reeb chords. A
straightforward calculation gives  
\[\mu(\Gamma_A)= \mu(B) + \sum_{r=1}^j\nu(a_r)-nj -\sum_{s=1}^k\nu(b_s).\]
Hence,
$$
\index=n+\mu(\Gamma_A)=\mu(B) + (1-j)n+\sum_{r=1}^j\nu(a_r)
-\sum_{s=1}^k\nu(b_s).
$$
\end{pf}

\subsection{The index and the Conley--Zehnder index at a self tangency}
\label{7CZST.section}
In this section we prove the analog of
Proposition~\ref{indexJtrick} for semi-admissible
submanifolds. First we need a definition of the Conley-Zehnder index
of a degenerate Reeb chord. Let $L\subset\R\times\C^n$ be a chord semi
generic Legendrian submanifold. Let $c$ be the Reeb chord of $L$ such
that $\Pi_\C(L)$ has a double point with self tangency along one
direction at $c^\ast$. Let $a$ and $b$ be the end points of $c$,
$z(a)>z(b)$. Let $V_0=d\Pi_\C(T_a L)$ and $V_1=d\Pi_\C(T_b L)$. Then
$V_0$ and $V_1$ are Lagrangian subspaces of $\C^n$ such that
$\dim_\R(V_0\cap V_1)=1$. Let $W\subset\C^n$ be the $1$-dimensional
complex linear subspace containing $V_0\cap V_1$ and let $\C^{n-1}$ be
the Hermitian orthogonal complement of $W$. Then
$V_0'=V_0\cap\C^{n-1}$ and $V_1'=V_1\cap\C^{n-1}$ are transverse
Lagrangian subspaces in $\C^{n-1}$. Pick a complex structure $I'$ on
$\C^{n-1}$ compatible with $\omega|\C^{n-1}$ such that
$I'V_1=V_0$. Define $\lambda(V_1,V_0)$ to be the path of Lagrangian planes 
$s\mapsto V_0\cap V_1\times e^{sI'}V_1'$. Also pick a capping path
$\gamma\colon[0,1]\to L$ with $\gamma(0)=a$ and $\gamma(1)=b$. Then
$\gamma$ induces a path $\Gamma$ of Lagrangian subspaces of
$\C^n$. Define the Conley--Zehnder index of $c$ as
$$
\nu_\gamma(c)=\mu(\Gamma\ast\lambda(V_1,V_0)).
$$

Let $0<\epsilon<\theta$, where $\theta$ is the smallest non-zero
complex angle of $L$ at $c$. Let $\cand_{2,\epsilon}({\mathbf
c};\kappa)$ denote the space of maps with boundary conditions
constructed from the Sobolev space with weight $\epsilon$ at each
puncture mapping to $c$ and define
$\widetilde{\cand}_{2,\epsilon}({\mathbf c};\kappa)$ as in Section
\ref{fcanI}. If $a$ is a Reeb chord of $L$ then let $\delta(a,c)=0$ if
$a\ne c$ and $\delta(c,c)=1$. Again we suppress capping paths from the
notation.    

\begin{prp}\label{indexJtrick1}
Let $(u,f)\in\cand_{2,\epsilon}({\mathbf c},\kappa;B)$ and 
$(v,g)\in\widetilde{\cand}_{2,\epsilon}({\mathbf c},\kappa;B)$ be
holomorphic disks. If
${\mathbf c}=(a; b_1,\dots,b_m)$ where $a\ne c$ then the index of  
$d\Gamma_{(u,f)}$ equals
$$
\mu(B)+\nu(a)-\sum_{r=1}^k(\nu(b_j)+\delta(b_j,c)),
$$
and the index of $d\Gamma_{(v,g)}$ equals
$$
\mu(B)+\nu(a)-\sum_{r=1}^k\nu(b_j).
$$
If ${\mathbf c}=(c;b_1,\dots,b_m)$ then the index of 
$d\Gamma_{(u,f)}$ equals
$$
\mu(B)+\nu(c)-\sum_{r=1}^k(\nu(b_j)+\delta(b_j,c)),
$$
and the index of $d\Gamma_{(v,g)}$ equals
$$
\mu(B)+(\nu(c)+1)-\sum_{r=1}^k\nu(b_j).
$$
\end{prp}

\begin{rmk}
Note again that the index computations are independent of the choices of
capping paths.
\end{rmk}

\begin{pf}
The proof is similar to the proof of
Proposition~\ref{indexJtrick}. Consider first the $\bar\pa$-operator
with boundary conditions determined by $(u,f)$ and acting on a Sobolev
space with small negative weight. Again we
need to compute the Maslov index contributions from the paths
${\mathbf f}_i$ in the loop $\Gamma_A$, where ${\mathbf f}_i$ fixes
the common direction in the tangent spaces at a self tangency double
point. Note that at a positive puncture $c$ the contribution is now the
contribution of $\lambda(V_1,V_0)$ minus $(n-1)$. At a negative puncture it
is again minus the contribution of $\lambda(V_1,V_0)$. Applying
Proposition \ref{7degenindex2.prop} the 
first and third index calculations above follow. Noting that the
tangent space of $\widetilde{\cand}_{2,\epsilon}({\mathbf c},\kappa,B)$ is obtained from that of
$\cand_{2,\epsilon}({\mathbf c},\kappa;B)$ by adding one $\R$-direction for each puncture
mapping to $c$ the other index formulas follow as well.
\end{pf}

\section{Transversality}
\label{8transversality.section}
In this section we show how to achieve transversality (or ``surjectivity'')
for the linearized $
\dbar$ equation by perturbing the Lagrangian
boundary condition. 
When proving transversality for some Floer-type theory, it is customary
to show that solution-maps are ``somewhere injective'' 
(see \cite{ms, FHS}, for example).
One then constructs a small perturbation, usually of the almost
complex structure or the Hamiltonian term, which is supported near
points where the map is injective. With a partial integration
argument, these perturbations eliminate non-zero elements of the
cokernel of $\dbar.$ 

For our set-up, we perturb the Lagrangian boundary condition.
In Sections \ref{8perturb.section} through 
\ref{8degenperturb.section}, we describe the space of perturbations
for the chord generic, one-parameter chord generic, and chord
semi-generic cases.
Although we do not have an injective (boundary) point, we exploit
the fact that there is only one positive puncture, and hence,
by Lemma \ref{1height}, the corresponding double point can represent
a corner only once. 
Of course other parts of the boundary can map to this corner elsewhere, but
not at other boundary punctures.
With this observation, we prove transversality in Sections 
\ref{8nonexceptionaltrans.section} and \ref{8degentrans.section}
first for the open set of non-exceptional maps, defined in Section 
\ref{8exceptional.section} and from this for all maps provided the
expected kernel has sufficiently low dimension.
We also prove some results in Sections 
\ref{8enhancedtrans.section} and \ref{8auxspace.section}
which will be useful later for the degenerate gluing of Section~\ref{10glue.section}.

\subsection{Perturbations of admissible Legendrian submanifolds }
\label{8perturb.section}
Let $L\subset\C^n\times\R$ be an admissible Legendrian
submanifold. Let $a(L)$ denote the minimal distance between the images
under $\Pi_\C$ of two distinct Reeb chords of $L$ and let $A(L)$ be such that 
$\Pi_\C(L)$ is contained in the  ball $B(0,A(L))\subset\C^n$.
Fix $\delta>0$ and $R>0$ such that $\delta\ll a(L)$ and
such that $R\gg A(L)$. 

\begin{dfn}
Let $\ham(L,\delta,R)$ be the linear space of smooth functions
$h\colon\C^n\to\R$ with support in $B(0,R)$ and satisfying the
following two conditions for any Reeb chord $c$. 
\begin{itemize}
\item[{\rm (i)}]
The restriction of $h$ to $B(c^\ast,\delta)$ is real analytic,
\item[{\rm (ii)}]
The differential of $h$ satisfies $Dh(c^\ast)=0$ and also
$h(c^\ast)=0$.
\end{itemize}
\end{dfn}

We are going to use Hamiltonian vector fields of elements in
$\ham(L,\delta,R)$ to perturb $L$. 
Condition (i) ensures that $L$ stays
admissible, and (ii) that the set of Reeb chords
$\{c_0,\dots,c_m\}$  of $L$ remains fixed. 

\begin{lma}\label{Banachsp}
The space $\ham(L,\delta,R)$ with the $C^\infty$-norm is a Banach
space.
\end{lma}
\begin{pf}
Using the characterization of real analytic functions as smooth
functions the derivatives of which satisfy certain uniform growth
restrictions one sees that the limit of a $C^\infty$-convergent
sequence of real analytic functions on an open set is real analytic.
\end{pf}

\begin{lma}\label{admpres}
If $L$ is admissible and $h\in\ham(L,\delta,R)$ then 
$\tilde\Phi_h(L)$ (see Section \ref{5Isotopies.section}) is admissible.  
\end{lma}
\begin{pf}
For each Reeb chord $c$, the Hamiltonian vector field is 
real analytic in $B(c^\ast,\delta)$. Also, $\Phi_h(c^\ast)=c^\ast$ and
hence there exists a neighborhood $W$ of $c^\ast$ such that 
$\Phi_h^t(W)\subset B(c^\ast,\delta)$ for $0\le t\le 1$. A well-known
ODE-result implies that the flow of a real analytic vector field depends
in a real analytic way on its initial data. This shows that
$\tilde\Phi_h(L)$ is admissible. 
\end{pf}

\subsection{Perturbations of $1$-parameter families of 
admissible submanifolds}
\label{8oneparamperturb.section}
Let $L_t$, $t\in[0,1]$ be an admissible $1$-parameter family of Legendrian
submanifolds without self-tangencies. 
Let $a=\min_{0\le t\le 1}a(L_t)$ and $A=\max_{0\le t\le 1}A(L_t)$.
Fix $\delta>0$ and $R>0$ such that $\delta\ll a$ and $R\gg A$. 

We define a continuous family of isomorphisms
$\ham(\delta,R,L_{0})\to\ham(\delta,R,L_t)$, 
$0\le t\le 1$. Let $(c_1(t),\dots,c_m(t))$ be the Reeb chords of
$L_t$. Then $(c^\ast_1(t),\dots,c^\ast_m(t))$, $0\le t\le 1$ is a continuous
curve in $(\C^n)^m$. Let $\psi^t\colon B(0,R)\to B(0,R)$ be a
continuous family of compactly supported diffeomorphisms which when
restricted to $B(c^\ast_j(0),\delta)$, $j=1,\dots,m$ 
agree with the map
$$
z\mapsto z+(c^\ast_j(t)-c^\ast_j(0)).
$$ 
Composition with $\psi^t$ can be used to give the space 
$$
\pham(L_t,\delta,R)=\bigcup_{0\le t\le 1}\ham(L_t,\delta,R)
$$ 
the structure of a Banach manifold which is a trivial bundle over
$[0,1]$. We note that if $(h,t)$ in $\pham(L_t,\delta,R)$ then Lemma
\ref{admpres} implies that $\tilde\Phi_a(L_t)$ is admissible.

\subsection{Bundles over perturbations}
\label{8oneparambundle.section}
Let $L\subset\C^n\times\R$ be an admissible
chord generic Legendrian submanifold. Above we constructed a
smooth map of the Banach space $\ham(L,\delta,R)$ into the space of
admissible chord generic Legendrian embeddings of $L$ into $\C^n\times\R$.

Let ${\bf c}=(c_0,c_1,\dots,c_m)$ be Reeb chords of $L$ and let
$\epsilon\in [0,\infty)^m$ and
consider as in Section \ref{fcanA} the space
\begin{equation}
\cand_{2,\epsilon,\ham(L,\delta,R)}({\mathbf c})
\end{equation}
and its tangent space
\begin{align}
T_{(w,f,\kappa,a)}\cand_{2,\epsilon,\ham(L,\delta,R)}\approx
T_{(w,f)}\cand_{2,\epsilon}\oplus T_\kappa\conf_m\oplus \ham(L,\delta,R).
\end{align}

In a similar way we consider for a $1$-parameter family $L_t$ the
space
\begin{equation}
\cand_{2,\epsilon,\pham(L_t,\delta,R)}({\mathbf c})
\end{equation}
and its tangent space.

For $\Lambda=\ham(L,\delta,R)$ 
or $\Lambda=\pham(L_t,\delta,R)$,
consider also the bundle map
$\Gamma\colon\cand_{2,\epsilon,\Lambda}({\bf c}) \to \sblv_{1,\epsilon,\Lambda}[0]
({T^\ast}^{0,1}D_m\otimes\C^n)$ here we are thinking of the spaces as bundles over 
$\ham(L,\delta,R)$ and we denote projection onto this space by $\proj.$ 
To emphasize this we will be write $(\Gamma, \proj)$ instead of just $\Gamma$ in the sequel. 
The differential $d\Gamma$ was calculated in
Lemma~\ref{lem:fcan.gamma}.


\subsection{Perturbations in the semi-admissible case}
\label{8degenperturb.section}

Let $L\subset\C^n\times\R$ be a semi-admissible
Legendrian submanifold. Let $(c_0,\dots,c_m)$ be the Reeb chords of
$L$. Assume that the self tangency Reeb chord is $c_0$, that
$c_0^\ast=0$, and that $L$ has standard form in a neighborhood of $0$,
see Definition \ref{5admissile_param.dfn}.

Let $a(L)$ denote the minimal distance between the images under $\Pi_\C$
of two distinct Reeb chords of $L$. Fix $\delta>0$ such that
$\delta\ll a(L)$. For $r>0$ let $C(r)=\C\times B'(0,r)\subset\C^n$,
where $B'(0,r)$ is 
the $r$-ball in $\C^{n-1}\approx\{z_1=0\}$, where as always
$(z_1,\dots,z_n)=(x_1+iy_1,\dots,x_n+iy_n)$ are coordinates on $\C^n$.

\begin{dfn}\label{dfnham_0}
Let $\ham_0(L,\delta)$  be the linear space of smooth functions
$h\colon\C^n\to\R$ with support in 
$C(10\delta)\cup\bigcup_{j\ge 1}B(c_j^\ast,10\delta)$
and satisfying the following conditions. 
\begin{itemize}
\item[{\rm (i)}]
The restriction of $h$ to $B(c^\ast_j,\delta)$ $1\le j\le m$ is real
analytic,  
\item[{\rm (ii)}]
In $C(10\delta)$, $\frac{\pa h}{\pa x_1}=0=\frac{\pa h}{\pa y_1}$ and 
the restriction of $h$ to $C(\delta)$ is real analytic.  
\item[{\rm (iii)}]
The differential of $h$ satisfies $Dh(c_j^\ast)=0$ and also
$h(c_j^\ast)=0$, for all $j$.
\end{itemize}
\end{dfn}

\begin{lma} \label{8D2.lma}
The space $\ham_0(L,\delta)$ with the $C^\infty$-norm is a Banach
space.
\end{lma}

\begin{pf}
See Lemma \ref{Banachsp} and note that the restriction of $h$ to
$C(10\delta)$ can be 
identified with a function of $(n-1)$-complex variables supported in
$B'(0,10\delta)$. 
\end{pf}

Let $\tilde\Phi_h$ be the Legendrian isotopy which is defined by using
the flow of $h$ locally around the Reeb chords of $L$. This is
well-defined for $h$ sufficiently small. Let $\ham_0(L,\delta,s)$
denote the $s$-ball around $0$ in $\ham_0(L,\delta)$.

\begin{lma} \label{8D3.lma}
There exists $s>0$ such that for $h\in\ham_0(L,\delta,s)$, $\tilde
\Phi_h(L)$ is an admissible chord semi-generic Legendrian 
submanifold.
\end{lma}

\begin{pf}
Note that the product structure in $C(10\delta)$ is preserved since $h$ does
not depend on $(x_1,y_1)$. Moreover, the isotopy is fixed in the
region $B(0,2+\epsilon)\setminus B(0,2)$ for $s$ and $\delta$ sufficiently
small. 
\end{pf}

We have defined a smooth map of $\ham_0(L,\delta,s)$ into the space of
admissible chord semi-generic Legendrian submanifolds and this maps
fulfills the conditions on $\Lambda$ in Section~\ref{fcanI}. We can therefore
construct the spaces
\begin{equation}
\cand_{2,\epsilon,\ham_0(L,\delta,s)},\quad{ \mbox{and} }\quad
{\widetilde \cand}_{2,\epsilon,\ham_0(L,\delta)},
\end{equation}
see Sections \ref{fcanA} and \ref{fcanI}, respectively. Moreover, as there we
will consider the $\bar\pa$-map and its linearization.

\subsection{Consequences of real analytic boundary conditions}
\label{8.realanalytic.section}
For $r >0$, let $E_+=\{z\in\C\colon |z|< r, \Im(z)\ge 0\}$.
If $w\colon(E_+,\pa E_+)\to(\C^n,M)$ where $M$ is a real analytic
Lagrangian submanifold and $w$ is holomorphic in the interior and
continuous on the boundary, then by Schwartz-reflection principle, $w$
extends in a unique way to a holomorphic map $w^d\colon E\to\C^n$ mapping
$\Im(z)=0$ to $M$, where $E=\{z\in\C\colon |\zeta|< r\}$. 
We call $w^d$ the double of $w$. 

Let $L\subset\C^n\times\R$ be a chord (semi-)generic Legendrian
submanifold. 
\begin{lma}\label{lmaTaylor}
Let $p$ be a point in $U\subset L$ such that
$\Pi_\C(U)$ is real analytic, where $U\subset L$ is a neighborhood of 
$p$ on which $\Pi_\C$ is injective. Assume that 
$$
w\colon (E_+,\pa E_+,0)\to(\C^n, \Pi_\C(U), p),
$$
is holomorphic. Then there is a holomorphic function $u$ with Taylor
expansion at $0$, 
$$
u(z)=a_0 + a_1 z +\dots, a_0\ne 0
$$
such that $w(z)=p+ z^k u(z)$ for some integer $k>0$.
\end{lma}

\begin{pf}
The double $w^d$ has a Taylor expansion.
\end{pf}

\begin{lma}\label{lmafinite}
Let $p, U$ and $L$ satisfy the conditions of Lemma~\ref{lmaTaylor}.
Assume that 
$w \colon D_m \to \C^n$ is holomorphic with boundary
on $\Pi_\C(L)$.  
Then $w^{-1}(p)  \cap \pa D_m$ is a finite set. 
\end{lma}

\begin{pf}
Using Lemma \ref{SalRob},
we may find $M>0$ such that there are no preimages of $p$ in 
$\cup_j E_{p_j}[M]$. Since the complement of $\cup_j E_{p_j}[M]$ is
compact, the lemma now follows from Lemma ~\ref{lmaTaylor}.
\end{pf}

\begin{lma}\label{lmamult}
Let $p\in L$ satisfy the conditions of Lemma \ref{lmaTaylor}, and let
\begin{align}
w_1, w_2\colon (E_+,\pa E_+,0)\to(\C^n, \Pi_\C(L), p),
\end{align}
be holomorphic maps such that $w_2$ maps one of the components $I$ of 
$\pa E_+\setminus\{0\}$ to $w_1(I)$. Then there exists a map $\hat
w\colon E\to\C^n$ and integers  
$k_j\ge 1$ such that $w_j^d(z)=\hat w(z^{k_j})$, $j=1,2$.   
\end{lma}
\begin{pf}
As above we may reduce to the case when
$\Pi_\C(L)=\R^n\subset\C^n$. The images $C_j=w^d_j(E)$,
$j=1,2$ are analytic subvarieties of complex dimension $1$ which
intersects in a set of real dimension $1$. Hence they
agree. Projection of $C=C_1=C_2$ onto a generic complex line through
$p$ identifies $C$ (locally) with the standard cover of the disk
possibly branched at $0$. This gives the map $\hat w$.
\end{pf}


\subsection{Exceptional holomorphic maps}
\label{8exceptional.section}
Let $\Lambda$ be one of the spaces $\ham(L,\delta,R)$,
$\pham(L,\delta,R)$, or $\ham_0(L,\delta,s)$. Let
$(w,f,\lambda)\in\cand_{2,\epsilon,\Lambda}({\mathbf c})$ (or
${\widetilde\cand}_{2,\epsilon,\Lambda}({\mathbf c})$) be a
holomorphic disk and let $q$ 
be a point on $\pa D_m$ such that $w(q)$ lies in a region where
$\Pi_\C(L_\lambda)$  is real analytic. Assume that $dw(q)=0$.
Since $w$ has a Taylor expansion around $q$ in this case we know there
exists a half-disk neighborhood $E$ of $q$ in $D_m$ such that $q$ is the only
critical point of $w$ in $E$. The boundary $\pa E$ is  subdivided
by $q$ into two arcs $\pa E\setminus\{q\}=I_+\cup I_-$. 
We say that $q$ is {\em an exceptional point of $(w,f)$} if there exists a
neighborhood $E$ as above such that $w(I_+)=w(I_-)$.  

\begin{dfn} \label{8exceptional.dfn}
Let $(w,f,\lambda)\in\cand_{2,\epsilon,\Lambda}({\bf c})$, where 
${\bf c}=(c_0(\lambda),c_1(\lambda),\dots,c_m(\lambda))$  and
$c_0(\lambda)$ is the Reeb chord 
on $L_\lambda$ of the positive puncture of $D_{m+1}$. Let
$B_1(\lambda)$ and $B_2(\lambda)$ be the two local branches of
$\Pi_\C(L_\lambda)$ at 
$c^\ast_0(\lambda)$. Then $(w,f)$ is {\em exceptional holomorphic} if it has
two exceptional points $q_1$ and $q_2$ with $w(q_1)=w(q_2)=c_0^\ast(\lambda)$
and if a neighborhood in $\pa D_m$ of $q_j$ maps to $B_j(\lambda)$, $j=1,2$.
\end{dfn}

\begin{dfn} \label{8nonexceptional.dfn}
Let $\cand'_{2,\epsilon,\Lambda}({\mathbf c})$
(${\widetilde\cand}'_{2,\epsilon,\Lambda}({\mathbf c})$) denote the
complement of 
the closure of the set of all exceptional holomorphic maps in
$\cand_{2,\epsilon,\Lambda}({\mathbf c})$
(${\widetilde\cand}_{2,\epsilon,\Lambda}({\mathbf c})$).
\end{dfn}

We note that $\cand'_{2,\epsilon,\Lambda}({\mathbf c})$ is an open
subspace of a Banach manifold and hence a Banach manifold itself.

\subsection{Transversality on the complement of exceptional
holomorphic maps in the admissible case}
\label{8nonexceptionaltrans.section}

\begin{lma}\label{lmatv'}
For $L$ admissible (respectively $L_t$ a $1$-parameter family of
admissible submanifolds) the bundle map, see Section \ref{fcanH}  
$$
(\Gamma,\proj)\colon\cand_{2,\epsilon,\Lambda}'({\bf c})\to
\sblv_{1,\epsilon,\Lambda}[0](D_m,T^{*}D_m\otimes\C^n), 
$$ 
where $\Lambda=\ham(L,\delta,R)$ (respectively $\Lambda=\pham(L_t,\delta,R)$)
is transverse to the $0$-section.  
\end{lma}

\begin{pf}
The proof for $1$-parameter families $L_t$ is only notationally more
difficult. We give the proof in the stationary case.
We must show that if $w\colon D_m\to\C^n$ is a (non-exceptional)
holomorphic map (in the conformal structure $\kappa$ on $D_m$)
which represents a holomorphic disk $(w,f)$ with boundary on
$L=L_\lambda$ (without loss of generality we take $\lambda=0$
below) then
$$
d\Gamma\Bigl( T_{((w,f),\kappa,0)}\cand_{2,\epsilon,\Lambda}({\bf c})\Bigr)=
\sblv_{1,\epsilon}(D_m,T^\ast D_m\otimes \C^n),
$$  
i.e., $d\Gamma$ is surjective.
To show this it is enough to show that 
$$
\left\{d\Gamma \Bigl(T_{((w,f),\kappa,0)}\cand_{2,\epsilon,\Lambda}({\bf
c})\Bigl)\right\}^\perp =\{0\},  
$$
where $V^\perp$ denotes the annihilator with respect to the
$L^2$-pairing of 
$V\subset\sblv_{1,\epsilon}[0](D_m,T^\ast D_m\otimes\C^n)$
in its dual space. 

An element $u$ in this annihilator satisfies
\begin{equation}
\label{8annihilator.eqn}
\int_{D_m}\la\bar\pa v,u\ra\,dA=0,
\end{equation}
for all $v\in T_{w}\B_{2,\epsilon}(0,r)$. 
Here $dA$ is the area form on $D_m.$
Lemma~\ref{lmacokerreg} implies that $u$ can 
be represented by a $C^2$-function which is anti-holomorphic.

We note that integrals of the form
\begin{equation}\label{2confinvtint} 
\int_{D_m}\la\phi,\psi\ra\,dA,
\end{equation}
where $\la\,,\ra$ is the inner product on $T^\ast D_m$ and where $\phi$ and
$\psi$ are sections, are conformally
invariant. We may therefore compute integrals of this form in any
conformal coordinate system on the disk $D_m$.

Restrict attention to a
small neighborhood of the image of the  
positive corner at $c_0^\ast$. Recall that $w$ is
assumed non-exceptional and consider a branch $B$ of $\Pi_\C(L)$ at
$c_0^\ast$ such that $w$ does not have an exceptional point mapping to
$c_0^\ast\in B$. Since $B$ is real analytic we may biholomorphically identify
$(\C^n,B,c^\ast_0)$ with $(\C^n,\R^n,0)$.   

Let $p$ be the positive puncture on $D_m.$
For $M$ large enough, by Lemma \ref{SalRob}, the image of the component of
$\pa E_{p}[M]$  which lies in $B$ is a regular oriented curve. Denote 
it by $\gamma$. For simplicity we assume that the component mapping to
$\gamma$ is $[M,\infty)\times \{0\}\subset E_{p_0}[M]$ and we let
$E_0=[M,\infty)\times[0,\frac12)$.  

Let $p_1,\dots,p_r$ be the preimages under $w$ of $c^\ast_0$ with
the property that one of the components of a punctured
neighborhood of $p_j$ in $\pa D_m$ maps to $\gamma$.
Note that $r<\infty$ by Lemma \ref{lmafinite} and that by shrinking $\gamma$ we
may assume that all these images are exactly $\gamma$.

We say that a point $p_j$ is positive if close to $p_j$, $w$
and the natural orientation on the boundary of $\pa D_m$ induce the
positive orientation on $\gamma$ otherwise we say it is negative.

The image of the other half of the punctured neighborhood of $p_1$ in
$\pa D_m$ maps to a curve $\gamma'$ under $w$. Our assumption that
$w$ is non-exceptional guarantees that $\gamma$ and $\gamma'$ are
distinct.  

Let $w_j$ denote the restriction of $w$ to a small neighborhood of $p_j$.
Let $E=\{z\in\C\colon |z|<r\}$, let $E_+=\{z\in E\colon \Im(z)\ge
0\}$, and let $E_-=\{z\in E\colon \Im(z)\le 0\}$.
Lemma \ref{lmamult} implies that we can find a map 
$\hat w\colon E\to\C^n$ and coordinate neighborhoods 
$(E_\pm(j),\pa E_\pm(j))$ of $p_j$ (where the sign $\pm$ is that of
$p_j$) such that $w_j^d(z)=\hat w(z^{k_j})$ for each $j$. Note that
$w$ non-exceptional implies all $k_j$ are odd. 

Let $k=k_1k_2\dots k_r$ and let $\hat k_j=\frac{k}{k_j}$. 
Let $\phi_j\colon E\to E(j)$ be the map $z\mapsto z^{\hat k_j}$.
Consider the restrictions $u_j$ of the anti-holomorphic map 
$u$ to the neighborhoods $(E_\pm(j),\pa E_\pm(j))$.  
Because of the real analytic boundary conditions (recall that $(B,\C^n)$ is
biholomorphically identified with $(\R^n,\C^n)$),
these maps can be doubled using Schwartz reflection principle. 
Use $\phi_j$ to pull-back the maps $u_j$ and $w_j$ to $E$.

Let $a\colon\C^n\to\R$ be any smooth function with support in a
small ball around a point $q'\in\gamma'$, where $q'$ is chosen so that
no point outside $\bigcup_j E_\pm(j)$ in $\partial D_m$ maps to $q'$. (There
exists such a point because of the asymptotics of $w$ at punctures and 
Lemma \ref{lmaTaylor}.) 
Let $Y_a$ is the Hamiltonian vector field associated with $a,$ see Section~\ref{2Ham.iso}.

If $v$ is a smooth function with support in $\bigcup_j E_\pm(j)$ which
is real and holomorphic on $\bigcup_j \pa E_\pm(j)$, if
$\xi+i\eta$ are coordinates on $E_\pm(j)$, and if the support of $a$
is sufficiently small then
\begin{align}\label{gamma'1}
0=&\int_{D_m}\la\bar\pa (Y_a+v), u\ra\,dA\\\label{gamma'2}
=&-\sum_j \int_{E_\pm(j)}\la Y_a+v, \pa u \ra\,d\xi\wedge d\eta
+\sum_j \int_{\pa E_\pm(j)}\la -i(Y_a+v),u\ra\,d\xi\\\label{gamma'3}
=&\sum_j \int_{\pa E_\pm(j)}\la -i(Y_a+v),u\ra\,d\xi.
\end{align}
The equality in \eqref{gamma'1} follows since $u$ is
an element of the annihilator and 
since $a$ can be arbitrarily well $C^2$-approximated by elements in
$\ham(L,\delta,R)$. The equality in \eqref{gamma'2} follows by partial
integration and the restrictions on the supports of $a$ and $v.$
The equality in \eqref{gamma'3} follows
from $\pa u=0$. Taking $a=0$ we see, since we are free to choose
$v$, that $u$ must be real valued on $\pa E_\pm(j)$ for every $j$.

We then take $v=0$ and express the integral in \eqref{gamma'3} as an
integral over $I_+=\{x+0i\colon x>0\}\subset E$. Note that if $\xi+i\eta$
are coordinates on $E(j)$ then under the identification by $\phi_j$, $d\xi=d x^{\hat k_j}=
\hat k_jx^{\hat k_j-1}\,dx$ and 
\begin{align}
&\sum_j \int_{\pa E_\pm(j)}\la -iY_a,u\ra\,d\xi
&=\int_{I_+}\la -iY_a(\hat w(z^k)),\sum_j \sigma(j)
\hat k_j\bar z^{\hat k_j-1}u_j(z^{\hat k_j})\ra\,dx, 
\end{align}
where $\sigma(j)=\pm 1$ equals the sign of $p_j$.
Thus, if $\alpha(z)=\sum_j \sigma(j)
\hat k_j\bar z^{\hat k_j-1}u_j(z^{\hat k_j})$ then $\alpha$ is
antiholomorphic and by varying $a$ we see that $\alpha$ vanishes
in the $i\R^n$-directions along
an arc in $I_+$.
Therefore $\alpha$ vanishes identically on $E$. 

Pick now instead $a$ supported in a small ball around $q$ in
$\gamma$.
 With the same arguments as above we find
\begin{align}\notag
0=&\int_{D_m}\la(\bar\pa Y_a+v), u\ra\,dA\\\notag
=&-\int_{E_0}\la Y_a+v, \pa u \ra\,d\tau\wedge dt
-\sum_j \int_{E_\pm(j)}\la Y_a+v, \pa u \ra\,d\xi\wedge d\eta\\\notag
&+\int_{[M,\infty)}\la -i(Y_a+v), u \ra\,d\tau
+\sum_j \int_{\pa E_\pm(j)}\la -i(Y_a+v),u\ra\,d\xi\\\label{gamma}
=&\int_{[M,\infty)}\la -i(Y_a+v), u \ra\,d\tau+
\sum_j \int_{\pa E_\pm(j)}\la -i(Y_a+v),u\ra\,d\xi.
\end{align}
and conclude that $u(\tau,0)\in\R^n$ for $\tau\in[M,\infty)$ as well. 

Again taking $v=0$ we get for the last integral in \eqref{gamma}
\begin{equation}\label{gamma0}
\sum_j \int_{\pa E_\pm(j)}\la -iY_a,u\ra\,d\xi=
\int_{I_-}\la -iY_a(\hat w(z^k)),\alpha(z))\ra\,dx=0,
\end{equation}   
where $I_-=\{x+0i\colon x<0\}\subset E$, and where the last
equality follows since $\alpha=0$. Equations \eqref{gamma0} and
\eqref{gamma} 
together implies (by varying $a$) that $u$ must vanish along an arc in
$[M,\infty)$. Since $u$ is antiholomorphic it must then vanish 
everywhere. This proves the annihilator is $0$ and the lemma follows.
\end{pf}

\begin{rmk}\label{brnchpts}
In the case that $w$ has an injective point on the boundary, the above
argument can be shortened. 
Namely, under this condition there is an arc $A$ on the
boundary of $D_m$ where $w$ is 
injective and varying $v$ and $a$ there we see that $u$ must vanish
along $A$ and therefore everywhere.
Oh achieves transversality using boundary perturbations assuming
an injective point \cite{oh2}.
\end{rmk}

\begin{cor}\label{cormfd'}
Let ${\mathbf c}=ab_1\dots b_m$. For a Baire set of
$h\in\ham(L,\delta,R)=\Lambda$, 
$\Gamma^{-1}(0)\cap \proj^{-1}(h)\cap\cand_{2,\epsilon,\Lambda}'({\bf c};A)$ 
is a finite dimensional smooth manifold of dimension  
$$
\mu(A)+\nu_\gamma(a)-\sum_j \nu_\gamma(b_j)+\max(0,m-2).
$$
For a Baire set of
$h\in\pham(L,\delta,R)=\Lambda$ 
$\Gamma^{-1}(0)\cap \proj^{-1}(h)\cap\cand_{2,\epsilon,\Lambda}'({\bf c};A)$ 
is a finite dimensional smooth manifold of dimension  
$$
\mu(A)+\nu_\gamma(a)-\sum_j \nu_\gamma(b_j)+\max(0,m-2)+1.
$$
\end{cor}

\begin{pf}
Let $Z\subset\cand_{2,\epsilon,\Lambda}'({\mathbf c};A)$ denote the inverse 
image of the $0$-section in
$\sblv_{1,\epsilon,\Lambda}[0]({T^\ast}^{0,1}D_m\otimes\C^n)$ under 
$(\Gamma,\proj)$. By the implicit function theorem and Lemma \ref{lmatv'},
$Z$ is a  submanifold. Consider the restriction of the projection 
$\pi\colon Z\to \Lambda$. Then $\pi$ is a
Fredholm map of index equal to the index of the Fredholm section
$\Gamma$. An application of the Sard-Smale theorem shows that
for generic $\lambda\in\Lambda$, $\pi^{-1}(\lambda)$ is a submanifold
of dimension given by the Fredholm index of $\Gamma$. Note that in
the first case, the restriction of $d\Gamma$ to the complement of the
$\max(0,m-2)$-dimensional subspace 
$T\conf_m\subset T\cand_{2,\epsilon}({\mathbf c};A)$ is an operator of
the type considered in Proposition~\ref{indexJtrick}. 
Thus, the proposition follows in this first case. 
In the second case, we restrict to a
$(\max(0,m-2)+1)$-codimensional subspace instead.  
\end{pf}

\subsection{General transversality in the admissible case}
\label{8gentransv.section}

If ${\mathbf c}$ is a collection of Reeb chords we define $l({\mathbf
c})$ as the number of elements in ${\mathbf c}$. 
We note that if $(f,w)$ is a holomorphic disk with boundary on
$L$ with $r$ punctures, then, if $r\le 2$, the kernel of $d\Gamma$ at
$(f,w)$ is at least $(3-r)$-dimensional. This is a consequence of the
existence of conformal reparameterizations in this case.  

\begin{thm}\label{thmmfd}
For a dense open set of $h\in\ham(L,\delta,R)$ ($h\in\pham(L,\delta,R)$),
$\Gamma^{-1}(0)\cap\proj^{-1}(h)
\subset\cand_{2,\epsilon,\Lambda}({\mathbf c})$ is a finite dimensional
$C^1$-smooth manifold of dimension as in 
Corollary~\ref{cormfd'}, provided this dimension is $\le 1$ if 
$l({\mathbf c})\ge 3$ and
$\le 1+(3-l({\mathbf c}))$ otherwise.
\end{thm}
\begin{pf}
After Corollary~\ref{cormfd'} we need only exclude holomorphic
disks in the closure of exceptional holomorphic disks.
Let $a\in\ham$ ($a\in\pham$) be such that $\Gamma^{-1}(0)$ is
regular. Then the same is true for $\tilde a$ in a neighborhood of
$a$. Now assume there exists a holomorphic disk in the closure of
exceptional holomorphic disks at $a$. Then there must exists
an exceptional holomorphic disk for some $\tilde a$ in the
neighborhood. However, such a disk $w$ has $k\ge 2$ points mapping to
the image of the positive puncture and with $w(I_+)=w(I_-)$.  
It is then easy to construct (by ``moving the branch point'') a 
$k$-parameter ($k+(3-l({\mathbf c}))$-parameter if $l({\mathbf c})\le
2$) family of distinct (since the location of the branch
point changes) non-exceptional holomorphic disks with
boundary on $L(\tilde a)$. This contradicts the fact that the
dimension of $\Gamma^{-1}(0)$ is $<k$ ($<k+(3-l({\mathbf c})$) for
every $\tilde a$ in the neighborhood.   
\end{pf}

\begin{proof}[Proof of Proposition~\ref{1manifold}]
If the number of punctures is $\ge 3$ the proposition is just
Theorem \ref{thmmfd}. The case of fewer punctures can be reduced to
that of many punctures as in Section \ref{ssecmarked}.
\end{proof}

\begin{cor}\label{8decay.cor}
For chord generic admissible Legendrian submanifolds in a Baire set of such 
manifolds, no rigid holomorphic disk with boundary on $L$
decays faster than $e^{-(\theta+\delta)|\tau|}$ close to any of its punctures
mapping to a Reeb chord $c.$ 
Here $\theta$ is the smallest complex
angle of the Reeb chord $c$, $\delta>0$ is arbitrary, and $\tau+it$
are coordinates near the puncture.
\end{cor}

\begin{pf}
Such a holomorphic disk would lie in
$\cand_{2,\epsilon}({\bf c})$, where the component of $\epsilon$
corresponding to the puncture mapping to the Reeb chord $c$ is larger
than $\theta$. By Proposition~\ref{7degenindex2.prop} this change of
weight lowers the 
Fredholm index of $d\Gamma$ by at least $1$. Since the Fredholm index
of $d\Gamma$ with smaller weight (e.g. $0$-weight) is the minimal
which allows for existence of disks the lemma follows from Theorem
\ref{thmmfd}.  
\end{pf}

\begin{proof}[Proof of Proposition~\ref{1param-manifold}]
The first statement in the proposition follows exactly as above. To
see that handle slides appear at distinct times, let 
$(a_1{\mathbf b}_1;A_1)$ and $(a_2{\mathbf b}_2;A_2)$ be such that 
$$
\mu(A_1)+|a_1|-|{\mathbf b}_1|=\mu(A_2)+|a_2|-|{\mathbf b}_2|=0
$$ 
and consider the bundle $\cand_{2,\Lambda}(a_1{\mathbf
b}_1;A_1)\tilde\times\cand_{2,\Lambda}(a_2{\mathbf b}_2;A_2).$ Here
$\tilde\times$ denotes the fiberwise product where, in the fibers, the
deformation coordinates $(t_1,t_2)$ are restricted to lie in the
diagonal: $t_1=t_2=t$. This is a bundle over $\Lambda$, and $\Gamma$
induces a bundle map to the bundle
$\sblv_{1,\Lambda}(D_{m_1},\C^n)\times\sblv_{1,\Lambda}(D_{m_2},\C^n)$,
where $\times$ denotes fiberwise product. It is then easy to check that
$\Gamma$ is a Fredholm section of index $-1$. As in Theorem \ref{thmmfd} we
see that $d\Gamma$ is surjective and that the inverse image of the $0$-section
intersected with $\proj^{-1}(h)$ is empty for generic $h$. This shows
that the handle slides appear at distinct times.

The statement about all rigid disks being transversely cut out at a
handle slide instant can be proved in a similar way: let $(a_1{\mathbf
b_1};A_1)$ be as above and let $(a_3{\mathbf b_3};A_3)$ be such that 
$$
\mu(A_3)+|a_3|-|{\mathbf b}_3|=1.
$$
Consider the bundle 
$$
\cand_{2,\Lambda}(a_3{\mathbf b_3};A_3)\tilde\times
\cand_{2,\Lambda}(a_3{\mathbf b_3};A_3)\tilde\times
\cand_{2,\Lambda}(a_1{\mathbf b_1};A_1),
$$
and the bundle map $\Gamma$ defined in the natural way with target 
$$ 
\sblv_{1,\Lambda}(D_{m_3},\C^n)\times
\sblv_{1,\Lambda}(D_{m_3},\C^n)\times
\sblv_{1,\Lambda}(D_{m_1},\C^n).
$$
Then the map
$\Gamma$ has Fredholm index $0$ and as above we see $d\Gamma$ is
surjective. Hence $\Gamma^{-1}(0)\cap\proj^{-1}(h)$ is a transversely
cut out $0$-manifold for generic $h$. We show that this implies that
if $t$ is such that 
$\M^t_{A_1}(a_1;{\mathbf b}_1)=\{(v,g)\}\ne\emptyset$ then
$\M^t_{A_3}(a_3;{\mathbf b}_3)$ is transversally cut out. Let 
$(u,f)\in\M^t_{A_3}(a_3;{\mathbf b}_3)$ and assume the differential
$d\Gamma_{(u,f)}^t$, which is a Fredholm operator of index $0$ is not
surjective. Then it has a cokernel of dimension $d>0$. Furthermore,
the image of the tangent space to the fiber under the 
differential $d\Gamma$ at the point
$\Bigl(((u,f),(u,f),(v,g)),h\Bigr)$ is contained in a subspace of
codimension $\ge 2d-1$ in the tangent space to the fiber of the target
space. This contradicts $\Gamma^{-1}(0)\cap\proj^{-1}(h)$ being
transversely cut out.       
\end{proof}

\subsection{Transversality in the semi-admissible case}
\label{8degentrans.section}

\begin{lma}\label{lmasgtv'}
Suppose $L$ is admissible chord-semi-generic and 
$\Lambda=\ham_0(L,\delta,s)$, then the bundle maps 
\begin{align}
\notag
(\Gamma,\proj)\colon\cand'_{2,\epsilon,\Lambda}({\bf c})\to
\sblv_{1,\epsilon,\Lambda}(D_m,T^{0,1}D_m\otimes\C^n),\\ 
\notag
(\Gamma,\proj)\colon{\widetilde\cand}'_{2,\epsilon,\Lambda}({\bf
c})\to \sblv_{1,\epsilon,\Lambda}(D_m,T^{0,1}D_m\otimes\C^n)
\notag 
\end{align}
are transverse to the $0$-section.
\end{lma}

\begin{pf}
We proceed as in the proof of Lemma~\ref{lmatv'}.
Let $u$ be an element in the annihilator.
The argument of Lemma~\ref{lmatv'} still applies up to the
point where we conclude $\alpha|I_+$ equals $0$. 
In the present setup not all Hamiltonian vector fields are allowed
(see Definition \ref{dfnham_0}). However, the ones that are allowed can be
used exactly as in the proof of Lemma \ref{lmatv'} to conclude the last
$(n-1)$ components of $u$ must vanish identically. 

Since $D_m$ is conformally equivalent to the unit disk $\Delta_m$ with
$m$ punctures on the boundary and since integrals as in 
\eqref{2confinvtint} are
conformally invariant, we have for any smooth compactly supported $v$
with appropriate boundary conditions  
\begin{equation}\label{intsg}
0=\int_{\Delta_m}\la \bar\pa v,u\ra\, dA=
\int_{\Delta_m}\la v,\pa u\ra\, dA
+\int_{\pa\Delta_m}\la u,e^{-i\theta}v\ra\,d\theta.
\end{equation}
As usual the first term in \eqref{intsg} vanishes and we find that $u$ is
orthogonal to $e^{i\theta}T_{w(e^{i\theta})}\Pi_\C(L)$. 

Now the boundary of the holomorphic disk must
cross the region $X=B(0,2+\epsilon)\setminus B(0,2)$, and the inverse image
of this region contains an arc $A$ in the boundary. The
intersection between the tangent
plane of $T_p\Pi_\C(L)$, $p\in X$ and the $z_1$-line equals $0$ and
the $z_1$-line is invariant under multiplication by $e^{i\theta}$. 
Hence the orthogonal complement of
$e^{i\theta}T_{w(e^{i\theta})}\Pi_\C(L)$ intersects 
the $z_1$-line trivially as well (for $\theta\in A$). 
We conclude that the first
component of $u$ must vanish identically along $A$ and by anti-analytic 
continuation vanish identically. It follows that $u$ is
identically zero. 
\end{pf}

In analogy with Corollary \ref{cormfd'} we get (with $c$ denoting the
degenerate Reeb chord of $L$)

\begin{cor} \label{8degenexpotrans.cor}
For a dense open set of $h\in\ham_0(L,\delta,s)$,
$\Gamma^{-1}(0)\cap\proj^{-1}(h)\subset\cand_{2,\epsilon,\Lambda}'({\mathbf
c};A)$ and 
$\Gamma^{-1}(0)\cap\proj^{-1}(h)\subset{\widetilde\cand}_{2,\epsilon,\Lambda}' 
({\mathbf c};A)$
are finite dimensional manifolds. If ${\mathbf c}=a b_1\dots b_m$ with
$a\ne c$ then the dimensions are
\begin{align*}
&\mu(A)+\nu(a)-\sum_{r=1}^m(\nu(b_j)+\delta(b_j,c))+\max(0,m-2)\text{
and }\\
&\mu(A)+\nu(a)-\sum_{r=1}^m(\nu(b_j))+\max(0,m-2),\text{ respectively.}
\end{align*}
If ${\mathbf c}=c b_1\dots b_k$ then the dimensions are
\begin{align*}
&\mu(B)+\nu(c)-\sum_{r=1}^m(\nu(b_j))+\max(0,m-2)\text{
and }\\
&\mu(B)+\nu(c)+1-\sum_{r=1}^m(\nu(b_j))+\max(0,m-2),\text{ respectively.}
\end{align*}
\end{cor}

The same argument as in the proof of Theorem \ref{thmmfd} gives

\begin{thm} \label{8degentrans.thm}
For a dense open set of $h\in\ham_0(L,\delta,s)$,
$\Gamma^{-1}(0)\cap\proj^{-1}(h)\cap\cand_{2,\epsilon,\Lambda}({\mathbf
c})$ and
$\Gamma^{-1}(0)\cap\proj^{-1}(h)\cap{\widetilde\cand}_{2,\epsilon,\Lambda}
({\mathbf c})$
are finite dimensional manifolds of dimensions given by the
dimension formula in Corollary~\ref{8degenexpotrans.cor},
provided this dimension is $\le 1$ if $l({\mathbf c})\ge 3$ and
$\le 1+(3-l({\mathbf c}))$ otherwise.  
\end{thm}

\begin{rmk}\label{8nodegenexpodisks.rmk}
Note that the expected dimension of the set of disks with dimension
count in $\tilde\cand_{2,\epsilon}$ equal to $1$ in 
$\cand_{2,\epsilon}$ is
equal to $-k$ (or $-k+(3-l({\mathbf c}))$ if $l({\mathbf c}\le 2)$),
where $k$ is the number of punctures mapping to the self-tangency
Reeb-chord. Therefore for a dense open set in the space of chord
semi-generic Legendrian submanifolds this space is empty. Since
any disk with exponential decay at the self-tangency point has a
neighborhood in $\cand_{2,\epsilon}$, we see that generically such
disks do not exist, provided their dimension count in
$\tilde\cand_{2,\epsilon}$ is as above. 
\end{rmk}

\subsection{Enhanced transversality}
\label{8enhancedtrans.section}

Let $L$ be a (semi-)admissible
submanifold. If $q\in L$ and $\zeta_0\in \pa D_m$ then define 
\begin{equation} \notag
\cand_{2,\epsilon}({\bf c},\zeta_0,p)=
\{(w,f)\in\cand_{2,\epsilon}({\bf c})\colon (w,f)(\zeta_0)=p\}
\end{equation}
and in the semi-admissible case also 
${\widetilde\cand}_{2,\epsilon}({\bf c},\zeta_0,p)$
in a similar way.

If $\ev_{\zeta_0}\colon\cand_{2,\epsilon}({\bf c})\to L$ denotes the map
$\ev_{\zeta_0}(w,f)=(w,f)(\zeta_0)$. Then $\ev_{\zeta_0}$ is 
smooth and transverse to $p$ (as is seen by using local coordinates
on $\cand_{2,\epsilon}({\mathbf c})$). Moreover,
$\ev_{\zeta_0}^{-1}(p)=\cand_{2,\epsilon}({\bf c},p)$ and hence 
$\cand_{2,\epsilon}({\bf c}, p)$ is a closed submanifold of
$\cand_{2,\epsilon}({\bf c})$ of codimension $\dim(L)$. 
Note that the tangent space
$T_{(w,f)}\cand_{2,\epsilon}({\mathbf c},p,\zeta_0)$ is
the closed subspace of elements $(v,\gamma)$ in the tangent space 
$T_{(w,f)}\cand_{2,\epsilon}({\mathbf c})$ which are such that
$v\colon D_m\to\C^n$ satisfies $v(\zeta_0)=0$.   

We consider  
\begin{equation} \notag
\cand_{2,\epsilon}({\mathbf c},p)=
\bigcup_{\zeta\in\pa D_m}\cand_{2,\epsilon}({\mathbf c},\zeta,p)
\end{equation}
as a locally trivial bundle over $\pa D_m$. Local
trivializations are given compositions with suitable diffeomorphisms
which move the boundary point $\zeta$ a little. 

We define perturbation spaces as the closed subspaces 
$\ham^p(L,\delta,R)\subset\ham(L,\delta,R)$  
and $\ham_0^p(L,\delta)\subset\ham_0(L,\delta)$ of functions $h$ 
such that $h(p)=0$ and $Dh(p)=0$. 
Thus, $\tilde\Phi_h$ fixes $p$. (Note
that if $p$ is the projection of a Reeb chord this is no {\em
additional} restriction.) If $\Lambda$ denotes one of these
perturbation spaces we form the 
bundles  
\begin{align}
\notag
\cand_{2,\epsilon,\Lambda}({\mathbf c},p)=
\bigcup_{L_\lambda,\lambda\in\Lambda}\cand_{2,\epsilon}({\mathbf
c},p),\\
\notag
{\widetilde\cand}_{2,\epsilon,\Lambda}({\mathbf c},p)=
\bigcup_{L_\lambda,\lambda\in\Lambda}\cand_{2,\epsilon}({\mathbf
c},p)
\end{align}
with local coordinates as before.

As before let $'$ denote exclusion of exceptional holomorphic maps.

\begin{lma}\label{lmatv'p}
Assume that $p\in L$ has a neighborhood $U$ such that  $\Pi_\C(U)$ is real
analytic. Then the bundle maps 
\begin{align}\notag 
(\Gamma,\proj)\colon\cand_{2,\epsilon,\Lambda}'({\mathbf c},p)\to
\sblv_{1,\epsilon,\Lambda}(D_m,{T^{0,1}}^\ast D_m\otimes\C^n)\\\label{Gp}
(\Gamma,\proj)\colon{\widetilde\cand}_{2,\epsilon,\Lambda}'({\mathbf c},p)\to
\sblv_{1,\epsilon,\Lambda}(D_m,{T^{0,1}}^\ast D_m\otimes\C^n)
\end{align}
are transverse to the $0$-section.    
\end{lma}

\begin{pf}
The proof is the same as the proof of Lemma \ref{lmatv'} in the admissible
case and the same as
that of Lemma \ref{lmasgtv'} in the semi-admissible case provided the arcs
$\gamma$ and $\gamma'$ used there do not contain the special point
$p$. On the other hand, if one of these arcs does contain $p$ we may
shorten it until it does not. 
(The key point is that the condition that the Hamiltonian vanishes at
a point does not destroy the approximation properties of the elements
in the perturbation space for smooth functions supported away from
this point). 
\end{pf}

\begin{cor}\label{lmagood}
Let $n>1.$ For $L$ in a Baire subset of the space of (semi-)admissible 
Legendrian $n$-submanifolds no rigid holomorphic disk passes
through the end point of any Reeb chord of $L$.
\end{cor}

Note, when $n=1$ this corollary is not true.

\begin{pf}
The proof of Theorem \ref{thmmfd} shows that for a Baire set
there are no exceptional holomorphic disks. The Sard-Smale theorem in
combination 
with Lemma \ref{lmatv'p} implies that for a Baire subset of this Baire
set the dimension of the space of rigid holomorphic disks with
some point mapping to the end point of a specific Reeb chord is given
by the Fredholm index of the operator $d\Gamma$ corresponding to
$\Gamma$ in \eqref{Gp}. Since the
source space of this operator is the sum of a copy of $\R$ (from the
movement of $\zeta$ on the boundary) and a closed codimension $\dim L$
subspace of the source space of $d\Gamma$ in Lemma \ref{lmatv'} which
has minimal index for disks to appear generically, we see the index in
the present case is too small. This 
implies that the subset is generically empty. Taking the intersection
of these Baire subsets for the finite collection of Reeb chord
endpoints of $L$ we get a Baire subset with the required properties.
\end{pf}

\begin{cor} \label{8injectivedegen.cor}
If $L$ is as in Corollary~\ref{lmagood} then 
there are no rigid holomorphic disks with boundary on $L$ which are
nowhere injective on the boundary. 
\end{cor}
\begin{pf}
Let $w\colon D_{m+1}\to\C^n$ represent a holomorphic disk with boundary on
$L$. By Corollary ~\ref{lmagood} we may assume that no point in the
boundary of $\pa D_m$ maps to an intersection point of $\Pi_\C(L)$.

Assume that $w$ has no injective point on the boundary and let the
punctures of $D_{m+1}$ map to the Reeb chords $(c_0,\dots,c_m)$ where
the positive puncture maps to $c_0$. 
Let $C$ be the holomorphic chain which is the closure of image $w(D_m)$ of
$w$ with local multiplicity $1$ everywhere.  
Then 
\begin{equation}\label{a(C)<a(w)}
\area(C)<\area(w)
\end{equation}
since close 
to the point in $C$ most distant from the origin in $\C^n$, $w$ has
multiplicity at least two. 

The corners of $C$ is a subset $S$ of ${c_0^\ast,\dots,c_m^\ast}$ and by
integrating $\sum_j y_jdx_j$ along the boundary $\pa C$ of $C$ which
lies in the exact Lagrangian $\Pi_\C(L)$ we find
\begin{equation}
\area(C)=\action(c_0)-\sum_{c_j^\ast\in S, j>0}\action(c_j),
\end{equation} 
where the first term must be present (i.e. $C$ must have a corner at
$c_0^\ast$) since otherwise the area of $C$ would be negative
contradicting the fact that $C$ is holomorphic. On the other hand
\begin{equation}
\area(w)=\action(c_0)-\sum_{j>0}\action(c_j).
\end{equation}
Hence
\begin{equation}
\area(C)\ge\area(w),
\end{equation}
which contradicts \eqref{a(C)<a(w)}. This contradiction finishes the
proof. 
\end{pf}

\subsection{Transversality in a split problem}

In this section we discuss transversality for disks, with one or two
punctures, lying entirely in one complex coordinate plane. Let
$L\subset\C^n\times\R$  be an admissible Legendrian submanifold. Let
$\Delta\subset\R^2$ denote the standard simplex. Let $\Delta_1$
($\Delta_2$) be the subsets of $\R^2$ which is bounded by $\pa\Delta$,
smoothened at two (one) of its corners. Let $(z_1,\dots,z_n)$ be
coordinates on $\C^n$. Let $\pi_i\colon\C^n\to\C$ denote projection to
the $i$-th coordinate and let $\hat\pi_i\colon\C^n\to\C^{n-1}$ denote
projection to the Hermitian complement of the $z_i$-line. Finally, if
$\gamma(t)$, $t\in I\subset\R$ is a one parameter family of lines then
we let $\int_{\gamma}d\theta$ denote the (signed) angle $\gamma(t)$
turns as $t$ ranges over $I$.

\begin{lem}\label{splittran1}
Let $(u,h)\in\cand_{2}(ab;A)$, $\mu(A)+|a|-|b|=1$, be a holomorphic disk
with boundary on $L$ such that 
$\hat\pi_1\circ u$ is constant and such that $\pi_1\circ u=f\circ g$,
where $g\colon\Delta_2\to D_2$ is a diffeomorphism and
$f\colon\Delta_2\to\C$ is an immersion. Furthermore, if $t_1,t_2$ are
coordinates along components of $\pa D_2$, assume that the paths 
$\Gamma(t)=d\Pi_\C(T_{(u,h)(t_j)}L)$ of Lagrangian subspaces are
split: $\Gamma(t_j)=\gamma(t_j)\times \hat V_j$, where
$\gamma(t)\subset\C$ is a (real line) and $\hat V_j\subset\C^{n-1}$,
$j=1,2$, are transverse Lagrangian subspaces. Then $d\Gamma_{(u,h)}$ is
surjective. (In other words, $(u,h)$ is transversely cut out).
\end{lem}

\begin{proof}
The Fredholm index of $d\Gamma$ at $(u,h)$ equals $1$. If $v$ is the
vector field on $D_2$ which generates the $1$-parameter family of
conformal automorphisms of $D_2$ (the vector field $\pa_\tau$ in
coordinates $\tau+it\in\R\times[0,1]$ on $D_2$) then $\xi=du\cdot v$ lies in
the kernel of $d\Gamma$ and $d\hat\pi\cdot\xi=0$. 

Since the boundary conditions are split we may consider them
separately. It follows from Section \ref{7modelproblem2.section}
that the $\hat\pi_1 d\Gamma$ with 
boundary conditions given by the two transverse Lagrangian subspaces
$\hat V_1$ and $\hat V_2$ has index $0$, no kernel and no cokernel.  

Let $\theta_1$ and $\theta_2$ be the interior angles at the corners of
the immersion $f$. Since $f(\pa\Delta_2)$ bounds an immersed disk we have
$$
\int_{\gamma_1}d\theta+\int_{\gamma_2}d\theta+
(\pi-\theta_1)+(\pi-\theta_2)=2\pi.  
$$
If $\eta_1=\pi_1\circ\eta$, where $\eta$ is in the kernel of $d\Gamma$
then, thinking of $D_2$ as $\R\times[0,1]$, we find that,
asymptotically, for some integers $n_1\ge 0$ and $n_2\ge 0$
\begin{equation*}
\eta_1(\tau+it)=\begin{cases}
c_1e^{-(\theta_1+n_1\pi)(\tau+it)},& \text{ for }\tau\to+\infty,\\
c_2e^{(\theta_2+n_2\pi)(\tau+it)}, &\text{ for }\tau\to-\infty, 
\end{cases}
\end{equation*}
where $c_1$ and $c_2$ are real constants. Cutting $D_2$ off at
$|\tau|=M$  for some sufficiently large $M$ we thus find a solution
of the classical Riemann-Hilbert problem with Maslov-class
$$
\frac{1}{\pi}\Bigl(\theta_1+\theta_2-\theta_1-\theta_2-(n_1+n_2)\pi\Bigr).
$$
Since the classical Riemann-Hilbert problem has no solution if the
Maslov class is negative and exactly one if it is $0$ we see that
the solution $\xi=\xi_1$ produced above is unique up to multiplication
with real constants.
\end{proof}

\begin{lem}\label{splittran2}
Let $(u,h)\in\cand_{2}(a;A)$, $\mu(A)+|a|=1$, be a holomorphic disk
with boundary on $L$ such that 
$\hat\pi_1\circ u$ is constant and such that $\pi_1\circ u=f\circ g$,
where $g\colon\Delta_1\to D_2$ is a diffeomorphism and
$f\colon\Delta_1\to\C$ is an immersion. Furthermore, if $t$ is a
coordinate along $\pa D_1$, assume that the path 
$\Gamma(t)=d\Pi_\C(T_{(u,h)(t_j)}L)$ of Lagrangian subspaces is 
split:
$\Gamma(t)=\gamma_1(t_j)\times\gamma_2(t)\times\dots\times\gamma_n(t)$,
where $\gamma_j(t)\subset\C$ is a (real line) such that 
$$
\int_{\gamma_j}d\theta < 0,\text{ for }2\le j\le n.
$$ 
Then $d\Gamma_{(u,h)}$ is surjective.
\end{lem}

\begin{proof}
The proof is similar to the one just given. Using asymptotics and the
classical Riemann-Hilbert problem it follows that the kernel of
$d\Gamma$ is spanned by two linearly independent solutions $\xi^j$,
$j=1,2$, with $\hat\pi_1\xi^j=0$.  
\end{proof}

\subsection{Auxiliary tangent-like spaces in the semi-admissible
case}
\label{8auxspace.section}

Let $L$ be a chord semi-admissible Legendrian
submanifold and assume that $L$ lies in the open subset of such
manifolds where the moduli-space of rigid holomorphic disks 
with corners at ${\mathbf c}$ is
$0$-dimensional (and compact by Theorem \ref{9cpt.thm}). 
Now if $(w,f)$ is a
holomorphic disk with boundary on $L$ then by Lemma~\ref{lmasgtv'} we know that  
the operator 
\begin{equation}\label{truedG}
d\Gamma\colon T_{(w,f)}{\widetilde\cand}_{2,\epsilon}({\mathbf c})\to
\sblv_{1,\epsilon}(D_m,T^\ast D_m\otimes\C^n)
\end{equation}
is surjective.

For any $(w,f)$ with $m+1$ punctures which maps the punctures
$p_1,\dots p_k$ to the self tangency Reeb 
chord of $L_h$ let $\hat\epsilon\in
[0,\infty)^{m+1-k}\times(-\delta,0)^k$, where $\delta>0$ is small 
compared to the complex angle of the self tangency Reeb chord and the
components of $\hat\epsilon$ which are negative correspond to the
punctures $p_1,\dots,p_k$. Define
the {\em tangent-like space}  
\begin{equation} \notag
T_{(w,f,h)}\cand_{2,\hat\epsilon}({\mathbf c})
\end{equation} 
as the linear space of elements $(v,\gamma)$ where 
$\gamma\in T_\kappa\conf_{m+1}$ and where
$v\in\sblv_{2,\hat\epsilon}(D_{m+1},\C^n)$ satisfies
\begin{align}
\notag
v(\zeta)\in \Pi_\C(T_{(w,f)(\zeta)}L)&\text{ for all }\zeta\in\pa D_m,\\
\notag
\int_{\pa D_m}\la\bar\pa v,u\ra\,ds=0 &\text{ for all
}u\in\C^\infty_0(\pa D_m,\C^n).
\end{align} 
and consider the linear operator
\begin{equation}\label{auxdG}
d\hat\Gamma(v,\gamma)=\bar\pa_\kappa v+i\circ dw\circ\gamma.
\end{equation}
The index of this Fredholm operator equals that of the operator in
\eqref{truedG} and moreover by asymptotics of solutions to these
equations (close to the self-tangency Reeb chord we can use the same
change of coordinates in the first coordinate as in the non-linear
case, see Section~\ref{2asymp} to determine the behavior of solutions) we find
that the kernels are canonically isomorphic. Thus,  
since the operator in \eqref{truedG} is surjective so is the operator
in \eqref{auxdG}.

\section{Gluing theorems}\label{10glue.section}

In this section we prove the gluing theorems used in Sections
\ref{1modulisection} and \ref{1Invariance.section}.
In Section~\ref{10statments}
we state the theorems.
Our general method of gluing curves is the standard one in
symplectic geometry. However, some of our specific gluings require
a significant amount of analysis.
We first ``preglue'' the pieces of the broken curves together. 
For the stationary case this is done in Section \ref{10statpreglue.section}, 
for the handle slide case in Section \ref{10hspreglue.section}, 
and for the self-tangency case in Sections \ref{10stpreshort.section} and 
\ref{10stpreglue.section}.
We then apply Picard's Theorem, stated in Section \ref{10Picard.section}.
Picard's Theorem requires a sequence of uniformly bounded right inverses
of the linearized $\dbar$ map.
We prove the bound for the stationary case in Section \ref{10unistat.section},
for the handle slide case in Section \ref{10unihs.section}, 
and for the self-tangency case in Sections \ref{10stshortuni.section} and 
\ref{10stuni.section}.
Picard's Theorem also requires a bound on the non-linear part
of the expansion of $\dbar$, which we do in Section \ref{10nonlinear.section}.
To handle disks with too few boundary punctures, we show
in Sections \ref{ssecmarked} through \ref{manymarked}, 
how by marking boundary points the disks
can be thought of as sitting inside a moduli space of disks
with many punctures.

Recall the following notation.
Bold-face letters will denote ordered collections of Reeb chords. If
${\mathbf c}$ denotes a non-empty ordered collection 
$(c_1,\dots,c_m)$ of Reeb chords then we say that the {\em length} of
${\mathbf c}$ is $m$. We say that the length of the empty ordered
collection is $0$.
Let ${\mathbf c^1},\dots,{\mathbf c^r}$ be an ordered collection of
ordered collections of Reeb chords. Let the length of ${\mathbf c^j}$
be $l(j)$
and let ${\mathbf a}=(a_1,\dots,a_k)$ be an ordered collection of
Reeb chords of length $k>0$. Let $S=\{s_1,\dots,s_r\}$ be $r$ distinct
integers in $\{1,\dots,k\}$. Define the ordered collection
${\mathbf a}_S({\mathbf c^1},\dots,{\mathbf c^r})$ of Reeb chords of
length $k-r+\sum_{j=1}^rl(j)$ as 
follows. For each index $s_j\in S$ remove $a_{s_j}$ from the ordered
collection ${\mathbf a}$ and insert at its place the ordered
collection ${\mathbf c^j}$. 

Recall that if $a$ is a Reeb chord and ${\mathbf b}$ is a collection of Reeb
chords of a Legendrian submanifold, then $\M_A(a;{\mathbf b})$ denotes
the moduli space of holomorphic disks with boundary on
$L$, punctures mapping to $(a,{\mathbf b})$, and boundary in $L$ which
after adding the chosen capping paths represents the homology class
$A\in H_1(L)$. After Theorem~\ref{thmmfd}
we know that if the 
length of ${\mathbf b}$ is at least $2$ then $\M_A(a;{\mathbf b})$ is
identified  
with the inverse image of the regular value $0$ of the $\bar\pa$-map
$\Gamma$ in Section \ref{fcanH}. If the length of ${\mathbf b}$ is $0$
or $1$ then $\M_A(a;{\mathbf b})$ is identified with the quotient of
$\Gamma^{-1}(0)$ under the group of conformal reparameterizations of
the source of the holomorphic disk.
 
Similarly if $L_\lambda$, $\lambda\in\Lambda$ is a $1$-parameter
family of chord generic Legendrian submanifolds we write
$\M^\Lambda_A(a;{\mathbf b})$ for the parameterized moduli space of
rigid holomorphic disks with boundary in 
$L_\lambda$, and punctures at $(a(\lambda),{\mathbf b}(\lambda))$,
$\lambda\in\Lambda$. We also write $\M^\lambda_A(a,{\mathbf b})$ to
denote the moduli space for a fixed $L_\lambda$, $\lambda\in\Lambda$.

Finally if $K\subset \C^n$ and $\delta>0$ then $B(K,\delta)$ denotes
the subset of all points in $\C^n$ of distance less than $\delta$ from
$K$. 

\subsection{The Gluing Theorems}\label{10statments}
In this section we state the various gluing theorems.

\subsubsection{Stationary gluing}\label{10statglue.section}

Let $L$ be an admissible Legendrian submanifold. Recall
that a holomorphic 
disk with boundary on $L$ is defined as a pair of functions $(u,f)$,
where $u\colon D_m\to\C^n$ and $f\colon\pa D_m\to\R$. Below we will
often drop the function $f$ from the notation and speak of the
holomorphic disk $u$. Let
$\M_A(a;{\mathbf b})$ and $\M_C(c;{\mathbf d})$ be moduli spaces of 
rigid holomorphic disks, where ${\mathbf b}$ has length $m$, $1\le
j\le m$, and ${\mathbf d}$ has length $l$.

\begin{thm}\label{glud^2=0}
Assume that the $j$-th Reeb chord in ${\mathbf b}$ equals $c$.
Then there exists $\delta>0$, $\rho_0>0$ and an embedding 
\begin{align}\notag
\M_A(a;{\mathbf b})
\times\M_C(c;{\mathbf d})
\times[\rho_0,\infty)&\to
\M_{A+C}(a;{\mathbf b}_{\{j\}}({\mathbf d}));\\\notag
(u,v,\rho)&\mapsto u\,\sharp_\rho v,
\end{align} 
such that if $u\in\M_A(a;{\mathbf b})$ and 
$v\in \M_C(c;{\mathbf d})$ and the image of
$w\in\M_{A+C}(a;{\mathbf b}_j({\mathbf d}))$ lies
inside $B(u(D_{m+1})\cup v(D_{l+1});\delta)$ then 
$w=u\,\sharp_\rho v$ for some $\rho\in[\rho_0,\infty)$.
\end{thm}

\begin{proof}
The theorem follows from Lemmas \ref{lmaprestat}, \ref{statunif}, and 
\ref{statnl} and Proposition \ref{prpFP}.
\end{proof}

\subsubsection{Handle slide gluing}\label{sec12}

Let $L_\lambda$, $\lambda\in(-1,1)=\Lambda$ be a $1$-parameter family
of admissible Legendrian submanifolds such that 
\begin{equation}\label{eqhanddi}
\M^\Lambda_A(a;{\mathbf b})=\M^0_A(a;{\mathbf b})
\end{equation}
is a transversely cut out handle slide disk, represented by a map
$u\colon D_{m+1}\to \C^n$ (the length of
${\mathbf b}$ is $m$) and such that for all $\lambda\in\Lambda$, all 
moduli-spaces of rigid disks with boundary on $L_\lambda$ are
transversely cut out. 

We formulate two gluing theorems in this case. They differ in the
following way. In Theorem \ref{gluhandslide1} we consider what happens
when the positive punctures of rigid disks are glued to negative
punctures of the handle slide disk. 
In Theorem
\ref{gluhandslide2} we consider what happens when the positive
puncture of the handle-slide disk is glued to negative punctures of
a rigid disk.

\begin{thm}\label{gluhandslide1}
Assume that ${\mathbf b}$ in \eqref{eqhanddi} has positive length and
has $c$ in its $j$-th position.
Let $\M^0_C(c;{\mathbf d})$ be a moduli space
of rigid holomorphic disks, where the length of ${\mathbf d}$ is
$k$. Then there exist $\delta>0$, $\rho_0>0$, and an embedding
\begin{align}\notag
\M^0_A(a;{\mathbf b})\times\M^0_C(c;{\mathbf d})\times[\rho_0,\infty)
&\to
\M^\Lambda_{A+C}(a;{\mathbf b}_{\{j\}}({\mathbf d})),\\\notag
(u,v,\rho) &\mapsto u\,\sharp_\rho v,
\end{align}
such that if
$v\in\M^0_C(c;{\mathbf d})$ and the image of 
$w\in\M^\Lambda_{A+C}(a;{\mathbf b}_{\{j\}}({\mathbf d}))$
lies inside $B(u(D_{m+1})\cup v(D_{k+1});\delta)$ then 
$w=u\,\sharp_\rho v$ for some $\rho\in[\rho_0,\infty)$.
\end{thm}

\begin{proof}
The theorem follows from Lemmas \ref{lmaprehsl1},
\ref{handsl1unif}, and \ref{handslnl} 
and Proposition \ref{prpFP}.
\end{proof}

\begin{thm}\label{gluhandslide2}
Let $\M^0_C(c;{\mathbf d})$ be a moduli space
of rigid holomorphic disks, where ${\mathbf d}=(d_1,\dots,d_l)$. Let
$S=\{s_1,\dots,s_r\}\subset\{1,\dots,l\}$ be such that $d_{s_j} = a$ 
for all $s_j \in S.$
Then there exist $\rho_0>0$, $\delta>0$, and an embedding 
\begin{align}\notag
\M^0_C(c;{\mathbf d})\times\Pi_S\M^0_A(a;{\mathbf b})
\times[\rho_0,\infty)&\to  
\M^\Lambda_{C+r\cdot A}(c;{\mathbf d}_{S}({\mathbf b},\dots,{\mathbf
b})),\\\notag 
(v,\underbrace{u,\dots,u}_r,\rho)&\mapsto v\,\sharp_\rho^S u,
\end{align}
such that if
$v\in \M^\Lambda_{C}(c;{\mathbf d})$ 
and the image of
$w\in\M^\Lambda_{C+r\cdot A}(c;{\mathbf d}_{S}({\mathbf
b},\dots,{\mathbf b}))$ lies inside  
$B(v(D_{l+1})\cup u(D_{m+1}) \cup \ldots \cup u(D_{m+1}));\delta)$ then 
$w=v\,\sharp_\rho^S u$ for some $\rho\in[\rho_0,\infty)$.
\end{thm}

\begin{proof}
The theorem follows from Lemmas \ref{lmaprehs2},
\ref{handsl2unif}, and \ref{handslnl} 
and Proposition \ref{prpFP}.
\end{proof}

\subsubsection{Self tangency shortening and self tangency
gluing}\label{10stglue.section} 

Let $L_\lambda$, $\lambda\in(-1,1)=\Lambda$ be an admissible 
$1$-parameter family of Legendrian submanifolds such that $L_0$ 
is semi-admissible with self-tangency Reeb chord $a$. For
simplicity (see Section \ref{5Isotopies.section})
we assume that all Reeb chords outside a
neighborhood of $a$ remain fixed under $\Lambda$. We take
$\Lambda$ so that if $\lambda>0$ then $L_{-\lambda}$ has two 
new-born Reeb-chords $a^+$ and $a^-$, where 
$\action (a^+)> \action(a^-)$.  
Assume that all moduli spaces of rigid holomorphic disks with boundary
on $L_\lambda$ are
transversely cut out for all fixed $\lambda\in\Lambda$, that for
all $\lambda\in\Lambda$, there are no disks with negative formal dimension, and 
that all rigid disks with a puncture at $a$ satisfy the non-decay
condition of Lemma \ref{SalRob} (see Remark \ref{8nodegenexpodisks.rmk}).

\begin{thm}\label{shoslft}
Let $\Lambda^-=(-1,0)$. Let $\M^0_A(a,{\mathbf b})$ be a moduli space
of rigid holomorphic disks where the length of ${\mathbf b}$ is $l$.
Then there exist $\rho_0>0$, $\delta>0$ and a local homeomorphism
\begin{align}\notag
\M^0_A(a;{\mathbf b})
\times[\rho_0,\infty)&\to
\M^{\Lambda^-}_A(a^+;{\mathbf b});\\\notag
(u,\rho)&\mapsto \sharp_\rho u,
\end{align}
such that if
$u\in\M^0_A(a;{\mathbf b})$ and the image of 
$w\in\M^{\Lambda^-}_A(a^+;{\mathbf b})$
lies inside $B(u(D_{l+1});\delta)$ then 
$w=\sharp_\rho u$ for some $\rho\in[\rho_0,\infty)$.

Let $\M^0_C(c,{\mathbf d})$ be a moduli 
space of rigid holomorphic disks where the length of ${\mathbf d}$ is
$m$. Let $S\subset\{1,\dots,m\}$ be the subset of positions of ${\mathbf d}$
where the Reeb chord $a$ appears (to avoid trivialities, 
assume $S\ne \emptyset$). Then there exists $\rho_0>0$ and $\delta>0$
and a local homeomorphism
\begin{align}\notag
\M^0_C(c,{\mathbf d})
\times[\rho_0,\infty)&\to
\M^{\Lambda^-}_C(c,{\mathbf d}_S(a^-));\\\notag
(u,\rho)&\mapsto \sharp_\rho u,
\end{align} 
such that if $u\in\M^0_C(c;{\mathbf d})$ and the image of 
$w\in\M^{\Lambda^-}_C(c;{\mathbf d}_S(a^-))$
lies inside $B(u(D_{m+1});\delta)$ then 
$w=\sharp_\rho u$ for some $\rho\in[\rho_0,\infty)$.
\end{thm}

\begin{proof}
Consider the first case, the second follows
by a similar argument. 
Applying Proposition \ref{prpFP} and Lemmas \ref{stpresh}, \ref{stshunif} and
\ref{stshnl} we find a homeomorphism 
$\M^0_A(a;{\bf b})\to\M^{\lambda_-}_A(a^+,{\mathbf b})$ for $\lambda^-<0$
small enough. The proof of Corollary \ref{cormfd'} implies
that $\M_{\Lambda^-}(a^+,{\mathbf b})$ is a $1$-dimensional manifold
homeomorphic to $\M_{\lambda_-}(a;{\mathbf b})\times\Lambda_-$, the
theorem follows. 
\end{proof}

\begin{thm}\label{gluslft}
Let $\Lambda^+=(0,1)$ and
let $\M^0_{A_1}(a;{\mathbf b^1}),\dots,\M^0_{A_r}(a;{\mathbf b^r})$ and
$\M^0_C(c;{\mathbf d})$ be a moduli spaces of
rigid holomorphic disks where the length of ${\mathbf b^j}$ is $l(j)$,
and the length of ${\mathbf d}$ is $m$. Let $S\subset\{1,\dots,m\}$ be
the subset of positions of ${\mathbf d}$ where the Reeb chord $a$
appears and assume that $S$ contains $r$ elements. Then there exists
$\delta>0$, $\rho_0>0$ 
and an embedding 
\begin{align}\notag
\M_C^0(c;{\mathbf d})\times
\Pi_{j=1}^r\M^0_{A_j}(a;{\mathbf b^j})
\times[\rho_0,\infty)&\to
\M^{\Lambda^+}_{C+\sum_{j}A_j}(c;{\mathbf d}_S({\mathbf b^1},\dots,{\mathbf b^r}));\\\notag
(v,u_1,\dots,u_r,\rho)&\mapsto v\,\sharp_\rho(u_1,\dots, u_r), 
\end{align} 
such that if $v\in\M^0_C(c;{\mathbf d})$ and 
$u_j\in\M^0_{A_j}(a;{\mathbf b^j})$, $j=1,\dots,r$ and  
the image of
$w\in\M^{\Lambda^+}_{C+\sum_j A_j}(c;{\mathbf d}_{S}({\mathbf b^1},\dots,{\mathbf
b^r}))$ lies inside 
$B(v(D_{m+1})\cup u_1(D_{l(1)+1})\cup\dots\cup u_r(D_{l(r)+1}));\delta)$ 
then $w=v\,\sharp_\rho (u_1,\dots, u_r)$ for some 
$\rho\in[\rho_0,\infty)$.
\end{thm}

\begin{proof}
Apply Proposition \ref{prpFP} and Lemmas \ref{stpreglu},
\ref{stgluunif}, and 
\ref{stglunl} and reason as above.
\end{proof}

\subsection{Floer's Picard lemma}\label{10Picard.section}

The proofs of the theorems stated in the preceding subsections are
all based on the following.
\begin{prp}\label{prpFP}
Let $f\colon B_1\to B_2$ be a smooth map between Banach spaces which
satisfies
\begin{equation}\notag
f(v)=f(0)+df(0)v + N(v),
\end{equation}
where $df(0)$ is Fredholm and has a right inverse $G$ satisfying 
\begin{equation}\notag
\|GN(u)-GN(v)\|\le C(\|u\|+\|v\|)\|u-v\|,
\end{equation}
for some constant $C$. Let $B(0,\epsilon)$ denote the $\epsilon$-ball
centered at $0\in B_1$ and assume that 
\begin{equation}\notag
\|Gf(0)\|\le\frac{1}{8C}.
\end{equation}
Then for $\epsilon<\frac{1}{4C}$, the zero-set of 
$f^{-1}(0)\cap B(0,\epsilon)$ is a smooth submanifold of dimension
$\dim(\krn(df(0)))$ diffeomorphic to the $\epsilon$-ball
in $\krn(df(0))$.
\end{prp}
\begin{pf}
See \cite{Floer87}.
\end{pf}

In our applications of Proposition \ref{prpFP} the map
$f$ will be the $\bar\pa$-map, see Section \ref{fcanH}. 

\subsection{Notation and cut-off functions}\label{sectcut-off}
To simplify notation, we deviate slightly from our standard
notation for holomorphic disks. We use the convention that the neighborhood
$E_{p_0}$ of the positive puncture $p_0$ in the source $D_m$ of a
holomorphic disk $(u,f)$ will be 
parameterized by $[1,\infty)\times[0,1]$ and that neighborhoods of
negative punctures $E_{p_j}$, $j\ge 1$ are parameterized by
$(-\infty,-1]\times[0,1]$.
 
In the constructions and proofs below we will use certain cut-off
functions repeatedly. Here we explain how to construct them.
Let $K>0$, $a<b$, and let $\phi\colon [a,b+K+1]\to[0,1]$ be a smooth
function which equals $1$ on $[a,b]$ and equals $0$ in
$[b+K,b+K+1]$. It is easy to see there exists such functions with
$|D^k\phi|=\Ordo(K^{-k})$ for $k=1,2$.
Let $\epsilon>0$ be small. Let $\psi\colon[0,1]\to\R$ be a smooth
function such that 
$\psi(0)=\psi(1)=0$, $\psi'(0)=\psi'(1)=1$, with $|\psi|\le\epsilon$.  
We will use cut-off functions $\alpha\colon[a,b+K+1]\times[0,1]\to\C$
of the form
\begin{equation}\notag
\alpha(\tau+it)=\phi(\tau)+i\psi(t)\phi'(\tau).
\end{equation}
Note that $\alpha|\pa([a,b+K+1]\times[0,1])$ is real-valued and
$\bar\pa\alpha=0$ on $\pa([a,b+K+1]\times[0,1])$. Also, 
$|D^k\alpha|=\Ordo(K^{-1})$ for $k=1,2$. 

\subsection{A gluing operation}\label{10glueop.section}

Let $L$ be a chord generic Legendrian submanifold.
Let $(u,f)\in\cand_{2,\epsilon}(a,{\mathbf b})$ where 
${\mathbf b}=(b_1,\dots,b_m)$ and let
$(v_j,h_j)\in\cand_{2,\epsilon}(b_j,{\mathbf c}^j)$, $j\in
S\subset\{1,\dots,m\}$. Denote the punctures on $D_{m+1}$ by $p_j$,
$j=0,\dots,m$ and the positive puncture on $D_{l(j)+1}$ by $q_j$.   
 
Let $g^\sigma$, $\sigma\in[0,1]$ be a $1$-parameter family of metrics
as in Section \ref{fcanC}. Then for $M>0$ large enough there exists
unique functions  
\begin{align}\notag
&\xi\colon E_{p_j}[-M]\to T_{b_j^\ast}\C^n\\\notag
&\eta_j\colon E_{q_j}[M]\to T_{b_j^\ast}\C^n,
\end{align}
such that
\begin{align}\notag
\exp_{b_j^\ast}^{t}(\xi(\tau+it)) &= u(\tau+it),\\\notag
\exp_{b_j^\ast}^{t}(\eta_j(\tau+it)) &= v_j(\tau+it),
\end{align}
where $\exp^\sigma$ denotes the exponential map of the metric
$g^\sigma$. Note that by our special choice of metrics the functions
$\xi$ and $\eta$ are tangent to $\Pi_\C(L)$ and holomorphic on the
boundary. 

For large $\rho>0$, let $D_r^S(\rho)$, $r=1+m+\sum_{j\in S}(l(j)-1)$ be
the disk obtained by gluing to the end of  
\begin{equation}\notag
D_{m+1}\setminus\bigcup_{j\in S} E_{p_j}[-\rho]
\end{equation}
corresponding to $p_j$ a copy of 
\begin{equation}\notag
D_{l(j)+1}-E_{q_j}[\rho]
\end{equation}
by identifying $\rho\times[0,1]\subset E_{p_j}$ with
$-\rho\times[0,1]\subset E_{q_j}$, for each $j\in S$. Note that
the metrics (and the complex structures $\kappa_1$ and $\kappa_2(j)$) on
$D_{m+1}$ and $D_{l(j)+1}$ 
glue together to a unique metric (and complex structure $\kappa_\rho$)
on $D_{r}^S(\rho)$. 
We consider $D_{m+1}\setminus \bigcup_{j\in S} E_{p_j}[-\rho]$ and
$D_{l(j)+1}\setminus E_{q_j}[\rho]$ as subsets of $D_r^S(\rho)$.

For $j \in S$, let $\Omega_j\subset D_r^S(\rho)$ denote the subset
\begin{equation}\notag
E_{q_j}[\rho-2,\rho]\cup E_{p_j}[-\rho,-\rho+2] \approx [-2,2]\times[0,1]
\end{equation} 
of $D_r^S(\rho)$. Let $z=\tau+it$ be a complex coordinate on $\Omega_j$ and let
$\alpha^\pm\colon\Omega_j\to\C$ be cut-off functions which are real valued
and holomorphic on the boundary and with 
$\alpha^+=1$ on $[-2,-1]\times[0,1]$, $\alpha^+=0$ on 
$[0,2]\times[0,1]$, $\alpha^-=1$ on $[1,2]\times[0,1]$, and
$\alpha^-=0$ on $[-2,0]\times[0,1]$.
Define the function $\Sigma_\rho^S(u,v_1,\dots,v_r)\colon
D_r\to\C^n$ as 
\begin{equation}\notag
\Sigma_\rho^S(u,v_1,\dots,v_r)(\zeta)=
\begin{cases}
&v_j(\zeta),\quad\quad \zeta\in D_{l(j)+1}\setminus E_{q_j}[\rho-2],\\
&u(\zeta),\quad\quad \zeta\in D_{m+1}\setminus \bigcup_{j\in S}
E_{p_j}[-\rho+2],\\ 
&\exp_{b_j^\ast}^{t}(\alpha^-(z)\xi_j(z)+\alpha^+(z)\eta_j(z)),
\quad\quad z=\tau+it\in\Omega_j. 
\end{cases}
\end{equation}

\subsection{Stationary pregluing}\label{10statpreglue.section}
Let $L\subset\C^n\times\R$ be an admissible Legendrian
submanifold. Let $u\colon D_{m+1}\to\C^n$ be a holomorphic
disk with its $j$-th negative puncture $p$ mapping to $c$,
$(u,f)\in\cand_2(a,{\mathbf b})$, and let
$v\colon D_{l+1}\to\C^n$ be a holomorphic disk with the positive 
puncture $q$ mapping to $c$, 
$(v,h)\in\cand_{2}(c, {\mathbf d})$. Define
\begin{equation}
w_\rho=\Sigma^{\{j\}}_\rho(u,v).
\end{equation}

\begin{lma}\label{lmaprestat}
The function $w_\rho$ satisfies 
$w_\rho\in\cand_2(a,{\mathbf b}_{\{j\}}({\mathbf d}))$ and  
\begin{equation}\label{statdbarw}
\|\bar\pa w_\rho\|_1=\Ordo(e^{-\theta\rho}), 
\end{equation} 
where $\theta$ is the smallest complex angle at the Reeb chord $c$. In
particular, $\|\bar\pa w_\rho\|_1\to 0$ as $\rho\to\infty$.  
\end{lma}

\begin{pf}
The first statement is trivial. 
Outside $\Omega_j$, $w_\rho$ agrees with $u$ or $v$ which are
holomorphic. Thus it is sufficient to consider  
the restriction of $w_\rho$ to $\Omega_j$. 
To derive the necessary estimates we Taylor expand the
exponential map at $c^\ast$. To simplify notation we let
$c^\ast=0\in\C^n$ and let $\xi\in\R^{2n}$ be coordinates in $T_0\C^n$ and
$x\in\R^{2n}$ coordinates around $0\in\C^n$. Then  
\begin{equation}\label{Taylorexp}
\exp^{t}_0(\xi)=\xi-\Gamma^k_{ij}(t)\xi^i\xi^k\pa_k
+\Ordo(|\xi|^3). 
\end{equation}
This implies the inverse of the exponential map has Taylor expansion
\begin{equation}\label{Taylorexp^-1}
\xi=x+\Gamma^k_{ij}(t)x^ix^k\pa_k+\Ordo(x^3).
\end{equation}
From \eqref{Taylorexp} and \eqref{Taylorexp^-1} we get
\begin{equation}\label{Taylortrue}
\exp_0^{t}(\alpha^+\xi_j)=
\alpha^+u+((\alpha^+)^2-\alpha^+)\Gamma^k_{ij}(t)u^iu^j\pa_k+\Ordo(|u|^3)
\end{equation}
and a similar expression for $\alpha^-\eta_j$ in terms of $v_j$.
Lemma \ref{SalRob} implies that $u$ and $Du$ are $\Ordo(e^{-\theta\rho})$ in
$E_{p_j}[\rho]$, which together with \eqref{Taylortrue} implies
\eqref{statdbarw}.      
\end{pf}

\subsection{Handle slide pregluing}\label{10hspreglue.section}

Let $L_\lambda$, $\lambda\in(-1,1)=\Lambda$ be a $1$-parameter family
of chord-generic Legendrian submanifolds such that 
\begin{equation}\notag
\M^\Lambda_A(a;{\mathbf b})=\M^0_A(a;{\mathbf b})
\end{equation}
is a transversally cut out handle slide disk, represented by a map
$u\colon D_{m+1}\to\C^n$ with punctures $q, p_1,\dots,p_m$.

We first consider the case studied in Theorem \ref{gluhandslide1}.
Let $v\colon D_{l+1}\to\C^n$ represent an element in $\M^0_C(c,{\mathbf
d})$. Let
\begin{equation}\notag
w_\rho^0=\Sigma^{\{j\}}_\rho(u,v),
\end{equation}
and define
\begin{equation}\notag
w_\rho=w^0_\rho[\rho], 
\end{equation}
see Section \ref{fcanF}. (We cut off $w_\rho^0$ close to its punctures in
order to be able to use local coordinates as in Lemma \ref{lem:fcan.gamma} in a
neighborhood of $w_\rho$.)

\begin{lma}\label{lmaprehsl1}
The function $w_\rho$ satisfies
$(w_{\rho},0)\in\cand_{2,\Lambda}(a,{\mathbf b}_j({\mathbf d}))$ and
\begin{equation}\label{hsdbarw1}
\|\bar\pa w_{\rho}\|_1=\Ordo(e^{-\theta\rho}), 
\end{equation}
where $\theta>0$ is the smallest complex angle at any Reeb chord of
$L_0$. In particular $\|\bar\pa w_\rho\|\to 0$ as $\rho\to 0$. 
\end{lma}
\begin{pf}
We must consider the gluing regions and the effect of making
$w_\rho^0$ constant close to all punctures. The argument 
in the proof of Lemma \ref{lmaprestat} gives the desired estimate.
\end{pf}

For Theorem \ref{gluhandslide2}, we need to handle a few
more punctures. Let $v\colon D_{l+1}\to\C^n$ be an element in
$\M^0_C(c,{\mathbf d})$ where the elements in $S=\{s_1,\dots,s_r\}$
are indices of punctures which map to $a$. Let
\begin{equation}\notag
w_\rho^{0,S}=\Sigma_\rho^S(v,u,\dots,u)
\end{equation} 
and define
\begin{equation}\notag
w_\rho^S=w^{0,S}_\rho[\rho]
\end{equation}

\begin{lma}\label{lmaprehs2}
The function $w_{\rho}^S$ satisfies
$(w_{\rho}^S,0)\in\cand_{2,\Lambda}(c,{\mathbf d}_S({\mathbf b}))$ and
\begin{equation}
\|\bar\pa w_{\rho}^S\|_1=\Ordo(e^{-\theta\rho}), 
\end{equation} 
where $\theta>0$ is the smallest complex
angle at any Reeb chord of 
$L_0$. In particular $\|\bar\pa w_\rho\|\to 0$ as $\rho\to 0$. 
\end{lma}
\begin{pf}
See the proof of Lemma \ref{lmaprehsl1}.
\end{pf}

\subsection{Marked points}\label{ssecmarked}
In order to treat disks with less than three punctures
(i.e., disks with conformal reparameterizations) in the same way as disks
with more than three punctures we introduce special points which we
call {\em marked points} on the boundary. When a disks
with few punctures and marked points are glued to a disk with many
punctures there arises a disk with many punctures and marked points
and we must study also that situation.  

\begin{rmk}\label{dropthings.rmk}
Below we will often write simply $\cand_{2,\epsilon}$ to denote spaces
like $\cand_{2,\epsilon}({\mathbf c})$, dropping the Reeb chords from
the notation.
\end{rmk}

Let $L\subset\C^n\times\R$ be a (semi-)admissible 
Legendrian submanifold and let $u\colon D_m\to\C^n$ represent
$(u,f)\in\cand_{2,\epsilon}(\kappa)$ where $\kappa$ is a fixed
conformal structure on $D_m$. Let
$U_r\subset\Pi_\C(L)$, $r=1,\dots,k$ be disjoint open subsets where
$\Pi_\C(L)$ is real 
analytic and let $q_r\in\pa D_m$ be points such that $u(q_r)\in U_r$ and
$du(q_r)\ne 0$. After
possibly shrinking $U_r$ we may biholomorphically identify
$(\C^n,U_r,u(q_r))$ with $(\C^n,V\subset\R^n,0)$. Let
$(x_1+iy_1,\dots,x_n+iy_n)$ be coordinates on $\C^n$ and assume these
coordinates are chosen so $du(q_r)\cdot v_0=\pa_1$, where $v_0\in
T_{q_r}D_m$ is a unit 
vector tangent to the boundary. Let $H_r\subset\R^n$ denote an open
neighborhood of $0$ in the submanifold $\{x_1=0\}$. 

Let $S$ denote the cyclically ordered set of points
$S=\{p_1,\dots,p_m,q_1,\dots,q_k\}$ where $p_i \in \partial D_m$ are
the punctures of $D_m$. 
Fix three points  
$s_1,s_2,s_3\in\{p_1,p_2,p_3,q_1,\dots,q_k\}$ then the positions of
the other points in $S$ parameterizes the
conformal structures on $\Delta_{m+k}$.
As in Section~\ref{fcanG}, we pick vector fields
$\tilde v_j$, $j=1,\dots,m+k-3$ supported around the non-fixed points
in $S$. Given a conformal structure on $\Delta_{m+k}$ we endow it with
the metric which makes a neighborhood of each puncture $p_j$ look like
the strip and denote disks with such metrics $\tilde D_{m,k}$. 

If $(u,f)$, $u\colon \tilde D_{m,k}\to\C^n$ and $f\colon\pa \tilde
D_{m,k}\to\R$ are 
maps and $\tilde\kappa$ is a conformal structure on $\tilde D_{m+k}$
then forgetting the marked points $q_1,\dots,q_k$ we may view the maps 
as defined on $D_m$ and the conformal structure $\tilde \kappa$
gives a conformal structure $\kappa$ on $D_m$. Note though that the
standard metrics on $D_m$ corresponding to $\kappa$ may be different
from the metric corresponding to $\tilde\kappa$ (this happens when one
of the punctures $q_j$ is very close to one of the punctures
$p_r$). However, the metrics differ only on a compact set and thus 
using this forgetful map we define for a
fixed conformal structure $\tilde\kappa$ on $\tilde D_{m,k}$ the space
\begin{equation}\notag
\cand_{2,\epsilon}^{S}(\tilde\kappa)\subset\cand_{2,\epsilon}(\kappa)
\end{equation}
as the subset of elements represented by maps $w\colon D_m\to\C^n$ such that
$w(q_r)\in H_r$ for $r=1,\dots,k$. Using local coordinates on
$\cand_{2,\epsilon}(\kappa)$ around 
$(u,f)$ we see that for some ball $B$ around $(u,f)$,
$\cand_{2,\epsilon}^{S}(\tilde\kappa)\cap B$ is a
codimension $k$ submanifold with
tangent space at $(w,g)$ the closed subset of
$T_{(w,g)}\cand_{2,\epsilon}$ consisting of $v\colon D_m\to\C^n$ with
$\la v(q_r),\pa_1\ra=0$. We call $\tilde D_{m,k}$ a disk with $m$
punctures and $k$ marked points.  

The diffeomorphisms $\tilde\phi_j^{\sigma_j}$, $\sigma_j\in\R$
generated by to $\tilde v_j$ 
gives local coordinates
$\sigma=(\sigma_1,\dots,\sigma_{m+k-3})\in\R^{m+k-3}$ on the 
space of conformal structures on 
$\tilde D_{m,k}$ and the structure of a locally trivial bundle to the
space
\begin{equation}\notag
\cand_{2,\epsilon}^S=
\bigcup_{\tilde\kappa\in\conf_{m+k}}\cand_{2,\epsilon}^S(\tilde\kappa).
\end{equation} 
The $\bar\pa$-map is defined in the natural way on this space and we
denote it $\tilde\Gamma$.

\subsubsection{Marked points on disks with few punctures}\label{fewmarked}
Let $L\subset\C^n\times\R$ be a (semi-)admissible
submanifold, let $m\le 2$ and consider a holomorphic disk $(u,f)$ with
boundary on $L$, represented by a map $u\colon D_m\to\C^n$ . We shall
put $3-m$ marked points on $D_m$.  

Pick $U_r\subset\Pi_\C(L)$, $1\le r\le 3-m$ as disjoint open subsets
in which $\Pi_\C(L)$ is real
analytic and let $q_r\in\pa D_m$ be points such that $u(q_r)\in U_r$ and
$du(q_r)\ne 0$. Such points exists by Lemma \ref{lmaTaylor}. As in
Section \ref{ssecmarked} we then consider the $q_r$ as marked points
and as there we use the notation $H_r$ for the submanifold into which
$q_r$ is mapped.

Then the class in the moduli space of
holomorphic disks of every holomorphic disk $(w,g)$ which is
sufficiently close to 
$(u,f)$ in $\cand_{2,\epsilon}$ has a unique representative 
$(\hat w,\hat g)\in \cand_{2,\epsilon}^S$. Namely, any such $(w,g)$ must 
intersect $H_r$ in a point $q_r'$ close to $q_r$, $1\le r\le 3-m$. If
$\psi$ denotes the unique conformal reparameterization of $D_m$ which  
takes $q_r$ to $q'_r$, $1\le r\le 3-m$  
then $\hat w(\zeta)=w(\psi(\zeta))$. Moreover, if 
\begin{equation}
d\Gamma_{(u,f)}\colon
T_{(u,f)}\cand_{2,\epsilon}\to\sblv_{1,\epsilon}[0](D_m,{T^\ast}^{0,1}D_m\otimes\C^n), 
\end{equation}
has index $k$ (note $k\ge 3-m$ since the space of conformal
reparameterizations of $D_m$ is $(3-m)$-dimensional) then the
restriction of $d\Gamma_{(u,f)}$ to $T_{(u,f)}\cand_{2,\epsilon}^S$
has index $k-(3-m)$. In particular if $d\Gamma_{(u,f)}$ is surjective so
is its restriction. 
 
We conclude from the above that to study the moduli space of holomorphic
disks in a neighborhood of a given holomorphic disk we may (and will) use 
a neighborhood of that disk in $\cand_{2,\epsilon}^S$ and
$\tilde\Gamma$ rather then a 
neighborhood in the bigger space $\cand_{2,\epsilon}$ and $\Gamma$.

\subsubsection{Marked points on disks with many punctures}\label{manymarked}

Let $L\subset\C^n\times\R$ be as above, let $m\ge 3$ and consider a
holomorphic disk $(u,f)$ with 
boundary on $L$, represented by a map $u\colon D_m\to\C^n$ . We shall
put $k$ marked points on $D_m$.  

Pick $U_r\subset\Pi_\C(L)$, $1\le r\le k$ as disjoint open subsets
in which $\Pi_\C(L)$ is real
analytic and let $q_r\in\pa D_m$ be points such that $u(q_r)\in U_r$ and
$du(q_r)\ne0$. As in
Section \ref{ssecmarked} we then consider the $q_r$ as marked points
and as there we use the notation $H_r$ for the submanifold into which
$q_r$ is mapped.

Note that $(u,f)$ lies in $\cand_{2,\epsilon}$ as well as in
$\cand_{2,\epsilon}^S$. We define a map 
\begin{equation}\notag
\Omega\colon U\subset\cand_{2,\epsilon}^S\to\cand_{2,\epsilon};\quad
\Omega((w,g,\tilde\phi^s))=(\hat w,\hat g,\tilde\phi^t)
\end{equation} 
where $U$ is a neighborhood of $((u,f),\tilde\kappa)$ in as follows.

Consider the local coordinates $\omega\in\R^{m+k-3}$ on $\conf_{m+k}$
around $\tilde\kappa$ and the product structure 
\begin{equation}\notag
\R^{m+k-3}=\R^{m-3}\times\R^j\times\R^{k-j},
\end{equation}
where $\R^{m-3}$ is identified with the diffeomorphisms 
\begin{equation}\notag
\phi^{\tau}=
\tilde\phi_{p_4}^{\tau_1}\circ
\dots\circ\tilde\phi_{p_m}^{\tau_{m-3}},\quad
\tau=(\tau_1,\dots,\tau_{m-3})\in\R^{m-3},
\end{equation}
where $j$ is the number of elements in
$\{s_1,s_2,s_3\}\setminus\{p_1,p_2,p_3\}$, and where 
$\R^{k-j}$ is identified with the diffeomorphisms
\begin{equation}\notag
\phi^\sigma=\tilde\phi_{\hat s_1}^{\sigma_1}\circ\dots\circ
\tilde\phi_{\hat s_{k-j}}^{\sigma_{k-j}},\quad
\sigma=(\sigma_1,\dots,\sigma_{k-j})\in\R^{k-j} 
\end{equation}
where $\{\hat s_1,\dots,\hat
s_{k-j}\}=S\setminus(\{p_4,\dots,p_m\}\cup\{s_1,s_2,s_3\})$.     

For $\tilde\theta$ near $\tilde\kappa$, let
$\{s_1',\dots,s'_{m+k-3}\}$ denote the 
corresponding positions of punctures and marked points in $\pa\Delta$
and let $\psi\colon\Delta\to\Delta$ be the
unique conformal reparameterization such that 
$\psi(p_j)=p_j'$ for $j=1,2,3$ and note
that we may view $\psi$ as a map from $D_m$
to $\tilde D_{m,k}$. Let $s_l''=\psi^{-1}(s_l')$ for $3\le l \le k+m-3$
and $s_l \ne p_i$, $i=1,2,3$ let $(\tau,\sigma)\in\R^{m+k-3-j}$ be the unique
element such that $\phi^\tau\circ \phi^\sigma(s_l)=s''_l$. Define
\begin{equation}\notag
\Omega(w,\tilde\theta)=(w\circ\psi\circ\phi^\sigma,(\phi^\tau)^{-1}),
\end{equation}
where $(\phi^\tau)^{-1}$ is interpreted as a conformal structure on $D_m$ in
a neighborhood of $\kappa$ in local coordinates given by $\phi^\tau$,
$\tau\in\R^{m-3}$ and where we drop the boundary function from the
notation since it is uniquely determined by the $\C^n$-function component
of $\Omega(w,\tilde\theta)$ and $g$.   

\begin{lma}\label{lmamanymarks}
The map $\Omega$ maps $U\cap\tilde\Gamma^{-1}(0)$ into
$\Gamma^{-1}(0)$. Moreover, $\Omega$ is a $C^1$-diffeomorphism on a
neighborhood of $(u,f)$. 
\end{lma}
\begin{pf}
Assume that $(w,g)\in\tilde\Gamma^{-1}(0)$. Then $w$ is holomorphic in the
conformal structure $\tilde\theta$. Since $\psi$ is a conformal
equivalence and since the conformal structure $\tilde\theta$ is
obtained from $\tilde\kappa$ by action of the inverses of
$\phi^\tau\circ\phi^\sigma$ this implies 
\begin{equation}\notag
0=dw\circ d\psi+ i\circ (dw\circ d\psi)\circ (d\phi^\sigma)
\circ (d\phi^\tau)\circ
j_\kappa\circ(d\phi^\tau)^{-1}\circ(d\phi^\sigma)^{-1}.
\end{equation} 
Thus 
\begin{equation}\notag
0= \Bigl(dw\circ d\psi\circ d\phi^\sigma + 
i\circ (dw\circ d\psi\circ d\phi^\sigma)
\circ (d\phi^\tau)\circ j_\kappa\circ(d\phi^\tau)^{-1}\Bigr)
\circ (d\phi^\sigma)^{-1},
\end{equation}
and $w\circ\psi\circ\phi^\sigma$ is holomorphic in the conformal
structure $d\phi^\tau j_\kappa(d\phi^\tau)^{-1}$ as required.

For the last statement we use the inverse function theorem. It is
clear that the map
$\Omega$ is $C^1$ and that the
differential of $\Omega$ at $(u,f)$ is a Fredholm operator. In
fact, on the complement of all conformal variations on 
$\tilde D_{m,k}$ not supported around any of $p_3,\dots,p_m$ it is just an
inclusion into a subspace of codimension $k$, which consists of
elements $v$ which vanish at $q_1,\dots,q_k$. Since $du(q_r)\ne 0$ for
all $r$ it follows easily that the image of the remaining $k$
directions in $T_{(u,f)}\cand_{2,\epsilon}^S$ spans the complement of
this subspace.     
\end{pf}

It is a consequence of Lemma \ref{lmamanymarks} that if $(u,f)$ is a
holomorphic disk with boundary on $L$ and more than $3$ punctures then
we may view a neighborhood of $(u,f)$ in the moduli space of such
disks either as a submanifold in $\cand_{2,\epsilon}^S$
or in $\cand_{2,\epsilon}$ in a neighborhood of $(u,f)$.

\begin{rmk}
Below we extend the use of the notion $\cand_{2,\epsilon}$ to include
also spaces $\cand_{2,\epsilon}^S$, when this is convenient. The point
being that after Sections \ref{fewmarked} and \ref{manymarked} we
may always assume 
the number of marked points and punctures is $\ge 3$, so that the 
moduli space of holomorphic disks (locally) may be viewed as a
submanifold of $\cand_{2,\epsilon}$.   
\end{rmk}

\subsection{Uniform invertibility of the differential in the
stationary case}\label{10unistat.section}
Let 
\begin{equation}\notag
\Gamma\colon\cand_2\to
\sblv_1[0](D_m,{T^\ast}^{0,1} D_m\otimes\C^n) 
\end{equation}
be the $\bar\pa$-map defined in Section \ref{fcanH} (see Remark
\ref{dropthings.rmk} for notation). 
Let $u\colon D_{m+1}\to\C^n$ and $v\colon D_{l+1}\to\C^n$ be as in
Section \ref{10statglue.section} and consider the
differential $d\Gamma_\rho$ at $(w_{\rho},\kappa_\rho)$, where $w_\rho$ is
as in Lemma \ref{lmaprestat} and $\kappa_\rho$ is the natural metric
(complex structure) on $D_r(\rho)$, $r=m+l$. 
After Sections \ref{fewmarked} and \ref{manymarked} we know that
after adding $3-m$ or $3-l$
marked points on holomorphic disks with $\le 2$ punctures, we may assume that
$m\ge 2$ and $l\ge 2$ below.

\begin{lma}\label{statunif}
There exist constants $C$ and $\rho_0$ such that 
if $\rho>\rho_0$ then there are continuous right inverses 
\begin{equation}\notag
G_\rho\colon\sblv_1[0]({T^\ast}^{0,1}D_r(\rho)\otimes\C^n)\to
T_{(w_\rho,\kappa_\rho)}\cand_2
\end{equation}
of $d\Gamma_\rho$ with 
\begin{equation}\notag
\|G_\rho(\xi)\|\le C\|\xi\|_1.
\end{equation}
\end{lma}

\begin{pf}
The kernels 
\begin{align*}
&\ker(d\Gamma_{(u,\kappa_1)})\subset T_{u}\cand_2\oplus
T_{\kappa_1}\conf_{m+1},\\ 
&\ker(d\Gamma_{(v,\kappa_2)})\subset T_{v}\cand_2\oplus
T_{\kappa_2}\conf_{l+1}
\end{align*}
are both $0$-dimensional.  
As in Section \ref{ssecmarked}, we view elements 
$\gamma_1\in T_{\kappa_1}\conf_{m+1}$  
($\gamma_2\in T_{\kappa_2}\conf_{l+1}$) as linear combinations of
sections of $\End(T D_{m+1})$ ($\End(T D_{l+1})$) supported in 
compact annular regions close to all punctures and marked points,
except at three. Since these annular regions are disjoint from the regions
affected by the gluing of $D_{m+1}$ and $D_{l+1}$, we get an embedding 
$$
T_{\kappa_1}\conf_{m+1}\oplus T_{\kappa_2}\conf_{l+1}\to 
T_{\kappa_\rho}\conf_{m+r},
$$
In fact, using this embedding,
$$
T_{\kappa_\rho}\conf_{r}=
T_{\kappa_1}\conf_{m+1}\oplus T_{\kappa_2}\conf_{l+1}\oplus\R, 
$$
where the last summand can be taken to be generated by a section
$\gamma_0$ of 
$\End(T D_r(\rho))$ supported in an annular region around a puncture
(marked point)
in $D_r(\rho)$ where there was no conformal variation before the
gluing. Then $\gamma_0$ spans a subspace of dimension $1$ in
$T_{(w_\rho,\kappa_\rho)}\cand_2$. Let $a$ be a
coordinate along this $1$-dimensional subspace.
We prove that for  $(\xi,\gamma)\in W_\rho=\{a=0\}$ we have the
estimate 
\begin{equation}\label{stat(i)}
\|(\xi,\gamma)\|\le
C\|d\Gamma_\rho(\xi,\gamma)\|_{1,\epsilon},
\end{equation}
for all sufficiently large $\rho$. Since the Fredholm-index of
$d\Gamma_\rho$ equals $1$, this shows $d\Gamma_\rho$ are surjective and
with uniformly bounded inverses $G_\rho$ as claimed and thus
finishes the proof. 

Assume \eqref{stat(i)} is not true. Then there exists a sequence 
of elements $(\xi_N,\gamma_N)\in W_{\rho(N)}$, 
$\rho(N)\to\infty$ as $N\to\infty$ with
\begin{align}\label{stat(i)norm}
&\|(\xi_N,\gamma_N)\|=1,\\\label{stat(i)to0} 
&\|\bar\pa_{\kappa_{\rho(N)}}\xi_N+i\circ 
dw_{\rho(N)}\circ \gamma_N\|_1\to 0.
\end{align} 
As in Section \ref{10glueop.section}, we glue a negative puncture at $p$ to
a positive one at $q.$
Note that on the strip 
\begin{equation}
\Theta_\rho=
(E_p[-1]\setminus E_{p}[-\rho])\cup (E_q[1]\setminus E_q[\rho])\approx
[-\rho,\rho]\times[0,1]\subset D_r(\rho)
\end{equation}
the conformal structure $\kappa_\rho$ is the standard one and therefore
$\bar\pa_{\kappa_\rho}$ is just the standard $\bar\pa$ operator. 
Also, since $\gamma_N$ does not have support in $\Theta_\rho$
the second term in \eqref{stat(i)to0} equals $0$ when restricted to
$\Theta_\rho$.   

Let $\alpha_\rho\colon\Theta_\rho\to\C$ be cut-off functions which are
real and holomorphic on the boundary,
equal $1$ on $[-2,2]\times[0,1]$, equal $0$ outside
$[-\frac12\rho,\frac12\rho]\times[0,1]$, and satisfy
$|D^k\alpha_\rho|=\Ordo(\rho^{-1})$, $k=1,2$. 

Then $\alpha_{\rho(N)}\xi_N$ is a sequence of functions on
$\R\times[0,1]$ which satisfy boundary conditions converging to two
transverse Lagrangian subspaces. 
Just as we prove in Lemma \ref{lmadef} that the (continuous) index 
is preserved under small
perturbations, we conclude that the (upper semi-continuous) dimension
of the kernel stay zero for large enough $\rho(N)$; thus,
there exists a constant $C$ such
that 
\begin{equation}\label{stat(i)inner'}
\|\xi_N|[-2,2]\times[0,1]\|_2\le
\|\alpha\xi_N\|_2\le C\Bigl(\|\alpha_{\rho(N)}(\bar\pa\xi_N)\|_1
+\|(\bar\pa\alpha_{\rho(N)})\xi_N\|_1\Bigr).
\end{equation}
As $N\to\infty$ both terms on the right hand side in
\eqref{stat(i)inner'} approaches $0$. Hence, 
\begin{equation}\label{stat(i)inner}
\|\xi_N|[-2,2]\times[0,1]\|_2\to 0,\text{ as }N\to\infty.
\end{equation}

Pick cut-off functions $\beta_N^+$ and $\beta_N^-$ on $D_r(\rho)$ which are
real valued and holomorphic on the boundary and have the following
properties. On $D_{m+1}\setminus E_{p}[-\rho+1]$, 
$\beta_N^+=1$ and on $D_{l+1}\setminus E_{q}[\rho]$, $\beta_N^+=0$.
On $D_{l+1}\setminus E_{q}[\rho-1]$,  
$\beta_N^-=1$ and on $D_{m+1}\setminus E_{p}[-\rho]$, 
$\beta_N^-=0$. Let
$(\xi_N,\gamma_N)^\pm=(\beta_N^{\pm}\xi_N,\beta_N^{\pm}\gamma_N)$. 
Using the invertibility of $d\Gamma_+=d\Gamma_{(u,\kappa_1)}$
and $d\Gamma_-=d\Gamma_{(v,\kappa_2)}$ we find a constant $C$ such
that
\begin{align}\notag
\|(\xi_N,\Gamma_N)^\pm\|&\le
C\|d\Gamma_\pm(\xi_N,\gamma)^\pm\|_{1}\\\label{stat(i)outer}
&\le C\Bigl(\|\beta_N^\pm d\Gamma_\rho(\xi_N,\gamma_N)\|_{1}
+\|(\bar\pa\beta_N^\pm)\xi_N\|_{1}\Bigr). 
\end{align}  
The first term in the last line of \eqref{stat(i)outer} tends to $0$
as $N\to\infty$ by \eqref{stat(i)to0}, the second term
tends to $0$ by \eqref{stat(i)inner}. Hence, the left hand
side of \eqref{stat(i)outer} tends to $0$ as $N\to\infty$.
Thus, \eqref{stat(i)inner} and \eqref{stat(i)outer} contradict
\eqref{stat(i)norm} and we conclude \eqref{stat(i)} holds. 
\end{pf}

\subsection{Uniform invertibility in the handle slide case}\label{10unihs.section}

Two handle-slide gluing theorems 
(Theorem \ref{gluhandslide1} and Theorem \ref{gluhandslide2}) were
formulated. 
Here the handle-slide analogs of Lemma \ref{statunif} corresponding to
these two theorems will be proved. The lemma corresponding to Theorem
\ref{gluhandslide1} is proved first, using a straightforward extension of the
proof of Lemma \ref{statunif}. Second the more difficult lemma
corresponding to Theorem \ref{gluhandslide2} is dealt with.  
Compare with the Appendix in \cite{Sullivan}.

Let $L_\lambda$ be as in Section \ref{sec12} and let
$u\colon D_{m+1}\to\C^n$, be the handle slide disk and let $p$ denote
one of its negative punctures mapping to the Reeb chord $b$. Also, let
$v\colon D_{l+1}\to\C^n$ be 
a rigid disk with positive puncture $q$ mapping to $b$.
Consider the differential $d\Gamma_\rho$ at
$(w_{\rho},\kappa_\rho,0)$, where 
$w_\rho$ is as in Lemma \ref{lmaprehsl1}, $\kappa_\rho$ is the natural
metric (complex 
structure) on $D_r(\rho)$, $r=m+l$, and
$0\in\Lambda\approx(-1,1)$. Again we assume $m\ge 2$ and $l\ge 2$.

\begin{lma}\label{handsl1unif}
There exist constants $C$ and $\rho_0$ such that 
if $\rho>\rho_0$ then there is a continuous right inverse 
$G_\rho$ of $d\Gamma_\rho$ 
$$
G_\rho\colon\sblv_{1}[0]({T^\ast}^{0,1}D_r(\rho)\otimes\C^n)\to
T_{(w_\rho,\kappa_\rho,0)}\cand_{2,\Lambda}
$$ 
such that
$$
\|G_\rho(\xi)\|\le C\|\xi\|_{1}.
$$
\end{lma}

\begin{pf}
The kernels 
\begin{align*}
&\ker(d\Gamma_{(u,\kappa_1,0)})\subset
T_{(u,0)}\cand_{2,\Lambda}\oplus 
T_{\kappa_1}\conf_{m+1},\\ 
&\ker(d\Gamma_{(v,\kappa_2)})\subset
T_{(v,\kappa_2)}\cand_{2}\oplus T_{\kappa_2}\conf_{l+1}
\end{align*}
are $0$-dimensional 
(note that $\Lambda$ is not involved in the second line).  
As in the proof of Lemma \ref{statunif}, we consider the embedding 
\begin{equation}\notag
T_{\kappa_1}\conf_{m+1}\oplus T_{\kappa_2}\conf_{l+1}\to 
T_{\kappa_\rho}\conf_{m+l},
\end{equation}
and use this inclusion to get the isomorphism
\begin{equation}\notag
T_{\kappa_\rho}\conf_{m+l}=
T_{\kappa_1}\conf_{m+1}\oplus T_{\kappa_2}\conf_{l+1}\oplus\R, 
\end{equation}
where the last summand can be taken to be generated by a section of 
$\End(T D_r(\rho))$ supported around a puncture (marked point) where
no conformal variation was supported before the gluing.
Let $\gamma_0$ denote a tangent vector in this direction. 
Then $\gamma_0$ spans a subspace of dimension $1$ in
$T_{(w_\rho,\kappa_\rho,0)}\cand_{2,\epsilon,\Lambda}$. Let $a$ be a
coordinate along this $1$-dimensional subspace.
We prove that for  $(\xi,\gamma,\lambda)\in W_\rho=\{a=0\}$ we have the
estimate 
\begin{equation}\label{handsl(i)}
\|(\xi,\gamma,\lambda)\|\le
C\|d\Gamma_\rho(\xi,\gamma,\lambda)\|_{1},
\end{equation}
for all sufficiently large $\rho$. Since the Fredholm-index of
$d\Gamma_\rho$ equals $1$, this finishes the proof. 

Assume \eqref{handsl(i)} is not true. Then there exists a sequence 
of elements $(\xi_N,\gamma_N,\lambda_N)\in W_{\rho(N)}$, 
such that, $\rho(N)\to\infty$ as $N\to\infty$ and 
\begin{align}\label{handsl(i)norm}
&\|(\xi_N,\gamma_N,\lambda_N)\|=1,\\\label{handsl(i)to0} 
&\|\bar\pa_{\kappa_{\rho(N)}}\xi_N+i\circ 
dw_{\rho(N)}\circ
\gamma_N+\lambda_N\bar\pa_{\kappa_\rho(N)}Y(w_{\rho(N)})\|_{1}\to
0. 
\end{align} 

Let $\alpha_\rho\colon D_r(\rho)\to\C$ be smooth functions which
equal $1$ on $(D_{m+1}\setminus E_p[-\frac12\rho])$, equal $0$ on the
complement of 
$D_{m+1}\setminus E_p[-\rho+3]$, are holomorphic and real
valued on the boundary, and satisfy
$|D^k\alpha_\rho|=\Ordo(\rho^{-1})$, $k=1,2$. Consider the support of
$\alpha_\rho$ as a subset of $D_{m+1}$. Then
$(\xi_N^+,\gamma_N^+)=(\alpha_{\rho(N)}\xi_N,\alpha_{\rho(N)}\gamma_N)$
is a sequence 
of elements in
$T_{(u,\kappa_1,0)}\cand_{2,\epsilon,\Lambda}\oplus T_{\kappa_1}\conf_{m+1}$
The invertibility of $d\Gamma_{(u,\kappa_1,0)}$ implies
there exists a constant $C$ such that 
\begin{equation}\label{lambdato0}
\|(\xi_N^+,\gamma_N^+,\lambda_N)\|\le 
C\|\bar\pa_{\kappa_1}\xi_N^+ +i\circ du\circ\gamma_N^++\lambda_N\bar\pa
Y(u)\|_{1}. 
\end{equation}

On the other hand, noting that $\bar\pa_{\kappa_1}$ is the standard
$\bar\pa$-operator in a neighborhood of $q$ and that
$\bar\pa_{\kappa_\rho}$ is the standard $\bar\pa$-operator in the
gluing region we find that   
\begin{align}\notag
\|\bar\pa_{\kappa_1}\xi_N^+ +i\circ du\circ\gamma_N^++\lambda\bar\pa
Y(u)\|_{1}\le&
C\Bigl(\|\alpha_\rho(\bar\pa_{\kappa_\rho}\xi_N +i\circ dw\circ\gamma_N+
\lambda\bar\pa Y(w_\rho))\|_1\\\notag
&+\|(\bar\pa\alpha_\rho)\xi_N\|_1+\|(1-\alpha_\rho)\bar\pa
Y(u)\|_1\Bigr) . 
\end{align}
The first term in the right hand side goes to zero as $\rho\to\infty$ by
\eqref{handsl(i)to0}. The second term goes to $0$ since
$|\bar\pa\alpha_\rho|=\Ordo(\rho^{-1})$. Finally, the last term goes to
$0$ since $\bar\pa
Y(u)\in\sblv_1(D_{m+1},{T^\ast}^{0,1}D_{m+1}\otimes\C^n)$. Hence, 
the right hand side goes to $0$ as $\rho = \rho(N) \to\infty$. We conclude from
\eqref{lambdato0} that $\|(\xi_N^+,\gamma_N^+,\lambda_N)\|\to 0$ as
$\rho\to\infty$. In particular, it follows that $\lambda_N\to 0$ as
$\rho\to\infty$ and once this has been established we can repeat the
argument from the proof of Lemma \ref{statunif} to conclude that
\eqref{handsl(i)} holds true.    
\end{pf}

We now turn to the second case. Let $L_\lambda$ be as in Section
\ref{sec12} and let again
$u\colon D_{m+1}\to\C^n$, be the handle slide disk and let $q$ denote
its positive puncture mapping to the Reeb chord $a$. Also, let
$v\colon D_{l+1}\to\C^n$ be 
a rigid disk with a set $S$ of negative punctures
$S=\{p_1,\dots,p_k\}$ mapping to $a$. 
Again we assume $m\ge 2$ and $l\ge 2$.

Consider the differential $d\Gamma^S_\rho$ at
$(w_{\rho}^S,\kappa_\rho^S,0)$, where 
$w_\rho^S$ is as in Lemma \ref{lmaprehs2}, $\kappa_\rho^S$ is the
natural metric (complex 
structure) on $D_r^S(\rho)$, $r=l+km$, and $0\in\Lambda\approx(-1,1)$.
We let $D^j_{m+1}$, $j=1,\dots,k$ denote $k$ distinct copies of
$D_{m+1}$.  

\begin{lma} \label{handsl2unif}
There exist constants $C$ and $\rho_0$ such that 
if $\rho>\rho_0$ then there is a continuous right inverse 
$G_\rho^S$ of $d\Gamma_\rho^S$ 
$$
G_\rho^S\colon\sblv_{1}[0]({T^\ast}^{0,1}D_r(\rho)\otimes\C^n)\to
T_{(w_\rho^S,\kappa_\rho^S,0)}\cand_{2,\Lambda}
$$ 
such that
$$
\|G_\rho^S(\xi)\|\le C\|\xi\|_{1}.
$$
\end{lma}

\begin{pf} 
The kernels 
\begin{align*}
&\ker(d\Gamma_{(v, \kappa_2,0)})\subset
T_{(v,0)}\cand_{2,\Lambda}\oplus 
T_{\kappa_2}\conf_{l+1},\\ 
&\ker(d\Gamma_{(u,\kappa_1,0)})\subset T_{(u,0)}\cand_{2,\Lambda}\oplus
T_{\kappa_1}\conf_{m+1}
\end{align*}
are $1$-dimensional, respectively $0$-dimensional. Note that
$\ker(d\Gamma_{(v, \kappa_2, 0)})$ is 
{\em not} contained in the subspace $\{\lambda=0\}$. 
We consider first the special case $S=\{1\}$.
As in previous proofs we write,
$$
T_{\kappa_\rho}\conf_{l+m}=
T_{\kappa_2}\conf_{l+1}\oplus
T_{\kappa_1}\conf_{m+1}\oplus\R,  
$$
where the generators of $\R$ is a section $\gamma_0$,
of $\End(T D_{l+m}(\rho))$ which is supported around one of the
previously fixed punctures of $D_{m+1}$. As in the proof of Lemma
\ref{handsl1unif}, we prove the the estimate 
\begin{equation}
\|(\xi,\gamma,\lambda)\|\le
C\|d\Gamma^{\{1\}}_\rho(\xi,\gamma,\lambda)\|_{1}. 
\end{equation}
on the complement of the new conformal direction. This finishes the
proof in the special case when $k=1$. 

In order to move to higher $k$, we prove that the kernel of
$d\Gamma_\rho^{\{1\}}$ is {\em not} contained in $\{\lambda=0\}$ for
$\rho$ large enough. Assume this is not the case. Then there exists a
sequence $(\xi_N,\gamma_N,0)$ such that
\begin{equation}\label{handsl2lambda}
\|(\xi_N,\gamma_N,0)\|=1,\quad d\Gamma_\rho^{\{1\}}(\xi_N,\gamma_N,0)=0.
\end{equation} 
Let $\alpha$ be a smooth function, real valued and holomorphic on
the boundary of $D_r^{\{1\}}$, equal to $1$ on 
$(D_{l+1}\setminus E_{p_1}[-\frac12\rho])\cup (D_{m+1}\setminus
E_{q}[\frac12\rho])$ and equal to $0$ on the gluing region
$\Omega$. Using the uniform 
invertibility of $d\Gamma_{(v,\kappa_2,0)}$ and of $d\Gamma_{(u,\kappa_1,0)}$
on the complement of $\{\lambda=0\}$ we conclude that
$\|\alpha(\xi_N,\gamma_N)\|\to 0$. Finally, using the elliptic
estimate on the strip we also find that the $2$-norm of $\xi_N$
restricted to 
$(E_{p_1}[-\frac12\rho]\setminus E_{p_1}[-\rho])\cup
(E_{q}[\frac12\rho]\setminus E_{q}[\rho])$ goes to zero. This contradicts   
\eqref{handsl2lambda} and we find the kernel has
$\lambda$-component. 

Now consider the case $k>1$. We assume by induction that the desired
invertibility of $d\Gamma^S_\rho$ is proved for $S$ with $|S|<k$ and
that the kernel of $d\Gamma^S_\rho$ has $\lambda$ component.  
Letting
$T$ be the union of $S$ and the new puncture, we view $w_\rho^T$ as
obtained by gluing $w_\rho^S$ and $u$. Repeating the above
argument we prove the estimate on the complement of the new-born
conformal structure and also the fact that the kernel of
$d\Gamma_\rho^T$ has $\lambda$ component. Since $k$ is finite this
finishes the proof.
\end{pf}

\subsection{Self-tangencies, coordinates and genericity assumptions}
Let $z=x+iy=(z_1,\dots,z_n)=(x_1+iy_1,\dots,x_n+iy_n)$ be coordinates
on $\C^n$. Let $L\subset \C^n\times\R$ be a semi-admissible 
Legendrian submanifold with self-tangency double point at
$0$. We assume that the self-tangency point is standard, see
Section \ref{5Isotopies.section}.

Theorems~\ref{8degentrans.thm} and \ref{9cpt.thm} imply 
that the moduli-space of rigid holomorphic disks
with boundary on $L$ is a $0$-dimensional compact manifold. Moreover,
because of the enhanced transversality discussion in
Section \ref{8enhancedtrans.section}, 
we may assume that there exists $r_0>0$ such that for
all $0<r<r_0$, if $u\colon D_m\to\C^n$ is a rigid holomorphic disk
with boundary on $L$ then $\pa D_m\cap u^{-1}(B(0,r))$ is a disjoint
union of subintervals of $\pa E_{p_j}[\pm M]$, for some $M>0$ and some
punctures $p_j$ on $\pa D_m$ mapping to $0$.

By Lemma \ref{SalRob},
if $u\colon D_m\to\C^n$ is a
rigid holomorphic disk with $q^+$ a positive ($q^-$ a negative)
puncture mapping to $0$ then there exists $c\in\R$ such
that for $\zeta=\tau+it\in E_{q^{\pm}}[\pm M]$ 
$$
u(\zeta)=\Bigl(-2(\zeta+c)^{-1},0,\dots,0\Bigr)+\Ordo(e^{-\theta|\tau|}),
$$    
for some $\theta>0$. For simplicity, we assume below that
coordinates on $E_{q^\pm}[M]$ are chosen in such a way that $c=0$
above. 
 
\subsection{Perturbations for self-tangency shortening}\label{10perturbstshort.section}

For $0<a<1$, with $a$ very close to $1$ and $R>0$ with $R^{-1}\ll r_0$, let
$b_R\colon[0,\infty)\to\R$ be a smooth non-increasing function with
support in $[0,R^{-1})$ and with the following properties
\begin{align}\notag
&b_R(r)=(R+R^a)^{-2}\text{ for } r \in \left[0,(R+\frac12 R^a)^{-1}\right],\\\notag
&|Db_R(r)|=\Ordo(R^{-a}),\\\notag
&|D^2 b_R(r)|=\Ordo(R^{2-2a}),\\\notag
&|D^3 b_R(r)|=\Ordo(R^{4-3a}),\\\label{b_Rbounds}
&|D^4 b_R(r)|=\Ordo(R^{6-4a}).
\end{align}
The existence of such a function is easily established using the fact
that the length of the interval where $Db_R$ is supported equals
$$
R^{-1}-\left(R+\frac12 R^{a}\right)^{-1}=\frac12 R^{a-2} +\Ordo(R^{2(a-2)}).
$$

Let 
\begin{equation}  
h_R(z)=-x_1(z)b_R(|z|).
\end{equation}
Let $L^1$ and $L^2$ be the two branches of the local Lagrangian
projection near the self-tangency, see Section \ref{5Isotopies.section}
or Figure \ref{1fig:TypeB}.
For $s>0$, let $\Psi_R^s$ denote the time $s$ Hamiltonian flow
of $h_R$ and let $L_R(s)$ denote the Legendrian
submanifold which results when $\Psi^s_R$ is lifted to a contact
flow on $\C^n\times\R$ (see Section \ref{5Isotopies.section})
which is used to move $L^2$. Let 
$L^2_R(s)=\Psi_R^s(L^2)$. Let
$g(R,s,\sigma)$ be a $3$-parameter family of metrics on $\C^n$
such that $L^1$ is totally geodesic for $g(R,s,0)$,
$L^2_R(s)$ is totally geodesic for $g(R,s,1)$ and such
that $g(R,s,0)$ and $g(R,s,1)$ have properties as the metrics
constructed in Section~\ref{fcanC}.

Note that $L^2_R(1)\cap L^1$ consists of exactly two points with
coordinates $(\pm (R+R^a)^{-1}+\Ordo(R^{-3}),0,\dots,0).$  

We will use $\Psi_R^s$ to deform holomorphic disks below. 
It will be
important for us to know they remain almost holomorphic in a rather
strong sense, for which we need to derive some estimates on
the flow $\Psi^s_R$ and its derivatives. Let $X_R$ denote the
Hamiltonian vector field of $h_R$. 
Then if $D$ denotes derivative with respect to the variables in $\C^n$
and $\cdot$ denotes contraction of tensors
\begin{align}\label{sh0flow}
\frac{d}{ds}\Psi^s_R =& X_R;\quad \Psi_\mu^0=\id,
\\\label{sh1flow} 
\frac{d}{ds}D\Psi^s_R =& DX_R\cdot D\Psi^s_R;\quad
D\Psi_\mu^0=\id,\\\label{sh2flow}
\frac{d}{ds}D^2\Psi^s_R  =& D^2 X_R\cdot D\Psi^s_R\cdot D\Psi^s_R+
D X_R\cdot D^2\Psi^s_R;\quad D^2\Psi_R^0=0,\\\notag
\frac{d}{ds}D^3\Psi^s_R  =&  D^3 X_R\cdot D\Psi^s_R\cdot 
D\Psi^s_R\cdot D\Psi^s_R+2D^2 X_R\cdot D^2\Psi^s_R\cdot 
D\Psi^s_R\\\label{shflow3} 
&+  D^2 X_R\cdot D\Psi^s_R\cdot D^2\Psi^s_R+
DX_R\cdot D^3\Psi^s_R
;\quad D^3\Psi_R^0=0. 
\end{align}

Since $X_R=i\cdot Dh_R$ and 
$x_1(z)=\Ordo(R^{-1})$ for $|z|$ in the support of $b_R$, 
\eqref{b_Rbounds} implies
\begin{align}\label{X_R0}
|X_R|=\Ordo(R^{-(1+a)}),\\\label{X_R1}
|DX_R|=\Ordo(R^{(1-2a)}),\\\label{X_R2}
|D^2 X_R|=\Ordo(R^{(3-3a)}),\\\label{X_R3}
|D^2 X_R|=\Ordo(R^{(5-4a)}).
\end{align}

If $0\le s\le 1$ then
\begin{itemize}
\item[{\bf F0}] \eqref{sh0flow} and \eqref{X_R0} imply
$|\Psi^s_R-\id|=\Ordo(R^{-(1+a)})$.
\item[{\bf F1}] \eqref{sh1flow} and \eqref{X_R1} first give
$|D\Psi_R^s|=\Ordo(1)$. This together with \eqref{X_R1}  
imply $|D\Psi^s_R-\id|=\Ordo(R^{1-2a})$.
\item[{\bf F2}] \eqref{sh2flow}, \eqref{X_R1}, \eqref{X_R2}, {\bf F1},
and Duhamel's principle imply
$|D^2\Psi_R^s|=\Ordo(R^{3-3a})$.
\item[{\bf F3}] In a similar way as in {\bf F2} we derive
$|D^3\Psi_R^s|=\Ordo(R^{5-4a})$. 
\end{itemize}

Let $u\colon\R\times[0,1]\to\C^n$ be a holomorphic function and
and let $\omega\colon[0,1]\to[0,1]$ be a
smooth non-decreasing surjective
approximation of the identity which is constant in a
$\delta$-neighborhood of the ends of the interval.
Consider the function $u_R(\tau+it)=\Psi_R^{\omega(t)}(u(\tau+it))$.   
We want estimates for $u_R$, $\bar\pa u_R$ and $D\bar\pa u_R$ and
$\pa_\tau D\bar\pa u_R$. 

{\bf F0} implies
\begin{equation}\label{u_R0}
u_R=u+\Ordo(R^{-(1+a)}).
\end{equation}
For the estimates on $\bar\pa u_R$ and its derivatives we note
\begin{equation}\label{dbaru_R}
\bar\pa u_R=D\Psi_R^{\omega(t)}\frac{\pa u}{\pa\tau}+
i\Bigl(D\Psi_R^{\omega(t)}\frac{\pa u}{\pa t}+
\frac{d\omega}{dt}X_R(u)\Bigr)
\end{equation}

By \eqref{X_R0}, \eqref{X_R1}, {\bf F1}, and the holomorphicity of $u$, 
\begin{equation}\label{dbaru_Rest}
|\bar\pa u_R|=\Ordo(R^{1-2a})|Du|+\Ordo(R^{-(1+a)}).
\end{equation}
 
Taking derivatives of \eqref{dbaru_R} with respect to $\tau$ and $t$
we find (using {\bf F0-3} and \eqref{X_R0}-\eqref{X_R3})
\begin{align}\notag
|D\bar\pa u_R|=& \Ordo(R^{1-2a})|D^2u|
+\Ordo(R^{3-3a})|Du|^2\\\label{Ddbaru_R}
&+\Ordo(R^{1-2a})|Du|
+\Ordo(R^{-(1+a)})\\\notag
|\pa_\tau D\bar\pa u_R|=& \Ordo(R^{1-2a})|D^3 u|
+\Ordo(R^{3-3a})|Du||D^2u|
+\Ordo(R^{5-4a})|Du|^3\\\label{tauDdbaru_R}
&+\Ordo(R^{1-2a})|D^2u|
+\Ordo(R^{3-3a})|Du|^2
+\Ordo(R^{-(1+a)})|Du|.
\end{align}

Finally, let $\theta\colon[0,1]\to\R$ be a smooth function supported
in a $\frac12\delta$-neighborhood of the endpoints of the interval
with $\theta'(0)=\theta'(1)=1$. Define
\begin{equation}
\hat u_R(\tau+it)=u_R(\tau+it)+i\theta(t)\bar\pa u_R(\tau+it).
\end{equation}
Then $u_R=\hat u_R$ on the boundary and $\hat u_R$ is holomorphic on
the boundary. Also for some constant $C$
\begin{align}\label{hatu_R0}
|\hat u_R|&\le C(|u_R|+|\bar\pa u_R|),\\\label{hatu_R1}
|\bar\pa\hat u_R|&\le C(|\bar\pa u_R|+|D\bar\pa u_R|),\\\label{hatu_R2}
|D\bar\pa\hat u_R|&\le C(|\bar\pa u_R|+|D\bar\pa u_R|+
|\pa_\tau D\bar\pa u_R|). 
\end{align}

\subsection{Self-tangency preshortening}\label{10stpreshort.section}

Let $u\colon D_{m+1}\to\C^n$ be a rigid holomorphic disk with boundary
on $L$ and with negative punctures $p_1,\dots,p_k$ mapping to
$0$. (The case of one positive puncture mapping to $0$ is completely
analogous to the case of one negative puncture so for simplicity we
consider only the case of negative punctures.)

For large $\rho>0$ let $R=R(\rho)$ be such that the intersection
points of $L^1$ and 
$L^2_R(1)$ are $a^\pm=(\pm (\rho+\rho^a)^{-1},0,\dots,0)$. Then 
$R(\rho)=\rho+\Ordo(\rho^{-1})$. 
Define
\begin{equation}\notag
u_\rho(\zeta)=
\begin{cases}
u(\zeta)&\text{ for }\zeta\in 
D_{m+1}\setminus\Bigl(\bigcup_{j=1}^k E_{p_j}[-\frac12\rho]\Bigr),\\
\hat u_{R(\rho)}(\tau+it)&\text{ for }
\zeta=\tau+it\in E_{p_j}[-\frac12\rho].
\end{cases}
\end{equation}

Then there exist unique functions 
\begin{equation}\notag
\xi_R(j)\colon E_{p_j}[-\rho]\to T_{a^-}\C^n
\end{equation}
such that
\begin{equation}\notag
\exp^{R,t}(\xi_R(j)(\zeta))= u_\rho(\zeta),\quad\zeta\in
E_{p_j}[-\rho], 
\end{equation} 
where $\exp^{R,t}$ denotes the exponential map in the metric
$g(R,\omega(t),t)$ at $a^-$.

Let $\alpha_\rho\colon(-\infty,-\rho]\times[0,1]\to\C$ be a
smooth cut-off function, real valued and holomorphic on the
boundary and such 
that $\alpha_\rho(\tau+it)=1$ for $\tau$ in a small neighborhood of
$-\rho$, $\alpha_\rho(\tau+it)=0$ 
for $\tau\le -\rho-\frac12\rho^a$, and $|D^k\alpha_\rho|=\Ordo(\rho^{-a})$,
$k=1,2$. Define $w_\rho\colon D_{m+1}\to\C^n$ as 
\begin{equation}\label{wrhostsh}
w_\rho(\zeta)=
\begin{cases}
u_\rho(\zeta)&\text{ for }\zeta\in D_{m+1}\setminus\Bigl(\bigcup_j
E_{p_j}[-\rho]\Bigr),\\
\exp^{R,t}(\alpha_\rho(\zeta)\xi_R(j)(\zeta)) &\text{ for }
\zeta=\tau+it\in E_{p_j}[-\rho],\quad
j=1,\dots,k,\\
a_-&\text{ for }\zeta\in \bigcup_jE_{p_j}[-\rho-\frac12\rho^a].
\end{cases}
\end{equation}

\subsection{Weight functions for shortened disks}\label{10shortweight.section}

Let $u\colon D_{m+1}\to\C^n$ be a rigid holomorphic disk with boundary
on $L$. Let $\epsilon>0$ be small and let $e_\rho\colon D_{m+1}\to\R$
be a function which equals $e^{-\epsilon|\tau|}$ for $\tau+it\in
E_{p_j}\setminus E_{p_j}[-\rho]$ and is constantly equal to
$e^{-\epsilon\rho}$ for $\tau+it\in E_{p_j}[-\rho]$. Define
$\cand_{2,-\epsilon,\rho}$ just as in Section \ref{fcanI} but replacing
the weight function $e_\epsilon$ with the new weight function
$e_\rho$. The corresponding weighted norms will be denoted 
$\|\cdot\|_{2,-\epsilon,\rho}$. 
We also write
$\sblv_{1,-\epsilon,\rho}[0](D_{m+1},{T^\ast}^{0,1}\otimes\C^n)$ to
denote the subspace of elements in the Sobolev space with the weight
function $e_\rho$ which vanishes on the boundary. 

\subsection{Estimates for self-tangency preshortened disks}

\begin{lma}\label{stpresh}
The function $w_\rho$ in
\eqref{wrhostsh} lies in $\cand_{2,-\epsilon,\rho}$ (see Remark
\ref{dropthings.rmk} for notation) and there exists a
constant $C$ such that 
\begin{equation}\notag
\|\bar\pa w_\rho\|_{1,-\epsilon,\rho}\le
Ce^{-\epsilon\rho}\rho^{-1-\frac12 a}.
\end{equation}
\end{lma}
\begin{pf}
The first statement is obvious. Consider the second. In
$D_{m+1}\setminus\Bigl(E_{p_j}[-\rho]\Bigr)$, $w_\rho$ equals
$u$ which is holomorphic. It thus remains to check
$E_{p_j}[-\rho]\approx (-\infty,-\rho]\times[0,1]$.

Taylor expansion of $\exp^{R,t}$ gives
\begin{equation}
\exp^{R,t}\xi=
\xi-\Gamma_{ij}^k(R,t)\xi^i\xi^j\pa_k+\Ordo(|\xi|^3).
\end{equation}
The Taylor expansion of the inverse then gives
\begin{equation}
\xi_R=\hat u_R+\Gamma_{ij}^k(R,t)\hat u_R^i\hat u_R^j\pa_k+\Ordo(|\hat u_R|^3).
\end{equation} 
Thus in $(-\infty,-\rho]\times[0,1]$ we have
\begin{equation} \label{10wrho.eq}
w_\rho=\alpha_\rho
\hat u_{R}+
(\alpha_\rho-\alpha^2_\rho)\Gamma_{ij}^k(R,t)\hat u_R^i\hat u_R^j\pa_k
+\Ordo(|\hat u_R|^3).  
\end{equation}
Now $R=\rho+\Ordo(\rho^{-1})$ from Section \ref{10stpreshort.section}, 
$|D^k\alpha_\rho|=\Ordo(\rho^{-a})$ for all cut-off functions, 
and by Lemma \ref{SalRob} $|D^k u|=\Ordo(\rho^{-(1+k)})$,
in $(-\infty,-\rho]\times[0,1]$; thus, applying (\ref{hatu_R0}) through
(\ref{hatu_R2})
to (\ref{10wrho.eq}) we get
\begin{equation}\notag
|\bar\pa w_\rho|+|D\bar\pa w_\rho|=\Ordo(\rho^{-(1+a)}).
\end{equation}

Noting that $\bar\pa w_\rho$ is supported on an interval of length
$\frac12\rho^a$, so multiplying with the weight function we find
\begin{equation}\notag
\|\bar\pa w_\rho\|_{1,-\epsilon, \rho}\le
Ce^{-\epsilon\rho}\rho^{-1-\frac12 a}.
\end{equation}
\end{pf}

\subsection{Controlled invertibility for self-tangency shortening}\label{10stshortuni.section}

Let $d\Gamma_\rho$ denote the differential of the map
\begin{equation}\notag
\Gamma_\rho\colon\cand_{2,-\epsilon,\rho}\to\sblv_{1,-\epsilon,\rho}[0](D_{m+1},
{T^\ast}^{0,1}\otimes\C^n). 
\end{equation}
Referring to Sections \ref{fewmarked} and
\ref{manymarked}, we assume that $m\ge 2$ and 
$l(j)\ge 2$ for each $j$.

\begin{lma}\label{stshunif} 
There exist constants $C$ and $\rho_0$ such that 
if $\rho>\rho_0$ then there is a continuous right inverse 
$G_\rho$ of $d\Gamma_\rho$ 
\begin{equation}\notag
G_\rho\colon\sblv_{1,-\epsilon,\rho}[0]({T^\ast}^{0,1}D_r(\rho)\otimes\C^n)\to
T_{(w_\rho,\kappa_\rho,0)}\cand_{2,-\epsilon,\rho}
\end{equation} 
such that for any $\delta>0$
\begin{equation}\label{10stshortuni.eqn}
\|G_\rho(\xi)\|\le C\rho^{1+\delta}\|\xi\|_{1,-\epsilon,\rho}.
\end{equation}
\end{lma}

\begin{pf}
The kernel
\begin{equation}\notag
\ker(d\Gamma_{u})\subset T_u\cand_{2,-\epsilon}\oplus T_\kappa\conf_{m+1}
\end{equation}
has dimension $0$. By the invertibility
of $d\Gamma_u$ we conclude there is a constant $C$ such that for
$\xi\in T_u\cand_{\epsilon,2,\rho}$ we have 
\begin{equation} \label{10stshortuni2.eqn}
\|\xi\|\le C\|d\Gamma_{u,\rho}\xi\|_{1,-\epsilon}.
\end{equation}    
Assume that \eqref{10stshortuni.eqn} is not true. Then there exists a sequence
$\xi_N\in T_{w_\rho}\cand_{2,-\epsilon,\rho(N)}$ with
$\rho(N)\to\infty$ as $N\to\infty$ such that 
\begin{align}\label{sh=1}
&\|\xi_N\|=1,\\\label{shto0}
&\|d\Gamma_\rho\xi_N\|_{1,-\epsilon,\rho(N)}\le C\rho^{-1-\frac{\delta}{2}}.
\end{align}

Let $\alpha\colon D_{m+1}\to\C$ be a smooth function which equals $0$
on $E_{p_j}[-\rho-\frac14\rho^a]$ and equals $1$ on
$D_{m+1}\setminus\Bigl(\bigcup E_{p_j}[-\rho-10]\Bigr)$, which is real
valued and holomorphic on the boundary and with
$|D^k\alpha|=\Ordo(\rho^{-a})$, $k=1,2$. Then (\ref{10stshortuni2.eqn})
implies
\begin{equation}\label{extsh}
\|\alpha\xi_N\|\le C(\|(\bar\pa\alpha)\xi_N\|_{1,-\epsilon}
+\|\alpha d\Gamma_{u,\rho}\xi_N\|_{1,-\epsilon})=\Ordo(\rho^{-a})
\end{equation}

Finally, we let $\hat \phi\colon (-\infty,-\rho+\rho^a]\to\C$ be the function
which equals $\theta(\rho)-\theta(\tau)$, where $\theta(\tau)$ denotes
the angle between the tangent line of $L^2_\rho(1)$ intersected with
the $z_1$-plane (the plane of the first coordinate in $\C^n$) at
$u(\tau+i)$ and the real line in that plane. From Lemma \ref{SalRob}
we calculate that 
$|D^k\hat\phi|=\Ordo(\rho^{a-2})$, $0\le k\le 2$. Using the same procedure
as for cut-off functions we extend it to a function 
$\phi\colon (-\infty,-\rho+\rho^a)\times[0,1]$ which is holomorphic on
the boundary, which equals $\hat\phi$ on 
$(-\infty,-\rho+\rho^a)\times \{1\}$ and which equals $0$ on
$(-\infty,-\rho+\rho^a)\times\{0\}$ and with the same derivative estimates. 
Let ${\mathbf M}=\diag(\phi,1,\dots,1)$.

Let $\alpha$ be a cut-off function which is $0$ in 
$D_{m+1}\setminus E_{p_j}[-\rho+\rho^a]$ and $1$ on 
$E_{p_j}[-\rho].$ Having frozen the angle away from $0$, we can use
Lemmas \ref{linfred} and \ref{lmadef} (assuming
that $\epsilon$ is smaller than the smallest component of the
complex angle) to get
\begin{equation}
\|e^{-\epsilon\rho}\alpha {\mathbf M}\xi_N\|\le
C\rho(\|e^{-\epsilon\rho}(\bar\pa\alpha {\mathbf M})\xi_N\|+
\|e^{-\epsilon\rho}\alpha{\mathbf M}d\Gamma_\rho\xi_N\|)
\end{equation}
The first term on the right hand side inside the brackets is
$\Ordo(\rho^{a-2})+\Ordo(\rho^{-2a})$ the second term is
$\Ordo(\rho^{-1-\delta})$. Hence as $\rho\to\infty$ the right hand
side goes to $0$. This together with \eqref{extsh} contradicts
\eqref{sh=1} 
and we conclude the lemma holds. 
\end{pf}

\subsection{Perturbations for self-tangency gluing}
For $R>0$ with $R^{-1}<< r_0$, let
$a_R\colon[0,\infty)\to\R$ be a smooth non-increasing function with
support in $[0,\frac12R^{-1})$ and with the following properties
\begin{align}\notag
a_R(r)&=R^{-1}\text{ for }r\in [0,R^{-2}],\\\notag
|Da_R(r)|&=\Ordo(1),\\\label{Da_Rbounds}
|D^2 a_R(r)|&=\Ordo(R).
\end{align}
The existence of such functions is easily established. Let
$h_R\colon\C^n\to\C^n$ be given by 
\begin{equation}  
h_R(z)=x_1(z)a_R(|z_1|).
\end{equation}
For $s>0$, let $\Phi_R^s$ denote the 
time $s$ Hamiltonian flow of $h_R$ and let $L_R(s)$ denote
the Legendrian 
submanifold which results when $\Phi^s_R$ is lifted
to a local contact flow on $\C^n\times\R$ which is used to move
$L^2$. (Note that $\Phi_R^s$ fixes the last $n-1$ 
coordinates and has small support in the $z_1$-direction and so its
lift can be extended to the identity outside $L^2(s)$.) 
Let $L^2_R(s)=\Phi_R^s(L^2)$. We pick $a_R$ so that $L^2_R(s)\cap
L^1=\emptyset$, for $0<s\le (KR)^{-1}$ for some fixed $K>4$.   

As in Section~\ref{10perturbstshort.section} we derive the estimates
\begin{align}\label{Flow0}
&|\Phi^s_R-\id|\le \Ordo(R^{-2}),\\\label{Flow1}
&|D\Phi^s_R-\id|=\Ordo(R^{-1}),\\\label{Flow2}
&|D^2\Phi^s_R|=\Ordo(1).
\end{align}

For convenient notation we write
$\gamma^2_R(s)$ for the curve in which $L^2_R(s)$ intersects the
$z_1$-line in a neighborhood of $0$.

\subsection{Self-tangency pregluing}\label{10stpreglue.section}
Let $u\colon D_m\to\C^n$ be a rigid holomorphic disk with boundary on
$L$ and with negative punctures $p_1,\dots,p_k$
(as above we write $S=\{p_1,\dots,p_k\}$) mapping to $0$. Let 
$v_j\colon D_{l(j)+1}\to\C_n$ be rigid holomorphic disks with positive
punctures $q_j$ mapping to $0$. 

For $0<\rho<\infty$ let $R=\rho$, $s=(K\rho)^{-1}$ and let $L_\rho$ be
the Legendrian submanifold which results when $\Phi^s_R$ is applied.
Consider the region $\Xi_\rho$ in the $z_1$-line bounded by the curves
$\gamma^2_R(s)$, $\gamma^1_R(s)$, $u^1(\rho+it)$, $0\le t\le 1$, and
$v(j)^1(\rho+it)$, $0\le t\le 1$. By the Riemann mapping theorem
there exists a holomorphic map from a rectangle 
$\phi_\rho\colon [-A(\rho),A(\rho)]\times[0,1]\to\C$ which parameterizes
this region in such a way that $[-A(\rho),A(\rho)]\times\{j-1\}$ maps
to $\gamma_j(s)$, $j=1,2$. Moreover, since $\Xi_\rho$ is symmetric with
respect to reflections in the $\Im z_1=y_1$-axis we have
$\phi_\rho(0+i[0,1])\subset\{\Re z_1=x_1=0\}$.

\begin{lma}\label{lmashape}
The shape of the rectangle depends on $\rho$. More precisely, there
exists constants $0<K_1<K_2<\infty$ such that
$K_1\rho\le A(\rho)\le K_2\rho$ for all $\rho$.
\end{lma}

\begin{pf}
Identify the $z_1$-line with $\C$. Consider the region $\Theta_\rho$
bounded by the circles of radius $1$ and 
$1+4\rho^{-2}$ both centered at $i\in\C$, and the lines through $i$ which
intersects the $x_1$-axis in the points $\pm 2(\rho)^{-1}$. Mark the
straight line segments of its boundary. The
conformal modulus of this region is easily seen to be
$\rho+\Ordo(\rho^{-1})$. 

On the other hand, using \eqref{Flow0} and \eqref{Flow1} one 
constructs a $(K+\Ordo(\rho^{-1}))$-quasi conformal map from
$\Theta_\rho$ to $\Xi_\rho$, for some $K>0$ independent of
$\rho$. This implies the conformal modulus $m_\rho$ of $\Theta_\rho$
satisfies  
\begin{equation}
(K+\Ordo(\rho^{-1}))^{-1}(\rho+\Ordo(\rho^{-1}))\le
m_\rho\le
(K+\Ordo(\rho^{-1}))(\rho+\Ordo(\rho^{-1}))
\end{equation}   
and the lemma follows.
\end{pf}

Let $u^1$ and $v_j^1$ denote the $z_1$-components of the maps $u$ and
$v_j$. Since $\Phi^s_R$ fixes $\gamma_2$ outside $|x_1|\le(2\rho)^{-1}$
we note that 
\begin{equation}\label{extension+}
\text{$u^1$ maps the region 
$E_{p_j}[-\rho]\setminus E_{p_j}[-2\rho]$ 
into $\Theta_\rho\setminus \phi_\rho(0\times[0,1])$.}  
\end{equation}
and that
\begin{equation}\label{extension-}
\text{$v^1_j$ maps the region 
$E_{p_j}[\rho]\setminus E_{p_j}[2\rho]$ 
into $\Theta_\rho\setminus \phi_\rho(0\times[0,1])$.}  
\end{equation}

Fix $0<a<\frac14$. Using $u^1$, $v_j^1$, the conformal map
$\phi_\rho$ and their inverses, we construct a complex $1$-dimensional
manifold $D_r(\rho)$ by gluing  
$\Omega_j(\rho)=[-A(\rho),A(\rho)]\times[0,1]$ to 
$D_{m+1}\setminus \Bigl(\bigcup_j E_{p_j}[-(1+a)\rho]\Bigr)$ and
$D_{l(j)+1}E_q{[(1+a)\rho]}$. Note that, by construction, $D_r(\rho)$  
comes equipped with a holomorphic function 
\begin{equation}\label{w^1}
w_\rho^1\colon D_{r}(\rho)\to\C,
\end{equation}
which equals $u^1$ on $D_{m+1}\setminus\bigcup_j E_{p_j}[-(1+a)\rho]$,
which equals $v_j^1$ on $D_{l(j)+1}\setminus E_{q_j}[(1+a)\rho]$, and
which equals $\phi_\rho$ on $\Omega_j$, for all $j$.

We next exploit the product structure of $\Pi_\C L$ in a neighborhood
of $0$. If $u'$ and $v'_j$ denotes the remaining components of $u$ and
$v_j$ so that $u=(u^1,u')$ and $v_j=(v_j^1,v_j')$ then in some
neighborhood of the punctures $q_j$ and $p_j$, $v'_j$ and $u'_j$ are
holomorphic functions with boundary on the two transverse Lagrangian
manifolds $P_1$ and $P_2$, see Section~\ref{5Isotopies.section}. As in Section \ref{fcanC} we
find a $1$-parameter family $g(\sigma)$ of metrics on
$\C^{n-1}\approx\{z_1=0\}$. Then,
for $M$ sufficiently large,
there exist unique vector valued functions $\xi'_j$ and $\eta'_j$ such
that 
\begin{align}
\exp^t_0\xi_j'(\tau+it)=u'(\tau+it),\quad \tau+it\in
E_{p_j}[-M],\\
\exp^t_0\eta_j'(\tau+it)=v_j'(\tau+it),\quad \tau+it\in
E_{q_j}[M].
\end{align}
Now pick a cut-off function $\alpha^+$ which equals $1$ on 
$D_{m+1}\setminus \bigcup_jE_{p_j}[-\rho+5]$ and $0$ on
$E_{p_j}[-\rho+3]$. Pick similar cut-off functions $\alpha^-$ on
$D_{l(j)+1}$. Define $w_\rho'\colon D_{r}(\rho)\to\C^{n-1}$ by
\begin{equation}\label{w'}  
w_\rho'(\zeta)=\begin{cases}
u'(\zeta), &\quad \zeta\in D_{m+1}\setminus \bigcup E_{p_j}[-\rho+5],\\
v_j'(\zeta), &\quad \zeta\in D_{l(j)+1}\setminus E_{q_j}[\rho-5],\\
\exp^t_0(\alpha^+(\zeta)\xi_j(\zeta)),&\quad
\zeta\in E_{p_j}[-\rho+5]\setminus E_{p_j}[-\rho],\\
\exp^t_0(\alpha^-(\zeta)\eta_j(\zeta)),&\quad
\zeta\in E_{q_j}[\rho-5]\setminus E_{p_j}[\rho],\\
0,&\quad \zeta\in\Omega_j.
\end{cases}
\end{equation}
Finally, combining \eqref{w^1} and \eqref{w'}, we define 
\begin{equation}\label{wrhostglu}
w_\rho=(w^1_\rho,w'_\rho)
\end{equation}

\subsection{Weight functions and Sobolev norms for self tangency gluing}
Consider $D_r(\rho)$ from the previous section,
$\epsilon>0$, and a smooth function
$f\colon D_r(\rho)\to\C^n.$ 
Let 
\begin{itemize}
\item
$f^+$ denote the restriction of $f$ to 
\begin{equation*}
\inr\Bigl(D_{m+1}\setminus\bigcup_j E_{p_j}[-(1+a)\rho]\Bigr),
\end{equation*}
which we consider as a subset of $D_{m+1}$.
\item
$f^-$ denote the restriction of $f$ 
\begin{equation*}
\bigcup_j \inr\Bigl(D_{l(j)+1}\setminus E_{q_j}[(1+a)\rho]\Bigr),
\end{equation*}
which we consider as subset of the disjoint union
$\bigcup_j D_{l(j)+1}$   
\item
$f^0$ denote the restriction of $f$ to the disjoint union
$\bigcup_j\inr(\Omega_j(\rho))$. 
\end{itemize}

For $\epsilon>0$, let $e_\epsilon^-$ denote the weight function on
$\bigcup_j D_{l(j)+1}$ which equals $1$ on  $D_{l(j)+1}\setminus
E_{q_j}$ and equals $e^{\delta|\tau|}$ in $E_{q_j}$, each $j$. 
Let $\|\cdot\|_{k,\epsilon,-}$ denote the Sobolev norm with weight 
$e_\delta^-$. Let $e_\epsilon^0$ denote the weight function on
$\Omega_j$ which   
equals $e^{\epsilon(A(\rho)+\rho+\tau)}$ and 
$\|\cdot\|_{k,\epsilon,0}$ denote the
Sobolev norm with this weight. Finally, let $e_\epsilon^+$ be the
function on $D_{m+1}$ which 
equals $e^{2\epsilon (A(\rho)+\rho)}$ on $D_{m+1}\setminus \bigcup_j
E_{p_j}]$ and equals $e^{2\epsilon(A(\rho)+\rho)-\epsilon|\tau|}$ in
$E_{p_j}$. Let $\|\cdot\|_{k,\epsilon,+}$ denote the corresponding norm.

Define
\begin{equation}
\|f\|_{k,\epsilon,\rho}=\|f_+\|_{k,\epsilon,+}+\|f_0\|_{k,\epsilon,0}
+\|f_-\|_{k,\epsilon,-}. 
\end{equation}
 
Using this norm we define as in the shortening case the
spaces $\cand_{2,\epsilon,\rho}$ and
$\sblv_{1,\epsilon,\rho}[0](D_{m+1},{T^\ast}^{0,1}\otimes\C^n)$.
The $\bar\pa$-map 
$\Gamma\colon\cand_{2,\epsilon,\rho}\to\sblv_{1,\epsilon,\rho}[0]({T^\ast}^{0,1}D_m\otimes\C^n)$
is defined in the natural way.

\subsection{Estimates for self tangency glued disks}
\begin{lma}\label{stpreglu}
The function $w_\rho$ in
\eqref{wrhostglu} lies in $\cand_{2,\epsilon,\rho}$ and there exists a
constant $C$ such that 
\begin{equation}\notag
\|\bar\pa w_\rho\|_{1,\epsilon,\rho}\le
Ce^{(-\theta+2K_2\epsilon)\rho},
\end{equation}
where $\theta\gg\epsilon$ is the smallest non-zero complex angle at the
self tangency point $0$ and where $K_2$ is as in Lemma \ref{lmashape}.
\end{lma}

\begin{pf}
Note that the first coordinate of $w_\rho$ is holomorphic and that the
support of $\bar\pa w_\rho$ is disjoint from $\Omega_j$. 
Using the asymptotics
of $u'$ and $v'_j$ the proof of Lemma \ref{lmaprestat} applies to give the
desired estimate once we note that the weight function is
$\Ordo(e^{K_2\epsilon \rho})$ by Lemma \ref{lmashape}.
\end{pf}

\subsection{Estimates for real boundary conditions}
In order to prove the counterpart of Lemma \ref{stshunif} in the self tangency
gluing case we study an auxiliary non-compact counterpart of the
gluing region. 

Let $\Omega(\rho)=[-A(\rho),A(\rho)]\times [0,1]$ and let $M_\rho$ be
the complex 
manifold which results when $(-\infty,-(1-a)\rho]\times[0,1]$ and
$[(1-a)\rho,\infty)\times[0,1]$ are glued to $\Omega(\rho)$ with the
holomorphic gluing maps $u^1\circ(\phi_\rho)^{-1}$ and
$v_j^{1}\circ(\phi_\rho)^{-1}$, respectively.
(That is, the maps which were used to construct $D_r(\rho)$.) 
We consider Sobolev norms on $M_\rho$ similar to those used above. 

For $\epsilon>0$, let 
\begin{itemize}
\item
$e^0_\epsilon\colon\colon[-A(\rho),A(\rho)]\times[0,1]\to\R$ 
be the 
function $e^0_\epsilon(\tau+it)=e^{\epsilon\tau}$, 
\item
$e^-_\epsilon\colon (-\infty,-(1-a)\rho]\times[0,1]\to\R$ be the
function $e^-_\epsilon(\tau+it)=e^{\epsilon(\rho-A(\rho)+\tau)}$       
\item
$e^+_\epsilon\colon [(1-a)\rho,\infty)\times[0,1]\to\R$ be the
function $e^+_\epsilon(\tau+it)=e^{\epsilon(-\rho+A(\rho)+\tau)}$       
\end{itemize}

If $f\colon M_\rho\to\C$ is function we let as above $f^-, f^0, f^+$
denote the restrictions of $f$ to the interiors of the pieces from
which $M_\rho$ was constructed and define the Sobolev norm
\begin{equation}
\|f\|_{k,\rho,\epsilon}=\|f^-\|_{k,\epsilon}+\|f^0\|_{k,\epsilon}
+\|f^+\|_{k,\epsilon}. 
\end{equation} 
 
\begin{lma}\label{lmaauxest}
There are constants $C$ and $\rho_0$ if $\rho>0$ and if 
$f\colon M_\rho\to\C$ is function which is real valued and holomorphic
on the boundary and has $\|f\|_{k,\rho,\epsilon}\le \infty$ then 
\begin{equation}\label{auxRstg}
\|f\|_{k,\rho,\epsilon}\le C\|\bar\pa f\|_{k-1,\rho,\epsilon},
\end{equation}
for $k=1,2$.
\end{lma}

\begin{pf}
To prove the lemma we first study the gluing functions. 
Let $\psi\colon [-\rho,-(1-a)\rho)\times[0,1]\to
[-\A(\rho),0]\times[0,1]$ be the function
$u^1\circ(\phi_\rho)^{-1}$. Note that $\psi$ is holomorphic and that by
\eqref{extension+} has a holomorphic extension (still denoted $\psi$)
to $[-\rho, 0)\times [0,1]$. 

To simplify notation we change coordinates and
think of the source $[-\rho,0)\times[0,1]$ as $[0,\rho)\times[0,1]$
and of the target $[-A(\rho),0]\times[0,1]$ as
$[0,A(\rho)]\times [0,1]$. Thus 
\begin{equation}
\psi\colon [0,\rho]\times[0,1]\to [0,A(\rho)]\times[0,1]
\end{equation}
is a holomorphic map. Consider the complex derivative
$\frac{\pa\psi}{\pa z}$. This is again a holomorphic function which is
real on the boundary of $[0,\rho)\times[0,1]$. In analogy with Lemma
\ref{lmaralocsol} we conclude that
\begin{equation}
\frac{\pa\psi}{\pa z}=\sum_{n\in Z}c'_n e^{n\pi z},
\end{equation}
for some real constants $c'_n$. Integrating this and using $\psi(0)=0$,
we find
\begin{equation}
\psi(z)=c_0 z +\sum c_n e^{n\pi z}, 
\end{equation}
for some real constants $c_n$. Then
\begin{equation}
i=\psi(i)=c_0 i +\sum c_n e^{n\pi i}
\end{equation}  
and we conclude $c_0=1$. Moreover, if $\psi^d$ denotes the double of
$\psi$ (which has the same Fourier expansion) then since $\psi^d(it)$
is purely imaginary for $0\le t\le 2$, we find that $c_n=-c_n$ for all
$n\ne 0$. Thus 
\begin{equation}\label{explform}
\psi(z)=z+\sum_n c_n(e^{n\pi z}-e^{-n\pi z}).
\end{equation} 
The area of the image of $\psi^d$ is $\Ordo(\rho)$ by Lemma
\ref{lmashape}. Since this area equals the $L^2$-norm of the derivative of
$\psi^d$ we conclude that
\begin{equation}
2\int_0^\rho 1^2\,d\tau + \sum_{n\in Z}\int_0^\rho n^2\pi^2 |c_n|^2
e^{2n\pi\tau}\,d\tau =\Ordo(\rho).
\end{equation}
Integrating we find there exists a constant $K$ and $0<\delta\ll 1$
such that  
\begin{equation}
|c_n|\le K \rho (n)^{-\frac12} e^{-n\pi\rho}\le K e^{-n(\pi-\delta)\rho},
\end{equation} 
for each $n\ne 0$. Thus, in the gluing region $[0,a\rho)$ we find
\begin{equation}
|\psi(z)-z|\le K \sum_{n>0} e^{-n(\pi-\delta-a)\rho}\le 
K' e^{-(\pi-2(\delta+a))\rho}=K'e^{-\eta\rho},
\end{equation}
where $\eta>0$. Similarly one shows 
$|D\psi-\id|\le K e^{-\frac12\eta\rho}$ and 
$|D^2\psi|\le K e^{-\frac12\eta\rho}$.  

Assume \eqref{auxRstg} is not true then there exists a sequence $f_j$
of functions on $M_{\rho(j)}$, $\rho(j)\to\infty$ as $j\to\infty$,
with 
\begin{align}\label{aux=1}
&\|f_j\|_{2,\rho,\epsilon}=1,\\\label{auxto0}
&\|\bar\pa f_j\|_{1,\rho,\epsilon}\to 0,\quad\text{ as } j\to\infty.
\end{align}

Let $\gamma\colon (-\infty,-(1-a)\rho]\times[0,1]$ be a cut-off
function which equals $1$ on $(-\infty,-(1-\frac14a)\rho]\times[0,1]$ which equals
$0$ on $[-(1-\frac12 a)\rho,-(1-a)\rho)$, has
$|D^k\gamma|=\Ordo(\rho^{-1})$, $k=1,2$, and is real valued and
holomorphic on the boundary. Then by uniform invertibility of the
$\bar\pa$-operator on the strip with constant weight $\epsilon$ we
find
\begin{equation}
\|\gamma f\|_{2,\epsilon}\le C(\|(\bar\pa\gamma)
f\|_{1,\epsilon}+\|\gamma\bar\pa f\|_{1,\epsilon}).
\end{equation}   
Here both terms on the right hand side goes to $0$ as
$\rho\to\infty$. In a similar way we conclude
\begin{equation}
\|\beta f\|_{2,\epsilon}\to 0,
\end{equation}
for $\beta$ a cut-off function on $[(1-a)\rho,\infty)$.

Now let $\alpha$ be a cut-off function on
$[-A(\rho),A(\rho)]\times[0,1]$ which equals $1$ on $[-A(\rho)+2,
A(\rho)-2]\times [0,1]$ and equals $0$ outside $[-A(\rho)+1,
A(\rho)-1]\times[0,1]$. We find
\begin{equation}
\|\alpha f\|_{2,\epsilon}\le C(\|(\bar\pa\alpha) f\|_{1,\epsilon}
+\|\alpha \bar\pa f\|_{1,\epsilon}).
\end{equation} 
Here the second term on the right hand side goes to $0$ as
$\rho\to\infty$ by \eqref{auxto0}. The first goes to $0$ as well since
$\|\gamma f\|\to 0$ and $\|\beta f\|\to 0$ and since the transition
functions are very close to the identity for $\rho$ large. 

In
conclusion we find
$\|f\|_{2,\rho,\epsilon}\to 0$, contradicting \eqref{aux=1}, and
\eqref{auxRstg} 
holds. 
\end{pf}

\subsection{Uniform invertibility for self tangency gluing}\label{10stuni.section}
Let $d\Gamma_\rho$ denote the differential of the map
\begin{equation}\notag
\Gamma\colon\cand_{2,\epsilon,\rho}\to\sblv_{1,\epsilon,\rho}[0]
(D_{m+1}, {T^\ast}^{0,1}\otimes\C^n), 
\end{equation}
at $w_\rho$. Referring to Sections \ref{fewmarked} and
\ref{manymarked}, we assume that $m\ge 2$ and 
$l(j)\ge 2$ for each $j$.

\begin{lma}\label{stgluunif} 
There exist constants $C$ and $\rho_0$ such that 
if $\rho>\rho_0$ and then there is a continuous right inverse 
$G_\rho$ of $d\Gamma_\rho$ 
\begin{equation}\notag
G_\rho\colon\sblv_{1,\epsilon,\rho}[0]({T^\ast}^{0,1}D_r(\rho)\otimes\C^n)\to
T_{(w_\rho,\kappa_\rho)}\cand_{2,\epsilon,\rho}
\end{equation} 
such that 
\begin{equation}\notag
\|G_\rho(\xi)\|\le C\|\xi\|_{1,\epsilon,\rho}.
\end{equation}
\end{lma}

\begin{pf}
Recall $0<\epsilon\ll\theta$, where $\theta>0$ is the smallest
non-zero complex angel at the self-tangency point.
Assume we glue $k$ disks $v_1,\dots,v_k$ to $u$.
The kernels 
\begin{align}\notag
d\Gamma_{(u,\kappa_1)}\subset
T_{(u,\kappa_1)}\cand_{2,-\epsilon},\\
d\Gamma_{(v_j,\kappa_2(j))}\subset
T_{(v_j,\kappa_2(j))}\widetilde{\cand}_{2,\epsilon},
\end{align}
are both of dimension $0$ and $d\Gamma_{(u,\kappa_1)}$ and
$d\Gamma_{(v_j,\kappa_2(j))}$ are invertible.

As usual we consider the embedding
\begin{equation}
T_{\kappa_1}\conf_{m+1}\oplus
\bigoplus_{j=1}^k T_{\kappa_2(j)}\conf_{l(j)+1}\to
T_{\kappa_\rho}\conf_r, 
\end{equation}
which identifies the left-hand side with a subspace of codimension $k$
in $T_{\kappa_\rho}\conf_r$. Let $W_\rho$  denote the complement of this
subspace in $T_{(w_\rho,\kappa_\rho)}\cand_{2,\epsilon,\rho}$. We show
that there exists a constant $C$ such that for $\rho$ large enough and
$(\xi,\gamma)\in W_\rho$ 
\begin{equation}\label{stglW}
\|(\xi,\gamma)\|\le C\|d\Gamma_\rho(\xi,\gamma)\|.
\end{equation}
Assume \eqref{stglW} is not true then there exists a sequence
$(\xi_N,\gamma_N)\in W_{\rho(N)}$, where $\rho(N)\to\infty$ as
$N\to\infty$ with
\begin{align}\label{stgl=1}
&\|(\xi_N,\gamma_N)\|=1,\\\label{stglto0}
&\|d\Gamma_{\rho(N)}(\xi_N,\gamma_N)\|\to 0,\quad\text{ as }N\to\infty.
\end{align} 

Let $\beta^0_{\rho}\colon D_r(\rho)\to\C$ be a cut-off function which
equals $1$ on $D_{m+1}\setminus (\bigcup_j E_{p_j}[-\frac12\rho])$,
equals $0$ outside $D_{m+1}\setminus (\bigcup_j E_{p_j}[-\frac34\rho])$,
with $|D^k\beta^0_{\rho}|=\Ordo(\rho^{-1})$, $k=1,2$. By the uniform
invertibility of $d\Gamma_{(u,\kappa_1)}$ we find
\begin{align}
\|\beta^0_{\rho(j)}(\xi_N,\gamma_N)\|_{2,\epsilon,\rho}&\le
C\|d\Gamma_{(u,\kappa_1)}\beta^0_{\rho(N)}
(\xi_N,\Gamma_N)\|_{1,\epsilon,\rho}\le\\
&C\Bigl(\|(\bar\pa \beta^0_{\rho(N)})\xi_N\|_{1,\epsilon,\rho}
+\|\beta^0_{\rho(N)}d\Gamma_{\rho}(\xi_N,\gamma_N)\|_{1,\epsilon,\rho}\Bigr).
\end{align} 
Both terms on the right hand side goes to $0$ as $N\to\infty$. Hence
\begin{equation}\label{out+}
\|\beta^0_{\rho(N)}(\xi_N,\gamma_N)\|_{2,\epsilon,\rho}\to
0,\quad\text{ as }N\to\infty.
\end{equation} 

Similarly, with $\beta^j_\rho\colon D_r(\rho)\to\C$
a cut-off function which
equals $1$ on $D_{l(j)+1}\setminus E_{q_j}[\frac12\rho]$,
equals $0$ outside $D_{l(j)+1}\setminus  E_{q_j}[\frac34\rho]$,
with $|D^k\beta^0_{\rho}|=\Ordo(\rho^{-1})$, $k=1,2$, we find, by the
uniform invertibility of $d\Gamma_{(v_j,\kappa_2(j))}$ that
\begin{equation}\label{out-}
\|\beta^j_{\rho(N)}(\xi_N,\gamma_N)\|_{2,\epsilon,\rho}\to
0,\quad\text{ as }N\to\infty\text{ for all }j.
\end{equation} 

For $1\le j\le k$ we consider the region 
\begin{align}
&\Theta_j(\rho)=\\
&\Bigl(E_{p_j}\setminus E_{p_j}[-(1+a)\rho]\Bigr)
\cup_{\bigl((\phi_{\rho})^{-1}\circ
u^1\bigr)}\Omega_j\cup_{\bigl((\phi_{\rho})^{-1}\circ v_j^1\bigr)} 
\Bigl(E_{q_j}\setminus E_{q_j}[(1+a)\rho]\Bigr).
\end{align}
Note that there is a natural inclusion $\Theta_j(\rho)\subset M_\rho$,
where $M_\rho$ is as in Lemma \ref{lmaauxest}. Also note that the boundary
conditions of the linearized equation over 
$\Omega_j(\rho)$ splits into a $1$-dimensional problem corresponding
to the first coordinate and an $(n-1)$-dimensional problem with
boundary conditions converging to two transverse Lagrangian subspaces
in the remaining $(n-1)$ coordinate directions. 

Let $\alpha^+_\rho$ be a cut-off
function on $\Theta_j(\rho)$ which equals $1$ on 
\begin{equation}
E_{p_j}[-\frac14\rho]\setminus E_{p_j}[-(1+\frac12a)\rho],
\end{equation}
equals $0$ outside
\begin{equation}
E_{p_j}[-\frac18\rho]\setminus E_{p_j}[-(1+\frac23a)\rho],
\end{equation}
with $|D^k\alpha^+|=\Ordo(\rho^{-1})$, $k=1,2$, and which is real
valued and holomorphic on the boundary. Note that over the region
where $\alpha^+$ is supported the boundary conditions of $w_\rho$
agrees with those of $u$. Thus the angle between the line giving the
boundary conditions of $w_\rho$ and the real line is $\Ordo(\rho^{-1})$ and
it is easy to construct a unitary diagonal matrix function 
${\mathbf M}$ on the support $\alpha^+$ with 
$|D^k {\mathbf M}|=\Ordo(\rho^{-1})$, $k=1,2$ with the property that
${\mathbf M}\xi_N$ has the boundary conditions of $w_\rho$ in the last
$(n-1)$ coordinates and has real boundary conditions in the first
coordinate. Thus Lemma \ref{lmaauxest} implies that
\begin{equation}\label{middle+}
\|\alpha^+\xi_N\|_{2,\rho,\epsilon}\le C\|{\mathbf M}\alpha^+\xi_N\|
\le C\Bigl(\|(\bar\pa {\alpha^+\mathbf
M})\xi\|_{1,\epsilon,\rho}
+\|{\mathbf M}\bar\pa \xi_N\|_{1,\epsilon,\rho}
\Bigr).
\end{equation}  
Here both terms in the right hand side goes to $0$ as $N\to\infty$.

In exactly the same way we show that 
\begin{equation}\label{middle-}
\|\alpha^-\xi_N\|\to 0\quad \text{ as }\rho\to \infty,
\end{equation}
for a cut-off function $\alpha^-$ with support on the other end of
$\Theta_\rho$. 

Let $\alpha^0$ be a cut-off function which equals $1$ on 
$[-A(\rho)+2,A(\rho)-2]\times[0,1]$ and equals $0$
outside
$[-A(\rho)+1,A(\rho)-1]\times[0,1]$
and with the usual properties. 
Then the function 
\begin{equation}
(\tau+it)\mapsto (d\Phi^{ts(N)}_{R(N)}(\phi_\rho(\tau+it)))^{-1}
\cdot\alpha^0(\tau+it)\xi_N(\tau+it)
\end{equation}
has the boundary conditions of $w_\rho$ in the last
$(n-1)$ coordinates (two transverse Lagrangian subspaces in this
region) and has real boundary conditions in the first coordinate.

Lemma \ref{lmaauxest} implies
\begin{align}\notag
\|&\alpha^0\xi_N\|_{1,\rho,\epsilon}\le C
\|(d\Phi^{ts(N)}_R(N))^{-1}\cdot 
\alpha^0\xi_N\|_{1,\rho,\epsilon}\le\\\label{rhs''}
&C\Bigl(\|\bar\pa (\alpha^0 d\Phi^{ts(N)}_R(N))^{-1})\cdot
\xi_N\|_{0,\rho,\epsilon}+ 
\|(\alpha^0d\Phi^{ts(N)}_R(N))^{-1})\cdot 
\bar\pa \xi_N\|_{0,\rho,\epsilon}\Bigr).
\end{align}
Using \eqref{Flow1} and \eqref{Flow2} in combination with
\eqref{middle+} and \eqref{middle-} 
we find that 
the first term in \eqref{rhs''} goes to $0$ as $N\to 0$. By
\eqref{stglto0}, so does the second. Hence
\begin{equation}
\|\alpha^0\xi_N\|_{1,\rho,\epsilon}\to 0.
\end{equation}
Applying the same argument to $\pa_\tau \xi_N$ and $i\pa_t \xi_N$ we
conclude that
\begin{equation}\label{middle}
\|\alpha^0\xi_N\|_{2,\rho,\epsilon}\to 0.
\end{equation}

Now \eqref{out-}, \eqref{out+}, \eqref{middle+}, \eqref{middle-}, and
\eqref{middle} contradict \eqref{stgl=1} and we find that
\eqref{stglW} holds.

To finish the proof we let $\mu_j=\bar\pa\frac{\pa\phi^{C_j}}{\pa C_j}$,
see Section \ref{8degentrans.section}. Then $\mu_j$ anti-commutes with
$j_{\kappa_\rho}$ and we consider the $\mu_j$ as newborn conformal
variations spanning the complement of $W_\rho$ in
$T_{(w_\rho,\kappa_\rho)}\cand_{2,\epsilon,\rho}$. 

The images of $\mu_j$, $j=1,\dots,k$ under $d\Gamma_\rho$ are clearly
linearly independent since they have mutually disjoint supports. We
show that their images stays a uniformly bounded distance away from
the subspace $d\Gamma_\rho(W_\rho)$. Assume not, then there exists a
sequence of elements $(\xi_\rho,\gamma_\rho)$ in $W_\rho$ with 
\begin{equation}
\|d\Gamma_\rho\bigl((\xi_\rho,\gamma_\rho)-\mu_j\bigr)\|_{1,\epsilon,\rho}\to
0 \quad \text{ as }\rho\to\infty.
\end{equation}
Since $\|d \Gamma_\rho\mu_j\|_{1,\epsilon,\rho}=\Ordo(1)$ we conclude
from \eqref{stglW} that
$\|(\xi_\rho,\gamma_\rho)\|_{2,\epsilon,\rho}=\Ordo(1)$. 
Then, with the cut-off function $\beta^j_\rho$ from above and notation
as in Section \ref{8degentrans.section} we find
\begin{align}
&\|d\Gamma_{(v_j,\kappa_2(j))}
(\beta^j_\rho(\xi_\rho,\gamma_\rho)-\hat C_j)\|_{1,\epsilon}= 
\|d\Gamma_\rho
(\beta^j_\rho(\xi_\rho,\gamma_\rho)-\mu_j)\|_{1,\epsilon,\rho}\le\\
&\|\beta^j_\rho(d\Gamma_\rho(\xi_\rho,\gamma_\rho)-\mu_j)\|_{1,\epsilon,\rho}
+\|(\bar\pa\beta^j_\rho)((\xi_\rho,\gamma_\rho)-\mu_j)\|_{1,\epsilon,\rho}.
\end{align}
The right hand side of the above equation goes to $0$ as
$\rho\to\infty$. Hence so does the left hand side. This however
contradicts the invertibility of $d\Gamma_{(v_j,\kappa_2(j))}$ and we
conclude $d\Gamma_\rho(W_\rho)$ stays a bounded distance away from
$d\Gamma_\rho(\mu_j)$. Thus, defining
$G_\rho(d\Gamma_{\rho}\mu_j)=\mu_j$ finishes the proof.
\end{pf}

\subsection{Estimates on the non-linear term}\label{10nonlinear.section}

In Section \ref{fcanH}, we linearized the map $\Gamma$ using local
coordinates $B$ around $(w,f)\in\cand_{2,\epsilon}$. To apply Floer's
Picard lemma, we must study also higher order
variations of $\Gamma$. 

For $(w,f,0)\in\cand_{2,\epsilon,\Lambda}$, $w\colon D_m\to\C^n$ and
conformal structure $\kappa$ on $D_m$, we take as in
Section \ref{fcanF} local 
coordinates $(v,\kappa)\in B_{2,\epsilon}\times\R^{m-3}\times\Lambda$ 
on $\cand_{2,\epsilon,\lambda}$ around $(w,f)$ and write (in these
coordinates)
\begin{equation}\notag
\Gamma(v,\lambda, \gamma)=\bar\pa_\kappa v+i\circ dw\circ\gamma
+\lambda\cdot \bar\pa_\kappa Y_0^\sigma +N(v,\lambda, \gamma).
\end{equation}  
We refer to $N(v,\lambda,\gamma)$ as the {\em non-linear term}. 
We consider $\Lambda$ to have dimension $0$ in the stationary and
the self tangency cases and have dimension $1$ in the handle
slide case. We first consider stationary gluing
\begin{lma}\label{statnl}
There exists a constant $C$ such that
the non-linear term $N(v,\gamma)$ of $\Gamma$ in a neighborhood
$w_\rho$, where $w_\rho$ is as in Section \ref{10statpreglue.section} 
satisfies
\begin{equation}\label{statnonl}
\|N(u,\beta)-N(v,\kappa)\|_1\le
C\Bigl(\|u\|_2+|\beta|+\|v\|_2+|\gamma|\Bigr)
\Bigl(\|u-v\|_2+|\beta-\gamma|\Bigr) 
\end{equation} 
\end{lma}
\begin{pf}
With notation as in Section \ref{fcanF} we have
\begin{equation}\notag
\Gamma(v,\gamma)=\bar\pa_{\kappa+\gamma}
\Bigl(\exp^{\sigma(\zeta)}_{w_\rho(\zeta)}v(\zeta)\Bigr).
\end{equation}
We prove \eqref{statnonl} first in the special case $\gamma=\beta=0$.
We perform our calculation in coordinates $x+iy$ on $D_r(\rho)$, which
agree with the standard coordinates on the ends and in the gluing
region on $D_r(\rho)$. On the remaining parts of the disk the
metric of these coordinates differs from the usual metric by a
conformal factor but since the remaining part is compact the estimates
are unaffected by this change of metric. In these coordinates we write
$\bar\pa_\kappa=\pa_x+i\pa_y$. Now, as in Lemma \ref{Bdryval} we find
\begin{equation}\notag
\pa_x\exp^{\sigma}_{w_\rho}v=
J[w_\rho,v,\pa_x w_\rho,\pa_x
v,\sigma](1)+\pa_\sigma(\exp^\sigma_{w_\rho}v)\cdot\pa_x\sigma, 
\end{equation}
where $J[x,\xi,x',\xi',\sigma]$ denotes the Jacobi field in the metric
$g(\sigma)$ along the geodesic $\exp_x^\sigma t\xi$ with initial
conditions $J(0)=x'$, $J'(0)=\xi'$. Of course a similar equation
holds for $\pa_y\exp_{w_\rho}^\sigma v$. 

Let $G\colon(\C^{n})^4\times[0,1]\times\R\to\C^n$ be the function 
\begin{equation}\notag
G(x,\xi,x',\xi',\sigma,\sigma')=
J[x,\xi,x',\xi',\sigma](1)-x'-\xi'+\pa_\sigma\exp^\sigma_x\xi\cdot\sigma'
\end{equation}
(unrelated to the earlier right inverses $G_\rho$) 
then with $w_\rho=w$,
\begin{equation}\notag
N(v)=G(w,v,\pa_xw,\pa_xv,\sigma,\pa_x\sigma)+iG(w,v,\pa_y w,\pa_y
v,\sigma,\pa_y\sigma). 
\end{equation}
Moreover, the function $G$ is smooth with uniformly bounded
derivatives, it is linear in $x',\xi',\sigma'$, and satisfies
\begin{align}\notag
&G(x,0,x',\xi',\sigma,\sigma') =0,\\\label{getgoing}
&D_2 G(x,0,x',\xi',\sigma,\sigma')=0,
\end{align}
where the last equation follows from Taylor expansion of the
exponential map and the Jacobi field.

We estimate the $1$-norm of 
\begin{equation}\notag
G(w,u,\pa_xw,\pa_xu,\sigma,\pa_x\sigma)- 
G(w,v,\pa_xw,\pa_xv,\sigma,\pa_x\sigma).
\end{equation}
For the $0$-norm, we write
\begin{align}\notag
&\Bigl|G(w,u,\pa_xw,\pa_xu,\sigma,\pa_x\sigma)- 
G(w,v,\pa_xw,\pa_xv,\sigma,\pa_x\sigma)\Bigr|\le\\\notag
&\Bigl|G(w,u,\pa_xw,\pa_xu,\sigma,\pa_x\sigma)- 
G(w,u,\pa_xw,\pa_xv,\sigma,\pa_x\sigma)\Bigr|\\\notag
&+\Bigl|G(w,u,\pa_xw,\pa_xv,\sigma,\pa_x\sigma)- 
G(w,v,\pa_xw,\pa_xv,\sigma,\pa_x\sigma)\Bigr|\le\\\notag
&C\Bigl(|u||Du-Dv|+|Dv||u-v|\Bigr)\le\\\label{0norm'}
&C\Bigl((|u|+|v|)|Du-Dv|+(|Du|+|Dv|)|u-v|\Bigr),
\end{align}
where we use \eqref{getgoing}. Noting that the $\|\cdot\|_2$-norm
controls the $\sup$-norm we square and integrate \eqref{0norm'} to
conclude
\begin{equation}\label{0norm}
\|N(u)-N(v)\|_0\le C\Bigl(\|u\|_2+\|v_2\|\Bigr)\|u-v\|_2.
\end{equation}

For the $1$-norm we find an $L^2$-estimate of
\begin{align}\notag
&\Bigl|
D_1G(w,u,\pa_xw,\pa_xu,\sigma,\pa_x\sigma)- 
D_1G(w,v,\pa_xw,\pa_xv,\sigma,\pa_x\sigma)
\Bigr||Dw|\\\notag
&+\Bigl|
D_2G(w,u,\pa_xw,\pa_xu,\sigma,\pa_x\sigma)\cdot Du- 
D_2G(w,v,\pa_xw,\pa_xv,\sigma,\pa_x\sigma)\cdot Dv
\Bigr|\\\notag
&+\Bigl|
D_3G(w,u,\pa_xw,\pa_xu,\sigma,\pa_x\sigma)\cdot D\pa_x w- 
D_3G(w,v,\pa_xw,\pa_xv,\sigma,\pa_x\sigma)\cdot D\pa_x w
\Bigr|\\\notag
&+\Bigl|
D_4G(w,u,\pa_xw,\pa_xu,\sigma,\pa_x\sigma)\cdot D\pa_x u- 
D_4G(w,v,\pa_xw,\pa_xv,\sigma,\pa_x\sigma)\cdot D\pa_x v
\Bigr|\\\notag
&+\Bigl|
D_5G(w,u,\pa_xw,\pa_xu,\sigma,\pa_x\sigma)- 
D_5G(w,v,\pa_xw,\pa_xv,\sigma,\pa_x\sigma)
\Bigr||D\sigma|\\\label{toolong}
&+\Bigl|
D_6G(w,u,\pa_xw,\pa_xu,\sigma,\pa_x\sigma)\cdot D\pa_x\sigma- 
D_6G(w,v,\pa_xw,\pa_xv,\sigma,\pa_x\sigma)\cdot D\pa_x\sigma
\Bigr|.
\end{align}
Using \eqref{getgoing} the first and fifth 
terms in \eqref{toolong} are estimated by
\begin{equation}\label{15} 
C(|u|+|v|)(|u-v|+|Du-Dv|)(|Dw|+|DF|),
\end{equation}
where $F$ is the extension of $f\colon \pa D_r(\rho)\to\R$ as in
Section \ref{fcanF}. 
The second term in \eqref{toolong} is estimated by
\begin{equation}\label{2}
C\Bigl((|u|+|v|)|Du-Dv|+(|Du|+|Dv|)|Du-Dv|
+|u-v|(|Du|+|Dv|)
\Bigr).
\end{equation}
For the remaining terms we use the linearity of $G$ in
$x',\xi',\sigma'$ to write them as
\begin{align}\notag
&\Bigl|
D_3G(w,u,D \pa_xw,\pa_xu,\sigma,\pa_x\sigma)- 
D_3G(w,v,D \pa_xw,\pa_xv,\sigma,\pa_x\sigma)
\Bigr|\\\notag
&+\Bigl|
D_4G(w,u,\pa_xw,D \pa_xu,\sigma,\pa_x\sigma)- 
D_4G(w,v,\pa_xw,D \pa_xv,\sigma,\pa_x\sigma)
\Bigr|\\\notag
&+\Bigl|
D_6G(w,u,\pa_xw,\pa_xu,\sigma,D \pa_x\sigma)- 
D_6G(w,v,\pa_xw,\pa_xv,\sigma,D \pa_x\sigma)
\Bigr|.
\end{align}
Thus the third and sixth terms in \eqref{toolong} are estimated as in
\eqref{0norm'} by
\begin{equation}\label{36}
C\Bigl((|u|+|v|)|Du-Dv|+(|Dv|+|Du|)|u-v|\Bigr).
\end{equation}
Finally, the fourth term in \eqref{toolong} is estimated by
\begin{equation}\label{4}
C\Bigl((|u|+|v|)|D^2u-D^2 v|+(|D^2 u|+|D^2 v|)|u-v|
\Bigr)
\end{equation}
To estimate the $L^2$-norms of these expressions we use the
$\sup$-norm bound of $|u|$ and $|v|$ and the fact that the
$\|\cdot\|_1$-norm controls the $L^p$-norm for all $p>2$. For 
example 
\begin{align}\notag
&\int_{D_m}(|Du|+|Dv|)^2|Du-Dv|^2\,dA\le\\\notag
&\Bigl(\int_{D_m}(|Du|+|Dv|)^4\Bigr)^{\frac12}
\Bigl(\int_{D_m}(|Du-Dv|)^4\Bigr)^{\frac12}\le\\\label{L4ex}
&(\|u\|_2+\|v\|_2)^2(\|u-v\|_2)^2.
\end{align} 
We conclude that we have the estimate
\begin{equation}\notag
\|N(u)-N(v)\|_1\le C\Bigl(\|u\|_2+\|v\|_2\Bigr)\|u-v\|_2,
\end{equation}
as desired.

Finally in the case when also the conformal structures changes we note
that if $j_\kappa$ in coordinates $x+iy$ is represented by the matrix 
\begin{equation}\notag
\left(\begin{matrix}
0 & -1\\
1 & 0
\end{matrix}\right)
\end{equation}
then $j_{\kappa+\gamma}$ is represented by the matrix
\begin{equation}\notag
\left(\begin{matrix}
\phi_\gamma & -1\\
1+\phi_\gamma^2 & -\phi_\gamma
\end{matrix}\right),
\end{equation}
where $\phi_\gamma\colon D_m\to\R$ is a compactly supported function. The
extra term which enters in the non-linear term is then 
\begin{equation}\notag
(\phi_\beta^2G(w,u,\pa_yw,\pa_yu,\sigma,\pa_y\sigma)
-\phi_\gamma^2G(w,v,\pa_y w,\pa_y v,\sigma,\pa_y\sigma)),  
\end{equation}
which is easily estimated using the techniques above.
\end{pf} 

For handle slide gluing the corresponding lemma reads.
\begin{lma}\label{handslnl}
There exists a constant $C$ such that
the non-linear term $N(v,\lambda, \gamma)$ of $\Gamma$ in a neighborhood
of $w_\rho$, where $w_\rho$ is as in Section \ref{10hspreglue.section} 
satisfies
\begin{equation}\notag
\|N(u,\beta,\lambda)-N(v,\kappa,\mu)\|_1\le
C\Bigl(\|u\|_2+|\beta|+|\lambda|+\|v\|_2+|\gamma|+|\mu|\Bigr)
\Bigl(\|u-v\|_2+|\beta-\gamma|+|\lambda-\mu|\Bigr) 
\end{equation} 
\end{lma}
\begin{pf}
The proof is similar to the proof of Lemma \ref{statnl}. The main
difference arises 
since the local coordinate map also depends on
$\lambda\in\Lambda$. Replacing the function $G$ in the proof of Lemma
\ref{statnl} with the function (notation as in Section \ref{fcanF})
\begin{align}\notag
&H(x,\xi,x',\xi',\sigma,\sigma',\lambda)=\\\notag
&J[x,\xi,x',\xi',\sigma,\lambda](1)-x'-\xi'-\lambda\cdot
DY^\sigma_0(x)\cdot x' +
\pa_\sigma(\exp_{\psi_\lambda^\sigma(x)}^{\lambda,\sigma}
A_\sigma^\lambda \xi)\cdot\sigma',
\end{align}
where $J[x,\xi,x',\xi',\sigma,\lambda]$ is the Jacobi field along the
geodesic $\exp_{\psi_\lambda^\sigma(x)}^{\lambda,\sigma}A_\lambda^\sigma
\xi$ with initial conditions $J(0)=d\psi_\lambda^\sigma\cdot
x'+\pa_\sigma\psi_\lambda^\sigma\cdot\sigma'$,  
$J'(0)=A_\lambda^\sigma\xi'+dA_\lambda^\sigma\cdot x'\cdot
\xi+\pa_\sigma A^\sigma_\lambda\cdot\sigma'\cdot\xi$ and repeating the
argument given there proves the lemma.
\end{pf} 

In the self-tangency shortening case the estimate is somewhat changed
since we work in Sobolev spaces with negative exponential weights in the
gluing region. Here we have
\begin{lma}\label{stshnl}
There exists a constant $C$ such that
the non-linear term $N(v,\gamma,\lambda)$ of $\Gamma$ in a neighborhood
of $w_\rho$, where $w_\rho$ is as in Section \ref{10stpreshort.section} 
satisfies
\begin{align}\notag
\|N(u,\beta)-N(v,\kappa)\|_{1,-\epsilon,\rho}\le 
C&e^{\epsilon\rho}\Bigl(
\|u\|_{2,-\epsilon,\rho}+|\beta|+\|v\|_{2,-\epsilon,\rho}+|\gamma|
\Bigr)\\\notag
&\times\Bigl(\|u-v\|_{2,-\epsilon,\rho}+|\beta-\gamma|\Bigr) 
\end{align} 
\end{lma}
\begin{pf}
The proof is exactly the same as the proof of Lemma \ref{statnl}. We must
only take into account in what way the weights affect the
estimates. Starting with \eqref{0norm}, we see that the norm
$\|\cdot\|_{2,\rho,-\epsilon}$ does not control the $\sup$-norm
uniformly in $\rho$. But it does control $e^{-\epsilon\rho}$ times the
$\sup$-norm. Thus we conclude instead of \eqref{0norm} 
\begin{equation}
\|N(u)-N(v)\|\le
Ce^{\epsilon\rho}(\|u\|_{2,-\epsilon,\rho}+\|v\|_{2,-\epsilon,\rho})
\|u-v\|_{2,-\epsilon,\rho}. 
\end{equation}
Similarly, we loose this factor in the other estimates where we use
the $\sup$-norm. 
Let $e_\rho$ denote the weight function from Section 
\ref{10shortweight.section}.
When we use the $L^4$-estimate we have instead of
\eqref{L4ex} the following 
\begin{align}\notag
&\int_{D_m}(|Du|+|Dv|)^2|Du-Dv|^2e_\rho^2\,dA\le\\\notag
e^{2\epsilon\rho}&\int_{D_m}(|Du|+|Dv|)^2|Du-Dv|^2e_\rho^4\,dA\le\\\notag
e^{2\epsilon\rho}&\Bigl(\int_{D_m}(|Du|+|Dv|)^4e_\rho^4\Bigr)^{\frac12}
\Bigl(\int_{D_m}(|Du-Dv|)^4e_\rho^4\Bigr)^{\frac12}\le\\\notag
Ce^{2\epsilon\rho}&(\|u\|_{2,-\epsilon,\rho}+\|v\|_{2,-\epsilon,\rho})^2
(\|u-v\|_{2,-\epsilon,\rho})^2.  
\end{align}  
We conclude finally 
\begin{equation}\notag
\|N(u)-N(v)\|_{1,-\epsilon,\rho}\le
Ce^{\epsilon\rho}(\|u\|_{2,-\epsilon,\rho}+\|v\|_{2,-\epsilon,\rho})
\|u-v\|_{2,-\epsilon,\rho}. 
\end{equation}
The same argument as in Lemma \ref{statnl} then takes care of the conformal
structures and the lemma follows. 
\end{pf}

Finally, we consider self-tangency gluing, where we have a large weight 
function which does not interfere (destructively) with the $\sup$-norm
and the $L^4$ estimates. 
\begin{lma}\label{stglunl}
There exists a constant $C$ such that
the non-linear term $N(v,\gamma)$ of $\Gamma$ in a neighborhood
of $w_\rho$, where $w_\rho$ is as Section \ref{10stpreglue.section} satisfies
\begin{align}\notag
\|N(u,\beta)-N(v,\kappa)\|_1\le&
C\Bigl(\|u\|_{2,\epsilon,\rho}+|\beta|+
\|v\|_{2,\epsilon,\rho}+|\gamma|\Bigr)\\\notag
&\times\Bigl(\|u-v\|_{2,\epsilon,\rho}+|\beta-\gamma|\Bigr) 
\end{align} 
\end{lma}
\begin{pf}
See the proof of Lemma \ref{statnl}
\end{pf}


%

\section{Gromov compactness} \label{9compactness.section}

In this section we prove a version of the Gromov compactness theorem.
In Section \ref{9conformalcompactness.section}, we discuss 
the compactification of the space of conformal
structures which is done is \cite{fo}.
In Section \ref{9statement.section}, we 
translate the notions of convergence and (limiting) broken curves
from \cite{ms} to our setting.
There are two notions of convergence we must prove: a strong local
convergence and a weak global convergence.
In Sections \ref{9strongI.section} and \ref{9strongII.section}, 
we adopt Floer's original approach, \cite{Floer88a}, to prove
the strong local convergence.
In proving local convergence, we in fact prove that our holomorphic
disks, away from the punctures, are smooth up to and including the
boundaries, see Remark \ref{9bootstrap.rem}.
To prove the weak global convergence in Section \ref{9weak.section}, 
we analyze where the area (or energy)
of a sequence of disks accumulates, and construct an appropriate
sequence of reparameterizations of the domain to recover this area.

We note that although our holomorphic curves map to a non-compact
space, $\C^n$, the set of curves we consider
lives in a compact subset.
This follows because $\C^n$ is a symplectic manifold with
``finite geometry at infinity'': a holomorphic curve with
a non-compact image must contain infinite area.
And the area of any disk we consider is bounded above
by the action of the chords mapped to at its corners.
Thus, we can prove the Gromov compactness theorem in this
non-compact set-up.
For a review of finite geometry at infinity (also known as ``tame''), 
see \cite{al} Chapter 5, as well as \cite{eg,g,Sullivan}.

\subsection{Notation and conventions for this section}
\label{9notation.section}

Unlike in the other sections, we need to consider Sobolev spaces
with derivatives in $L^p$ for $p \ne 2$.
We define in the obvious way the spaces
$\W^{p, \rm{loc}}_k(\Delta_m, \C^n)$ to indicate
$\C^n$-valued functions on $\Delta_m$ whose first $k$
derivatives are locally $L^p$-integrable.
For this section only, we denote the corresponding norm
by $\| \cdot \|_{k,p}$.

In order to define broken curves in the next subsection, we will
need to extend the disk continuously to the boundary punctures.
Of course the extra Legendrian boundary condition, $h$,  does not
extend continuously.
For this reason, we will only extend $u$ to $\Disk_m$, the closure
of $\Delta_m$; thus, $u\colon \Disk_m\to \C^n.$ 
Note that $\|u \|_{p,k}$ might still blow up at these punctures.
We sometimes only consider $u$ and $u| \partial \Disk_m$ in which case we write
$u: ( \Disk_m, \partial \Disk_m) \rightarrow ( \C^n, \projL(L)).$
For $X\subset\Disk_m$, let $\| u \|_{k,p:X} = \| u | X \|_{k,p}$, and 
$\| u \|_{k,p:\epsilon}$ denote the norm restricted to some
disk (or half-disk) of radius $\epsilon$.

Because we sometimes change the number of boundary punctures,
we will denote by $D$ the unit disk in $\C^n$.

\subsection{Compactification of space of conformal structures}
\label{9conformalcompactness.section}

Recall $\conf_m$ is the space of conformal structures (modulo conformal
reparameterizations) on the unit disk in $\C$ with $m$ boundary
punctures.

When $m \ge 3$, we define a {\it stable cusp disk representative} with
$m$ {\it marked boundary points}, 
$(\Sigma; p_1, \ldots, p_m)$, to be a connected,
simply-connected union of unit disks in $\C$ where pairs of disks may overlap
at isolated boundary points (which we call {\it double points} of $\Sigma$)
and each disk in $\Sigma$ has at least 3 points, called {\it
marked points}, which correspond either to double points or 
the original boundary marked points.
When $m =1$ or 2, the stable cusp disk representative shall be a single disk.
Two stable cusp disk representatives are {\it equivalent} if there exist a conformal
reparameterization of the disks taking one set of marked points to the other.
We define a {\it stable (cusp) disk} with $m$ marked points
to be an equivalence class of stable disk representatives with $m$ 
marked points.

In Section 10 of \cite{fo}, Fukaya and Oh prove that $\bar{\conf}_m$, the
compactification of $\conf_m$, is the space of stable disks
with $m$ marked points.

\subsection{The statement}
\label{9statement.section}

A {\it broken curve} 
$(u,h) = ((u^1,h^1), \ldots (u^N,h^N))$ is a connected
union of holomorphic disks, 
$(u^j,h^j),$ (recall $u^j$ is extended to $\Disk_{m_j}$)
where each $u^j$ has exactly one positive puncture and except for one disk, say $u^1,$ 
the positive puncture of $u^j$ agrees with the negative puncture of some other $u^{j'}.$ 
One may easily check that 
a broken curve can be parameterized by a single smooth 
$v : (D_m, \partial D) \rightarrow (\C^n, \projL(L)),$ such that
$v^{-1}$ is finite except 
at points where two punctures were identified, here $v^{-1}$ is an arc
in $\Delta_m.$  
%

\begin{dfn}
\label{9conv.dfn} \rm
A sequence of holomorphic disks $(u_\alpha, h_\alpha)$ {\em converges}
to a broken curve $(u,h) = ((u^1, h^1), \ldots, (u^N, h^N))$ if the
following holds

\begin{enumerate}

\item
({\em{Strong local convergence}}) For every $j \le N$, there exists
a sequence $\phi_\alpha^j:D \rar D$ of linear fractional transformations
and a finite set $X^j \subset D$ such that $u_\alpha \circ \phi^j_\alpha$
converges to $u^j$ uniformly with all derivatives on compact subsets
of $D \setminus X^j$

\item
({\em{Weak global convergence}}) There exists a sequence of 
orientation-preserving diffeomorphisms $f_\alpha:D \rar D$ such that
$u_\alpha \circ f_\alpha$ converges in the $C^0$-topology to a parameterization
of $u$.

\end{enumerate}

\end{dfn}

Henceforth, to simplify notation when passing to a subsequence, we will not
change the indexing.

\begin{thm}
\label{9cpt.thm}

Assume $(u_\alpha, h_\alpha) \in \M(a; b_1, \ldots, b_m)$
is a sequence of holomorphic disks
with $L_\alpha$ Legendrian boundary condition.
Let $\kappa_\alpha \in \conf_{m+1}$ denote the conformal
structure on the domain of $u_\alpha.$
Assume $L_\alpha$ converges to an embedded Legendrian $L$ in the
$C^\infty$-topology.
Then there exists a subsequence $(u_\alpha, h_\alpha, \kappa_\alpha)$ 
such that $\kappa_\alpha$ converges to $\kappa \in \bar{\conf}_{m+1}$ and 
$(u_\alpha, h_\alpha)$ converges to a broken curve $(u,h)$
whose domain is a stable disk representative of $\kappa.$

\end{thm}

Note that using the strong local convergence property {\it a posteriori},
this compactness result proves that all derivatives of a holomorphic
disk $(u,h)$ are locally integrable away from the finite set of points.
In particular, such disks are smooth at the boundary away from these
points. See Remark \ref{9bootstrap.rem}.

We also remark that, the appropriately modified, Theorem \ref{9cpt.thm} holds if the disks have more
than one positive puncture. 


\subsection{Area of a disk}
\label{9area.section}
For holomorphic $u:(D, \partial D) \rar (\C^n, \projL(L))$, 
recall that $\area(u) = \int u^*\omega$, where $\omega=\sum_i
dx_i\wedge dy_i$, denotes its (signed) area.

\begin{lma}
\label{9hbar.lma}
Consider an admissible Legendrian isotopy
parameterized by $\lambda \in \Lambda$. 
We assume $\Lambda \subset \R$ is compact.
Denote by $L_\lambda$ the moving Legendrian submanifold.
There exists a positive upper
semi-continuous function $\hbar: \Lambda \rar \R^+$
such that for any non-constant holomorphic map
$u:(D, \partial D) \rar (\C^n, \projL(L_\lambda))$,
$\area(u) \ge \hbar (\lambda)$.
\end{lma}

\begin{proof}
We need the following statement from Proposition 4.3.1 (ii) of Sikorav
in \cite{al}: There are constants $r_1, k$ (depending 
only on $\C^n$) such that if $r\in(0,r_1]$ and $u\colon \Sigma\to
B(x,r)$ is a holomorphic map of a Riemann surface containing $x$ in
its image and with $u(\pa\Sigma)\subset\pa B(x,r)$ then  
$\text{Area}(\Sigma)\geq k r^2.$

Since $u$ is non-constant, Stokes Theorem implies $u$ must
have boundary punctures.
Choose $r >0$, an upper semi-continuous function of $\lambda$, such that:

\begin{itemize}
\item
for all Reeb chords $c,$ 
$\projL(L_\lambda) \cap B(c^*,r)$ is real analytic and diffeomorphic
either to $\R^n \times \{0\} \cup \{0\} \times \R^n$ or 
the local picture of the singular moment in a standard self-tangency
move (see Definition~\ref{5admissile_param.dfn}).
\item
for all distinct Reeb chords $c_1, c_2,$ 
$B(c_1^*,r) \cap B(c_2^*,r) = \emptyset$ and
\item $r<r_1.$
\end{itemize}

Let $\theta_\lambda$ be the smallest angle among all the complex angles associated to all the tranverse double points
of $\Pi_{\C}(L_\lambda).$ Now set
\begin{equation*}
\hbar(\lambda) = \min \left\{ \min_{c \in \mathcal{C}(L_\lambda)} \action(c), \frac{k r^2 \cos^2 \theta_\lambda}{8}
\right\} > 0.
\end{equation*}

Suppose $u$ maps all $n$ of 
its punctures to the same double point $c^*,$ 
then by
(\ref{1action.eqn}) 
\begin{equation*}
\area(u) \ge  \action(c) \ge \hbar.
\end{equation*}
(Note the number of positive punctures of $u$ must be larger than the number of negative ones since $u$ is not constant.)

Otherwise, assume $u$ maps boundary punctures to 
at least two distinct double points $c_1^*, c_2^*$ where $c_1^*$ is a non-degenerate double  point.
Then $c_2^* \notin \bar{B}(c_1^*, r)$ implies that there exists a point $x\in u(D)\cap \Pi_{\C}(L)\cap \partial B(c_1^*,\frac{r}{2}).$
Moreover, $B(x,\frac{r\cos\theta_\lambda}{2})\subset B(x,r)$ intersects $\Pi_{\C}(L)$ in only one sheet. 
Using the real-analyticity of the boundary, we double $u(D)\cap B(x,\frac{r\cos\theta_\lambda}{2})$ and apply
the proposition of Sikorav to conclude
\[\area(u)\geq \area (u(D)\cap B(x,\frac{r\cos\theta_\lambda}{2}))\geq \frac {k r^2 \cos^2 \theta_\lambda}{8}\geq \hbar.\]

\end{proof}

We introduce one more area-related notion, again borrowed from \cite{ms}.
Given a sequence of holomorphic maps $u_\alpha$ 
we say $z \in D$ is a {\it{point mass}} of $\{u_\alpha\}$ 
with {\it{mass}} $m$ if there exists a sequence $z_\alpha \in D$ 
converging to $z \in D$ such that
\begin{equation*}
\lim _{\epsilon \rar 0}
\lim_{\alpha \rar \infty} \area \left( u \left| 
B_\epsilon (z_\alpha) \cap D \right. \right) = m.
\end{equation*}

\subsection{Strong local convergence I: bootstrapping}
\label{9strongI.section}

In this subsection we prove the following ``bootstrap''
elliptic estimate: if we know a holomorphic curve lies locally
in $W^p_k$ with $p > 2, k \ge 1$, then the $\| \cdot\|_{p',k}$-(local)-norm
controls the $\| \cdot \|_{p', k+1}$-(local)-norm for $p' \in [2, p)$.


This proof first appeared as Lemma 2.3 in \cite{Floer88a} and later
corrected as Proposition 3.1 \cite{oh}.
Floer and Oh both prove the $k = 1$ case and state the general case.
Although there are no new techniques here, we present for the
reader the general case in more detail.

Let $A \subset \C$ denote the open disk or half-disk with boundary on the real
line.
Let $\W^p_k(A, \C^n)$ denote the closure, under the 
$\| \cdot \|_{k,p}$-norm, of the set of all smooth compactly supported
functions from $A$ to $\C^n$.
\begin{lma}
\label{9hormander.lma}
For every $l > k$, $ l - 2/q > k - 2/p$, there exists a constant $C$
such that if $\xi \in \W^p_k(A, \C^n)$ is compactly supported, 
$\xi| \partial A \subset \R^n$, and $\dbar \xi \in \W^q_{l-1}(A, \C^n)$
then
\begin{equation}
\label{9hormander.eq}
\| \xi \|_{l,q} \le C \| \dbar \xi \|_{l-1,q}.
\end{equation}
\end{lma}

This is stated as Lemma 2.2 of \cite{Floer88a} and Lemma 3.2 of \cite{oh}.
Floer attributes this result to Theorem 20.1.2 of \cite{ho}.
However, we were unable to deduce Lemma \ref{9hormander.lma}
for $k > 1$ from H{\"o}rmander's theorem. 
Alternatively, one can use the Seeley extension theorem (see \cite{me},
section 1.4 for example)
to extend the map to the full disk (in the case of the
half disk) and then use the well-known full disk version of Lemma
\ref{9hormander.lma}. 

We can now state and prove this subsection's main theorem.

\begin{thm}
\label{9bootstrap.thm}
Fix $k \ge 1$ and (not necessarily small) 
$\delta_{k-1} > \delta_{k} \ge 0$.
For any compact $K \subset A$, there exists a ``constant'' 
$C_1 = C_1(\|u \|_{k, 2+\delta_{k-1}})$ depending
continuously on $\| u\|_{k, 2+\delta_{k-1}}$ such that for all
holomorphic maps $u \in \W^{2+ \delta_{k-1}}_k(A, \C^n)$ with
$u(\partial A) \subset \projL(L)$, we have
\begin{equation}
\label{9bootstrap.eq}
\| u \|_{k+1, 2 +\delta_k: K} \le C_1 \| u\|_{k, 2+ \delta_k: A}.
\end{equation}
Moreover, if $u_\alpha$ is a sequence of holomorphic maps in 
$\W^{2+\delta_{k-1}}_k(A, \C^n)$ such that
$u_\alpha(\partial A) \subset \projL(L)$ and
$\|u_\alpha \|_{k, 2+\delta_{k-1}}$ 
is uniformly bounded,
then there exists a subsequence $u_\alpha$ converging in
$\W^{2+\delta_{k}}_k(K, \C^n)$ to some holomorphic map
$u:K \rar \C^n$.
\end{thm}

\begin{rem} \rm
\label{9bootstrap.rem}
Note how we can use the Sobolev embedding theorem to conclude that
all derivatives of the curve lie in $L^2$ locally, assuming we have
a finite local $\| \cdot \|_{1, 2 + \delta_0}$ norm to begin with.
In particular, a holomorphic disk $(h,u)$ with boundary
punctures becomes smooth at the boundary away from the punctures.
We did not have to assume this smoothness {\it{a priori}}.
\end{rem}

\begin{proof}
We shall only prove the first statement.
The second one easily follows from the first and the Sobolev embedding theorem.

Our goal is to prove (\ref{9bootstrap.eq}) in some small $\epsilon$
ball in $K$. 
The claim will then follow from the compactness of $K$.
 
Because $\W^{2+\delta_{k-1}}_{k}$ compactly sits in $C^0$ for $k \ge 1$,
we can choose small $\epsilon$ (continuously in 
$\|u \|_{k, 2+\delta_{k-1}})$ 
and $\epsilon'$ such that given $z_0 \in K$,
$u(B(z_0,\epsilon)) \subset B(u(z_0,\epsilon'))$.
We fix $\epsilon$ and $\epsilon'$ at the end of the proof.
Assume $\ep$ is small enough such that $B(z_0,\ep)$ does not
contain any boundary punctures of $u$.

Choose a diffeomorphism $\phi$ of $\C^n$ so that
$\phi(\projL(L)) \cap B_{\epsilon'}( \phi \circ u(z_0))$
corresponds to a piece of $\R^n \subset \C^n$ if $z_0 \in \partial A$.
Assume $z_0 = 0 \in \partial A$ as we will not consider the
easier interior estimate.
Locally near $0$ define $v = \phi \circ u;$ thus,
$v$ has $\R^n$ boundary conditions and
$\dbar_{\phi^*i} v = 0$ where $\dbar_{\phi^*i}$ uses the pull-back
(almost) complex structure.

Choose a compactly supported smooth function $\gamma: A \rar \R$ such that
$\gamma(z) = 1$ for $|z| \le \frac{1}{2}$ and set 
$\gamma_\epsilon(z) = \gamma(z/\epsilon)$.
Note that we can choose $\gamma$ such that the $C^k$-norm 
$\|\gamma_\epsilon\|_{C^k}$ is of order $\frac{1}{\epsilon^k}$.
By Lemma \ref{9hormander.lma}, there exists $C_2$ such that for any
$\epsilon$, 

\begin{eqnarray}
\label{9boot1.eq}
\| \gamma_\ep v \|_{k+1,2+\delta_k} & \le &
C_2 \| \dbar( \gamma_\ep v) \|_{k,2+\delta_k} \\
\notag & \le &
C_2 ( \| (\dbar \gamma_\ep) v \|_{k,2+\delta_k} + \|\gamma_\ep
\dbar v \|_{k,2+\delta_k} ) \\
\notag & \le &
C_2\left( \| \gamma_\ep\|_{C^k}  \| v \|_{k, 2+\delta_k: \ep}
+ \left\| \gamma_\ep( \phi^* i- i) \frac{\partial v}{\partial y}
\right\|_{k,2+\delta_k}\right)
\end{eqnarray}
where $x+iy$ is the complex coordinate in $A$.

To simplify notation, let $J(v) = \phi^*i - i$.
We can assume that all derivatives of $J(\cdot)$ are uniformly
bounded on $B(\phi\circ u(z_0),\ep')$.
We consider the last term of (\ref{9boot1.eq}):
\begin{eqnarray}
\label{9boot2.eq}
\left\| \gamma_\ep J(v) \frac{\pa v}{\pa y} \right\|_{k,2+\delta_k} & \le & 
C_3 \sum_{\{i,j\,\,|\,\,i+j = k\}} 
\left\| D^i (J (v)) D^j \left(\gamma_\ep
\frac{\pa v}{\pa y} \right)\right\|_{0, 2+\delta_k} \\
\notag & \le &
C_4 C_3 \| J(v) \|_{0, \infty: \epsilon} 
\left\| D^k \left( \gamma_\ep \frac{\pa v}{\pa y}
\right) \right\|_{0, 2+\delta_k} + \\
\notag
&& 
C_3 \sum_{\{i,j\,\,|\,\,i+j = k, j < k\}} 
\left\|
D^i (J (v)) \right\|_{0,2+\delta_{k-1}: \ep} 
\left\|D^j \left(\gamma_\ep \frac{\pa v}{\pa y} \right)\right\|_{0, p: \ep}
\end{eqnarray}
where $\frac{1}{p} = \frac{1}{2+ \delta_{k}} - \frac{1}{2+\delta_{k-1}}$.
(Here we use $\delta_{k-1} > \delta_k.$)

Choose $\ep$ small enough such that $|J(v)| \le \frac{1}{3C_4C_3 C_2};$
thus, combining (\ref{9boot1.eq}) and (\ref{9boot2.eq}), we get
\begin{eqnarray}
\label{9boot3.eq}
\frac{2}{3} \left\| 
\gamma_\ep v \right\|_{k+1, 2 +\delta_k} & \le &
C_2 \| \gamma_\ep\|_{C^k}  \| v \|_{k, 2+\delta_k: \ep} + \\
\notag &&
C_2 C_3 \sum_{\{i,j\,\,|\,\,i+j = k, j < k\}} 
\| D^i (J (v)) \|_{0,2+\delta_{k-1}: \ep} 
\left\|D^j \left(\gamma_\ep \frac{\pa v}{\pa y} \right) \right\|_{0, p}.
\end{eqnarray}

Since, $i \cdot (2+\delta_{k-1}) > 2$ for $1 \le i \le k $, we get 
(see \cite{ms} Proposition B.1.7)
\begin{equation}
\label{9boot4.eq}
\left \| D^i ( J(v)) \right \|_{0, 2 + \delta_{k-1}: \ep}
\le
\left\| J(v) \right\|_{i, 2+\delta_{k-1}: \epsilon}
\le
C_5 \left( \|J\|_{C^i} + 1 \right)
\left\| v \right\|_{i, 2+\delta_{k-1}: \epsilon}
\le
C_6,
\end{equation}
where $C_6 = C_6
\left(\left\| v \right\|_{k, 2+\delta_{k-1}: \epsilon}\right).$


Fix any $s >p$. 
Then for any $\kappa >0$, if we set
\begin{equation}
\notag
\mu = \frac{ \frac{1}{2 + \de_k} - \frac{1}{p}}{\frac{1}{p} - \frac{1}{s}}
\end{equation}
we get
\begin{eqnarray}
\label{9boot5.eq}
\left\|D^j \left(\gamma_\ep \frac{\pa v}{\pa y} \right) \right\|_{0, p}
& \le &
\kappa 
\left\|D^j \left(\gamma_\ep \frac{\pa v}{\pa y} \right) \right\|_{0, s}
+\kappa^{-\mu}
\left\|D^j \left(\gamma_\ep \frac{\pa v}{\pa y} \right) \right\|_{0, 
2+ \delta_k} \\
\notag & \le &
\kappa C_7
\left\|D^j \left(\gamma_\ep \frac{\pa v}{\pa y} \right) \right\|_{1, 2+\de_k}
+\kappa^{-\mu}
\left\|D^j \left(\gamma_\ep \frac{\pa v}{\pa y} \right) \right\|_{0, 
2+ \delta_k} \\
\notag & \le &
\kappa C_8
\left\|\gamma_\ep v \right\|_{j+2, 2+\de_k}
+\kappa^{-\mu}
\left\|D^j \left(\gamma_\ep \frac{\pa v}{\pa y} \right) \right\|_{0, 
2+ \delta_k}.
\end{eqnarray}
The first inequality uses an interpolation result and the second the embedding
theorem.

Choose $\kappa$ such that
\begin{eqnarray}
\label{9boot6.eq}
C_2 C_3 C_6 \kappa C_8
& 
\le & \frac{1}{3}.
\end{eqnarray}

Using (\ref{9boot4.eq}) to (\ref{9boot6.eq}) we bound the first term
of the sum on the left hand side of (\ref{9boot3.eq}):
\begin{eqnarray}
\label{9boot7.eq}
\lefteqn{
C_2 C_3 \left\| D^1 J(v) \right\|_{0,2+\de_{k-1}: \ep}
\left\|D^{k-1} \left(\gamma_\ep \frac{\pa v}{\pa y} \right) \right\|_{0, p}
} &&\\
\notag
& \le &
\frac{1}{3} \left\| \gamma_\ep v \right\|_{k+1, 2+\delta_k}
+ C_2 C_3 C_6
\cdot \kappa^{-\mu} \|v\|_{k, 2+\delta_k : \ep}
\end{eqnarray}

Combining (\ref{9boot3.eq}) with (\ref{9boot7.eq}),
there exists
$C_9 = C_9 \left(\|v\|_{k, 2+ \delta_{k-1} : \ep}\right)$
and 
$C_{10} = C_{10} \left(\|v\|_{k, 2+ \delta_{k-1} : \ep}\right)$
such that
\begin{eqnarray*}
\frac{1}{3} \left\| 
\gamma_\ep v\right\|_{k+1, 2 +\delta_k} & \le &
\left( C_2 \| \gamma_\ep\|_{C^k} + C_9 \right)
\| v \|_{k, 2+\delta_k: \ep} + \\
\notag &&
C_2 C_3  \sum_{\{i,j\,\,|\,\,i+j = k, j < k-1\}} 
C_6
\left\|D^j \left(\gamma_\ep \frac{\pa v}{\pa y} \right) \right\|_{0, p} \\
\notag
& \le &
C_{10} \| v \|_{k, 2+\delta_k: \ep}.
\end{eqnarray*}
This last line follows since $\|v\|_{k, 2+ \delta_{k}: \ep}$ controls
$\|v\|_{j+1, p: \ep}$ for $j < k-1$.

\end{proof}

\subsection{Strong local convergence II: uniformly bounding higher Sobolev 
norms}
\label{9strongII.section}

In order to apply Theorem \ref{9bootstrap.thm}, we need a uniform bound on
the $\| \cdot \|_{k, 2+ \delta}$-norm where $\delta >0$ might be large.
Our holomorphic disk only come with a bound on the $\| \cdot \|_{1,2}$-norm
in terms of the action.
In this subsection we indicate how the latter norm controls the former.

\begin{thm}
\label{9unibound.thm}
Consider the sequence of holomorphic disks 
$(u_\alpha, h_\alpha) \in \M(a; b_1, \ldots, b_m)$.
There exists a finite number of points 
$z_1, \ldots, z_l \in \partial \Delta_m$ 
and a ``constant'' $C_{11} = C_{11}(K, p,k)$
such that for any
positive integer $k$, for any $p \in \R$ with $k > \frac{2}{p}$,
and for any compact set 
$K \subset \Delta_m \setminus \{ z_1, \ldots, z_l\}$,
\begin{equation}
\notag
\| D^k u_\alpha \|_{0, p:K} \le C_{11}.
\end{equation}
\end{thm}

This result in the special case when $m=2$ and $k=1$
was proved as Theorem 2 in \cite{Floer88a} and then reproved as Proposition 3.3
in \cite{oh}.
Thus, the first part of this proof uses some ideas of \cite{Floer88a}.
For the cross-referencing inclined reader, we borrow the notation
from \cite{oh}.

\begin{proof}

Let $a = \area(u_\alpha)$ which is independent of $\alpha.$
Let
\begin{eqnarray*}
\eta_\alpha(K) & = & \inf\{ \eta >0  :  
\mbox{ there exists } z \in K \mbox{ such that } \\
&& \quad
\| D^k u_\alpha \|_{0,p:B(z,\eta)} \ge \eta ^{\frac{2}{p} - k} \}.
\end{eqnarray*}
Assume $\eta_\alpha(K) \rightarrow 0$ for some subsequence, 
otherwise our theorem holds.
Choose $z_\alpha \in K$ such that
\[
\left\|D^k u_\alpha \right\|^p_{0,p:B(z_\alpha, \eta_\alpha)} \ge
\frac{1}{2} \eta_\alpha^{2 - pk}.
\]
Define
$r_\alpha = \eta^{-1}_\alpha \mbox{dist}(z_\alpha, \partial \Delta_m)$.
There are two cases to consider. We must derive a contradiction
for both.

\medskip

\noindent
Case 1: $r_\alpha \rightarrow \infty.$

\smallskip

Pass to a subsequence and assume $z_\alpha$ converges to some 
$z_0 \in \Disk_m$.

Let $A \subset \R^2$ be a compact subset and use variables $(s,t)$ on $\R^2.$
Let $f: A \rightarrow \R^{2n}.$
We remark that when rescaling 
variables $(s,t) \mapsto (\beta s, \beta t)$, the $L^p$-norm
changes like:
\begin{equation}
\label{9scale.eq}
\left\|D^k f \right\|_p \rightarrow 
\left((\beta^{-k})^p (\beta^2) \right)^{\frac{1}{p}}
\left\| D^k f \right\|_p = \beta^{-k + \frac{2}{p}} \left \|D^k f \right\|_p.
\end{equation}
This remark allows us to use the Floer technique for our more general case:
$(k,p)$, $k > \frac{2}{p}$ versus $(1,p)$, $p>2$.

Since $r_\alpha \rightarrow \infty$,
for every $R >0$ and large enough $\alpha$ we have holomorphic maps
\begin{equation*}
v_\alpha: \C \supset B(0,R) \rightarrow \C^n,
\quad
v_\alpha(z) = u_\alpha(\eta_\alpha(z - z_\alpha)),
\end{equation*}
satisfying
\begin{eqnarray}
\label{9uni1.eq}
\left \| Dv_\alpha \right\|_{0,2} & \le & a \\
\label{9uni2.eq}
\left \| D^kv_\alpha \right\|_{0,p:B(0,1)} & \ge & \frac{1}{2} \\
\label{9uni3.eq}
\left \| D^kv_\alpha \right\|_{0,p:B(z,1)} & \le & 1, \text{ for all } z
\in B(0,R-1).
\end{eqnarray}
Equations (\ref{9uni1.eq}) to (\ref{9uni3.eq}) follow from (\ref{9scale.eq})
and the exponent in the definition of $\eta_\alpha$.
We use (\ref{9uni3.eq}) and
the last statement of Theorem \ref{9bootstrap.thm} to find a subsequence, 
$v_\alpha$ converging in $\W^p_k(B(0,R), \C^n)$
to some holomorphic map $w_R: B(0,R) \rightarrow \C^n$.
(Technically, the convergence is in $\W^{p'}_k$ for some $p' < p$;
however, since we still have $p' > 2$ we will ignore this.) 
Repeating this procedure for all positive integers $R$ and choosing converging 
subsequences, we obtain a holomorphic map $w: \C \rightarrow \C^n$ satisfying
\begin{eqnarray}
\label{9uni4.eq}
\left \| Dw \right\|_{0,2} & \le & a \\
\label{9uni5.eq}
\left \| D^k w \right\|_{0,p:B(0,1)} & \ge & \frac{1}{2}.
\end{eqnarray}
By (\ref{9uni5.eq}), $w$ is non-constant; hence,
$\left \| Dw \right\|_{0,2} = a' >0.$

To derive a contradiction, consider the sequence of annuli
\[
\C \supset A_\rho := \{ re^{i \theta} \,\, : \,\, \rho  \le r < \rho +1 \}.
\]
(\ref{9uni4.eq}) implies that 
$\left \| Dw \right\|_{0,2: A_\rho} \rightarrow 0$.
Thus, for some circle $C_{\rho'} = \{\rho' e^{i \theta} \}$
where $\rho' \in (\rho, \rho+1)$, we have
$\left \| Dw \right\|_{0,2: C_{\rho'}} \rightarrow 0$.
By Sobolev's Theorem, this implies
$\|w\|_{C^0: C_\rho'} \rightarrow 0.$

Thus, for sufficiently large $\rho$, $w|C_{\rho'}$ spans a small disk
$D_{\rho'} \subset \C^n$ 
whose (absolute) area is bounded by $\frac{a'}{3}$.
Let $B_{\rho'} \subset \C$ be the disk spanned by $C_{\rho'}$.
Choose $\rho$ large enough such that
$\|Dw \|_{0,2; B_{\rho'}} > \frac{2a'}{3}.$
Let $\tilde{w}: S^2 \rightarrow \C^n$ be the (not necessarily holomorphic)
map $w|B_{\rho'}$ capped off with $D_{\rho'}$.
Then, with $\omega=\sum_i dx_i\wedge dy_i$ we find, since $\pi_2(\C^n) = 0$, 
\[
0 = \int_{S^2} \tilde{w}^* (\omega)
\ge \|Dw\|_{0,2: B_{\rho'}} - 
\int_{D_{\rho'}} \omega
\ge \frac{a'}{3} >0
\]

\medskip

\noindent
Case 2: $r_\alpha \rightarrow r < \infty.$

\smallskip

In this case, $z_\alpha \rightarrow z_0 \in \partial \Delta_m$.
We proceed as before to construct a limiting holomorphic map $w$, where
this time
\[
w:\C_r = \{ z \in \C \,\, : \,\, \mbox{Im}(z) \ge -r \} \rightarrow \C^n
\]
with $w(\partial \C_r) \subset \projL(L)$ and satisfying
(\ref{9uni4.eq}) and (\ref{9uni5.eq}).
Here $B_1(0)$ is the unit ball in $\C_r$.

At this point our proof deviates from \cite{Floer88a} and \cite{oh}
because of multiple boundary punctures.

Repeat the second part of the discussion of Case 1, where 
$A_\rho \subset \C_r$
is a ``partial'' annulus and $C_{\rho'} \subset A_\rho$ a ``partial'' circle.
Instead of constructing a smooth (but not necessarily holomorphic)
map $\tilde{w}$ as before, we use the convergence
$\|w\|_{C^0: C_\rho'} \rightarrow 0$ to conclude that
(after precomposing with a linear fractional transformation from
the unit disk $D \subset \C$ to $\C_r$ which takes $-1$ to $\infty$)
\[
w: (D \setminus \{-1\}, \partial D \setminus \{-1\})
\rightarrow (\C^n, \projL(L))
\]
can be continuously extended to $-1$.

By Lemma \ref{9hbar.lma},
$\left\| w \right\|_{0,2} \ge \hbar.$

From the construction of $w$ and $v_\alpha$, it is easy to see that for large
enough $\alpha$,
\[
\area \left(u_\alpha \left| \Disk_m \setminus B(z_\alpha, \eta_\alpha) 
\right. \right)
\le a - \frac{\hbar}{2}.
\]

Note that $z_0$ is an example of a point mass for the sequence
$u_\alpha$.
We will relabel it $z_1$ to coincide with the statement of the theorem.

Now suppose there is a subsequence 
$z_\alpha \rightarrow z_2 \ne z_1$ such that
\[
\left\|D^k u_{\alpha} \right\|_{0,p:B(z_{\alpha}, \eta_{\alpha})} \ge
\frac{1}{2} \eta_{\alpha}^{2 - pk}.
\]
Repeating Case 1, we again conclude that $z_2 \in \partial \Delta_m$.
Furthermore,
\[
\area \left(u_\alpha \left| \Disk_m \setminus 
\left(
B(z_\alpha, \eta_\alpha)  \cup B(z_{\alpha, \eta_\alpha}) \right)\right.
\right)
\le a - 2\frac{\hbar}{2}.
\]
Assume $\alpha$ is large enough such that
$B(z_\alpha, \eta_\alpha)  \cup B(z_{\alpha, \eta_\alpha})
= \emptyset.$
Since $a$ is finite, this can only happen a finite number of times:
$z_1, \ldots z_l$.
So as long as $K \subset \Delta_m \setminus \{z_1, \ldots
z_l\}$, the bound of the theorem holds.

\end{proof}

\subsection{Recovering the bubbles}
\label{9weak.section}

The goal of this subsection is to construct a (not
necessarily conformal) reparameterization
of $\Disk_m$ which recovers all disks which bubble off.
This reparameterization implies the second convergence in
Definition \ref{9conv.dfn}.

Consider a sequence $(u_\alpha, h_\alpha)$ which converges strongly
on any compact $K \subset \Delta_m \setminus \{z_1, \ldots, z_l\}$.
By the proof of Theorem \ref{9unibound.thm},
we can assume that $z_1$ is a point mass with mass $m_1>0$.

Let $\C_+ \subset \C$ denote the upper-half plane.
Let $B_r = \{z \in \C_+ : \| z\| < r\}$
and $C_r = \partial B_r$.
Define the conformal map
\[
\psi_\alpha: \C_+ \rightarrow \Disk_m, \quad
\psi_\alpha(z) = \frac{-z + i R_\alpha^2}{z + i R_\alpha^2} \cdot z_1
\]
where $R_\alpha \in \R$ is such that
\[
\area \left( u_\alpha 
\left| \psi_\alpha \left(B_{R_\alpha}\right) \right. \right)
 = m_1.
\]

Pass to a subsequence and assume $\alpha < \alpha'$ implies
$R_\alpha < R_{\alpha'}$, which can be done since
by the definition of point mass, $\lim_{\alpha \rightarrow \infty}
R_\alpha = \infty.$
Note that 
\begin{equation}
\label{9global1.eq}
\lim_{\alpha \rightarrow \infty} \psi_\alpha \left( B_{R_\alpha} \right)
 = \lim_{\alpha \rightarrow \infty}
\psi_\alpha \left(B_{R_\alpha^{{3/2}}} \right) = z_1.
\end{equation}
Assume $\alpha$ is large enough so that 
$\psi_\alpha \left(B_{R^{{3/2}}_\alpha}\right)$ contains
no other point masses of the sequence $u_\alpha$.
However, $\psi_\alpha \left(B_{R^{{3/2}}_\alpha}\right)$ might contain
boundary punctures.

After passing to a subsequence, we can use Theorems \ref{9bootstrap.thm}
and \ref{9unibound.thm} to
assume that 
$u_\alpha$ converges to some $u$ on any compact
set in 
$\Delta_m \setminus\left( \{z_2, \ldots, z_l\} \cup
\psi_\alpha \left(B_{R^{3/2}_\alpha} \right) \right).$

The definition of $R_\alpha$ and (\ref{9global1.eq}) imply
\begin{equation}
\label{9global2.eq}
 \lim_{\alpha \rightarrow \infty}
\area \left( u_\alpha \left| 
\psi_\alpha \left(B_{R^{{3/2}}_\alpha} \setminus B_{R_\alpha} \right)
\right. \right)  = 0.
\end{equation}

Use (\ref{9global2.eq}) and argue as
in the previous subsection to find some half circle 
$C_{R'_\alpha} \subset \C_+$,
with $R'_\alpha \in ( R^{{3/2}}_\alpha-1, R^{{3/2}}_\alpha]$
such that
\[
\|u_\alpha \circ \psi_\alpha\|_{C^0: C_{R'_\alpha}} \rightarrow 0.
\]

Define the center of mass of $u_\alpha \circ \psi_\alpha$ to be
\begin{equation*}
z_\alpha = x_\alpha + i y_\alpha
 = \frac{1}{m_1} \int_{B_{R_\alpha}} \left| D(u_\alpha \circ \psi_\alpha)
\right|^2 (x+iy) \,dx\wedge dy
\, \in \, B_{R_\alpha},
\end{equation*}
where $x+iy$ are coordinates on $\C_+$.
Define the conformal map $\phi_\alpha$ which sends $i$ to $z_\alpha$:
\begin{equation*}
\phi_\alpha: \C_+ \rightarrow \C_+, \quad
\phi_\alpha(z) = y_\alpha z + x_\alpha.
\end{equation*}
Note that although
$\phi^{-1}_\alpha \left( C_{R_\alpha} \right)$ might remain bounded,
$\phi^{-1}_\alpha \left( C_{R'_\alpha} \right)$
converges to $\infty$ because
\[
\left| \phi^{-1}_\alpha \left(R'_\alpha e^{i \theta} \right) \right|
\ge
\frac{(R_\alpha^{{3/2}} - 1)}{\left| y_\alpha \right|}
\max \{|\cos \theta|, |\sin \theta| \}
\]
and $y_\alpha < R_\alpha.$

Define the conformal map
\begin{equation*}
\Psi: D \rightarrow \C_+, \quad
\Psi(z) = \frac{z-1}{iz+i}
\end{equation*}
where $D \subset \C$ is the unit disk.
Note that $\Psi^{-1} \phi^{-1}_\alpha \left( C_{R'_\alpha} \right)
\rightarrow -1$ and that 
$u_\alpha \circ \psi_\alpha \circ \phi_\alpha \circ \Psi$ all have
center of mass at $0 \in D$.
(Recall that the center of mass uses the Euclidean metric on $\C_+$,
not on $D$.)

Since
\[
\left\|u_\alpha \circ \psi_\alpha \circ \phi_\alpha \circ \Psi
\right\|_{C^0: \Psi^{-1} \circ \phi^{-1}_\alpha \left(C_{R'_\alpha}\right)} 
\rightarrow 0,
\]
pass to a subsequence as before and conclude that
$u_\alpha \circ \psi_\alpha \circ \phi_\alpha \circ \Psi$ 
converges to some holomorphic $w$ on compact sets outside of some
boundary point masses and punctures, as well as $-1$
(since 
$u_\alpha \circ \psi_\alpha \circ \phi_\alpha \circ \Psi$ 
is not defined at $-1$).

As before, $w$ can be continuously extended to $-1$.
We claim that under this reparameterization, $-1$ is not a point mass
of
$u_\alpha \circ \psi_\alpha \circ \phi_\alpha \circ \Psi.$
Otherwise, in the $\C_+$ set-up, as some mass escaped to $\infty$, the
center of mass would have to go to $\infty$ as well, contradicting
the fact that it is fixed at $i \in \C_+$.

Because $u_\alpha$ converges to $u$ outside of $\psi_\alpha \left(
B_{R^{3/2}_\alpha } \right)$, and because no area
is ``unaccounted'' for by (\ref{9global2.eq}), we can continuously
extend $u$ to $z_1$ so that
$u(z_1) = w(-1).$ 
Considering how $u$ and $w$ were obtained from 
$u_\alpha$ it is easy to see that the sign of the punctures ($z_1$ for $u$ and $-1$ for $w$)
will be opposite. Thus since each of $u$ and $w$ must have a positive puncture each
will have exactly one.
Repeat the above argument at all the other point masses $z_j.$  Then repeat for
any
new point masses in the sequences defining the holomorphic disks $w_j$ associated
to $z_j.$ Continuing until all point masses have been dealt with 
we see no holomorphic curves were overlooked in the reparameterization.

\subsection{Proof of Theorem \ref{9cpt.thm}}
Let $\projL(L)$ denote the limiting Lagrangian boundary condition.
Let $\hbar = \hbar \left( \projL(L) \right)$ be the minimal area
of non-constant maps defined in Section \ref{9area.section}.
Use the discussion in Section \ref{9conformalcompactness.section}
to pass to a subsequence whose conformal structures converge to a stable disk.

We wish to apply Theorem \ref{9bootstrap.thm} to derive strong
local convergence.
To achieve the required uniform bound on 
$\|u_\alpha \|_{k, 2+\delta_{k-1}: K}$ for some compact set
$K \subset \Delta_m$ which lies away from point masses,
we apply Theorem \ref{9unibound.thm} $k$ times to bound
$\|u_\alpha \|_{i, 2+\delta_{k-1}: K}$ for $i = 1, \ldots, k$.
The reparameterizations $\phi_\alpha^j$ in Definition \ref{9conv.dfn}
come from the discussion in Section \ref{9weak.section}.

The weak global convergence follows readily from Section
\ref{9weak.section}.\qed


\end{document}